\documentclass[10pt,leqno]{amsart}
\title{Finiteness obstructions and Euler characteristics of categories}
\author{Thomas M. Fiore}
\author{Wolfgang L\"uck}
\author{Roman Sauer}
\address{Thomas M. Fiore (Corresponding author) \\ Department of Mathematics\\
University of Chicago\\
5734 South University\\
Chicago, IL 60637\\
U.S.A. }
\address{Thomas M. Fiore (Present Address) \\
Department of Mathematics and Statistics \\
University of Michigan-Dearborn \\  4901 Evergreen Road \\
Dearborn, Michigan 48128 \\ U.S.A.} \email{tmfiore@umd.umich.edu}
\urladdr{http://www-personal.umd.umich.edu/~tmfiore/}
\address{Wolfgang L\"uck\\ Fachbereich Mathematik\\ Universit\"at M\"unster\\
Einsteinstr.~62\\ 48149 M\"unster\\Germany\\
FAX: 49 251 8338370}
\email{lueck@math.uni-muenster.de}
\urladdr{http://www.math.uni-muenster.de/u/lueck}
\address{Roman Sauer \\ Fachbereich Mathematik\\ Universit\"at M\"unster\\
Einsteinstr.~62\\ 48149 M\"unster\\Germany\\
FAX: 49 251 8338370} \email{sauerr@uni-muenster.de}
\urladdr{http://www.math.uni-muenster.de/u/sauerr/}

\usepackage{color}

\usepackage{amsmath}
\usepackage{amsthm}
\usepackage{amssymb}
\usepackage{amscd}
\usepackage{amsfonts}  
\usepackage{hyperref}
\usepackage[all]{xy} \SelectTips{cm}{}
\usepackage{graphicx}    
\usepackage{multicol}    

\DeclareMathAlphabet\EuR{U}{eur}{m}{n}
\SetMathAlphabet\EuR{bold}{U}{eur}{b}{n}

\makeindex             

\theoremstyle{plain}
\newtheorem{theorem}{Theorem}[section]
\newtheorem{lemma}[theorem]{Lemma}

\newtheorem{corollary}[theorem]{Corollary}

\theoremstyle{definition}
\newtheorem{definition}[theorem]{Definition}
\newtheorem{example}[theorem]{Example}
\newtheorem{condition}[theorem]{Condition}
\newtheorem{remark}[theorem]{Remark}
\newtheorem{notation}[theorem]{Notation}
\newtheorem{convention}[theorem]{Convention}

{\catcode`@=11\global\let\c@equation=\c@theorem}

\newcommand{\comsquare}[8]                   
{\begin{CD}
#1 @>#2>> #3\\
@V{#4}VV @V{#5}VV\\
#6 @>#7>> #8
\end{CD}
}

\newcommand{\xycomsquareminus}[8]                   
{\xymatrix
{#1 \ar[r]^-{#2} \ar[d]_{#4} &
#3 \ar[d]^{#5}  \\
#6\ar[r]^-{#7} &
#8
}
}

\newcommand{\calfin}{\mathcal{FIN}}

\newcommand{\calb}{{\mathcal B}}
\newcommand{\calc}{{\mathcal C}}

\newcommand{\cale}{{\mathcal E}}
\newcommand{\calf}{{\mathcal F}}
\newcommand{\calg}{{\mathcal G}}

\newcommand{\calj}{{\mathcal J}}

\newcommand{\caln}{{\mathcal N}}

\newcommand{\IA}{{\mathbb A}}

\newcommand{\IC}{{\mathbb C}}

\newcommand{\IN}{{\mathbb N}}

\newcommand{\IQ}{{\mathbb Q}}
\newcommand{\IR}{{\mathbb R}}

\newcommand{\IU}{{\mathbb U}}

\newcommand{\IZ}{{\mathbb Z}}

\newcommand{\curs}{\EuR}
\newcommand{\ABELIANGROUPS}{\curs{ABELIAN}\text{-}\curs{GROUPS}}

\newcommand{\CATS}{\curs{CAT}}

\newcommand{\DIRFINCAT}{\curs{DIR.FIN.}\text{-}\curs{CAT}}
\newcommand{\GROUPOIDS}{\curs{GROUPOIDS}}

\newcommand{\MOD}{\curs{MOD}}
\newcommand{\Or}{\curs{Or}}

\newcommand{\SETS}{\curs{SETS}}

\newcommand{\aut}{\operatorname{aut}}

\newcommand{\CAT}{\operatorname{CAT}}

\newcommand{\ch}{\operatorname{ch}}

\newcommand{\coker}{\operatorname{coker}}
\newcommand{\colim}{\operatorname{colim}}

\newcommand{\hocolim}{\operatorname{hocolim}}
\newcommand{\HS}{\operatorname{HS}}
\newcommand{\id}{\operatorname{id}}

\newcommand{\im}{\operatorname{im}}

\newcommand{\ind}{\operatorname{ind}}

\newcommand{\Irr}{\operatorname{Irr}}
\newcommand{\iso}{\operatorname{iso}}

\newcommand{\map}{\operatorname{map}}
\newcommand{\mor}{\operatorname{mor}}

\newcommand{\ob}{\operatorname{ob}}
\newcommand{\op}{\operatorname{op}}

\newcommand{\pr}{\operatorname{pr}}

\newcommand{\res}{\operatorname{res}}
\newcommand{\Res}{\operatorname{Res}}
\newcommand{\rk}{\operatorname{rk}}

\newcommand{\Split}{\operatorname{Split}}

\newcommand{\tr}{\operatorname{tr}}

\newcommand{\pt}{\{\bullet\}}

\newcommand{\eub}[1]{\underline{E}#1}              
\newcommand{\OrGF}[2]{\Or_{#2}(#1)}                
\newcommand{\uor}[1]{\underline{\Or}(#1)}
\newcommand{\higherlim}[3]{{\setbox1=\hbox{\rm lim}
        \setbox2=\hbox to \wd1{\leftarrowfill} \ht2=0pt \dp2=-1pt
        \mathop{\vtop{\baselineskip=5pt\box1\box2}}
        _{#1}}^{#2}#3}

\newcommand{\version}[1]                       
{\begin{center} last edited on #1\\
last compiled on \today\\
name of tex-file: \jobname
\end{center}
}

\newcounter{commentcounter}

 \newcommand{\squarematrix}[4]{\left( \begin{array}{cc} #1 & #2 \\ #3 &
 #4
 \end{array} \right)
 }

\begin{document}

\typeout{----------------------------  eulcat.tex  ----------------------------}

\maketitle


\typeout{-----------------------  Abstract  ------------------------}

\begin{abstract}
We introduce notions of finiteness obstruction, Euler
characteristic, $L^2$-Euler characteristic, and M\"obius inversion
for wide classes of categories. The finiteness obstruction of a
category $\Gamma$ of type (FP$_R$) is a class in the projective
class group $K_0(R\Gamma)$; the functorial Euler characteristic and
functorial $L^2$-Euler characteristic are respectively its
$R\Gamma$-rank and $L^2$-rank. We
also extend the second author's $K$-theoretic M\"obius inversion
from finite categories to quasi-finite categories. Our main example
is the proper orbit category, for which these invariants are
established notions in the geometry and topology of classifying
spaces for proper group actions. Baez--Dolan's groupoid cardinality
and Leinster's Euler characteristic are special cases of the
$L^2$-Euler characteristic. Some of Leinster's results on
M\"obius-Rota inversion are special cases of the $K$-theoretic
M\"obius inversion.
  \\[2mm]
\end{abstract}


\typeout{--------------------   Section 0: Introduction
--------------------------}

\setcounter{section}{-1}
\section{Introduction and statement of results}

The Euler characteristic is one the earliest and most elementary
homotopy invariants. Though purely combinatorially defined for
finite simplicial complexes as the alternating sum of the numbers of
simplices in each dimension, the Euler characteristic has remarkable
connections to geometry. For example, for closed connected
orientable surfaces $M$, the Euler characteristic determines the
genus: $g=1-\frac{1}{2}\chi(M)$. For such $M$, if $\chi(M)$ is
negative, then $M$ admits a hyperbolic metric. More substantially,
the celebrated Gauss--Bonnet Theorem computes the Euler
characteristic in terms of curvature. A further example of geometry
in the Euler characteristic is provided by the Hopf-Singer
conjecture.

Of course, Euler characteristics are not only defined for finite
simplicial complexes or manifolds, but also for a great variety of
objects, such as equivariant spaces, orbifolds, or finite posets.
Baez--Dolan considered in~\cite{Baez-Dolan(2001)} an Euler
characteristic (groupoid cardinality) for finite groupoids and
certain infinite ones, such as the groupoid of finite sets. Leinster
and Berger--Leinster have considered Euler characteristics not just
of finite posets and groupoids, but more generally of finite
categories
in~\cite{Leinster(2008)}~and~\cite{Berger+Leinster(2008)}. If a
finite category admits both a weighting and coweighting, then it
admits an Euler characteristic in the sense of Leinster.

In the present paper, we define Euler characteristics for wide
classes of categories, provide a unified conceptual framework in
terms of finiteness obstructions and projective class groups, and
extract geometric and algebraic information from our invariants in
certain cases. This obstruction-theoretic framework works well for
both finite and infinite categories. Our main example is the proper
orbit category of a group $G$. In this case, our invariants are
established geometric invariants of the classifying space for proper
$G$-actions. We also extend the second author's $K$-theoretic
M\"obius inversion from finite EI-categories to quasi-finite
EI-categories (a category $\Gamma$ is said to be EI if each
endomorphism in $\Gamma$ is an isomorphism). The $K$-theoretic
M\"obius inversion does not require the categories in question to be
skeletal, unlike the M\"obius inversion of
Leinster~\cite{Leinster(2008)}. Several of the results
of~\cite{Leinster(2008)} are special cases.

Our point of departure is the theory of projective modules over a
category and the associated projective class group. Let $\Gamma$ be
a small category, and $R$ an associative commutative ring with
identity. An \emph{$R\Gamma$-module} is a functor from
$\Gamma^{\op}$ to the abelian category of left $R$-modules. If
$\Gamma$ is a group $G$ viewed as a one-object category, then an
$R\Gamma$-module is nothing more than a right $RG$-module.  The
category $\MOD\text{-}R\Gamma$ of $R\Gamma$-modules is an abelian
category, and therefore we automatically have the notions of
projective $R\Gamma$-module, chain complexes of $R\Gamma$-modules,
and resolutions of $R\Gamma$-modules. The finiteness obstruction,
whenever it exists, lives in the \emph{projective class group
$K_0(R\Gamma)$}, which is the free abelian group on the isomorphism
classes of finitely generated projective $R\Gamma$-modules modulo
short exact sequences. We say that $\Gamma$ is \emph{of type
(FP$_R$)} if the constant $R\Gamma$-module $\underline{R}\colon
\Gamma^{\op} \to R\text{-}\MOD$ admits a resolution by finitely
generated projective $R\Gamma$-modules in which only finitely many
of the $R\Gamma$-modules are nonzero. If $\Gamma$ is of type
(FP$_R$), the \emph{finiteness obstruction} $o(\Gamma;R) \in
K_0(R\Gamma)$ is the alternating sum of the classes of modules
appearing in a finite projective resolution of $\underline{R}$. For
example, if $\Gamma$ is a finite group of order invertible in $R$,
then $\underline{R}$ is itself a projective $R \Gamma$-module,
$\underline{R}$ provides us with a finite projective resolution of
$\underline{R}$, and $[\underline{R}]$ is the finiteness obstruction
$o(\Gamma;R)$. Further examples of categories of type (FP$_R$) are
provided by any finite EI-category such that $\vert \aut(x) \vert$
is invertible in $R$ for each object $x$, and any category $\Gamma$
which admits a finite $\Gamma$-$CW$-model for
$E\Gamma$. The basics of $R\Gamma$-modules and finiteness obstructions are
discussed in
Sections~\ref{sec:basics_about_modules_over_a_category}~and~\ref{sec:The_finiteness_obstruction_of_a_category}.

To obtain the Euler characteristic and the $L^2$-Euler
characteristic from the finiteness obstruction, we use L\"uck's
Splitting of $K_0$~\cite[Theorem~10.34 on page~196]{Lueck(1989)},
and two notions of rank for $R\Gamma$-modules: the $R\Gamma$-rank
$\rk_{R\Gamma}$ and the $L^2$-rank $\rk_\Gamma^{(2)}$. In the case
that every endomorphism in $\Gamma$ is an isomorphism, that is,
$\Gamma$ is an EI-category, L\"uck constructed in~\cite{Lueck(1989)}
the natural \emph{splitting} isomorphism \[S \colon K_0(R\Gamma) \to
\Split K_0(R\Gamma) := \bigoplus_{\overline{x} \in \iso(\Gamma)}
K_0(R\aut(x))\] and its natural inverse $E$, called
\emph{extension}. In
Section~\ref{sec:Splitting_the_projective_class_group} we recall the
splitting $(S,E)$, and prove that $S$ remains a left inverse to $E$
in the more general case of directly finite $\Gamma$. Let $S_x$
denote the $\overline{x}$-component of $S$ and let $U(\Gamma)$
denote the free abelian group on the isomorphism classes of objects
of $\Gamma$. The \emph{$R\Gamma$-rank} of a finitely generated
$R\Gamma$-module $M$ is the element $\rk_{R\Gamma}M \in U(\Gamma)$
which is $\rk_R\left(S_xM \otimes_{R\aut(x)} R\right)$ at
$\overline{x} \in \iso(\Gamma)$. This induces a homomorphism
$\rk_{R\Gamma} \colon K_0(R\Gamma) \to U(\Gamma)$. If $\Gamma$ is of
type (FP$_R$), we define the \emph{functorial Euler characteristic}
$\chi_f(\Gamma;R)$ to be the image of the finiteness obstruction
$o(\Gamma;R)$ under $\rk_{R\Gamma}$. The sum of the components of
$\chi_f(\Gamma;R)$ is called the \emph{Euler characteristic of
$\Gamma$}, denoted $\chi(\Gamma;R)$. Indeed,
if $R$ is Noetherian, and $\Gamma$ is \emph{directly finite} in
addition to type (FP$_R$), then $\chi(\Gamma;R)$ coincides with the
\emph{topological Euler characteristic} $\chi(B\Gamma;R)$. For
example, if $\Gamma$ is a finite group, then $\chi_f(\Gamma;\IQ)$ is
1, and so is the rational Euler characteristic. In
Section~\ref{sec:The_(functorial)_Euler_characteristic_of_a_category}
we treat the topological Euler characteristic $\chi(B\Gamma;R)$, the
$R\Gamma$-rank $\rk_{R\Gamma}$, the functorial Euler characteristic
$\chi_f(\Gamma;R)$, and the Euler characteristic $\chi(\Gamma;R)$.

To obtain the $L^2$-Euler characteristic from the finiteness
obstruction using the splitting functor $S_x$ and the $L^2$-rank
$\rk_\Gamma^{(2)}$, we need some elementary theory of finite von
Neumann algebras. For a group $G$, the \emph{group von Neumann
algebra of $G$} is the algebra of $G$-equivariant bounded operators
$\ell^2(G) \to \ell^2(G)$, which we denote by $\caln(G)$. If $G$ is
a finite group, $\caln(G)$ is simply the group ring $\IC G$. The
\emph{von Neumann dimension for $\caln(G)$-modules} is the unique
function $\dim_{\caln(G)}$ satisfying Hattori-Stallings rank,
additivity, cofinality, and continuity as recalled in
Theorem~\ref{the:properties_of_the_dimension_function}. In the case
of a finite group $G$, the von Neumann dimension of a $\IC G$-module
is the complex dimension divided by $|G|$. The \emph{$L^2$-rank} of
a finitely generated $\IC\Gamma$-module $M$ is the element
$\rk_{\Gamma}^{(2)}M \in U(\Gamma) \otimes_\IZ \IR$ which is
$\dim_{\caln(\aut(x))}\left(S_x M \otimes_{\IC\aut(x)}
\caln(\aut(x))\right)$ at $\overline{x} \in \iso(\Gamma)$. This
induces a homomorphism $\rk_{\Gamma}^{(2)} \colon K_0(R\Gamma) \to
U(\Gamma)\otimes_\IZ \IR$. If $\Gamma$ is of type (FP$_\IC$), the
\emph{functorial $L^2$-Euler characteristic} $\chi_f^{(2)}(\Gamma)$
is the image of the finiteness obstruction $o(\Gamma;\IC)$ under
$\rk_{\Gamma}^{(2)}$. The \emph{$L^2$-Euler characteristic}
$\chi^{(2)}(\Gamma)$ is the sum of the components of
$\chi_f^{(2)}(\Gamma)$. For example, if $\Gamma$ is a finite
groupoid of type (FP$_\IC$), its functorial $L^2$-Euler
characteristic has at $\overline{x}$ the value $1/|\aut(x)|$, and
the $L^2$-Euler characteristic is the sum of these. This agrees with
the \emph{groupoid cardinality} of
Baez--Dolan~\cite{Baez-Dolan(2001)} and also Leinster's Euler
characteristic in the case of finite groupoids. If $\Gamma$ is
directly finite and of type (FF$_\IZ)$, and $R$ is Noetherian, then
$\chi(B\Gamma;R)=\chi(\Gamma;R)=\chi^{(2)}(\Gamma)$. In
Section~\ref{sec:The_(functorial)_L2_Euler_characteristic_anf_L2-Betti_numbers_of_a_category}
we review the necessary prerequisites from the theory of finite von
Neumann algebras, and introduce the $L^2$-rank $\rk_\Gamma^{(2)}$,
the functorial $L^2$-Euler characteristic $\chi_f^{(2)}(\Gamma)$,
and the $L^2$-Euler characteristic $\chi^{(2)}(\Gamma)$. These are
defined for categories of type ($L^2$), a slightly weaker
requirement than type (FP$_\IC$).

The invariants we introduce in this paper have many desirable
properties. The finiteness obstruction, functorial Euler
characteristic, Euler characteristic, functorial $L^2$-Euler
characteristic, and $L^2$-Euler characteristic are all invariant
under equivalence of categories and are compatible with finite
products, finite coproducts, and homotopy colimits (see
Fiore--L\"uck--Sauer~\cite{FioreLueckSauerHoColim(2009)} for the
compatibility with homotopy colimits). Moreover, the $L^2$-Euler
characteristic is compatible with isofibrations and coverings
between finite groupoids (see
Subsection~\ref{subsec:Compatibility_with_Coverings}). The
$L^2$-Euler characteristic coincides with the classical $L^2$-Euler
characteristic in the case of a group, for finite groups this is
$\chi^{(2)}(G)=\frac{1}{|G|}$. Another advantage of the $L^2$-Euler
characteristic is that it is closely related to the geometry and
topology of the classifying space for proper $G$-actions, a topic to
which we return in Section~\ref{sec:The_proper_orbit_category}.

After this treatment of finiteness obstructions and various Euler
characteristics, we turn in Section~\ref{sec:Moebius-inversion} to
our next main result: the generalization of the second author's
$K$-theoretic M\"obius inversion to quasi-finite EI-categories. We
introduce the \emph{restriction-inclusion splitting} $\Res \colon
K_0(R\Gamma) \rightleftarrows \Split K_0(R\Gamma) \colon I$ in
Subsection~\ref{subsec:A_second_Splitting}. The \emph{$K$-theoretic
M\"obius inversion} \[\mu \colon \Split K_0(R\Gamma)
\rightleftarrows \Split K_0(R\Gamma) \colon \omega\] compares the
splitting $(\Res,I)$ with the splitting $(S,E)$ in
Theorem~\ref{the:Two_splittings_and_the_K-theoretic_Moebius_inversion}.
See Subsection~\ref{subsec:The-K-theoretic-Moebius-Inversion} for
the definition of $(\mu,\omega)$ in terms of chains in $\Gamma$ and
hom-sets of $\Gamma$. A computationally useful byproduct of the
comparison via M\"obius inversion is the equation
\[S\left(o(\Gamma;R)\right)=\mu\left(\left(o(\widehat{\aut(x)};R)\right)_{\overline{x}
\in \iso(\Gamma)}\right)\] for $\Gamma$ of type (FP$_R$). For
example, this enables us to compute in
Theorem~\ref{the:The_finiteness_obstruction_and_Moebius_inversion_finite_EI-categories}
the finiteness obstruction and Euler characteristics of a finite
EI-category in terms of chains. The $K$-theoretic M\"obius inversion
is also compatible with the $L^2$-rank $\rk_{\Gamma}^{(2)}$ and the
pair $(\overline{\mu}^{(2)}, \overline{\omega}^{(2)})$ as in
Subsection~\ref{subsec:The-K-theoretic-Moebius-Inversion_and_the_L2-rank}.
All of these splittings and homomorphisms are illustrated explicitly
for $G$-$H$-bisets in
Subsection~\ref{subsec:The_example_of_a_biset}. The rest of
Section~\ref{sec:Moebius-inversion} compares and contrasts the
invariants for $\Gamma$ and $\Gamma^{\op}$, which can generally be
quite different. Important special cases are \emph{M\"obius-Rota
inversion} for a finite partially ordered set
(Example~\ref{Finite_partially_ordered_set}), \emph{Leinster's
M\"obius inversion} for a finite skeletal category with trivial
endomorphisms
(Example~\ref{exa:Moebius_Inversion_for_fin_sk_cat_with_triv_endos}),
and \emph{rational M\"obius inversion} for a finite, skeletal, free
EI-category
(Example~\ref{exa:Rational_Moebius_Inversion_for_fin_sk_cat_free_EI}).

In Section~\ref{sec:Comparison_with_Leinsters_invariant} we recall
the groupoid cardinality of Baez--Dolan~\cite{Baez-Dolan(2001)} and
the Euler characteristic of Leinster~\cite{Leinster(2008)} and make
comparisons. The groupoid cardinality coincides with the $L^2$-Euler
characteristic for finite groupoids. Leinster's Euler characteristic
coincides with the $L^2$-Euler characteristic for finite, free,
skeletal EI-categories. Here ``free'' is not meant in the usual
category-theoretic sense, but rather in the sense of group actions.
We say that a category $\Gamma$ is \emph{free} if the left
$\aut(y)$-action on $\mor(x,y)$ is free for every two objects $x,y
\in \ob(\Gamma)$. If $\Gamma$ is not free, then $\chi^{(2)}(\Gamma)$
could very well be different from Leinster's Euler characteristic of
$\Gamma$ (see Remark~\ref{rem_free_is_necessary}). Our invariants
are more sensitive than Leinster's Euler characterstic. For example,
Leinster's Euler characteristic for finite categories only depends
on the set of objects $\ob(\Gamma)$ and the orders
$|\mor_\Gamma(x,y)|$. As such, it cannot distinguish between the
group $\IZ/2\IZ$ and the two-element monoid consisting of the
identity and an idempotent. The finiteness obstruction and the
$L^2$-Euler characteristic can distinguish these. Leinster's Euler
characteristic cannot distinguish between $\Gamma$ and
$\Gamma^{\op}$, while the functorial Euler characteristic, the
functorial $L^2$-Euler characteristic, and the $L^2$-Euler
characteristic can. In
Section~\ref{sec:Comparison_with_Leinsters_invariant} we also
explain how to construct weightings in the sense of Leinster from
finite free resolutions of the constant $R\Gamma$-module
$\underline{R}$ as well as from finite $\Gamma$-$CW$-models for the
classifying $\Gamma$-space. Several of the weightings in Leinster's
article~\cite{Leinster(2008)} arise in this way.

As mentioned at the outset, Euler characteristics of spaces and
manifolds contain geometric information, such as genus, curvature,
or evidence of a hyperbolic metric. Similarly, the Euler
characteristics of certain categories contain geometric and
algebraic information. The topic of
Section~\ref{sec:The_proper_orbit_category} is our main example: the
\emph{proper orbit category of a group $G$}, denoted $\uor{G}$. Its
objects are the homogeneous sets $G/H$ for finite subgroups $H$ of
$G$, and its morphisms are the $G$-equivariant maps between such
homogeneous sets. The invariants of the category $\uor{G}$ are
closely related to the equivariant invariants of a model $\eub{G}$
for the classifying space for proper $G$-actions. Namely, if the
model $\eub{G}$ is a finitely dominated $G$-$CW$-complex, then our
category-theoretic finiteness obstruction $o(\uor{G}; \IZ)$ agrees
with the equivariant finiteness obstruction of $\eub{G}$. If the
model $\eub{G}$ is even a finite $G$-$CW$-complex, then both the
functorial Euler characteristic $\chi_f(\uor{G};\IZ)$ and the
functorial $L^2$-Euler characteristic $\chi_f^{(2)}(\uor{G})$ agree
with the equivariant Euler characteristic of $\eub{G}$. Examples of
groups $G$ with finite models $\eub{G}$ include hyperbolic groups,
groups that act simplicially cocompactly and properly by isometries
on a CAT(0)-space, mapping class groups, the group of outer
automorphisms of a finitely generated free group, finitely generated
one-relator groups, and cocompact lattices in connected Lie groups.

In addition to these geometric aspects of our invariants in the case
of the category $\uor{G}$, we also have interesting algebraic
consequences of the $K$-theoretic M\"obius inversion and its
compatibility with the $L^2$-rank. For example, if the category
$\uor{G}$ is of type (FP$_\IQ$) and satisfies condition (I) of
Condition \ref{con:condition_(I)}, then the functorial $L^2$-Euler
characteristic of $\uor{G}$ is the $L^2$-M\"obius inversion of the
$L^2$-Euler characteristics of Weyl groups associated to finite $H <
G$: \[\chi_f^{(2)}(\uor{G}) =
\overline{\mu}^{(2)}\biggl(\bigl(\chi^{(2)}(W_GH)\bigr)_{(H), |H| <
\infty}\biggr).\] More substantially, for finite $G$ we deduce the
\emph{Burnside ring congruences}, which distinguish the image of the
character map  \[\ch = \ch^G \colon U(\uor{G})  \to \bigoplus_{(H)}
\IZ.\] Here $U(\uor{G})$ is the free abelian group on the set of
isomorphism classes of objects in $\uor{G}$, we identify
$U(\uor{G})$ with the Burnside ring $A(G)$, and the direct sum of
$\IZ$'s is over the conjugacy classes $(H)$ of subgroups of the
finite group $G$. The character map counts $H$-fixed points, namely,
for any finite $G$-set $S$ we have $\ch(S)=\bigl(|S^H|\bigr)_{(H)}$.
An element $\xi$ lies in the image of $\ch$ if and only if the
integral congruence \[\nu(\xi)_{(H)} \equiv 0 \mod |W_GH|\] holds
for every conjugacy class $(H)$ of subgroups of $G$ (the matrix
$\nu$ is specified in
Subsection~\ref{subsec:The_Burnside_integrality_relations_and_the_classical_Burnside_congruences}).
We finish Section~\ref{sec:The_proper_orbit_category} by working out
everything explicitly for the infinite dihedral group.

The last two sections of the paper are explicit examples. In
Section~\ref{sec:An_example_of_a_finite_category_without_property_EI}
we consider a small example of a category which is not EI and
calculate its various $K$-theoretic morphisms: the splitting functor
$S$, the extension functor $E$, the restriction functor $\Res$, and
the homomorphism $\omega$. In
Section~\ref{sec:A_finite_category_without_property_(FP)} we
consider a category $\IA$ which does not satisfy property (FP$_R$).
Leinster considered this category in Example~1.11.d
of~\cite{Leinster(2008)} and proved that it does not admit a
weighting. We prove that $\IA$ does not satisfy property (FP$_R$),
classify the finitely generated projective $R\IA$-modules, and
compute the projective class group $K_0(R\IA)$, the Grothendieck
group of finitely generated $\IQ \IA$-modules $G_0(\IQ \IA)$, and
the homology $H_n(B\IA;R)=H_n(\IA;R)$.

\section*{Acknowledgements}
Thomas M. Fiore was supported
at the University of Chicago by NSF Grant DMS-0501208. At the
Universitat Aut\`{o}noma de Barcelona he was supported by Grant
SB2006-0085 of the Spanish Ministerio de Educaci\'{o}n y Ciencia
under the Programa Nacional de ayudas para la movilidad de
profesores de universidad e investigadores
espa$\tilde{\text{n}}$oles y extranjeros. Thomas M. Fiore also
thanks the Centre de Recerca Matem\`{a}tica in Bellaterra
(Barcelona) for its hospitality during the CRM Research Program on
Higher Categories and Homotopy Theory in 2007-2008, where Tom
Leinster spoke about Euler characteristics.

Wolfgang L\"{u}ck was financially supported by the
Sonderforschungsbereich 478 \--- Geometrische Strukturen in der
Mathematik \--- and his Max-Planck-Forschungspreis.

Roman Sauer is grateful for support from the \emph{Deutsche
Forschungsgemeinschaft}, made through Grant SA 1661/1-2.

The authors were also supported by the Leibniz award of the second author.

\smallskip

We thank the referee for helpful corrections.

\tableofcontents


\typeout{-------  Section 1: Basics about modules over a category ----------------------}
\section{Basics about modules over a category}
\label{sec:basics_about_modules_over_a_category}

Throughout this paper, let $\Gamma$ be a small category and let $R$
be an associative, commutative ring with identity. We explain
some basics about modules over a category. More details can be found
in L\"uck\cite[Section~9]{Lueck(1989)}. An \emph{$R\Gamma$-module}
is a functor from $\Gamma^{\op}$ into the abelian category of left
$R$-modules. This is a natural generalization of the notion of right
$RG$-module for a group $G$. The category of $R\Gamma$-modules forms
an abelian category $\MOD\text{-}R\Gamma$. An object of
$\MOD\text{-}R\Gamma$ is projective if and only if it is a direct
summand in an $R\Gamma$-module which is free on a collection of sets
indexed by $\ob(\Gamma)$. Given a functor $F \colon \Gamma_1 \to
\Gamma_2$, we have induction and restriction functors $\ind_F \colon
\MOD\text{-}R\Gamma_1 \rightleftarrows \MOD\text{-}R\Gamma_2 \colon
\res_F$, and these are adjoint. We also introduce in this section
the projective class group $K_0(R\Gamma)$, which provides a home for
the finiteness obstruction $o(\Gamma;R)$. The \emph{projective class
group $K_0(R\Gamma)$} is the free abelian group on the isomorphism
classes of finitely generated projective $R\Gamma$-modules modulo
short exact sequences. The induction functor induces a homomorphism
of projective class groups, as does the restriction functor,
provided $F$ is admissible.

\begin{definition}[Modules over a category]
  \label{def:modules_over_a-category}
  A \emph{(contravariant) $R\Gamma$-module} is a contravariant functor
  $\Gamma \to R\text{-}\MOD$ from $\Gamma$ to the abelian
  category of $R$-modules.  A \emph {morphism of $R\Gamma$-modules}
  is a natural transformation of such functors.
  We denote by $\MOD\text{-}R\Gamma$ the category of (contravariant)
  $R\Gamma$-modules.
\end{definition}

\begin{example}[Modules over group rings]
  \label{exa:modules_over_group_rings}
  Let $G$ be a group. Let $\widehat{G}$ be the associated groupoid
  with one object and $G$ as its set of morphisms with the obvious composition law.
  Then the category
  $\MOD\text{-}R\widehat{G}$ of contravariant $R\widehat{G}$-modules
  agrees with the category of right $RG$-modules, where $RG$ is the
  group ring of $G$ with coefficients in $R$.
\end{example}

\begin{example} \label{exa:diagrams} Let $\Gamma$ be the category
  having one object and the natural numbers $\IN = \{0,1,2,\ldots \}$
  as morphisms with the obvious composition law.  Then
  $\MOD\text{-}R\Gamma$ is the category whose objects are
  endomorphisms of $R$-modules and whose set of morphisms from an
  endomorphism $f$ to an endomorphism $g$ is given by the set of
  commutative diagrams
$$\xymatrix{M \ar[r]^f \ar[d]_u & M \ar[d]^u\\
  N \ar[r]_g & N}
$$
If one replaces $\IN$ by $\IZ$ and endomorphisms by automorphisms, the
corresponding statement holds.
\end{example}

The (standard) structure of an abelian category on $R\text{-}\MOD$
induces the structure of an abelian category on $\MOD\text{-}R\Gamma$
in the obvious way, namely objectwise.
In particular, the notion of a projective $R\Gamma$-module is
defined. Namely, an $R\Gamma$-module $P$ is \emph{projective} if for
every surjective $R\Gamma$-morphism $p \colon M \to N$ and
$R\Gamma$-morphism $f \colon P \to N$ there exists an $R\Gamma$-morphism
$\overline{f} \colon P \to M$ such that $p \circ \overline{f} = f$,
where $p$ is called \emph{surjective} if for any object $x \in \Gamma$ the
$R$-homomorphism $p(x) \colon M(x) \to N(x)$ is surjective.

Consider an object $x$ in $\Gamma$. For a set $C$
we denote by $RC$ the free module with $C$ as basis, i.e., the $R$-module of
maps with finite support from $C$ to $R$.  Denote by
\begin{eqnarray}
& R\mor(?,x) & \text{for}\; x \in \ob(\Gamma)
\label{Rhom(?,x)}
\end{eqnarray}
the $R\Gamma$-module which sends an object $y$ to the $R$-module
$R\mor(y,x)$, and a morphism $u \colon y \to z$ to the $R$-map induced
by the morphism of sets $\mor(z,x) \to \mor(y,x)$ that maps
$v \colon z \to x$ to $v \circ u \colon y \to x$.

\begin{lemma} \label{lem:elementary_property_of_rmor(?,x))} Let $M$ be
  any $R\Gamma$-module. Consider any element $\alpha \in M(x)$.  Then
  there is precisely one map of $R\Gamma$-modules
$$F_{\alpha} \colon R\mor(?,x) \to M$$
such that $F_{\alpha}(x) \colon R\mor(x,x) \to M(x)$ sends $\id_x$ to
$\alpha$.
\end{lemma}
\begin{proof}
  This is a direct application of the Yoneda Lemma. Given $u \colon y
  \to x$, define $F_{\alpha}(u) := M(u)(\alpha)$.
\end{proof}

Since $\Gamma$ is by assumption small, its objects form a set denoted
by $\ob(\Gamma)$. An \emph{$\ob(\Gamma)$-set} $C$ is a collection of
sets $C = \{C_x \mid x \in \ob(\Gamma)\}$ indexed by $\ob(\Gamma)$. A
\emph{morphism of $\ob(\Gamma)$-sets} $f \colon C \to D$ is a
collection of maps of sets $\{f_x \colon C_x \to D_x \mid x \in
\ob(\Gamma)\}$. Denote by $\ob(\Gamma)\text{-}\SETS$ the category of
$\ob(\Gamma)$-sets. We obtain an obvious \emph{forgetful functor}
$$F \colon \MOD\text{-}R\Gamma \to \ob(\Gamma)\text{-}\SETS.$$
Let
$$B \colon \ob(\Gamma)\text{-}\SETS \to  \MOD\text{-}R\Gamma$$
be the functor sending an $\ob(\Gamma)$-set $C$ to the
$R\Gamma$-module
\begin{equation} \label{equ:free_RGamma_module}
B(C) := \bigoplus_{x \in \ob(\Gamma)} \bigoplus_{C_x} R\mor(?,x).
\end{equation}
We call $B(C)$ the \emph{free $R\Gamma$-module with basis the
$\ob(\Gamma)$-set $C$}.  This name is justified by the following
consequence of Lemma~\ref{lem:elementary_property_of_rmor(?,x))} and
the universal property of the direct sum.

\begin{lemma} \label{lem:adjunction_B_and_F} We obtain a pair of
  adjoint functors by $(B,F)$.
\end{lemma}

Lemma~\ref{lem:adjunction_B_and_F} implies that the abelian category
$\MOD\text{-}R\Gamma$ has enough projectives. Namely, any free
$R\Gamma$-module is projective and for any  $R\Gamma$-module $M$
there is a surjective morphism of $R\Gamma$-modules $B(F(M)) \to M$,
given by the adjoint of $\id \colon F(M) \to F(M)$.  Therefore the
standard machinery of homological algebra applies to
$\MOD\text{-}R\Gamma$.  We also conclude that an $R\Gamma$-module is
projective if and only if it is a direct summand in a free
$R\Gamma$-module.

An $\ob(\Gamma)$-set $C$ is  \emph{finite} if the cardinality of
$\coprod_{x \in \ob(\Gamma)} C_x$ is finite. An $R\Gamma$-module $M$
is \emph{finitely generated} if and only if there is a finite
$\ob(\Gamma)$-set $C$ together with a surjective $R\Gamma$-morphism $B(C)
\to M$.  An $R\Gamma$ module is finitely generated projective if and only
if it is a direct summand in free $R\Gamma$-module
$B(C)$ for a finite $\ob(\Gamma)$-set $C$.

\begin{definition} \label{def:tensor_product}
If $M\colon \Gamma^{\op} \to R\text{-}\MOD$  and $N\colon \Gamma \to
R\text{-}\MOD$ are functors, then the \emph{tensor product} $M
\otimes_{R\Gamma} N$ is the quotient of the $R$-module
\begin{equation*}
 \bigoplus_{x \in \ob(\Gamma)} M(x) \otimes_R N(x)
\end{equation*}
by the $R$-submodule generated by elements of the form
\begin{equation*}
(M(f)m) \otimes n - m \otimes (N(f)n)
\end{equation*}
where $f:x \to y$ is a morphism in $\Gamma$, $m \in M(y)$, and $n
\in N(x)$. The tensor product is an $R$-module, not an
$R\Gamma$-module.
\end{definition}

\begin{definition}[Projective class group]
  \label{def:projective_class_group}
    The \emph{projective class group $K_0(R\Gamma)$}
    is the abelian group whose generators $[P]$ are
    isomorphism classes of finitely generated projective
    $R\Gamma$-modules and whose relations are given by expressions
    $[P_0] - [P_1] + [P_2] = 0$ for every exact sequence $0 \to P_0
    \to P_1 \to P_2 \to 0$ of finitely generated projective
    $R\Gamma$-modules.
  \end{definition}

Given a functor $F \colon \Gamma_1 \to \Gamma_2$, \emph{induction with
  $F$} is the functor
\begin{eqnarray}
  \ind_F \colon \MOD\text{-}R\Gamma_1 & \to & \MOD\text{-}R\Gamma_2
  \label{ind_F}
\end{eqnarray}
which sends a contravariant $R\Gamma_1$-module $M = M(?)$ to the
contravariant $R\Gamma_2$-module $M(?) \otimes_{R\Gamma_1}
R\mor_{\Gamma_2}(??,F(?))$ which is the tensor product over
$R\Gamma_1$ with the $R\Gamma_1\text{-}R\Gamma_2$-bimodule
$R\mor_{\Gamma_2}(??,F(?))$ (see~L\"uck~\cite[9.15 on
page~166]{Lueck(1989)} for more details).  The functor $\ind_F$
respects direct sums over arbitrary index sets and satisfies $\ind_F
(R\mor_{\Gamma_1}(?,x)) = R\mor_{\Gamma_2}(??,F(x))$ for every $x
\in \ob(\Gamma_1)$. Hence $\ind_F$ sends finitely generated
$R\Gamma_1$-modules to finitely generated $R\Gamma_2$-modules and
sends projective $R\Gamma_1$-modules to projective
$R\Gamma_2$-modules. The
functor $\ind_F$ induces a homomorphism
\begin{eqnarray} \label{equ:induction_on_K0}
  F_* \colon K_0(R\Gamma_1) & \to & K_0(R\Gamma_2),
  \label{F_ast_colon_K_0_to_K_0}
\end{eqnarray}
which depends only on the natural isomorphism class of $F$.  Given
functors $F_0 \colon \Gamma_0 \to \Gamma_1$ and $F_1 \colon \Gamma_1
\to \Gamma_2$, the functors of abelian categories
$\ind_{F_1 \circ  F_0}$ and $\ind_{F_1} \circ \ind_{F_0}$ are naturally isomorphic and
hence $(F_1 \circ F_0)_* = (F_1)_* \circ (F_0)_*$.

Given a functor $F \colon \Gamma_1 \to \Gamma_2$, \emph{restriction
  with F} is the functor of abelian categories
\begin{eqnarray}
  \res_F \colon \MOD\text{-}R\Gamma_2 & \to & \MOD\text{-}R\Gamma_1, \quad M \mapsto M \circ F.
  \label{res_F}
\end{eqnarray}
It is exact and sends the constant $R\Gamma_2$-module
$\underline{R}$ to the constant $R\Gamma_1$-module $\underline{R}$.
In general it does not send a finitely generated projective
$R\Gamma_2$-module to a finitely generated projective
$R\Gamma_1$-module.  We call $F$ \emph{admissible} if $\res_F$ sends
a finitely generated projective $R\Gamma_2$-module to a finitely
generated projective $R\Gamma_1$-module. The question when $F$ is
admissible is answered in~L\"uck~\cite[Proposition~10.16 on
page~187]{Lueck(1989)}.  If $F$ is admissible, it induces a
homomorphism
\begin{eqnarray}
  F^* \colon K_0(R\Gamma_2) & \to & K_0(R\Gamma_1).
  \label{F_upperast_colon_K_0_to_K_0}
\end{eqnarray}

The following is proved in~L\"uck~\cite[9.22 on
page~169]{Lueck(1989)} and is based on the fact that $\res_F$ is the
same as the functor $- \otimes_{R\Gamma_2}
R\mor_{\Gamma_2}(F(?),??)$.

\begin{lemma} \label{lem:adjunction_ind_F_and_res_F}
  Given a functor $F \colon \Gamma_0 \to \Gamma_1$, we obtain an
  adjoint pair of functors $(\ind_F,\res_F)$.
\end{lemma}


\typeout{-------  Section 2: The finiteness obstruction of a category}
\section{The finiteness obstruction of a category}
\label{sec:The_finiteness_obstruction_of_a_category}

After the introduction to $R\Gamma$-modules in Section
\ref{sec:basics_about_modules_over_a_category}, we can now define
the finiteness obstruction of a category in terms of chain complexes
and establish its basic properties. Since $\MOD\text{-}R\Gamma$ is
abelian, we can talk about $R\Gamma$-chain complexes.  In the sequel
all chain complexes $C_*$ will satisfy $C_n = 0$ for $n \le -1$.  A
\emph{finite projective
  $R\Gamma$-chain complex} $P_*$ is an $R\Gamma$-chain complex such
there exists a natural number $N$ with $P_n = 0$ for $n > N$ and
each $R\Gamma$-module $P_i$ is finitely generated projective. Let
$M$ be an $R\Gamma$-module. A \emph{finite projective
$R\Gamma$-resolution} of $M$ is a finite projective $R\Gamma$-chain
complex $P_*$ satisfying $H_n(P_*) = 0$ for $n \ge 1$ together with
an isomorphism of $R\Gamma$-modules $M \xrightarrow{\cong}
H_0(P_*)$. If $P_*$ can be chosen as a finite free $R\Gamma$-chain
complex, we call it a \emph{finite free $R\Gamma$-resolution}.

If the constant $R\Gamma$-module $\underline{R}\colon \Gamma^{\op}
\to R\text{-}\MOD$ with value $R$ admits a finite projective
$R\Gamma$-resolution or a finite free $R\Gamma$-resolution, we say
that $\Gamma$ is \emph{of type (FP$_R$)} or \emph{of type (FF$_R$)}
respectively. Examples of categories of type (FP$_R$) are: any
finite group of order invertible in $R$, and more generally, any
finite category in which every endomorphism is an isomorphism and
$\vert \aut_\Gamma(x)\vert$ is invertible in $R$ for each object
$x$. Any category $\Gamma$ which admits a finite $\Gamma$-$CW$-model
for $E\Gamma$ is of type (FF$_R$) and therefore of type (FP$_R$), in
particular any category with a terminal object is of type (FF$_R$)
and (FP$_R$).

If $\Gamma$ is of type (FP$_R$), we define the \emph{finiteness
obstruction $o(\Gamma;R) \in K_0(R\Gamma)$} to be the alternating
sum of the classes $[P_n]$ appearing in a finite projective
resolution of $\underline{R}$. If $G$ is a finitely presented group
of type (FP$_\IZ$), then the finiteness obstruction is the same as
Wall's finiteness obstruction $o(BG) \in K_0(\IZ G)$.

Type (FP$_R$) and the finiteness obstruction have all the properties
one could hope for. Any category equivalent to a category of type
(FP$_R$) is also of type (FP$_R$), and the induced map of an
equivalence preserves the finiteness obstruction. If $\Gamma_1$ and
$\Gamma_2$ are of type (FP$_R$), then so are $\Gamma_1 \times
\Gamma_2$ and $\Gamma_1 \amalg \Gamma_2$, and the finiteness
obstructions behave accordingly. Restriction along admissible
functors preserves type (FP$_R$) and finiteness obstructions, as
does induction along right adjoints.
In~\cite{FioreLueckSauerHoColim(2009)}, we prove that type
(FP$_R$), type (FF$_R$), and the finiteness obstruction are
compatible with homotopy colimits.

\begin{definition}[Finiteness obstruction of an $R\Gamma$-module]
  \label{def:finiteness_obstruction_of_a_module}
  Let $M$ be an $R\Gamma$-module which possesses a finite projective
  resolution. The \emph{finiteness obstruction of $M$} is
$$o(M) := \sum_{n \ge 0} (-1)^n \cdot [P_n] \; \in K_0(R\Gamma),$$
where $P_*$ is any choice of a finite projective
$R\Gamma$-resolution of $M$.
\end{definition}

This definition is a special case of~L\"uck~\cite[Definition~11.1 on
page~211]{Lueck(1989)}.  It is indeed independent of the choice of
finite projective resolution.  If $P$  is a finitely generated
projective $R\Gamma$-module, then of course $o(P) = [P]$. Given an
exact sequence $0 \to M_0 \to M_1 \to M_2 \to 0$ of
$R\Gamma$-modules such that two of them possess finite projective
resolutions, then all three possess finite projective resolutions
and we get in $K_0(R\Gamma)$
\begin{eqnarray}
  o(M_0) - o(M_1) + o(M_2) & = & 0.
  \label{o(M_0)_minus_o(M_1)_plus_o(M_2)_is_zero}
\end{eqnarray}
All this follows for instance
from~L\"uck~\cite[Chapter~11]{Lueck(1989)}.

\begin{definition}[Type (FP$_R$) and (FF$_R$) for categories]
  \label{def:Type_(FP)_and_(FF)_for_categories}
  We call a category $\Gamma$ \emph{of type (FP$_R$)}
  if the constant functor $\underline{R} \colon \Gamma^{\op} \to
  R\text{-}\MOD$ with value $R$ defines a contravariant
  $R\Gamma$-module which possesses a finite projective
  resolution.

  We call a category $\Gamma$ \emph{of type (FF$_R$)}
  if $\underline{R}$ possesses a finite free resolution.
\end{definition}

If $G$ is a group and $\widehat{G}$ is the groupoid with one object
and $G$ as automorphism group of this object, then the notions
(FP$_R$) and (FF$_R$) for $\widehat{G}$ of
Definition~\ref{def:Type_(FP)_and_(FF)_for_categories} agree with
the classical notions (FP$_R$) and (FF$_R$) for the group $G$
(see~Brown~\cite[page~199]{Brown(1982)}).

\begin{example}[Finite groups of invertible order are of type
  (FP$_R$)] \label{examp:fin_gps_FP} Let $G$ be a finite group whose order is
  invertible in the ring $R$. Then the $RG$-map $RG \to \underline{R}$,
  \begin{equation*}
    \sum_{g \in G}r_g g \mapsto \sum_{g \in G } r_g
  \end{equation*}
  admits a right inverse, namely $1 \mapsto \frac{1}{\vert G \vert} \sum_{g \in
    G} g$. The trivial $RG$-module $\underline{R}$ is then a direct summand of a
  free $RG$-module, and is therefore projective. A finite projective resolution
  of $\underline{R}$ is simply the identity $\underline{R} \to
  \underline{R}$. The group $G$ and category $\widehat{G}$ are of type (FP$_R$).
\end{example}

\begin{example}[Finite EI-categories with automorphism groups of invertible
  order are of type (FP$_R$)] \label{examp:fin_EI_FP} We may extend
  Example~\ref{examp:fin_gps_FP} to certain categories. If $\Gamma$ is a category in
  which every endomorphism is an automorphism, $\vert\aut(x) \vert$ is invertible in $R$
  for every object $x$, the category $\Gamma$ has only finitely many isomorphism
  classes of objects, and $\vert \mor_\Gamma(x,y) \vert$ is finite for all objects $x$ and
  $y$, then $\Gamma$ is of type (FP$_R$). This will follow from
  Lemma~\ref{lem:finite_homological_dimension}~\ref{lem:finite_homological_dimension:finite_Gamma}.
\end{example}

\begin{example}[Categories $\Gamma$ with a finite $\Gamma$-$CW$-model for
  $E\Gamma$ are of type (FF$_R$)] \label{exa:finite_model_implies_FFQ} If
  $\Gamma$ is a category which admits a finite $\Gamma$-$CW$-model $X$ for the
  classifying $\Gamma$-space $E\Gamma$,
  then the cellular $R$-chains of $X$
  form a finite free resolution of the constant $R\Gamma$-module
  $\underline{R}$. For example, the categories $\{1 \leftarrow 0 \rightarrow
  2\}$ and $\{a \rightrightarrows b\}$ admit finite models, as does the poset of
  non-empty subsets of $[q]=\{0,1, \dots, q\}$. Every category with a terminal
  object also admits a finite model. (Our paper~\cite{FioreLueckSauerHoColim(2009)} recalls the $\Gamma$-$CW$-complexes of~Davis--L\"uck~\cite{Davis-Lueck(1998)} in the
  context of Euler characteristics and homotopy colimits.)
\end{example}

\begin{definition}[Finiteness obstruction of a category]
  \label{def:finiteness_obstruction_of_a_category}
  The \emph{finiteness obstruction with coefficients in $R$ of a
  category $\Gamma$ of type (FP$_R$)} is
$$o(\Gamma;R) := o(\underline{R}) \; \in K_0(R\Gamma),$$
where $o(\underline{R})$ is the finiteness obstruction in
Definition~\ref{def:finiteness_obstruction_of_a_module} for the
constant $R\Gamma$-module $\underline{R}$ . We also use the notation
$[\underline{R}]$, or simply $[R]$, to denote the finiteness
obstruction $o(\Gamma;R)$.
\end{definition}

The notation $[\underline{R}]$ for the finiteness obstruction is quite natural,
for in Example~\ref{examp:fin_gps_FP} the module $\underline{R}$ is projective,
and the alternating sum of
Definition~\ref{def:finiteness_obstruction_of_a_module} is merely
$[\underline{R}]$. However, in general, the module $\underline{R}$ may not be
projective.

The homomorphism $F_*$ of \eqref{equ:induction_on_K0} depends only on the
natural isomorphism class of $F$.  Hence $F_*$ is bijective if $F$ is an
equivalence of categories.  In general $\ind_F$ is not exact and $\ind_F
\underline{R}$ is not isomorphic to $\underline{R}$.  However, this is the case
if $F$ is an equivalence of categories.  This implies

\begin{theorem}[Invariance of the finiteness obstruction under equivalence of categories]
  \label{the:invariance_under_equivalence}
  Let $\Gamma_1$ and $\Gamma_2$ be two categories such that there
  exists an equivalence $F \colon \Gamma_1 \to \Gamma_2$ of
  categories.

Then $\Gamma_1$ is of type (FP$_R$) if and only if $\Gamma_2$ is of
type (FP$_R$). In this case the isomorphism induced by $F$
$$F_* \colon K_0(R\Gamma_1) \xrightarrow{\cong} K_0(R\Gamma_2)$$
maps $o(\Gamma_1;R)$ to $o(\Gamma_2;R)$.

Moreover, $\Gamma_1$ is of type (FF$_R$) if and only if $\Gamma_2$
is of type (FF$_R$).
\end{theorem}

One easily checks

\begin{theorem}[Restriction]
  \label{the:o(Gamma;R)_and_restriction}
  Suppose that $F \colon \Gamma_1 \to \Gamma_2$ is an admissible
  functor and $\Gamma_2$ is of type (FP$_R$).

  Then $\Gamma_1$ is of type (FP$_R$) and the homomorphism
  $F^* \colon  K_0(R\Gamma_2) \to K_0(R\Gamma_1)$ sends $o(\Gamma_2;R)$
  to $o(\Gamma_1;R)$.
\end{theorem}

\begin{theorem}[Right adjoints and induction]\label{the:adjoint_functors}
  Suppose for the functors $F \colon \Gamma_1 \to \Gamma_2$ and $G
  \colon \Gamma_2 \to \Gamma_1$ that they form an adjoint pair
  $(G,F)$. Suppose that $\Gamma_1$ is of type (FP$_R$).

  Then $\Gamma_2$ is of type (FP$_R$) and
  \begin{eqnarray*}
    F_*(o(\Gamma_1;R)) & = & o(\Gamma_2;R).
  \end{eqnarray*}
\end{theorem}
\begin{proof}
  Recall that $\ind_F$ agrees with $- \otimes_{R\Gamma_1}
  R\mor_{\Gamma_2}(??,F(?))$ and $\res_G$ agrees with $-
  \otimes_{R\Gamma_1} R\mor_{\Gamma_1}(G(??),?)$. The adjunction
  $(G,F)$ (see Lemma~\ref{lem:adjunction_ind_F_and_res_F}) implies
  that $\res_G = \ind_F$. Hence $G$ is admissible. We conclude
  from Theorem~\ref{the:o(Gamma;R)_and_restriction}
\[F_*(o(\Gamma_1;R)) = G^*(o(\Gamma_1;R)) = o(\Gamma_2;R).\qedhere\]
\end{proof}

\begin{example}[Category with a terminal object]
  \label{exa:terminal_object}
  Suppose that $\Gamma$ has a terminal object $x$.  Then the constant
  $R\Gamma$-module $\underline{R}$ with value $R$ agrees with the free
  $R\Gamma$-module $R\mor(?,x)$. Hence $\Gamma$ is of type (FF$_R$) and
  the finiteness obstruction satisfies
  $$o(\Gamma;R) = [R\mor(?,x)] \in K_0(R\Gamma).$$
  Let $i\colon  \{\ast\} \to \Gamma$ be the inclusion of the trivial category
  which has precisely one morphism and sends the only object in
  $\{\ast\}$ to $x$.  Then the induced map
  $$i_* \colon K_0(R) =  K_0(R\{\ast\}) \to K_0(R\Gamma)$$
  sends $[R]$ to $o(\Gamma;R)$. This
  follows also from Theorem~\ref{the:adjoint_functors} taking $F = i$
  and $G$ to be the obvious projection.
\end{example}

\begin{example}[Wall's finiteness obstruction]
\label{exa:Walls_finiteness_obstruction} Let $G$ be a  group. Let
$\widehat{G}$ be the groupoid with one object and $G$ as morphism
set with the composition law coming from the group structure.
Because of Example~\ref{exa:modules_over_group_rings} the group $G$
is of type (FP$_R$) in the sense of homological algebra
(see~Brown~\cite[page~199]{Brown(1982)}) if and only if
$\widehat{G}$ is of type (FP$_R$) in the sense of
Definition~\ref{def:Type_(FP)_and_(FF)_for_categories}, and the
projective class group $K_0(\IZ G)$ of the group ring $\IZ G$ agrees
with $K_0(\IZ \widehat{G})$ introduced in
Definition~\ref{def:projective_class_group}.

Suppose that $G$ is of type (FP$_\IZ$) and finitely presented. Then
there is a model for $BG$ which is finitely dominated
(see~Brown~\cite[Theorem~7.1 in VIII.7 on page 205]{Brown(1982)})
and Wall (see~\cite{Wall(1965a)} and~\cite{Wall(1966)}) has defined
its finiteness obstruction
$$o(BG) \in K_0(\IZ G).$$
It agrees with the finiteness obstruction $o(\widehat{G};\IZ)$
of Definition~\ref{def:finiteness_obstruction_of_a_category}.
\end{example}

The elementary proof of the next result is left to the reader.

\begin{theorem}[Coproduct formula for the finiteness obstruction]\label{the:disjoint_union}
Let $\Gamma_1$ and $\Gamma_2$ be categories of type (FP$_R$). Then
their disjoint union $\Gamma_1 \amalg \Gamma_2$ has type (FP$_R$)
and the inclusions induce an isomorphism
$$K_0(R\Gamma_1) \oplus K_0(R\Gamma_2) \xrightarrow{\cong} K_0(R(\Gamma_1 \amalg \Gamma_2))$$
which sends $(o(\Gamma_1),o(\Gamma_2))$ to $o(\Gamma_1 \amalg \Gamma_2)$.
\end{theorem}

Let $x$ be any object of $\Gamma$. We denote by $\aut(x)$ the group
of automorphisms of $x$. We often abbreviate the associated group
ring by
\begin{eqnarray}
R[x] & := & R[\aut(x)].
\label{R[x]}
\end{eqnarray}

\begin{example}[The finiteness obstruction of a finite groupoid]
\label{exa:o(calg)_for_finite_groupoids} Let $\calg$ be a finite
groupoid, i.e., a (small) groupoid such that $\iso(\calg)$ and
$\aut_{\calg}(x)$ for any object $x \in \ob(\calg)$ are finite sets.
Then $\Gamma$ is of type (FP$_R$) if and only if for every object $x
\in \ob(\calg)$, $|\aut_{\calg}(x)| \cdot 1_R$ is a unit in $R$ (see
Lemma~\ref{lem:finite_homological_dimension}~\ref{lem:finite_homological_dimension:finite_Gamma}).

Suppose that $\calg$ is of type (FP$_R$). Then the trivial
$R[x]$-module $R$ is finitely generated projective  and defines a
class $[R]$ in $K_0(R[x])$ for every object $x \in \ob(\calg)$. We
obtain from Theorem~\ref{the:invariance_under_equivalence} and
Theorem~\ref{the:disjoint_union} a decomposition
$$K_0(R\calg) =
\bigoplus_{\overline{x} \in \iso(\Gamma)} K_0(R[x]).$$
The finiteness obstruction $o(\calg)$ has under the decomposition above
the entry $[R] \in K_0(R[x])$ for $x \in \iso(\Gamma)$.
\end{example}

Let $\Gamma_1$ and $\Gamma_2$ be two small categories. Then their
product $\Gamma_1 \times \Gamma_2$ is a small category.  Since $R$ is
commutative, the tensor product $\otimes_R$ defines a functor
$$\otimes_R \colon \MOD\text{-}R\Gamma_1 \times \MOD\text{-}R\Gamma_2 \to
\MOD\text{-}R(\Gamma_1 \times \Gamma_2).$$ Namely, put $(M \otimes_R
N)(x,y) = M(x) \otimes_R N(y)$.  Obviously
$$(M_1 \oplus M_2) \otimes_R (N_1 \oplus N_2) \;\cong \;
(M_1 \otimes_R N_1) \oplus (M_1 \otimes_R N_2) \oplus (M_2 \otimes_R
N_1) \oplus (M_2 \otimes_R N_2),$$ and for $x_1 \in \ob(\Gamma_1)$
and $x_2 \in \ob(\Gamma_2)$ we obtain isomorphisms of $R(\Gamma_1
\times \Gamma_2)$-modules
$$R\mor_{\Gamma_1}(?,x_1) \otimes_R R\mor_{\Gamma_2}(??,x_2)
\cong R\mor_{\Gamma_1 \times \Gamma_2}\bigl((?,??),(x_1,x_2)\bigr).$$
Hence we obtain a well-defined pairing
\begin{eqnarray}
  \quad && \otimes_R \colon K_0(R\Gamma_1) \otimes_{\IZ} K_0(R\Gamma_2)
  \to
  K_0(R(\Gamma_1 \times \Gamma_2)),
  \quad [P_1] \otimes [P_2] \to [P_1 \otimes_R P_2].
  \label{otimes_R_on_K_0}
\end{eqnarray}

\begin{theorem}[Product formula for the finiteness obstruction]
\label{the:product_formula_for_o} Let $\Gamma_1$ and $\Gamma_2$ be
categories of type (FP$_R$).

Then $\Gamma_1 \times \Gamma_2$ is of type (FP$_R$) and  we get
\begin{eqnarray*}
o(\Gamma_1 \times \Gamma_2;R)
& = &
o(\Gamma_1;R) \otimes_R o(\Gamma_2;R)
\end{eqnarray*}
under the pairing~\eqref{otimes_R_on_K_0}.
\end{theorem}
\begin{proof}
  Let $P^{i}_*$ be a finite projective resolution of $\underline{R}$
  over $\MOD\text{-}R\Gamma_i$ for $i = 1,2$. The evaluation of a
  projective $R\Gamma_i$-module at an object is projective and hence
  flat as $R$-module since this is obviously true for
  $R\mor_{\Gamma_i}(?,x)$ and every projective $R\Gamma_i$-module is a
  direct sum in a free one. Hence the $R(\Gamma_1 \times
  \Gamma_2)$-chain complex $P^1_* \otimes_R P^2_*$ is a projective  $R\Gamma_1
  \times R\Gamma_2$-resolution of $\underline{R}$.  Now an easy
  calculation (see~L\"uck~\cite[11.18 on page~227]{Lueck(1989)} shows
\[o(\Gamma_1 \times \Gamma_2;R) = o(P^1_* \otimes_R P^2_*)
= o(P^1_*) \otimes_R o(P^2_*) = o(\Gamma_1;R) \otimes_R
o(\Gamma_2;R).\qedhere\]
\end{proof}

\begin{example}\label{exa:finiteness_obstruction_of_projection_category}
  Let $\Gamma$ be the category which has precisely one object $x$ and
  two morphisms $\id_x \colon x \to x$ and $p \colon x \to x$ such
  that $p \circ p = p$.  Given an $R$-module $M$, let $I_i(M)$ for $i
  = 0,1$ be the contravariant $R\Gamma$-module which sends $p \colon x
  \to x$ to $i \cdot \id_M \colon M \to M$. Given any $R\Gamma$-module
  $N$, we obtain an isomorphism of $R\Gamma$-modules
  $$f \colon I_0\bigl(\ker(N(p))\bigr) \oplus I_1\bigl(\im(N(p))\bigr)
  \xrightarrow{\cong} N$$ from the inclusions of $\ker(N(p))$ and
  $\im(N(p))$ to $N(x)$.  This isomorphism is natural in $N$ and
  respects direct sums. If $N = R\mor(?,x)$, we have
  $\ker(N(p)) \cong \im(N(p)) \cong R$. Hence $I_i(R)$ is a finitely generated
  projective $R\Gamma$-module for $i = 0,1$.
  This implies that $N$ is a finitely generated
  projective $R\Gamma$-module if and only if $\ker(N(p))$ and
  $\im(N(p))$ are finitely generated projective $R$-modules. Hence we
  obtain an isomorphism
  $$K_0(R\Gamma) \xrightarrow{\cong} K_0(R) \oplus K_0(R), \quad
  [P] \mapsto \bigl([\ker(P(p))], [\im(P(p))]\bigr).$$ Its inverse sends
  $([P_0],[P_1])$ to $[I_0(P_0) \oplus I_1(P_1)]$. The constant
  $R\Gamma$-module $\underline{R}$ agrees with $I_1(R)$.  Hence the
  category $\Gamma$ is of type (FP$_R$) and the finiteness obstruction
  $o(\Gamma;R)$ is sent under the isomorphism above to the element
  $(0,[R])$.
\end{example}


\typeout{------------------  Section 3: Splitting the projective class -----------------}
\section{Splitting the projective class group}
\label{sec:Splitting_the_projective_class_group}

In this section we will investigate the projective class group
$K_0(R\Gamma)$. In the case that every endomorphism in $\Gamma$ is
an isomorphism, we construct the natural \emph{splitting}
isomorphism
$$S \colon K_0(R\Gamma) \to \Split K_0(R\Gamma) :=
\bigoplus_{\overline{x} \in \iso(\Gamma)} K_0(R\aut_\Gamma(x))$$ and
its natural inverse $E$, called \emph{extension}. This is L\"uck's
Splitting of $K_0(R\Gamma)$ in~\cite[Theorem~10.34 on
page~196]{Lueck(1989)}. If $\Gamma$ is merely directly finite rather
than EI, we still have $S \circ E=\id_{\Split K_0(R\Gamma)}$ and the
naturality of $S$, though $S$ is no longer bijective. The splitting
functor $S_x$ of \eqref{equ:splitting_functor} and the extension
functor $E_x$ of \eqref{equ:extension_functor} respect direct sums
and send epimorphisms to epimorphisms. The extension functor $E_x$
sends free $R\aut_\Gamma(x)$-modules to free $R\Gamma$-modules. If
$\Gamma$ is directly finite, the restriction functor $S_x$ sends
free $R\Gamma$-modules to free $R\aut_\Gamma(x)$-modules and
respects finitely generated and projective. The relationship between
EI-categories, directly finite categories, and Cauchy complete
categories is clarified in Lemma
\ref{lem:directly_finite_and_idempotents_and_EI}.

Recall that a ring is called \emph{directly finite} if for two
elements $r,s \in R$ we have the implication $rs = 1 \implies sr =
1$. Therefore we define

\begin{definition}[Directly finite category]
\label{def:directly_finite_category}
A category is called \emph{directly finite} if for any two objects $x$ and $y$ and morphisms
$u \colon x \to y$ and $v \colon y \to x$ the implication
$vu = \id_x \implies uv = \id_y$ holds.
\end{definition}

\begin{lemma}[Invariance of direct finiteness under equivalence of categories] \label{lem:direct_finiteness_under_equivalence}
Suppose $\Gamma_1$ and $\Gamma_2$ are equivalent categories. Then
$\Gamma_1$ is directly finite if and only if $\Gamma_2$ is directly
finite.
\end{lemma}
\begin{proof}
Suppose $F \colon \Gamma_1 \to \Gamma_2$ is fully faithful and
essentially surjective, that $\Gamma_1$ is directly finite, and
$vu=\id_x$ in $\Gamma_2$. Then we can extend to a commutative
diagram
$$\xymatrix{x \ar[r]_u \ar@/^1pc/[rr]^{\id_x} \ar[d]^\cong_m & y \ar[r]_v \ar[d]^\cong_n & x
\ar[d]_\cong^m \ar[r]_u & y \ar[d]^n_\cong
\\ F(a) \ar[r]_{F(f)} & F(b) \ar[r]_{F(g)} & F(a) \ar[r]_{F(f)} & F(b).}$$
Hence $F(g \circ f)=\id_{F(a)}$, and $g\circ f=\id_a$. The direct
finiteness of $\Gamma_1$ then implies $f \circ g=\id_b$. Together
with the commutativity of the two right squares above, this implies
$u \circ v=\id_y$, so that $\Gamma_2$ is also directly finite.
\end{proof}

Let $M$ be any $R\Gamma$-module and let $x$ be any object. We denote
by $\aut_{\Gamma}(x)$ (or $\aut(x)$ when $\Gamma$ is clear) the
group of automorphisms of $x$. As in~\ref{R[x]}, we abbreviate the
associated group ring by $R[x]:= R[\aut(x)]$. Define an $R$-module
$S_xM$ by the cokernel of the map of $R$-modules
\[S_xM := \coker\left(\bigoplus_{\substack{u \colon x \to y\\u \;
\text{is not an isomorphism}}} M(u)
  \colon \bigoplus_{\substack{u \colon x \to y\\u \; \text{is not an isomorphism
        }}} M(y) \; \to \;M(x)\right).\]
In other words, $S_xM$ is the quotient of the $R$-module $M(x)$ by
the $R$-submodule generated by all images of $R$-module
homomorphisms $M(u)\colon M(y) \to M(x)$ induced by all
non-invertible morphisms $u\colon x \to y$ in $\Gamma$. One easily
checks that the right $R[x]$-module structure on $M(x)$ coming from
functoriality induces a right $R[x]$-module structure on $S_xM$.
Thus we obtain a functor called \emph{splitting functor at $x \in
  \ob(\Gamma)$}
\begin{eqnarray} \label{equ:splitting_functor}
  & S_x \colon \MOD\text{-}R\Gamma \to \MOD\text{-}R[x], &
  \label{S_x}
\end{eqnarray}
where $\MOD\text{-}R[x]$ denotes the category of right $R[x]$-modules.
Define a functor, called \emph{extension functor at $x \in \ob(\Gamma)$,}
\begin{eqnarray} \label{equ:extension_functor}
  & E_x \colon  \MOD\text{-}R[x] \to  \MOD\text{-}R\Gamma&
  \label{E_x}
\end{eqnarray}
by sending an $R[x]$-module $N$ to the  $R\Gamma$-module
$N \otimes_{R[x]} R\mor(?,x)$.

\begin{lemma}[Extension/splitting, direct sums, and free/projective modules]\label{lem:basic_properties_of-E_x_and_S_x}\hfill
  \begin{enumerate}
  \item \label{lem:basic_properties_of-E_x_and_S_x:E_x} The functor
    $E_x$ respects direct sums. It sends epimorphisms to epimorphisms.
    It sends a free $R[x]$-module with the set $C$ as basis to the
    free $R\Gamma$-module with the $\ob(\Gamma)$-set $D$ as basis, where
    $D_x = C$ and $D_y = \emptyset$ for $y \not= x$. It respects
    finitely generated and projective;

  \item \label{lem:basic_properties_of-E_x_and_S_x:S_x_circE_x_of_free}
    We have $S_y \circ E_x = 0$, if $x$ and $y$ are not isomorphic.
    For every projective right $R[x]$-module $P$ we have a surjective map of
    $R[x]$-modules, natural in $P$ and compatible with direct sums
    $$\sigma_P \colon P \to S_x \circ E_x(P);$$

\item \label{lem:basic_properties_of-E_x_and_S_x:S_x} The functor
    $S_x$ respects direct sums. It sends epimorphisms to epimorphisms
    and sends finitely generated $R\Gamma$-modules to finitely generated $R[x]$-modules;

 \item \label{lem:basic_properties_of-E_x_and_S_x:directly_finite}
    Suppose that $\Gamma$ is directly finite. Then $S_x$
    sends a free $R\Gamma$-module with the $\ob(\Gamma)$-set $C$ as
    basis to the free $R[x]$-module with
    $\coprod_{y \in \ob(\Gamma), \overline{y} = \overline{x}} C_y$ as basis and respects
    finitely generated and projective. Further, $\sigma_P$ appearing in
    assertion~\ref{lem:basic_properties_of-E_x_and_S_x:S_x_circE_x_of_free}
    is bijective for every projective right $R[x]$-module $P$.

\end{enumerate}
\end{lemma}
\begin{proof}\ref{lem:basic_properties_of-E_x_and_S_x:E_x} Obviously $E_x$ is
  compatible with direct sums. It sends epimorphisms to epimorphisms
  since tensor products   are right exact.  We have
  $$E_x(R[x]) = R[x] \otimes_{R[x]} R\mor(?,x) =  R\mor(?,x).
  $$%
\ref{lem:basic_properties_of-E_x_and_S_x:S_x_circE_x_of_free}
  Suppose that $x$ and $y$ are not isomorphic. Let $P$ be an $R[x]$-module.
  Consider an element $p \otimes u \in E_xP(y) = P \otimes_{R[x]} R\mor(y,x)$.
  Since $x$ and $y$ are not isomorphic, $u$ is not an  isomorphism.
  The element $p \otimes u$ lies in the  image of the map induced
  by composition from the right with $u$
  $$P \otimes_{R[x]} R\mor(x,x) \to P \otimes_{R[x]} R\mor(y,x),$$
  a preimage is given by $p \otimes \id_x$. Hence $S_y \circ E_x(P) = 0$.

  Define an $R[x]$-map $P \to P \otimes_{R[x]} \mor(x,x)$ by sending
  $p \in P$ to $p \otimes_{R[x]}\id_x$. Its composition with the
  canonical projection $P \otimes_{R[x]} \mor(x,x) \to S_x \circ E_x(P)$
  yields an $R[x]$-map
  $$\sigma_P \colon P \to S_x \circ E_x(P).$$
  Obviously it is surjective, natural in $P$ and compatible with direct sums.
  \\[1mm]\ref{lem:basic_properties_of-E_x_and_S_x:S_x}
  This is obvious except that $S_x$ respects finitely generated.
  We know already that $S_yR\mor(?,x) = 0$ if $x$ and $y$ are not isomorphic
  and that there is an epimorphism $R[x] \to S_xR\mor(?,x)$.
  Hence $S_xR\mor(?,y)$ is a finitely generated $R\aut(x)$-module for all
  $y \in \ob(\Gamma)$ and the claim follows.
  \\[1mm]\ref{lem:basic_properties_of-E_x_and_S_x:directly_finite}
  Consider an endomorphism $u \colon x \to x$. It lies
  in the image of the map $\mor(x,x) \to
  \mor(x,x), \; v \mapsto v \circ u$, a preimage is $\id_x$.
  If $u$ is an isomorphism, then there exists no morphism $w \colon x \to y$
  such that $w$ is not an isomorphism and $u$ lies in the image
  of $\mor(y,x) \to \mor(x,x), \; v \mapsto v \circ w$, since $\Gamma$
  is directly finite. This implies that
  $$\sigma_{R[x]} \colon R[x] \xrightarrow{\cong} S_x \circ E_x(R[x]) = S_xR\mor(?,x)$$
  is an isomorphism.
  Now assertion~\ref{lem:basic_properties_of-E_x_and_S_x:directly_finite}
  follows from compatibility with direct sums and the facts that an
  $R\Gamma$-module is projective if and only if it is a direct summand
  in a free $R\Gamma$-module and that $S_x$ respects epimorphisms.
\end{proof}

We denote by $\iso(\Gamma)$ the set of isomorphism classes of
objects of $\Gamma$.  Choose for any class $\overline{x} \in
\iso(\Gamma)$ a representative $x \in \overline{x}$. Define
\begin{eqnarray}
  \Split K_0(R\Gamma) & := & \bigoplus_{\overline{x} \in \iso(\Gamma)} K_0(R[x]).
  \label{Split_K_0}
\end{eqnarray}
Provided that $\Gamma$ is directly finite, we obtain from
Lemma~\ref{lem:basic_properties_of-E_x_and_S_x} homomorphisms

\begin{eqnarray}
S \colon K_0(R\Gamma)  \to  \Split K_0(R\Gamma),
& & [P] \mapsto  \{[S_xP] \mid \overline{x} \in \iso(\Gamma)\};
\label{S_colon_SplitK_0_toK_0}
\\
E \colon \Split K_0(R\Gamma)  \to  K_0(R\Gamma),
& & \{[Q_x] \mid \overline{x} \in \iso(\Gamma)\}  \mapsto
\sum_{\overline{x} \in \iso(\Gamma)} [E_xQ_x],
\label{E_colon_SplitK_0_to_K_0}
\end{eqnarray}
and get
\begin{lemma} \label{S_circ_E_is_id}
Suppose that $\Gamma$ is directly finite.
The composite $S \circ E$ is the identity. In particular $S$ is split surjective.
\end{lemma}

The group $\Split K_0(R\Gamma)$ is easier to understand than
$K_0(R\Gamma)$ since its input are projective class groups over group rings.  We
will later explain that for an EI-category the maps $E$ and $S$ are
bijective (see Theorem~\ref{the_splitting_of_K-theory_for_EI_categories}).

\begin{definition} \label{def:EI-category}
A category is an \emph{EI-category} if every endomorphism is an
isomorphism.
\end{definition}

The EI-property is invariant under equivalence of categories.
\begin{lemma} \label{lem:EI_under_equivalence}
Suppose $\Gamma_1$ and $\Gamma_2$ are equivalent categories. Then
$\Gamma_1$ is an EI-category if and only if $\Gamma_2$ is an
EI-category.
\end{lemma}
\begin{proof}
Let $\Gamma_1$ be an EI-category, $F\colon \Gamma_1 \to \Gamma_2$ an
equivalence of categories, and $b\in \ob(\Gamma_2)$. Then $b\cong
F(a)$ for some $a\in\ob(\Gamma_1)$. We have isomorphisms of monoids
$$\mor_{\Gamma_1}(a,a)\cong \mor_{\Gamma_2}(F(a),F(a))\cong\mor_{\Gamma_2}(b,b).$$
The first monoid is a group, and hence so is the last.
\end{proof}

\begin{definition} [Cauchy complete category]
  \label{def:Cauchy_complete}
  A category $\Gamma$ is  \emph{Cauchy complete} if every
  idempotent splits, i.e., for every idempotent $p \colon x \to x$
  there exists morphisms $i \colon y \to x$ and $r \colon x \to y$
  with $r \circ i = \id_y$ and $i \circ r = p$.
\end{definition}

\begin{lemma}\label{lem:directly_finite_and_idempotents_and_EI}
Consider a category $\Gamma$. Consider the statements

\begin{enumerate}
\item \label{lem:directly_finite_and_idempotents_and_EI:EI}
$\Gamma$ is an EI-category;

\item \label{lem:directly_finite_and_idempotents_and_EI:idempotents}
Every idempotent $p \colon x \to x$ in $\Gamma$ satisfies $p = \id_x$;

\item \label{lem:directly_finite_and_idempotents_and_EI:directly_finite_and_Cauchy}
$\Gamma$ is directly finite and Cauchy complete.
\end{enumerate}
Then
$\ref{lem:directly_finite_and_idempotents_and_EI:EI}
\implies%
$\ref{lem:directly_finite_and_idempotents_and_EI:idempotents}
and
$\ref{lem:directly_finite_and_idempotents_and_EI:idempotents}
\Longleftrightarrow%
\ref{lem:directly_finite_and_idempotents_and_EI:directly_finite_and_Cauchy}$.

If $\mor(x,x)$ is finite for all $x \in \ob(\Gamma)$, then
$\ref{lem:directly_finite_and_idempotents_and_EI:EI}
\Longleftrightarrow%
\ref{lem:directly_finite_and_idempotents_and_EI:idempotents}
\Longleftrightarrow%
\ref{lem:directly_finite_and_idempotents_and_EI:directly_finite_and_Cauchy}$.

\end{lemma}
\begin{proof}%
\ref{lem:directly_finite_and_idempotents_and_EI:EI}
$\implies$%
\ref{lem:directly_finite_and_idempotents_and_EI:idempotents}
If $p \colon x \to x$ is an idempotent, it is an endomorphism and hence an isomorphism.
Hence $\id_x = p^{-1} \circ p = p^{-1} \circ p \circ p = \id_x \circ p = p$.
\\[1mm]%
\ref{lem:directly_finite_and_idempotents_and_EI:idempotents}
$\implies $%
\ref{lem:directly_finite_and_idempotents_and_EI:directly_finite_and_Cauchy}
Consider morphisms $u \colon x \to y$ and $v \colon y \to x$ with $vu = \id_x$.
Then $(uv)^2 = uvuv = u \circ \id_x\circ v = uv$ is an idempotent
and hence by assumption $uv=\id_y$. Obviously $\Gamma$ is Cauchy complete.
\\[1mm]%
\ref{lem:directly_finite_and_idempotents_and_EI:directly_finite_and_Cauchy}
$\implies $%
\ref{lem:directly_finite_and_idempotents_and_EI:idempotents}
Consider an idempotent $p \colon x \to x$. Since $\Gamma$ is Cauchy complete,
we can choose morphisms $i \colon y \to x$ and $r \colon x \to y$
  with $r \circ i = \id_y$ and $i \circ r = p$. Since $\Gamma$ is directly finite,
$p = i \circ r = \id_x$.
\\[1mm]%
It remains to show~%
\ref{lem:directly_finite_and_idempotents_and_EI:idempotents}
$\implies$%
\ref{lem:directly_finite_and_idempotents_and_EI:EI}
provided that $\mor(x,x)$ is finite for all objects $x \in \ob(\Gamma)$.
Consider an endomorphism $f \colon x \to x$. Since $\mor(x,x)$ is finite, there
exists integers $m,n \ge 1$ with $f^m = f^{m+n}$. This implies
$f^m = f^{m+kn}$ for all natural numbers $k \ge 1$. Hence we get  $f^m = f^{m+n}$ for
some $n \ge 1$ with $n - m \ge 0$. Then
$$f^{n} \circ f^{n}  =  f^{2n} = f^{m+n} \circ f^{n-m} = f^m \circ
f^{n-m}=  f^{n}.$$
Hence $f^{n}$ is an idempotent.
Since then $f^{n} = \id$ for some $n \ge 1$, the endomorphism
$f$ must be an isomorphism.
\end{proof}

The next result is~from L\"uck~\cite[Theorem~10.34 on
page~196]{Lueck(1989)}.

\begin{theorem}[Splitting of $K_0(R\Gamma)$ for EI-categories]
\label{the_splitting_of_K-theory_for_EI_categories} If $\Gamma$ is
an EI-category, the group homomorphisms
\begin{eqnarray*}
S \colon K_0(R\Gamma)  \to  \Split K_0(R\Gamma),
& & [P] \mapsto  \{[S_xP] \mid \overline{x} \in \iso(\Gamma)\};
\\
E \colon \Split K_0(R\Gamma)  \to  K_0(R\Gamma),
& & \{[Q_x] \mid \overline{x} \in \iso(\Gamma)\}  \mapsto
\sum_{\overline{x} \in \iso(\Gamma)} [E_xQ_x],
\end{eqnarray*}
of~\eqref{S_colon_SplitK_0_toK_0}
and~\eqref{E_colon_SplitK_0_to_K_0} are isomorphisms and inverse to
one another. They are covariantly natural with respect to functors
$F:\Gamma_1 \to \Gamma_2$ between EI-categories, that is
$$(\Split F_*) \circ S^{R\Gamma_1}=S^{R\Gamma_2} \circ F_*$$
$$\text{and}$$
$$F_* \circ E^{R\Gamma_1}=E^{R\Gamma_2} \circ (\Split F_*).$$ The
functor $\Split F_*$ is defined in more detail in Lemma~\ref{lem:S_and_splitting}.
Moreover, $S$ and $E$ are also contravariantly natural with respect to admissible functors
$F:\Gamma_1 \to \Gamma_2$ between EI-categories, that is $$
S^{R\Gamma_1} \circ F^*=\Split F^* \circ S^{R\Gamma_2}$$
$$\text{and}$$
$$E^{R\Gamma_1} \circ (\Split F^*) =F^*  \circ E^{R\Gamma_2}.$$
\end{theorem}

Example~\ref{exa:finiteness_obstruction_of_projection_category}
shows that the EI hypothesis on $\Gamma$ in
Theorem~\ref{the_splitting_of_K-theory_for_EI_categories} is
necessary for $S$ and $E$ to be bijections. Though the splitting
homomorphism $S$ is no longer an isomorphism in general, it is covariantly
natural in the more general setting of directly finite categories.
\begin{lemma} \label{lem:S_and_splitting}
Let $\Gamma_1$ and $\Gamma_2$ be directly finite categories and $F
\colon \Gamma_1 \to \Gamma_2$ be a functor.

 Then the following
diagram commutes
$$\xymatrix@C=3pc{K_0(R\Gamma_1) \ar[r]^{F_\ast} \ar[d]_{S^{R\Gamma_1}} & K_0(R\Gamma_2)
\ar[d]^{S^{R\Gamma_2}}
\\ \Split K_0(R\Gamma_1) \ar[r]_{\Split
  F_*} & \Split K_0(R\Gamma_2) }
$$
where the vertical maps have been defined
in~\eqref{S_colon_SplitK_0_toK_0}, the upper horizontal map is
induced by induction with $F$, and the lower horizontal arrow is
given by the matrix of homomorphisms
$$
\bigl((F_{\overline{x},\overline{y}})_*
\bigr)_{\overline{x} \in \iso(\Gamma_1),\overline{y} \in \iso(\Gamma_2)}
$$
where $(F_{\overline{x},\overline{y}})_*$ is
trivial if $\overline{F(x)} \not= \overline{y}$ and given by induction
with the group homomorphism  $F_x \colon \aut_{\Gamma_1}(x) \to
\aut_{\Gamma_2}(F(x)), \; f \mapsto F(f)$ for $\overline{y} = \overline{F(x)}$.

In particular, the commutativity of the diagram guarantees
$$S_{F(x)}^{R\Gamma_2} \circ F_* =F_x \circ S^{R\Gamma_1}_{x}.$$
\end{lemma}
\begin{proof} For two objects $x$ and $y$ in $\Gamma_1$, let
  $\mor^{\cong}(x,y)$ be the set of isomorphisms from $x$ to $y$. The
  covariant $R\Gamma_1$-module $R\mor^{\cong}(x,?)$ assigns to an
  object $x$ the trivial $R$-module $\{0\}$ if $\overline{x} \not=
  \overline{y}$ and $R\mor^{\cong}(x,y)$ if $\overline{x} =
  \overline{y}$. The evaluation of  $R\mor^{\cong}(x,?)$
  at a morphism  $f \colon y_1 \to y_2$ is given by
  $$R\mor^{\cong}(x,y_1) \to R\mor^{\cong}(x,y_2), \quad g \mapsto f
  \circ g$$
  if $f$ is an isomorphism and $\overline{x} = \overline{y}$,
  and by the trivial $R$-homomorphism otherwise. This definition makes
  sense since $\Gamma_1$ is directly finite. Obviously
  $R\mor^{\cong}(x,?)$ is an $R\Gamma_1$-$R[x]$-bimodule. Hence we
  obtain a functor
  $$\MOD\text{-}R\Gamma_1 \to \MOD\text{-}R[x],
  \quad P \mapsto P \otimes_{R\Gamma_1} R\mor^{\cong}(x,?).$$
  It is naturally isomorphic to the splitting functor $S_x$
  defined in~\eqref{S_x}. Namely, a natural isomorphism is given by the
  $R[x]$-isomorphisms which are inverse to one another
  $$S_xP \to P \otimes_{R\Gamma_1} R\mor^{\cong}(x,?),
  \quad \overline{p} \mapsto p \otimes \id_x.$$
  and
  $$P \otimes_{R\Gamma_1} R\mor^{\cong}(x,?) \to S_xP,
  \quad p \otimes f \mapsto \overline{P(f)(p)}.$$

  Consider a projective $R\Gamma_1$-module $P$.
  Then we obtain for $y \in \iso(\Gamma_2)$ a natural
  isomorphism of $R[y]$-modules
  \begin{eqnarray*}
  S_y \circ \ind_F P
  & \cong &
  P \otimes_{R\Gamma_1} R\mor_{\Gamma_2}(??,F(?)) \otimes_{R\Gamma_2} R\mor^{\cong}_{\Gamma_2}(y,??)
  \\
  & \cong &
  P \otimes_{R\Gamma_1} R\mor^{\cong}_{\Gamma_2}(y,F(?))
  \\
  & \cong &
  P \otimes_{R\Gamma_1} \bigoplus_{\overline{x} \in \iso(\Gamma_1), \overline{F(x)} = \overline{y}}
  R\mor^{\cong}_{\Gamma_1}(x,?)   \otimes_{R[x]} R\mor^{\cong}_{\Gamma_2}(y,F(x))
  \\
  & \cong &
  \bigoplus_{\overline{x} \in \iso(\Gamma_1), \overline{F(x)} = \overline{y}}
  P \otimes_{R\Gamma_1} R\mor^{\cong}_{\Gamma_1}(x,?)
  \otimes_{R[x]} R\mor^{\cong}_{\Gamma_2}(y,F(x))
  \\
  & \cong &
  \bigoplus_{\overline{x} \in \iso(\Gamma_1), \overline{F(x)} = \overline{y}}
  \ind_{F_x} \circ S_xP.
  \end{eqnarray*}
  This finishes the proof of Lemma~\ref{lem:S_and_splitting}.
\end{proof}


\typeout{------  Section 4: The (functorial) Euler characteristic of a category --------}
\section{The  (functorial) Euler characteristic of a category}
\label{sec:The_(functorial)_Euler_characteristic_of_a_category}

Perhaps the most naive notion of Euler characteristic for a category
$\Gamma$ is the \emph{topological Euler characteristic}, namely the
classical Euler characteristic of the classifying space $B\Gamma$.
However, even in the simplest cases, $\chi(B\Gamma;R)$ may not
exist, for example $\Gamma=\widehat{\IZ_2}$ and $R=\IZ_2$. We
propose better invariants using the homological algebra of
$R\Gamma$-modules and von Neumann dimension.

Depending on which notion of rank we choose for $R\Gamma$-modules,
$\rk_{R\Gamma}$ vs. $\rk_\Gamma^{(2)}$, there are two possible ways
to define (functorial) Euler characteristics. In this section, we
start with the topological Euler characteristic $\chi(B\Gamma;R)$,
and then treat the homological Euler characteristic $\chi(\Gamma;R)$
and its functorial counterpart $\chi_f(\Gamma;R)$, both of which
arise from $\rk_{R\Gamma}$.  In Section
\ref{sec:The_(functorial)_L2_Euler_characteristic_anf_L2-Betti_numbers_of_a_category}
we take $R=\IC$ and $\rk_\Gamma^{(2)}$ (defined in terms of the von
Neumann dimension) to treat the $L^2$-Euler characteristic
$\chi^{(2)}(\Gamma)$ and its functorial counterpart
$\chi^{(2)}_f(\Gamma)$.

To obtain the Euler characteristic, we use the splitting functor
$S_x$ as follows. The \emph{$R\Gamma$-rank} of a finitely generated
$R\Gamma$-module $M$ is an element of $U(\Gamma)$, the free abelian
group on the isomorphism classes of objects of $\Gamma$. At
$\overline{x} \in \iso(\Gamma)$, $\rk_{R\Gamma}M$ is $\rk_R(S_xM
\otimes_{R\aut(x)} R)$.  This induces a homomorphism $\rk_{R\Gamma}$
from $K_0(R\Gamma)$ to $U(\Gamma)$. If $\Gamma$ is of type (FP$_R$),
we define the \emph{functorial Euler characteristic}
$\chi_f(\Gamma;R)$ to be the image of the finiteness obstruction
$o(\Gamma;R)$ under $\rk_{R\Gamma}$. The functorial Euler
characteristic is compatible with equivalences between directly
finite categories of type (FP$_R$). The \emph{Euler characteristic}
$\chi(\Gamma;R)$ is the sum of the components of the functorial
Euler characteristic $\chi_f(\Gamma;R)$. If $\Gamma$ is a directly
finite category of type (FP$_R$) and $R$ is Noetherian, then the
Euler characteristic $\chi(\Gamma;R)$ is equal to the topological
Euler characteristic $\chi(B\Gamma;R)$. If $R$ is Noetherian and
$\Gamma$ is of type (FP$_R$), but not necessarily directly finite,
then the image of the finiteness obstruction under $\rk_R\pr_*$ in
\eqref{equ:relating_o_and_chi} is the topological Euler
characteristic $\chi(B\Gamma;R)$. If $R$ is Noetherian and  $\Gamma$
is directly finite and of type (FF$_\IZ$), then
$\chi(B\Gamma;R)=\chi(\Gamma;R)=\chi^{(2)}(\Gamma)$, see
Theorem~\ref{the:coincidence}.

Each notion of Euler characteristic ($\chi$ vs. $\chi^{(2)}$) has
its advantages. Both are invariant under equivalence of categories
(assuming directly finite) and are compatible with finite products,
finite coproducts, and homotopy colimits
(see~Fiore--L\"uck--Sauer~\cite{FioreLueckSauerHoColim(2009)} for
the compatibility with homotopy colimits). The $L^2$-Euler
characteristic is compatible with isofibrations and coverings
between connected finite groupoids (see Subsection
\ref{subsec:Compatibility_with_Coverings}). If the groupoids are
additionally of type (FF$_\IC$), then the Euler characteristic and
topological Euler characteristic agree with the $L^2$-Euler
characteristic, and are therefore compatible with the isofibrations
and coverings at hand. For a finite discrete category (a set), both
$\chi$ and $\chi^{(2)}$ return the cardinality. For a finite group
$G$, we have $\chi(\widehat{G};\IQ)=1$, while the $L^2$-Euler
characteristic is $\chi^{(2)}(\widehat{G})=\frac{1}{\vert G \vert}$.
The groupoid cardinality of Baez--Dolan~\cite{Baez-Dolan(2001)} and
the Euler characteristic of Leinster~\cite{Leinster(2008)} will
occur as an $L^2$-Euler characteristic, see Section
\ref{sec:Comparison_with_Leinsters_invariant} for the comparison.
The main advantages of our $K$-theoretic approach are: 1) it works
for infinite categories, and 2) it encompasses important examples,
such as the $L^2$-Euler characteristic of a group and the
equivariant Euler characteristic of the classifying space $\eub{G}$
for proper $G$-actions.

To begin the details of the topological Euler characteristic and the
Euler characteristic, suppose that we have specified the notion of a
rank\label{page:notion of rank}
\begin{eqnarray}
&\rk_R(N) \in \IZ&
\label{rk_R}
\end{eqnarray}
for every finitely generated  $R$-module such that $\rk_R(N_1) =
\rk_R(N_0) + \rk_R(N_2)$ for any sequence $0 \to N_0 \to N_1 \to N_2
\to 0$ of finitely generated $R$-modules and $\rk_R(R) = 1$. If $R$
is a commutative principal ideal domain, we will use $\rk_R(N) :=
\dim_F(F \otimes_R N)$ for $F$ the quotient field of $R$.

\begin{definition}[The topological Euler characteristic of a category $\Gamma$]
\label{def:topological_Euler_characteristic_of_a_category} Let
$\Gamma$ be a category. Let $B\Gamma$ be its classifying space,
i.e., the geometric realization of its nerve. Suppose that
$H_n(B\Gamma;R)$ is a finitely generated $R$-module for every $n \ge
0$ and that there exists a natural number $d$ with $H_n(B\Gamma;R) =
0$ for $n > d$. The \emph{topological Euler characteristic of
$\Gamma$} is
$$\chi(B\Gamma;R) = \sum_{n \ge 0} (-1)^n \cdot \rk_{R}(H_n(B\Gamma;R)) \in \IZ.$$
\end{definition}

\begin{example}[The topological Euler characteristic of a finite groupoid]
\label{exa:chi(calg)_for_finite_groupoids} Let $\calg$ be a finite
groupoid, i.e., a (small) groupoid such that $\iso(\calg)$ and
$\aut(x)$ for any object $x \in \ob(\calg)$ are finite. Consider $R
= \IQ$. Then the assumptions in
Definition~\ref{def:topological_Euler_characteristic_of_a_category}
are satisfied and
$$\chi(B\calg) = |\iso(\calg)|.$$
\end{example}

\begin{notation}[The abelian group $U(\Gamma)$ and the augmentation homomorphism
$\epsilon$] \label{U(Gamma)_and_augmentation} Let $\Gamma$ be a
category. We denote by $U(\Gamma)$ the free abelian group on the set
of isomorphism classes of objects in $\Gamma$, that is $$U(\Gamma)
:= \IZ\iso(\Gamma).$$ For a functor $F\colon \Gamma_1 \to \Gamma_2$,
the group homomorphism $U(F)\colon \Gamma_1 \to \Gamma_2$ maps the
basis element $\overline{x}$ to the basis element
$\overline{Fx}$. The
\emph{augmentation homomorphism} $\epsilon\colon U(\Gamma) \to \IZ$
sends every basis element of $\iso(\Gamma)$ to $1 \in \IZ$. The
augmentation homomorphism is a natural transformation from the
covariant functor $U\colon \CATS \to \ABELIANGROUPS$ to the constant
functor $\IZ$, that is, for any functor $F\colon \Gamma_1 \to
\Gamma_2$ the diagram
\begin{equation} \label{equ:augmentation_naturality_diagram}
\xymatrix@C=3pc{U(\Gamma_1) \ar[r]^{U(F)} \ar[d]_{\epsilon} & U(\Gamma_2) \ar[d]^{\epsilon}
    &
    \\
    \IZ \ar[r]_-{\id_{\IZ}} & \IZ}
\end{equation}
commutes.
\end{notation}

\begin{definition}[Rank of a finitely generated $R\Gamma$-module]
\label{def:rank_of_fin_gen_prof_RGamma-module} Let $M$ be a finitely
generated $R\Gamma$-module $M$, define its \emph{rank}
$$\rk_{R\Gamma}(M) :=
\bigl\{\rk_R(S_xM \otimes_{R[x]} R) \mid \overline{x} \in
\iso(\Gamma)\bigr\} \quad \in U(\Gamma).$$
\end{definition}

The rank $\rk_{R\Gamma}$ defines a homomorphism
\begin{eqnarray}
  \rk_{R\Gamma} \colon K_0(R\Gamma) \to U(\Gamma), \quad [P] \to \rk_{R\Gamma}(P).&&
  \label{rk_K_o_to_U}
\end{eqnarray}
It obviously factorizes over $S \colon K_0(R\Gamma) \to \Split
K_0(R\Gamma)$.  Define
\begin{eqnarray}
  && \iota \colon  U(\Gamma) \to K_0(R\Gamma),
  \quad (n_{\overline{x}})_{\overline{x} \in \iso(\Gamma)} \mapsto
  \sum_{\overline{x} \in \iso(\Gamma)} n_{\overline{x}} \cdot [R\mor(?,x)].
  \label{iota}
\end{eqnarray}
This is the same as the composite
$$U(\Gamma) = \bigoplus_{\overline{x} \in \iso(\Gamma)} \IZ
  \xrightarrow{\bigoplus_{\overline{x} \in \iso(\Gamma)} i_x}
  \bigoplus_{\overline{x} \in \iso(\Gamma)} K_0(R[x]) = \Split K_0(R\Gamma)
  \xrightarrow{E} K_0(R\Gamma),$$
where $i_x \colon \IZ \to K_0(R[x])$ sends $n$ to $n \cdot [R[x]]$ and
$E$ has been defined in~\eqref{E_colon_SplitK_0_to_K_0}.

\begin{lemma}[Naturality of $\rk_{R\Gamma}$]
\label{lem:naturality_of_RGamma_rank} The rank $\rk_{R\Gamma}$ is
natural for functors $F\colon\Gamma_1 \to \Gamma_2$ between directly
finite categories. In particular, we have a natural transformation
$$\rk_{R-} \colon K_0(R-) \to U(-)$$ between covariant functors
$$K_0(R-), U(-) \colon \DIRFINCAT \to \ABELIANGROUPS.$$
\end{lemma}
\begin{proof}
The proof by~L\"uck~\cite[Proposition~10.44 (b) on page
202]{Lueck(1989)} for functors between EI-categories also works for
functors between directly finite categories. The rank
$\rk_{R\Gamma}$ is equal to $r \circ S$ where $r\colon \Split
K_0(R\Gamma) \to U(\Gamma)$ is the direct sum of
$$K_0(R[x]) \to \IZ$$
$$[P] \mapsto \rk_R (P \otimes_{R[x]} R)$$
over $\overline{x} \in \iso(\Gamma)$. By Lemma~\ref{lem:S_and_splitting},
the functor $S$ is covariantly natural
with respect to functors between directly finite categories. The
functor $r$ is also natural for such functors $F$, for if $F_x\colon
\aut_{\Gamma_1}(x) \to \aut_{\Gamma_2}(Fx)$ is the restriction of
$F$ to $\aut_{\Gamma_1}(x)$ we have
\[P \otimes_{R[x]} R \cong \ind_{F_x}(P) \otimes_{R[Fx]} R.\qedhere\]
\end{proof}

\begin{lemma} \label{lem:rank_RGamma_and_free_modules} Let $\Gamma$ be
  a directly finite category.

  \begin{enumerate}

  \item \label{lem:rank_RGamma_and_free_modules:iota_splits} The
    composite
$$U(\Gamma) \xrightarrow{\iota} K_0(R\Gamma) \xrightarrow{\rk_{R\Gamma}} U(\Gamma)$$
of the homomorphisms defined in~\eqref{rk_K_o_to_U} and~\eqref{iota}
is the identity;

\item \label{lem:rank_RGamma_and_free_modules_rank_determines_f.f.-modules}
Let $F$ be a finitely generated free $R\Gamma$-module. Then
$$F \cong \bigoplus_{\overline{x} \in \iso(\Gamma)}
\bigoplus_{i = 1} ^{\rk_{R\Gamma}(F)_{\overline{x}}} R\mor(?;x).$$ In
particular two finitely generated free $R\Gamma$-modules $F_1$ and
$F_2$ are isomorphic if and only if $\rk_{R\Gamma}(F_1) =
\rk_{R\Gamma}(F_2)$;
\end{enumerate}
\end{lemma}
\begin{proof}\ref{lem:rank_RGamma_and_free_modules:iota_splits}
  This follows from Lemma~\ref{lem:basic_properties_of-E_x_and_S_x}.
  \\[1mm]\ref{lem:rank_RGamma_and_free_modules_rank_determines_f.f.-modules}
Let $F$ be a free
  $R\Gamma$-module. By definition it looks like
$$F = \bigoplus_{\overline{x} \in \iso(\Gamma)} \bigoplus_{I_x} R\mor(?,x)$$
for some index sets $I_x$. It is finitely generated if there exist
natural numbers $m_x$ and an epimorphism
$$f \colon \bigoplus_{\overline{x} \in \iso(\Gamma)} \bigoplus_{i=1}^{m_x} R\mor(?,x)
\to \bigoplus_{\overline{x} \in \iso(\Gamma)} \bigoplus_{I_x}
R\mor(?,x)
$$
such that only finitely many $m_x$ are different from zero.
Lemma~\ref{lem:basic_properties_of-E_x_and_S_x} implies that we obtain
for every $\overline{x} \in \iso(\Gamma)$ an epimorphism $S_xf \colon
\bigoplus_{i=1}^{m_x} R[x] \to \bigoplus_{I_x} R[x]$.  This implies that each set
$I_x$ is finite and only finitely many of the sets $I_x$ are not
empty.  Hence we can find for a finitely generated free
$R\Gamma$-module $F$ natural numbers $n_x$ such that
$$F \cong \bigoplus_{\overline{x} \in \iso(\Gamma)} \bigoplus_{i=1}^{n_x} R\mor(?,x)$$
and only finitely many $n_x$ are different from zero.
Lemma~\ref{lem:basic_properties_of-E_x_and_S_x} implies
$$\rk_{R\Gamma}(F)_{\overline{x}} = n_x.$$
In particular $\rk_{R\Gamma}(F)$ determines the isomorphism type of a
finitely generated free $R\Gamma$-module $F$.
\end{proof}

\begin{definition}[The functorial Euler characteristic of a category]
  \label{def:functorial_Euler_characteristic_of_a_category}
  Suppose that $\Gamma$ is of type (FP$_R$).
  The \emph{functorial Euler characteristic of $\Gamma$ with coefficients in $R$},
\[
\chi_f(\Gamma;R)\in U(\Gamma),
\]
  is the image of the finiteness obstruction $o(\Gamma;R) \in K_0(R\Gamma)$ in
  Definition~\ref{def:finiteness_obstruction_of_a_module} under the
  homomorphism $\rk_{R\Gamma}\colon K_0(R\Gamma) \to U(\Gamma)$ in~\eqref{rk_K_o_to_U}.
\end{definition}

The word functorial refers to the fact that the group, in which
$\chi_f$ takes values, depends in a functorial way on $\Gamma$.

\begin{example}[The functorial Euler characteristic of a finite groupoid]
\label{exa:chi_f(calg)_for_finite_groupoids}
Let $\calg$ be a finite groupoid, i.e., a (small) groupoid such that
$\iso(\calg)$ and $\aut(x)$ for any object $x \in \ob(\calg)$ are finite.
Consider $R = \IQ$. Then $U(\calg)$ is the  abelian group generated by $\iso(\calg)$
and $\chi_f(\calg) \in U(\calg)$ is given by the sum of the basis elements.
\end{example}

\begin{theorem}[Invariance of the functorial Euler characteristic under equivalence of categories]
\label{invariance_of_functorial_Euler_characteristic_under_equivalence}
Let $F\colon \Gamma_1 \to \Gamma_2$ be an equivalence of categories
and suppose that $\Gamma_1$ is of type (FP$_R$) and directly finite.
Then $\Gamma_2$ is of type (FP$_R$) and directly finite, and
$$U(F) (\chi_f(\Gamma_1;R))=\chi_f(\Gamma_2;R).$$
\end{theorem}
\begin{proof}
The category $\Gamma_2$ is of type (FP$_R$) and
$F_*(o(\Gamma_1;R))=o(\Gamma_2;R)$ by
Theorem~\ref{the:invariance_under_equivalence}. The category
$\Gamma_2$ is directly finite by
Lemma~\ref{lem:direct_finiteness_under_equivalence}. We have $U(F)
(\chi_f(\Gamma_1;R))=\chi_f(\Gamma_2;R)$ by the naturality of
$\rk_{R-}$ in Lemma~\ref{lem:naturality_of_RGamma_rank} and
$F_*(o(\Gamma_1;R))=o(\Gamma_2;R)$.
\end{proof}

\begin{lemma} \label{lem:rank_RGamma_and_free_modules:o_and_chi_for_type_(FF)}
  Let $\Gamma$ be a directly finite category. Suppose that $\Gamma$ is of type (FF$_R$) (see
  Definition~\ref{def:Type_(FP)_and_(FF)_for_categories}).  Then the
  finiteness obstruction $o(\Gamma;R) \in K_0(R\Gamma)$ is the image
  of $\chi_f(\Gamma;R)$ under the homomorphism $\iota$ of~\eqref{iota}.
\end{lemma}
\begin{proof}
This follows from the definitions in combination with Lemma~\ref{lem:rank_RGamma_and_free_modules}.
\end{proof}

Obviously the functorial Euler characteristic $\chi_f(\Gamma;R)$ and
the topological Euler characteristic $\chi(B\Gamma;R)$  are weaker
invariants than the finiteness obstruction and carry less
information, but they live in explicit abelian groups and are easier
to compute.

\begin{theorem}[The finiteness obstruction determines the topological Euler characteristic]
  \label{the:relating_o_and_chi}
  Let $\Gamma$ be a category of type (FP$_R$). Suppose that $R$ is
  Noetherian.  We denote by $\pr \colon \Gamma \to \{\ast\}$ the
  projection to the trivial category with precisely one morphism.

  Then the assumptions in
  Definition~\ref{def:topological_Euler_characteristic_of_a_category} are
  satisfied and the composite
\begin{equation}\label{equ:relating_o_and_chi}
K_0(R\Gamma) \xrightarrow{\pr_*} K_0(R\{\ast\}) = K_0(R)
\xrightarrow{\rk_R} \IZ
\end{equation}
  sends the finiteness obstruction $o(\Gamma;R)$ to the topological Euler
  characteristic $\chi(B\Gamma;R)$.
\end{theorem}
\begin{proof}
  Associated to a category $\Gamma$ there is a classifying contravariant
  $\Gamma$-space $E\Gamma$ which is a $\Gamma$-$CW$-complex with the
  property that $E\Gamma$ evaluated at any object $x \in \ob(\Gamma)$
  is contractible. We refer to Davis--L\"uck~\cite[Definition~1.2,
  Definition~3.2, Definition~3.8, and page~230]{Davis-Lueck(1998)} for the definition
  of a contravariant $\Gamma$-space, a $\Gamma$-$CW$-complex (which is
  called free $\Gamma$-$CW$-complex there), the classifying
  $\Gamma$-space $E\Gamma$, and the bar construction. The cellular $R\Gamma$-chain complex
  $C_*(X)$ with $R$ coefficients of a $\Gamma$-$CW$-complex $X$ is the
  composition of the functor given by $X$ with the functor cellular
  chain complex with coefficients in $R$ and has free $R\Gamma$-chain
  modules. The proof of the last fact is analogous to the proof
  of~L\"uck~\cite[Lemma~13.2 on page~260]{Lueck(1989)}.  Since the evaluation
  of $E\Gamma$ at any object $x\in \ob(\Gamma)$ is contractible, the
  $R\Gamma$-module $H_n(C_*(E\Gamma;R))$ is trivial for $n \not=0$ and
  isomorphic to the constant $R\Gamma$-module $\underline{R}$ for $n = 0$. In
  particular $C_*(E\Gamma;R)$ is a projective $R\Gamma$-resolution of
  the constant $R\Gamma$-module $\underline{R}$.
  By assumption there exists a finite projective
  $R\Gamma$-resolution $P_*$ of $\underline{R}$. By the fundamental lemma of
  homological algebra (see~L\"uck~\cite[Lemma~11.3 on page~212]{Lueck(1989)})
  there exists an $R\Gamma$-chain homotopy equivalence $f_* \colon
  C_*(E\Gamma;R) \to P_*$.  If $\pr \colon \Gamma \to \{\ast\}$ is the
  projection to the trivial category, we obtain an $R$-chain homotopy
  equivalence $\ind_{\pr}f_* \colon \ind_{\pr} C_*(E\Gamma;R) \to
  \ind_{\pr}P_*$. There is also the notion of an induction functor for
  contravariant $\Gamma$-spaces
  (see~Davis--L\"uck~\cite[Definition~1.8]{Davis-Lueck(1998)} and a natural
  isomorphism of $R$-chain complexes
  $\ind_{\pr} C_*(E\Gamma;R)  \xrightarrow{\cong} C_*(\ind_{\pr}E\Gamma;R)$.
  The $CW$-complex $\ind_{\pr}E\Gamma$ is a model for $B\Gamma$
  (see~\cite[Definition~3.10, page~225 and page~230]{Davis-Lueck(1998)}).
  Hence we obtain a chain homotopy equivalence
  $$C_*(B\Gamma;R) \xrightarrow{\simeq} \ind_{\pr}P_*$$
  and $\ind_{\pr}P_*$ is an $R$-chain complex such that every $R$-chain
  module is finitely generated projective and only finitely many
  $R$-chain modules are non-trivial. Since $R$ is Noetherian, this
  implies that $H_n(\ind_{\pr}P_*)$ is finitely generated as an $R$-module
  for every $n \ge 0$ and that there is a natural number $d$ with
  $H_n(\ind_{\pr}P_*) = 0$ for $n > d$.  This implies that the same is true for the
  homology $H_*(B\Gamma;R)$. Our assumptions on the rank function
  $\rk_R$ of~\eqref{rk_R} imply
  $$\aligned
  \sum_{n \ge 0} (-1)^n \cdot \rk_R(\ind_{\pr} P_n) &= \sum_{n \ge 0}
  (-1)^n \cdot \rk_R(H_n(\ind_{\pr} P_*))
  \\
  &= \sum_{n \ge 0} (-1)^n \cdot \rk_R(H_n(B\Gamma)) \\ &= \chi(B\Gamma;R).
  \endaligned$$
  Since the composite
  $$K_0(R\Gamma) \xrightarrow{\pr_*} K_0(R\{\ast\}) = K_0(R) \xrightarrow{\rk_R} \IZ$$
  sends $o(\Gamma;R) = \sum_{n \ge 0} (-1)^n \cdot [P_n]$ to
  $\sum_{n \ge 0} (-1)^n \cdot \rk_R(\ind_{\pr} P_n)$,
  Theorem~\ref{the:relating_o_and_chi} follows.
\end{proof}

\begin{example}\label{exa:chi_and_chi_f_of_projection_category}
  Let $\Gamma$ be the category appearing in
  Example~\ref{exa:finiteness_obstruction_of_projection_category}.  It
  contains idempotents different from the identity, is directly
  finite, and of type (FP$_R$). We have $U(\Gamma) = \IZ$ and
$\chi_f(\Gamma;R) =
  \chi(B\Gamma;R) = 1$.
\end{example}

\begin{definition}[The Euler characteristic of a category] \label{def:Euler_characteristic_of_a_category}
Suppose that $\Gamma$ is of type (FP$_R$). The \emph{Euler
characteristic of $\Gamma$ with coefficients in $R$}
  is the sum of the components of the functorial Euler characteristic, that is,
\[\chi(\Gamma;R):=\epsilon(\chi_f(\Gamma;R)).\]
\end{definition}

\begin{theorem}[Invariance of the Euler characteristic under equivalence of
categories]
\label{thm:Euler_characteristic_invariant_under_equivalence} Let
$F\colon \Gamma_1 \to \Gamma_2$ be an equivalence of categories and
suppose that $\Gamma_1$ is of type (FP$_R$) and directly finite.
Then $\Gamma_2$ is of type (FP$_R$) and directly finite, and
$\Gamma_1$ and $\Gamma_2$ have the same Euler characteristic, that
is,
$$\chi(\Gamma_1;R)=\chi(\Gamma_2;R).$$
\end{theorem}
\begin{proof}
This follows from
Theorem~\ref{invariance_of_functorial_Euler_characteristic_under_equivalence}
and the naturality of the augmentation homomorphism in diagram
\eqref{equ:augmentation_naturality_diagram}.
\end{proof}

As we have seen in Theorem~\ref{the:relating_o_and_chi}, the
topological Euler characteristic is determined by the finiteness
obstruction when $\Gamma$ is of type (FP$_R$) and $R$ is Noetherian.
If we additionally assume $\Gamma$ is directly finite, then the
topological Euler characteristic and Euler characteristic agree.

\begin{theorem}[The Euler characteristic and topological Euler characteristic]
\label{the:chi_f_determines_chi} Let $R$ be a Noetherian ring and
$\Gamma$ a directly finite category of type (FP$_R$). Then the Euler
characteristic and topological Euler characteristic of $\Gamma$
agree.  That is, $H_n(B\Gamma;R)$ is a finitely generated $R$-module
for every $n \geq 0$, there exists a natural number $d$ with
$H_n(B\Gamma;R) = 0$ for all $n > d$, and
$$\chi(\Gamma;R) = \chi(B\Gamma;R)= \sum_{n \geq 0} (-1)^n \cdot \rk_{R}(H_n(B\Gamma;R)) \in \IZ,$$
where $\chi(\Gamma;R)$ is defined in
Definition~\ref{def:Euler_characteristic_of_a_category}.
\end{theorem}
\begin{proof}
Because of Theorem~\ref{the:relating_o_and_chi}, it suffices to show
that the diagram
$$\xymatrix@C=4pc{K_0(R\Gamma) \ar[r]^{\rk_{R\Gamma}} \ar[d]_{pr_*} & U(\Gamma) \ar[d]^\epsilon
\\ K_0(R\{\ast\})=K_0(R) \ar[r]_-{\rk_R} & \IZ}$$
commutes. However, this is precisely the $\rk_{R-}$ naturality
diagram associated to the functor $\Gamma \to \{*\}$. This diagram
commutes by Lemma~\ref{lem:naturality_of_RGamma_rank} because
$\Gamma$ and $\{*\}$ are directly finite categories.
\end{proof}

Euler characteristics are compatible with finite products. There is
an obvious pairing coming from the natural bijection $\iso(\Gamma_1)
\times \iso(\Gamma_2) \xrightarrow{\cong} \iso(\Gamma_1 \times
\Gamma_2)$
\begin{eqnarray}
\otimes \colon U(\Gamma_1) \otimes_{\IZ} U(\Gamma_2)
& \to &
U(\Gamma_1 \times \Gamma_2)
\label{product_pairing_for_U}
\end{eqnarray}

\begin{theorem}[Product formula for $\chi_f$, $\chi$, and $\chi(B-)$]
  \label{the:product_formula_for_chi_f_and_chi}
  Let $\Gamma_1$ and $\Gamma_2$ be categories of type (FP$_R$). Suppose
  that the rank $\rk_R$ satisfies $\rk_R(M \otimes N) = \rk_R(M) \cdot \rk_R(N)$
  for all finitely generated $R$-modules $M$ and $N$.

  Then $\Gamma_1 \times \Gamma_2$ is of type (FP$_R$), the
  functorial Euler characteristic satisfies
  \begin{eqnarray*}
    \chi_f(\Gamma_1 \times \Gamma_2;R)
    & = &
    \chi_f(\Gamma_1;R) \otimes \chi_f(\Gamma_2;R)
  \end{eqnarray*}
  under the pairing~\eqref{product_pairing_for_U}, the Euler
  characteristic satisfies
$$\chi(\Gamma_1 \times \Gamma_2;R) = \chi(\Gamma_1;R) \cdot \chi(\Gamma_2;R),$$
and the topological Euler characteristic satisfies
$$\chi(B(\Gamma_1 \times \Gamma_2);R) = \chi(B\Gamma_1;R) \cdot \chi(B\Gamma_2;R).$$
\end{theorem}
\begin{proof}
The product $\Gamma_1 \times \Gamma_2$ is of type (FP$_R$) by
Theorem~\ref{the:product_formula_for_o}.

Consider the diagram below,
  $$\xymatrix{K_0(R\Gamma_1) \otimes K_0(R\Gamma_2) \ar[r]^-{\otimes_R}
  \ar[d]_{\bigl(\rk_{R\Gamma_1} \circ S_{R\Gamma_1}\bigr) \otimes \bigl(\rk_{R\Gamma_2}
    \circ S_{R\Gamma_2}\bigr)} & K_0(R(\Gamma_1 \times \Gamma_2))
  \ar[d]^{\rk_{R(\Gamma_1 \times \Gamma_2)} \circ S_{R(\Gamma_1 \times \Gamma_2)}}
  \\
  U(\Gamma_1) \otimes U(\Gamma_2) \ar[r]_-{\otimes} & U(\Gamma_1 \times
  \Gamma_2) }
  $$
  where the horizontal pairings have been introduced
  in~\eqref{otimes_R_on_K_0} and~\eqref{product_pairing_for_U}, the
  homomorphisms $S$ in~\eqref{S_colon_SplitK_0_toK_0}, and the
  homomorphism $\rk_{R\Gamma}$ in~\eqref{rk_K_o_to_U}. One easily
  checks that it commutes.  Now the claim follows for $\chi_f$
  from Theorem~\ref{the:product_formula_for_o}.

  The claim for $\chi$ follows from that for $\chi_f$ because the
  pairing~\eqref{product_pairing_for_U} is compatible with the
  augmentation homomorphism.

  The claim for the topological Euler characteristic follows from the fact $B\Gamma_1 \times
  B\Gamma_2 = B (\Gamma_1 \times \Gamma_2)$ and the K\"unneth formula.
\end{proof}


\typeout{---  Section 5: The (functorial) L^2-Euler characteristic of a category -------}
\section{The  (functorial) $L^2$-Euler characteristic and $L^2$-Betti numbers
of a category}
\label{sec:The_(functorial)_L2_Euler_characteristic_anf_L2-Betti_numbers_of_a_category}

In this section we introduce the (functorial) $L^2$-Euler
characteristic and $L^2$-Betti numbers of a category. This requires
some input from the theory of finite von Neumann algebras and their
dimension theory which we briefly record next. For more information
we refer for instance to~L\"uck~\cite{Lueck(2002)},
\cite{Lueck(2003h)}.

In Subsection~\ref{subsec:Group_von_Neumann_Algebras_and_dimension}
we recall the group von Neumann algebra $\caln(G)$ associated to a
group $G$, the von Neumann dimension $\dim_{\caln(G)}$ for right
$\caln(G)$-modules, its properties, and compatibility with induction
and restriction for modules over group von Neumann algebras. For a
finite group $G$, the von Neumann algebra $\caln(G)$ is $\IC G$ and
the von Neumann dimension of a $\IC G$-module is the complex
dimension divided by $|G|$. For general $G$, the von Neumann algebra
$\caln(G)$ is a $\IC G$-$\caln(G)$-bimodule.

In Subsection~\ref{subsec:L2_Euler_characteristic_of_chain_complex}
we recall the $L^2$-Euler characteristic $\chi^{(2)}(C_*)$ of an
$\caln(G)$-chain complex $C_*$ as the alternating sum of the von
Neumann dimensions of the homology groups, and discuss the relevant
properties.

In Subsection~\ref{subsec:L2_Euler_characteristic} we define the
$L^2$-Euler characteristic for categories of type ($L^2$) using the
splitting functor $S_x$. A category $\Gamma$ is \emph{of type
($L^2$)} if the constant $\IC \Gamma$-module $\underline{\IC}$
admits a (not necessarily finite) projective $\IC\Gamma$-resolution
$P_*$ such that the sum over all $\overline{x} \in \iso(\Gamma)$ of
all von Neumann dimensions of the homology groups of all
$\caln(\aut(x))$-chain complexes $S_x P_* \otimes_{\IC\aut(x)}
\caln(\aut(x))$ converges to a finite number. Any directly finite
category of type (FP$_\IC$) is of type ($L^2$). For example, finite
groupoids, finite posets, and more generally finite EI-categories
are of type ($L^2$).

Let $U^{(1)}(\Gamma)$ be the set of absolutely convergent sequences
on the index set $\iso(\Gamma)$. The \emph{functorial $L^2$-Euler
characteristic} $\chi_f^{(2)}(\Gamma) \in U^{(1)}(\Gamma)$ has at
index $\overline{x}$ the number $\chi^{(2)}\bigl(S_x P_*
\otimes_{\IC\aut(x)} \caln(\aut(x))\bigr)$, where $P_*$ is a
projective $\IC \Gamma$-resolution of $\underline{\IC}$. The
\emph{$L^2$-Euler characteristic} $\chi^{(2)}(\Gamma)\in \IR$ is the
sum of the sequence $\chi_f^{(2)}(\Gamma)$. For example, if $\Gamma$
is a finite groupoid, then $\chi_f^{(2)}(\Gamma)$ has at index
$\overline{x}$ the value $1/|\aut(x)|$, and the $L^2$-Euler
characteristic is the sum of these.

Like the topological Euler characteristic and the Euler
characteristic, the $L^2$-Euler characteristic comes from the
finiteness obstruction in certain cases. However, for the
$L^2$-Euler characteristic, we use the $L^2$-rank $\rk_\Gamma^{(2)}$
instead of the $R\Gamma$-rank $\rk_{R\Gamma}$. In
Subsection~\ref{subsec:finiteness_obstruction_and_L2_Euler_characteristic}
we define the $L^2$-rank and prove that
$\rk_\Gamma^{(2)}o(\Gamma;\IC)=\chi_f^{(2)}(\Gamma)$ whenever
$\Gamma$ is directly finite and of type (FP$_\IC$).

The $L^2$-Euler characteristic is compatible with covering maps and
isofibrations between connected finite groupoids, as we prove in
Subsection~\ref{subsec:Compatibility_with_Coverings}.

We now recall the prerequisites from the theory of finite von
Neumann algebras and motivate its use.

\subsection{Group von Neumann algebras and their dimension
theory}\label{subsec:Group_von_Neumann_Algebras_and_dimension}

The appearence of (group) von Neumann algebras and their dimension theory in our context stems
from the task to assign some sort of rational- or real-valued dimension to projective
modules over group rings (coming from automorphism groups in a category), which itself is needed to extract a number, namely the Euler characteristic, from the finiteness obstruction.

The well-known \emph{Hattori-Stallings rank} $\HS(M)$
in~Brown~\cite[Chapter IX, 2]{Brown(1982)} of a finitely generated
projective $R$-module $M$ over an arbitrary ring $R$ is a way to
assign a ``dimension'' to $M$. However, $\HS(M)$ is not a number but
an element in the quotient $R/[R,R]$ of $R$ by the additive subgroup
$[R,R]$ generated by all commutators $ab-ba$, $a,b\in R$. In order
to get, say, a $\IC$-valued invariant one needs an additive
homomorphism $t:R\to\IC$ satisfying the \emph{trace property}
$t(ab)=t(ba)$.

Consider the case of the complex group ring $R=\IC G$ of a group $G$. The map
    $\tr_{\caln(G)}\colon  \IC G \to \IC$, the notation of which already anticipates
    a more general setup, is defined by
\[
    \tr_{\caln(G)}\bigl(\sum_{g\in G}\lambda_g g\bigr)=\lambda_e
\]
and satisfies the trace property,
thus providing a notion of dimension for finitely generated projective $\IC G$-modules.
This dimension does not extend to arbitrary $\IC G$-modules, which is a major drawback
as we would like to define the dimension of certain homology groups of
projective resolutions that are \emph{not} projective anymore.
Next we explain work of the second
author that allows to define
a dimension for all modules -- if one works with the larger ring $\caln(G)$, the group von Neumann algebra of $G$, instead.

Let $l^2(G)$ be the Hilbert space with Hilbert basis $G$; it
consists of formal sums $\sum_{g \in G} \lambda_g \cdot g$ for
complex numbers $\lambda_g$ such that $\sum_{g \in G} |\lambda_g|^2
< \infty$. The complex group ring $\IC G$ is a dense subset of
$l^2(G)$. In fact, $l^2(G)$ is the Hilbert space completion of the
complex group ring $\IC G$ with respect to the pre-Hilbert structure
for which $G$ is an orthonormal basis. Left and right multiplication
with elements in $G$ induce respectively isometric left and right
$G$-actions on $l^2(G)$.

\begin{definition}[Group von Neumann algebra]
\label{exa:group_von_neumann_algebra}

The \emph{group von Neumann algebra} of the group $G$
\begin{equation*}
\caln(G) =  \calb(l^2(G))^G
\end{equation*}
is the algebra of bounded operators that are equivariant with
respect to the right $G$-action. The \emph{standard trace} on
$\caln(G)$ is defined by
$$\tr_{\caln(G)} \colon \caln(G) \to \IC,
\quad f \mapsto \langle f(e),e \rangle_{l^2(G)}.$$
\end{definition}

The standard trace extends the definition on $\IC G$ given earlier on.
From now on we view $\caln(G)$ simply as a ring, ignoring its functional-analytic
origin. The latter is only important for the proof of our `blackbox'
Theorem~\ref{the:properties_of_the_dimension_function} below.
Modules over $\caln(G)$ are understood in the purely algebraic sense.

Sending an element $g \in G$ to the isometric $G$-equivariant
operator $l^2(G) \to l^2(G)$ given by left multiplication with $g
\in G$ induces an embedding of $\IC G$ into $\caln(G)$ as a subring.
In particular, we can view $\caln(G)$ as a $\IC
G$-$\caln(G)$-bimodule.

\begin{theorem}[Properties of the dimension function]
  \label{the:properties_of_the_dimension_function}
There exists a dimension function $\dim_{\caln(G)}$ that assigns to
every right $\caln(G)$-module $M$ a number, possibly infinite,
\[\dim_{\caln(G)}(M)\in [0,\infty]=\IR_{\ge 0}\cup\{\infty\}\]
and satisfies the following properties:

\begin{enumerate}
\item
\label{the:properties_of_the_dimension_function:Hattori-Stallings-rank}
Hattori-Stallings rank\\[1mm]
If $M$ is a finitely generated projective $\caln(G)$-module, then
$$\dim_{\caln(G)}(M) = \sum_{i=1}^n \tr_{\caln(G)}(a_{i,i}) \in [0,\infty),$$
where $A = (a_{i,j})$ is any $(n,n)$-matrix over $\caln(G)$ with $A^2 = A$
such that the image of the $\caln(G)$-homomorphism $\caln(G)^n \to \caln(G)^n$
given by left multiplication with $A$ is $\caln(G)$-isomorphic to $M$;

\item
\label{the:properties_of_the_dimension_function:additivity}
Additivity\\[1mm]
If $0 \to M_0 \to M_1 \to M_2 \to 0$ is an
exact sequence of $\caln(G)$-modules, then
$$\dim_{\caln(G)}(M_1) =  \dim_{\caln(G)}(M_0) + \dim_{\caln(G)}(M_2),$$
where for $r,s \in [0,\infty]$ we define $r + s$  by the ordinary sum of two
real numbers if both $r$ and $s$ are not $\infty$, and by $\infty$ otherwise;

\item
\label{the:properties_of_the_dimension_function:cofinality}
Cofinality\\[1mm]
Let $\{M_i\mid i \in I\}$ be a cofinal system of submodules of $M$,
i.e., $M = \bigcup_{i \in I} M_i$ and for two indices $i$ and $j$
there is an index $k$ in $I$ satisfying $M_i,M_j \subseteq M_k$.
Then
$$\dim_{\caln(G)}(M)  =  \sup\{\dim_{\caln(G)}(M_i) \mid i \in I\}.$$
\end{enumerate}
\end{theorem}
\begin{proof} See~L\"uck~\cite[Theorem~6.5 and Theorem~6.7  on  page~239]{Lueck(2002)}.
\end{proof}

Let $i \colon H \to G$ be an injective group homomorphism.  Then the
induced injective ring homomorphism $i_* \colon\IC H \to \IC G$
extends to an injective ring homomorphism denoted in the same way
$i_* \colon \caln(H) \to \caln(G)$.

\begin{lemma} \label{lem:dim_caln(G)_and_induction_and_restriction}
Let $i \colon H \to G$ be an injective group homomorphism.

\begin{enumerate}

\item \label{lem:dim_caln(G)_and_induction_an_restriction:faithfully_flat}
The induction functor $\ind_{i_*} \colon \MOD\text{-}\caln(H) \to \MOD\text{-}\caln(G)$
sending $M$ to $M \otimes_{\caln(H)} \caln(G)$
is faithfully flat, i.e., a sequence of $\caln(H)$-modules
$M_1 \to M_2 \to M_3 $ is exact if and only if the induced
sequence of $\caln(G)$-modules $\ind_{i_*} M_1  \to \ind_{i_*} M_2 \to \ind_{i_*}M_3$
is exact;

\item \label{lem:dim_caln(G)_and_induction_and_restriction:dim_and_induction}
If $M$ is an $\caln(H)$-module, then
$$\dim_{\caln(G)}(\ind_{i_*}M) = \dim_{\caln(H)}(M);$$

\item \label{lem:dim_caln(G)_and_induction_and_restriction:dim_and_restriction}
Suppose that the index $[G:i(H)]$ of $i(H)$ in $G$ is finite. Then
we get for every $\caln(G)$-module $M$, if $\res_{i_*}$ denotes its
restriction to an $\caln(H)$-module by $i_*$
$$\dim_{\caln(H)}(\res_{i_*} M) = [G:i(H)] \cdot \dim_{\caln(G)}(M),$$
where $[G:i(H)] \cdot \infty$ is defined to be $\infty$.

\end{enumerate}
\end{lemma}
\begin{proof}
See~L\"uck~\cite[Theorem~6.29 on page~253 and Theorem~6.54 (6) on
page~266]{Lueck(2002)}.
\end{proof}

Here are some useful examples of the von Neumann dimension.

\begin{example}\label{exa:examples of the dimension}\
    \begin{enumerate}
    \item
  \label{exa:The_von_Neumann_dimension_for_finite_groups}
  (von Neumann dimension for finite groups){\bf .} Let $G$ be a finite group. Then $\caln(G) = \IC G$ and we get for a
  $\IC G$-module $M$
  $$\dim_{\caln(G)}(M) = \frac{1}{|G|} \cdot \dim_{\IC}(M);$$
  where $\dim_{\IC}$ is the dimension of $M$ viewed as a complex vector space.
    \item\label{exa:von_Neumann_dimension_and_permutation_modules}(von Neumann dimension and permutation
    modules){\bf .} Let $G$ be a (not necessarily finite) group and $S$ a cofinite $G$-set, i.e., $S$ is the disjoint union of
homogeneous $G$-spaces $\coprod_{i \in I} G/L_i$ for finite $I$.
By~L\"uck~\cite[Lemma~4.4]{Lueck(1998a)}, we have
$$\dim_{\caln(G)}\bigl(\IC S \otimes_{\IC G} \caln(G)\bigr) = \sum_{\substack{i \in I\\|L_i| < \infty}} \frac{1}{|L_i|}.$$
    \item\label{exa:The_von_Neumann_dimension_for_Z}(von Neumann dimension for
    $\IZ$){\bf .}
Let $G = \IZ$. Then $\caln(\IZ) = L^{\infty}(S^1)$ by Fourier transformation.
Under this identification we obtain that
$$\tr_{\caln(\IZ)} \colon \caln(\IZ) \to \IC, \quad f \mapsto \int_{S^1} f d\mu,$$
where $\mu$ is the probability Lebesgue measure on $S^1$.

Let $X \subseteq S^1$ be any measurable set and $\chi_X \in L^{\infty}(S^1)$
be its characteristic function. Since $\chi_X$ is an idempotent,
its image $P$ is a finitely generated projective
$\caln(\IZ)$-module, whose von Neumann dimension
$\dim_{\caln(\IZ)}(P)$ is the volume $\mu(X)$ of $X$.
In particular any non-negative real number occurs as
$\dim_{\caln(\IZ)}(P)$ for some finitely generated projective
$\caln(\IZ)$-module $P$.
\end{enumerate}
\end{example}

\subsection{The $L^2$-Euler characteristic and $L^2$-Betti
numbers}\label{subsec:L2_Euler_characteristic_of_chain_complex}

In this section we briefly recall some basic facts about $L^2$-Betti
numbers and $L^2$-Euler characteristics. For more information we
refer to~L\"uck~\cite[Section~6.6.1 on~page~277ff]{Lueck(2002)}.

\begin{definition}[$L^2$-Betti numbers]
\label{def:L2-Betti_number} Let $C_*$ be an $\caln(G)$-chain
complex. The \emph{$p$-th $L^2$-Betti number of $C_*$} is the von
Neumann dimension of the $\caln(G)$-module given by its $p$-th
homology, namely
$$b_p^{(2)}(C_*) := \dim_{\caln(G)}(H_p(C_*)) \quad \in [0,\infty].$$
\end{definition}

\begin{definition}[$L^2$-Euler characteristic]\label{def:L2-Euler_characteristic}
Let $C_*$ be an $\caln(G)$-chain complex. Define
$$h^{(2)}(C_*) := \sum_{p \ge 0} b_p^{(2)}(C_*) \quad \in [0,\infty].$$
If $h^{(2)}(C_*) < \infty$, the \emph{$L^2$-Euler characteristic of
$C_*$} is
$$\chi^{(2)}(C_*) := \sum_{p \ge 0} (-1)^p \cdot b_p^{(2)}(C_*) \quad \in \IR.$$
\end{definition}

Notice that $h^{(2)}(C_*)$ can be finite also in the case, where infinitely many
$L^2$-Betti numbers are different from zero.

\begin{lemma}
\label{lem:basic_properties_of_L2-Euler_characteristic}\
  \begin{enumerate}
  \item \label{lem:basic_properties_of_L2-Euler_characteristic:homology}
    Let $C_*$ be an $\caln(G)$-chain complex. Suppose that
   $\sum_{p \ge 0} \dim_{\caln(G)}(C_p)$ is finite.
    Then $h^{(2)}(C_*)$ is finite and
   $\sum_{p \ge 0} (-1)^p \cdot \dim_{\caln(G)}(C_p)=\chi^{(2)}(C_*)$;

  \item \label{lem:basic_properties_of_L2-Euler_characteristic:homotopy_invariance}
  Let $C_*$ and $D_*$ be $\caln(G)$-chain complexes which are $\caln(G)$-homotopy equivalent.
   Then we get $b_p^{(2)}(C_*) = b_p^{(2)}(D_*)$ and $h^{(2)}(C_*) = h^{(2)}(D_*)$ and,
   provided that $h^{(2)}(C_*)$ is finite, $\chi^{(2)}(C_*) = \chi^{(2)}(D_*)$;

  \item \label{lem:basic_properties_of_L2-Euler_characteristic:Additivity}
    Let $0 \to C_* \to D_* \to E_* \to 0$ be an exact sequence of
    $\caln(G)$-chain complexes.  Suppose that two of the elements
    $h^{(2)}(C_*)$, $h^{(2)}(D_*)$, and $h^{(2)}(E_*)$ in $[0,\infty]$
    are finite.
    Then this is true for all three and we obtain that
    $$\chi^{(2)}(C_*) - \chi^{(2)}(D_*) + \chi^{(2)}(E_*) = 0;$$

    \item \label{lem:basic_properties_of_L2-Euler_characteristic:induction}
    Let $i \colon H \to G$ be an injective group homomorphism and let
    $C_*$ be an $\caln(H)$-chain complex. Then $h^{(2)}(C_*) = h^{(2)}(\ind_{i_*} C_*)$
    and, provided that $h^{(2)}(C_*) < \infty$, we have
    $\chi^{(2)}(C_*) = \chi^{(2)}(\ind_{i_*} C_*)$;

    \item \label{lem:basic_properties_of_L2-Euler_characteristic:restriction}
    Let $i \colon H \to G$ be an injective group homomorphism with finite index
    \mbox{$[G:i(H)]$}. Let $C_*$ be an $\caln(G)$-chain complex.
    Then
    $$h^{(2)}(\res_{i_*} C_*) = [G:i(H)] \cdot  h^{(2)}(C_*)$$
    and, provided that $h^{(2)}(C_*) < \infty$, we have
    $$\chi^{(2)}(\res_{i_*}  C_*) = [G:i(H)] \cdot \chi^{(2)}( C_*).$$
\end{enumerate}
\end{lemma}
\begin{proof} ii) is obvious from the definition. The rest easily follows from Theorem~\ref{the:properties_of_the_dimension_function}
and Lemma~\ref{lem:dim_caln(G)_and_induction_and_restriction}.
\end{proof}

\subsection{The (functorial) $L^2$-Euler characteristic}\label{subsec:L2_Euler_characteristic}

In the following, $\Gamma$ is always a small category. For every $x
\in \ob(\Gamma)$ let
$$\caln(x):=\caln(\aut(x))$$
be the group von Neumann algebra of the automorphism group $\aut(x)$.

Recall that two projective $\caln(G)$-resolutions $P_*$ and $Q_*$
of the constant $\IC\Gamma$-module $\underline{\IC}$ are $\IC\Gamma$-chain
homotopy equivalent and hence the $\IC[x]$-chain complexes $S_xP_*$ and
$S_xQ_*$ and the $\IC[x]$-chain complexes $\Res_xP_*$ and
$\Res_x Q_*$ are $\IC[x]$-chain homotopy equivalent.
Therefore the following definitions will be independent of the choice of a
projective $\IC \Gamma$-resolution of $\underline{\IC}$.

\begin{definition}[Type ($L^2$)] \label{def:Type_(L^2)}
We call $\Gamma$ \emph{of type }($L^2$) if for some (and hence for
every) projective $\IC\Gamma$-resolution $P_*$ of the constant $\IC
\Gamma$-module $\underline{\IC}$ we have
$$\sum_{\overline{x} \in \iso{\Gamma}} h^{(2)}\bigl(S_x P_*\otimes_{\IC [x]}\caln(x)\bigr) < \infty.$$
\end{definition}
We shall see in Example \ref{exa:chi_f(calg)L2_for_finite_groupoids}
that any finite groupoid is of type ($L^2$). We shall see in
Theorem \ref{the:comparing_o_and_chi(2)} that any directly finite
category of type (FP$_\IC$) is of type ($L^2$).

\begin{definition}[The functorial $L^2$-Euler characteristic of a category]
  \label{def:functorial_L2-Euler_characteristic_of_a_category}
  Suppose that $\Gamma$ is of type ($L^2$) and let
  $$U^{(1)}(\Gamma) :=
  \biggl\{\sum_{\overline{x} \in \iso(\Gamma)} r_{\overline{x}} \cdot \overline{x}
  \;\bigg|\;
  r_{\overline{x}} \in \IR,
  \sum_{\overline{x} \in \iso(\Gamma)} |r_{\overline{x}}| < \infty\biggr\}\subseteq\prod_{\bar{x}\in\iso(\Gamma)}\IR.$$
  The \emph{functorial $L^2$-Euler  characteristic of  $\Gamma$} is
  \[ \chi_f^{(2)}(\Gamma):=\Bigl \{\chi^{(2)}\bigl(S_xP_*\otimes_{\IC [x]} \caln(x)\bigr)\mid \bar{x} \in \iso(\Gamma) \Bigr\} \in U^{(1)}(\Gamma),
  \]
  where $P_*$ is a projective $\IC \Gamma$-resolution of the constant $\IC\Gamma$-module
  $\underline{\IC}$.
\end{definition}

The word functorial refers to the fact that the group
$U^{(1)}(\Gamma)$, in which $\chi_f^{(2)}$ takes values, depends in
a functorial way on $\Gamma$.

We can also get a real-valued invariant as follows.

\begin{definition}[The $L^2$-Euler characteristic of a category]
  \label{def:L2-Euler_characteristic_of_a_category}
  Suppose that $\Gamma$ is of type ($L^2$).
  The \emph{$L^2$-Euler characteristic of $\Gamma$} is the sum over $\bar{x} \in
\iso(\Gamma)$ of the components of the functorial Euler
characteristic, that is,
\[
\chi^{(2)}(\Gamma) := \sum_{\overline{x} \in
\iso(\Gamma)} \chi^{(2)}\bigl(S_xP_*\otimes_{\IC [x]}
\caln(x)\bigr)\in\IR,
\]
where $P_*$ is a projective $\IC \Gamma$-resolution of the constant
$\IC\Gamma$-module $\underline{\IC}$.
\end{definition}

Notice that this definition makes sense since the condition ($L^2$)
ensures that the sum
$\sum_{\overline{x} \in \iso(\Gamma)} \chi^{(2)}\bigl(S_xP_* \otimes_{\IC [x]} \caln(x)\bigr)$ is
absolutely convergent.

\begin{remark} \label{rem:L2_Euler_char_defn_like_augmentation} In
  Definition~\ref{def:L2-Euler_characteristic_of_a_category}, the $L^2$-Euler
  characteristic is defined to be the sum of the components of the functorial
  $L^2$-Euler characteristic. This is analogous to the situation for the
  ordinary Euler characteristic in
  Definition~\ref{def:Euler_characteristic_of_a_category}.
\end{remark}

\begin{example}[The (functorial) $L^2$-Euler characteristic of groupoids]
  \label{exa:chi_f(calg)L2_for_finite_groupoids}
  Let $\calg$ be a (small) groupoid such that
  $\aut_\calg(x)$ is
  finite for any object $x \in \ob(\calg)$  and
  \begin{equation}\label{equ:finiteness of aut-sum}
  \sum_{\overline{x} \in \iso(\calg)} \frac{1}{|\aut_\calg(x)|}<\infty.
  \end{equation}
  Let $P_\ast$ be any projective $\IC\calg$-resolution of $\underline{\IC}$; a (not necessarily finite) projective
  resolution always exists.
  Since $\calg$ is a groupoid, for every $x\in\ob\calg$ and every $\IC\calg$-module $M$
  we have $S_x M=\Res_x M$.
  Thus $S_x$ is exact. By Lemma~\ref{lem:basic_properties_of-E_x_and_S_x}, $S_x$ respects
  projectives. Hence $S_xP_x$ is a projective $\IC [x]$-resolution of the trivial $\IC [x]$-module
  $\IC$. Since $\aut_\calg(x)$ is finite, $\IC$ is already a projective $\IC [x]$-module.
  This implies that
  \[H_p\bigl(S_xP_*\otimes_{\IC [x]} \caln(x)\bigr)=\begin{cases}
                                                      \IC\otimes_{\IC [x]}\caln(x)& p=0\\
                                                      0                           & p>0.
                                                    \end{cases}\]
  Example~\ref{exa:examples of the dimension}~\ref{exa:The_von_Neumann_dimension_for_finite_groups}
  and~\eqref{equ:finiteness of aut-sum} yield that
  $\calg$ is of type ($L^2$), the
  functorial $L^2$-Euler characteristic
  $\chi_f^{(2)}(\calg) \in \prod_{\overline{x} \in \iso(\calg)} \IR$
  has at $\overline{x} \in \iso(\calg)$ the value $1/|\aut_\calg(x)|$, and
  $$\chi^{(2)}(\calg) = \sum_{\overline{x} \in \iso(\calg)} \frac{1}{|\aut_\calg(x)|}.$$

  In particular, we can
  conclude that, for all groupoids such that \eqref{equ:finiteness of aut-sum} holds,
  the $L^2$-Euler characteristic coincides with the Baez--Dolan groupoid cardinality, and also with Leinster's Euler
  characteristic when the groupoid is finite.

A concrete case of a groupoid satisfying our conditions is a
skeleton $\calg$ of the groupoid of nonempty finite sets. This
groupoid has objects (isomorphic to) $\underline{1}=\{1\}$,
$\underline{2}=\{1,2\}$, $\underline{3}=\{1,2,3\}$, and so on. The
morphisms are the permutations. This example was studied by
Baez--Dolan~\cite{Baez-Dolan(2001)}. The groupoid $\calg$ is of type
($L^2$), and the functorial $L^2$-Euler characteristic has at the
object $\underline{n}$ the value
$1/|\aut_\calg(\underline{n})|=1/n!$. The $L^2$-Euler characteristic
is
$$\chi^{(2)}(\calg)=\sum_{n \geq 1} \frac{1}{|S_n|}=\sum_{n \geq 1} \frac{1}{n!}=e.$$
\end{example}

\begin{remark} \label{exa:Gamma_a_group} If $G$ is a group and
  $\widehat{G}$ denotes the groupoid with precisely one object and $G$ as
  automorphism group of this object, then $\chi^{(2)}(\widehat{G})$ in
  the sense of
  Definition~\ref{def:L2-Euler_characteristic_of_a_category} agrees
  with the classical definition of the $L^2$-Euler characteristic
  $\chi^{(2)}(G)$ of a group which has been intensively studied in the
  literature (see for instance~L\"uck~\cite[Chapter~7]{Lueck(2002)}).
\end{remark}

\begin{lemma}[Invariance of $L^2$-Euler characteristic under equivalence of categories]\label{lem:L2_under_equivalence}\
\begin{enumerate}
\item \label{lem:L2_under_equivalence:(L2)}
Suppose $\Gamma_1$ and $\Gamma_2$ are equivalent categories. Then
$\Gamma_1$ is both directly finite and of type ($L^2$) if and only if $\Gamma_2$ is
both directly finite and of type ($L^2$).

\item \label{lem:L2_under_equivalence:chi}
Let $F \colon \Gamma_1 \to \Gamma_2$ be an equivalence of categories.
Suppose that $\Gamma_i$ is both directly finite and of type  ($L^2$) for $i = 1,2$.

Then the bijection
$$U^{(1)}(F) \colon U^{(1)}(\Gamma_1) \xrightarrow{\cong} U^{(1)}(\Gamma_2)$$
induced by $F$ sends $\chi^{(2)}_f(\Gamma_1)$ to $\chi^{(2)}_f(\Gamma_2)$ and we have
$$\chi^{(2)}(\Gamma_1) = \chi^{(2)}(\Gamma_2).$$
\end{enumerate}
\end{lemma}
\begin{proof}
We have already shown that the property of being directly finite depends only on the
equivalence class of a category (see Lemma~\ref{lem:direct_finiteness_under_equivalence}).
So in the sequel we can assume that $\Gamma_1$ and $\Gamma_2$ are directly finite.

Let $F \colon \Gamma_1 \to \Gamma_2$ be an
equivalence of categories. It induces a bijection
$$F_* \colon \iso(\Gamma_1) \xrightarrow{\cong} \iso(\Gamma_2),
\quad \overline{x}\mapsto \overline{F(x)},$$
and thus a bijection
$$U^{(1)}(F) \colon U^{(1)}(\Gamma_1) \xrightarrow{\cong} U^{(1)}(\Gamma_2).$$
Recall from Section~\ref{sec:basics_about_modules_over_a_category}
that the induction functor $\ind_F$ associated to $F$ sends
projective $\IC\Gamma_1$-modules to projective
$\IC\Gamma_2$-modules. The equivalence $F$ induces for every object
$x$ in $\Gamma_1$ an isomorphism of groups
$$F_x \colon \aut_{\Gamma_1}(x) \xrightarrow{\cong} \aut_{\Gamma_2}(F(x)), \quad f \mapsto F(f).$$
From the proof of Lemma~\ref{lem:S_and_splitting}, we have for every
object $x$ in $\Gamma_1$ and projective $\IC\Gamma_1$-module $P$ a
natural isomorphism of $\IC[F(x)]$-modules
$$\alpha(P) \colon \ind_{F_x} \circ S_x P \xrightarrow{\cong} S_{F(x)} \circ \ind_F P$$
(the direct sum in the proof of Lemma~\ref{lem:S_and_splitting} has
only one summand because $F$ is an equivalence).

Fix an object $x$ in $\Gamma_1$.
The argument in the proof of Theorem~\ref{the:adjoint_functors} shows that the
induction functor $\ind_F$ associated to $F$ is an exact functor and sends
$\underline{\IC}$ to $\underline{\IC}$.  Let $P_*$ be a free
$\IC \Gamma_1$-resolution of $\underline{\IC}$.
Then $\ind_F P_*$ is a free $\IC \Gamma_2$-resolution of $\underline{\IC}$.
The various isomorphisms $\alpha(P_n)$ induce an isomorphism of
$\IC[F(x)]$-chain complexes
$$\alpha(P_*)\colon \ind_{F_x} \circ S_x P_* \xrightarrow{\cong} S_{F(x)} \circ \ind_F P_*.$$
We have for every $R[x]$-module $M$ a canonical
$\caln(F(x))$-isomorphism
$$\bigl(\ind_{F_x} M\bigr) \otimes_{\IC[F(x)]} \caln(F(x))
\xrightarrow{\cong} \ind_{F_x} \bigl(M \otimes_{\IC[x]} \caln(x)\bigr).$$
If we apply $-\otimes_{\IC[F(x)]} \caln(F(x))$ to $\alpha(P_*)$ and use the isomorphisms above
we obtain an isomorphism of $\caln(F(x))$-chain complexes
$$\alpha^{(2)}(P_*) \colon \ind_{F_x} \bigl(S_x P_* \otimes_{\IC[x]} \caln(x)\bigr)
\xrightarrow{\cong} \bigl(S_{F(x)} \circ \ind_F P_*\bigr) \otimes_{\IC[F(x)]} \caln(F(x)).$$
We conclude from
Lemma~\ref{lem:basic_properties_of_L2-Euler_characteristic}~%
\ref{lem:dim_caln(G)_and_induction_and_restriction:dim_and_induction}
\begin{eqnarray*}
h^{(2)}\bigl(S_x P_* \otimes_{\IC[x]} \caln(x) \bigr) & = &
h^{(2)}\left(\bigl(S_{F(x)} \circ \ind_F P_*\bigr)
\otimes_{\IC[F(x)]} \caln(F(x))\right)
\end{eqnarray*}
and, provided that $h^{(2)}\bigl(S_x P_* \otimes_{\IC[x]} \caln(x)
\bigr) < \infty$
\begin{eqnarray*}
\chi^{(2)}\bigl(S_x P_* \otimes_{\IC[x]} \caln(x) \bigr) & = &
\chi^{(2)}\left(\bigl(S_{F(x)} \circ \ind_F P_*\bigr)
\otimes_{\IC[F(x)]} \caln(F(x))\right).
\end{eqnarray*}
Now Lemma~\ref{lem:L2_under_equivalence} follows.
\end{proof}

Next we consider products of categories.
Since $\iso(\Gamma_1 \times \Gamma_2) = \iso(\Gamma_1) \times
\iso(\Gamma_2)$, we obtain a pairing
\begin{multline}
\otimes \colon U^{(1)}(\Gamma_1) \otimes  U^{(1)}(\Gamma_2)
\to  U^{(1)}(\Gamma_1 \times \Gamma_2),
\\
\sum_{\overline{x_1} \in \iso(\Gamma_1)} r_{\overline{x_1}} \cdot \overline{x_1} \; \otimes
\sum_{\overline{x_2} \in \iso(\Gamma_2)} s_{\overline{x_2}} \cdot \overline{x_2}
\mapsto
\sum_{\overline{(x_1,x_2)} \in \iso(\Gamma_1 \times \Gamma_2)}
r_{\overline{x_1}}  s_{\overline{x_2}} \cdot \overline{(x_1,x_2)}.
\label{pairing_for_prod_iso_R}
\end{multline}

\begin{theorem}[Product formula for $\chi_f^{(2)}$ and $\chi^{(2)}$]
  \label{the:product_formula_for_chi_f(2)_and_chi(2)}
  Let $\Gamma_1$ and $\Gamma_2$ be categories of type ($L^2$).

  Then $\Gamma_1 \times \Gamma_2$ is of type ($L^2$), we get for the
  functorial $L^2$-Euler characteristic
  \begin{eqnarray*}
    \chi_f^{(2)}(\Gamma_1 \times \Gamma_2)
    & = &
    \chi_f^{(2)}(\Gamma_1) \otimes \chi_f^{(2)}(\Gamma_2)
  \end{eqnarray*}
  under the pairing~\eqref{pairing_for_prod_iso_R}, and we get for the
  $L^2$-Euler characteristic
  $$\chi^{(2)}(\Gamma_1 \times \Gamma_2)
  = \chi^{(2)}(\Gamma_1) \cdot \chi^{(2)}(\Gamma_2).$$
\end{theorem}
\begin{proof}
  If $P_*$ is a projective $\IC \Gamma_1$-resolution of the constant
  $\IC\Gamma_1$-module $\underline{\IC}$ and $Q_*$ is a projective
  $\IC \Gamma_2$-resolution of the constant $\IC\Gamma_2$-module
  $\underline{\IC}$, then $P_* \otimes Q_*$ is a projective
  $\IC (\Gamma_1 \times \Gamma_2)$-resolution of the constant
  $\IC(\Gamma_1 \times \Gamma_2)$-module $\underline{\IC}$. Given
  $\overline{x} \in \iso(\Gamma_1)$ and $\overline{y} \in \iso(\Gamma_2)$, there
  is a canonical isomorphism of chain complexes over
  $\IC[(x,y)] = \IC[x] \otimes_{\IC} \IC[y]$
  $$S_xP_* \otimes_{\IC} S_yP_* = S_{(x,y)}(P_* \otimes_{\IC} Q_*).$$
  Since the Cauchy product of two absolutely convergent series of real
  numbers is again an absolutely convergent series, it suffices to show
  for two groups $H$ and $G$, a projective $\IC H$-chain complex $C_*$
  and a projective $\IC G$-chain complex $D_*$, that for the projective
  $\IC [G \times H]$-chain $C_* \otimes_{\IC} D_*$ we have
  \begin{eqnarray*}
  h^{(2)}(C_* \otimes_{\IC} D_*)  & <  & \infty;
  \\
  \chi^{(2)}(C_* \otimes_{\IC} D_*) & = & \chi^{(2)}(C_*) \cdot \chi^{(2)}(D_*)
  \end{eqnarray*}
  provided that $h^{(2)}(C_*)$ and $h^{(2)}(D_*)$ are finite.  The proof
  of this claim is the chain complex analogue of the proof
  of~L\"uck~\cite[Theorem~6.80~(6) on page~278]{Lueck(2002)}.
\end{proof}

\subsection{The finiteness obstruction and the (functorial) $L^2$-Euler
characteristic}\label{subsec:finiteness_obstruction_and_L2_Euler_characteristic}
Next we compare these definitions with the finiteness obstruction
and Euler characteristic.

\begin{definition}[$L^2$-rank of a finitely generated $\IC\Gamma$-module]
\label{def:L2rank_of_fin_gen_prof_RGamma-module} Let $M$ be a
finitely generated $\IC\Gamma$-module $M$.  The \emph{$L^2$-rank of
$M$ is}
\begin{equation*}
 \rk_{\Gamma}^{(2)}(M)  :=
\bigl\{\dim_{\caln(x)}(S_xM \otimes_{\IC[x]} \caln(x))\mid \bar{x}
\in \iso(\Gamma)\bigr\}
\in U(\Gamma) \otimes_{\IZ} \IR = \bigoplus_{\iso(\Gamma)} \IR.
\end{equation*}
\end{definition}

The rank $\rk_{\Gamma}^{(2)}$ defines a homomorphism
\begin{eqnarray}
\rk_{\Gamma}^{(2)} \colon K_0(\IC\Gamma) \to U(\Gamma) \otimes_{\IZ} \IR,
\quad [P] \to \rk_{\Gamma}^{(2)}(P)&&
\label{rk(2)_K_o_to_U}
\end{eqnarray}
since for a finitely generated $\IC\Gamma$-module $M$ the value of $S_xM$ is non-trivial
only for finitely many elements $\overline{x} \in \iso(\Gamma)$ and
the $\IC\aut(x)$-module $S_xM$ is finitely generated for every $x \in \ob(\Gamma)$
(see Lemma~\ref{lem:basic_properties_of-E_x_and_S_x}).

If $\Gamma$ is directly finite, then the map $\rk_{\Gamma}^{(2)}$
obviously factorizes over $S \colon
K_0(\IC\Gamma) \to \Split K_0(\IC\Gamma)$.

\begin{example} \label{exa:L2_rank_for_modules_over_finite_groups}
If $H$ is a subgroup of $G$ of finite index $[G:H]$, and $i$ denotes the inclusion, then the diagram
$$\xymatrix@C=3pc{K_0(\IC G) \ar[r]^-{\rk^{(2)}_G} \ar[d]_{i^\ast} & \IR \ar[d]^{[G:H]\cdot} \\
K_0(\IC H) \ar[r]_-{\rk^{(2)}_H} & \IR}$$
commutes.
\end{example}

\begin{proof}
It follows from existence of a $\IC H$-isomorphism $\IC
G^n=\oplus_{G/H}\IC H^n$ that the restriction $i^\ast P$ of a
finitely generated projective $\IC G$-module $P$ is a finitely
generated projective $\IC H$-module. So the left vertical map in the
above diagram is well defined. It directly follows from the proof
of~L\"uck~\cite[Theorem~6.54 (6) on page~266]{Lueck(2002)} that
\[
    (i^\ast P)\otimes_{\IC H}\caln(H)\cong\res_{i_\ast}\bigl(P\otimes_{\IC G}\caln(G)\bigr).
\]
Now the assertion follows from
Lemma~\ref{lem:dim_caln(G)_and_induction_and_restriction}~\ref{lem:dim_caln(G)_and_induction_and_restriction:dim_and_restriction}.
\end{proof}

\begin{remark}[$L^2$-rank of a finitely generated
  $R\Gamma$-module] \label{rem:L2_rank_wrt_subring_of_complexes} In
  Definition~\ref{def:L2rank_of_fin_gen_prof_RGamma-module} we have defined the $L^2$-rank
  of a finitely generated $\IC\Gamma$-module. If $R$ is a subring of $\IC$, we
  may analogously define the $L^2$-rank of a finitely generated $R\Gamma$-module
  $M$. Namely, we view $\caln(x)$ as an $R\aut(x)$-$\caln(x)$-bimodule via the
  embedding of rings $R\aut(x) \to \IC \aut(x) \to \caln(\aut(x))$ and then take
  $\dim_{\caln(x)}(S_xM \otimes_{R[x]} \caln(x))$ as the components of the
  $L^2$-rank of $M$. We will primarily be interested in the case $R=\IC$, so we
  omit $\IC$ from the notation $\rk_{\Gamma}^{(2)}$. Occasionally we will also
  consider $R=\IQ$.
\end{remark}

\begin{theorem}[Relating the finiteness obstruction and the
  $L^2$-Euler characteristic]
  \label{the:comparing_o_and_chi(2)}
  Suppose that $\Gamma$ is a directly finite category of type (FP$_{\IC}$).
  Then
  $\Gamma$ is of type ($L^2$) and the image of the finiteness
  obstruction $o(\Gamma;\IC)$ (see
  Definition~\ref{def:finiteness_obstruction_of_a_category})
  under the homomorphism
  $$\rk_{\Gamma}^{(2)} \colon K_0(\IC\Gamma) \to U(\Gamma) \otimes_{\IZ} \IR
  = \bigoplus_{\overline{x} \in \iso(\Gamma)} \IR$$ defined
  in~\eqref{rk(2)_K_o_to_U} is $\chi_f^{(2)}(\Gamma)$.
\end{theorem}
\begin{proof}
Since $\Gamma$ is of type (FP$_{\IC}$), we can  find a finite projective
$\IC \Gamma$-resolution $P_*$ of $\underline{\IC}$. Hence $S_xP_*$ is non-trivial
only for finitely many objects $x$ in $\Gamma$ and
a finite projective $\IC[x]$-chain complex for all objects $x$ in $\Gamma$
by Lemma~\ref{lem:basic_properties_of-E_x_and_S_x}.
Hence $\Gamma$ is of type ($L^2$). Now apply
Lemma~\ref{lem:basic_properties_of_L2-Euler_characteristic}~%
\ref{lem:basic_properties_of_L2-Euler_characteristic:homology}.
\end{proof}

\begin{example}
Finite EI-categories are of type ($L^2)$ by
Theorem~\ref{the:comparing_o_and_chi(2)},
Lemma~\ref{lem:directly_finite_and_idempotents_and_EI}, and
Lemma~\ref{lem:finite_homological_dimension}~\ref{lem:finite_homological_dimension:finite_Gamma}.
\end{example}

\begin{lemma}
  \label{lem:rank_is_rankL2_for_fin_free_modules}
  Suppose that $\Gamma$ is directly finite. Then:

  \begin{enumerate}
  \item \label{lem:rank_is_rankL2_for_fin_free_modules:rank_is_L2-rank}
    If $F$ is a finitely generated free $\IC\Gamma$-module, the rank
    $\rk_{\IC\Gamma}(F)$ of
    Definition~\ref{def:rank_of_fin_gen_prof_RGamma-module} and the
    rank $\rk_{\Gamma}^{(2)}(F)$ of
    Definition~\ref{def:L2rank_of_fin_gen_prof_RGamma-module} agree;

  \item \label{lem:rank_is_rankL2_for_fin_free_modules:iota_splits}
    The composite
    $$U(\Gamma) \xrightarrow{\iota} K_0(\IC\Gamma)
    \xrightarrow{\rk_{\Gamma}^{(2)}} U(\Gamma) \otimes_{\IZ} \IR$$ of the
    homomorphisms defined in~\eqref{iota} and~\eqref{rk(2)_K_o_to_U} is
    the obvious inclusion $U(\Gamma) \to U(\Gamma) \otimes_{\IZ} \IR$;
\end{enumerate}
\end{lemma}
\begin{proof}\ref{lem:rank_is_rankL2_for_fin_free_modules:rank_is_L2-rank}
  This follows from Lemma~\ref{lem:rank_RGamma_and_free_modules} since
  for $\overline{y} = \overline{x}$ we have
  \begin{multline*}
    \rk_{\Gamma}^{(2)}(\IC\mor(?,x))_{\overline{y}}
    = \dim_{\caln(x)}\bigl(S_x\IC\mor(?,x) \otimes_{\IC[x]} \caln(x)) \\
    = \dim_{\caln(x)}(\caln(x)\bigr) = 1 = \rk_\IC(S_x\IC\mor(?,x)
    \otimes_{\IC[x]} \IC) = \rk_{\IC\Gamma}(\IC\mor(?,x))_{\overline{y}}.
  \end{multline*}
  and for $\overline{y} \not= \overline{x}$ we get
$$\rk_{\Gamma}^{(2)}(\IC\mor(?,x))_{\overline{y}} = 0
= \rk_{\IC \Gamma}(\IC\mor(?,x))_{\overline{y}}.$$%
\ref{lem:rank_is_rankL2_for_fin_free_modules:iota_splits} This follows
from
assertion~\ref{lem:rank_is_rankL2_for_fin_free_modules:rank_is_L2-rank}
and
Lemma~\ref{lem:rank_RGamma_and_free_modules}~\ref{lem:rank_RGamma_and_free_modules:iota_splits}.
\end{proof}
%

\begin{theorem}[Invariants agree for directly finite and type (FF$_\IZ$)]\label{the:coincidence}
Suppose $\Gamma$ is directly finite and of type
(FF$_\IZ$). Then the
functorial $L^2$-Euler characteristic of
Definition~\ref{def:functorial_L2-Euler_characteristic_of_a_category}
coincides with the functorial Euler characteristic of
Definition~\ref{def:functorial_Euler_characteristic_of_a_category}
for any associative, commutative ring $R$ with identity
\[
\chi_f^{(2)}(\Gamma)=\chi_f(\Gamma;R) \in U(\Gamma) \subseteq
U^{(1)}(\Gamma),
\]
and thus $\chi^{(2)}(\Gamma)=\chi(\Gamma;R)$ in
Definition~\ref{def:L2-Euler_characteristic_of_a_category} and
Definition~\ref{def:Euler_characteristic_of_a_category}.

If $R$ is additionally Noetherian, then
\begin{equation} \label{equ:coincidence_for_type_FFZ_and_directly_finite}
\chi(B\Gamma;R)=\chi(\Gamma;R)=\chi^{(2)}(\Gamma).
\end{equation}
Moreover, if $\Gamma$ is merely of type (FF$_\IC$) rather than
(FF$_\IZ$), then equation
\eqref{equ:coincidence_for_type_FFZ_and_directly_finite} holds for
any Noetherian ring $R$ containing $\IC$.
\end{theorem}
\begin{proof}
If $\Gamma$ is of type (FF$_\IZ$), it is also of type (FF$_R$), since any
(augmented) resolution of $\IZ$ is contractible as a complex of $\IZ$-modules,
thus stays exact after applying $\_\otimes_\IZ R$.
Using
Lemma~\ref{lem:basic_properties_of-E_x_and_S_x}~\ref{lem:basic_properties_of-E_x_and_S_x:directly_finite},
we can show
\[
\rk_R \left(S_x \left(F_n\otimes_{\IZ} R  \right) \otimes_{R[x]} R
\right) = \rk_\IZ \left(S_x F_n \otimes_{\IZ[x]} \IZ \right).
\]
Consequently, $\chi_f(\Gamma;R)=\chi_f(\Gamma;\IZ)$ and
$\chi(\Gamma;R)=\chi(\Gamma;\IZ)$ for any ring $R$.

By
Lemma~\ref{lem:rank_is_rankL2_for_fin_free_modules}~\ref{lem:rank_is_rankL2_for_fin_free_modules:rank_is_L2-rank},
the $\IC\Gamma$-rank $\rk_{\IC\Gamma}$ coincides with the $L^2$-rank
$\rk_\Gamma^{(2)}$ for finitely generated free $\IC\Gamma$-modules,
and we have
$\chi^{(2)}_f(\Gamma)=\rk_\Gamma^{(2)}o(\Gamma;\IC)=\rk_{\IC\Gamma}o(\Gamma;\IC)=\chi_f(\Gamma;\IC)=\chi_f(\Gamma;R)$
by Theorem~\ref{the:comparing_o_and_chi(2)} and the above (here we
use a finite free resolution in $o(\Gamma;\IC)$). Summing up, we
have $\chi^{(2)}(\Gamma)=\chi(\Gamma;R)$.
If $R$ is additionally Noetherian, then
Theorem~\ref{the:chi_f_determines_chi} implies
$\chi(\Gamma;R)=\chi(B\Gamma;R)$.

The statement after\eqref{equ:coincidence_for_type_FFZ_and_directly_finite} follows
by a similar argument as above.
\end{proof}

We may contrast the assumptions of (FF$_\IZ$) and direct finiteness
in Theorem~\ref{the:coincidence} with the relaxed assumptions of
(FP$_R$) and direct finiteness. If we only assume type (FP$_R$) and
direct finiteness, then $\chi(\Gamma;R)$ and $\chi(B\Gamma;R)$
coincide by Theorem~\ref{the:chi_f_determines_chi}, but these may be
different from $\chi^{(2)}(\Gamma)$. For example, if $G$ is a
nontrivial finite group, then it is of type (FP$_\IC$) but not of
type (FF$_\IC$), and we have $\chi(B\Gamma;\IC) = \chi(\Gamma;\IC) =
1$, but $\chi^{(2)}(\Gamma)=\frac{1}{\vert G \vert}$.

\begin{corollary}
Suppose $\Gamma$ is directly finite and of type (FF$_\IZ$). We have
$$\iota\bigl(\chi_f^{(2)}(\Gamma;\IC)\bigr) = o(\Gamma;\IC)$$
for the homomorphism $\iota$ defined in equation~\eqref{iota}.
\end{corollary}
\begin{proof}
This follows from Theorem~\ref{the:coincidence} and
Lemma~\ref{lem:rank_RGamma_and_free_modules:o_and_chi_for_type_(FF)}.
\end{proof}

\begin{remark} \label{rem:L2_of_Bcalc} Recall that $\chi(B\calc;\IQ)$ is
  the Euler characteristic of $B\calc$.  However, it is not true that
  $\chi^{(2)}(\calc)$ is related to the $L^2$-Euler characteristic
  $\chi^{(2)}\bigl(\widetilde{B\calc};\caln(\pi_1(B\calc))\bigr)$ in
  the sense of~L\"uck~\cite[Definition~6.20]{Lueck(2002)}. We will compute
  $\chi^{(2)}(\uor{D_{\infty}}) = 0$ in
  Subsection~\ref{subsec:infinite_dihedral_group}.  On the other hand
  $B\uor{D_{\infty}} = D_{\infty}\backslash\eub{D_{\infty}}$ is contractible and hence
  $\chi^{(2)}\bigl(\widetilde{B\calc};\caln(\pi_1(B\calc))\bigr) =
  \chi(B\uor{D_{\infty}}) = 1$.
\end{remark}

\subsection{Compatibility of Euler characteristics with coverings and isofibrations}
\label{subsec:Compatibility_with_Coverings} Our next task is to show
that the $L^2$-Euler characteristic is compatible with covering maps
and isofibrations between connected finite groupoids.
In the context of groupoids, the role of a covering neighborhood is played by
the star of an object.  If $\cale$ is a small groupoid and $e$ is an object of
$\cale$, we denote by $St(e)$ the \emph{star of $e$}, namely the set of all
morphisms in $\cale$ with domain $e$.
\begin{definition}[Covering of a groupoid]
  A functor $p\colon \cale \to \calb$ between connected small groupoids is a
  \emph{covering} if it is surjective on objects and restricts to a bijection
$$\xymatrix{St(e) \ar[r] & St(p(e))}$$
for each object $e$ of $\cale$. We say that a covering $p$ is \emph{$n$-sheeted}
if $|\ob (p^{-1}(b))|=n$ for all objects $b$ of $\calb$.
\end{definition}

Recall that a small groupoid $\cale$ is \emph{finite} if
$\iso(\cale)$ is finite and for any object $e \in \ob(\cale)$ the
set $\aut(e)$ is finite.

\begin{theorem}[Compatibility of the $L^2$-Euler characteristic with coverings of finite groupoids] \label{the:compatibility_with_coverings}
Let $\cale$ and $\calb$ be connected finited groupoids. If $p\colon
\cale \to \calb$ is an $n$-sheeted covering, then
\begin{equation} \label{equ:compatibility_with_coverings}
\chi^{(2)}(\cale)=n\chi^{(2)}(\calb).
\end{equation}
\end{theorem}

\begin{proof}
  We present two proofs, one counting morphisms and the other using the
  technology of the finiteness obstruction.

  To prove the theorem by counting morphisms, we first reduce to the case where
  the base groupoid has only one object. If $b \in \calb$ and $\cale_b$ denotes
  the groupoid $p^{-1}(\widehat{\aut(b)})$, then the diagram
$$\xymatrix{\cale_b \ar@{^{(}->}[r] \ar[d]_{p\vert_{\cale_b}} & \cale \ar[d]^p \\ \widehat{\aut(b)} \ar@{^{(}->}[r] & \calb}$$
commutes and the horizontal functors are equivalences of categories.
The groupoid $\cale_b$ is connected; for if $e,e' \in \cale_b$, then
$f\colon e \cong e'$ in $\cale$, and $p(f) \in \aut(b)$, so $f \in
\mor(\cale_b)$. Moreover, $St_{\cale_b}(e) \subseteq St_{\cale}(e)$
for all $e \in \cale_b$, $St_{\aut(b)}(b) \subseteq St_{\calb}(b)$,
and $p\vert_{\cale_b}$ is an $n$-sheeted covering.  By
Theorem~\ref{the:invariance_under_equivalence},
Lemma~\ref{lem:directly_finite_and_idempotents_and_EI},
Theorem~\ref{the:comparing_o_and_chi(2)}, and
Definition~\ref{def:L2-Euler_characteristic_of_a_category}, the
groupoids $\cale_b$ and $\cale$ have the same $L^2$-Euler
characteristic, as do $\widehat{\aut(b)}$ and $\calb$.
Alternatively, we know from
Example~\ref{exa:chi_f(calg)L2_for_finite_groupoids} directly that
$$\chi^{(2)}(\cale_b)=\frac{1}{|\aut(e)|}=\chi^{(2)}(\cale)$$
$$\chi^{(2)}(\widehat{\aut(b)})=\frac{1}{|\aut(b)|}=\chi^{(2)}(\calb).$$
Thus, if the theorem holds in the case where the base groupoid has only one
object, it holds in general.

Suppose now that $\calb$ has only one object $b$, so that
$\calb=\widehat{\aut(b)}$. Then $\cale$ has only $n$ objects, say
$e_1, \dots, e_n$. Since $\cale$ is a connected finite groupoid, all
of its hom-sets have the same number of elements. Let $e \in \cale$.
We have
\begin{align} \label{equ:counting_morphisms}
  \vert \aut(b) \vert = \vert St(e) \vert
  = \vert \bigcup_{i=1}^n \mor_\cale(e,e_i) \vert
  = \sum_{i=1}^n \vert \mor_\cale(e,e_i) \vert
  &= \sum_{i=1}^n \vert \aut(e) \vert\notag \\
  &= n\vert\aut(e)\vert.
\end{align}
In conclusion, $\chi^{(2)}(\cale)=n\chi^{(2)}(\calb)$.

We may also prove Theorem~\ref{the:compatibility_with_coverings} on the level of
finiteness obstructions as follows, without reduction to the case of one object
in the base groupoid.

The covering $p\colon \cale \to \calb$ is admissible in
the sense that $\res_p$ sends a finitely generated projective
$R\calb$-module to a finitely generated projective $R\cale$-module
as a consequence of~L\"uck~\cite[Proposition~10.16 on
page~187]{Lueck(1989)} as follows. A morphism $h \colon p(x) \to y$
in $\calb$ is said to be \emph{irreducible} if for any factorization
$h=f \circ p(g)$ the morphism $g$ in $\cale$ is an isomorphism.
Clearly, the set $\Irr(x,y)$ of irreducible morphisms $p(x) \to y$
in $\calb$ is $\mor_\calb(p(x),y)$, since $\cale$ is a groupoid.
Since $\cale$ is finite, for a given $y \in \calb$, the set
$\Irr(x,y)$ is nonempty for only finitely many $\overline{x} \in
\iso(\cale)$. Since $\calb$ is finite, for each $x \in \cale$ the
right $\aut_\cale(x)$-set $\Irr(x,y)$ has only finitely many orbits.
The right action of $\aut_\cale(x)$ on $\Irr(x,y)$ is free because
$\calb$ is a groupoid and $p$ is a covering: if $h \in
\mor_\calb(p(x),y)$ and $h \circ pm=h \circ pn$, then $pm=pn$ and
$m=n$. Every morphism $h$ in $\mor_\calb(p(x),y)$ is irreducible, so
clearly we have a factorization $f \circ p(g) = h$ with $f$
irreducible, namely $f=h$ and $g=\id_x$. Any two factorizations $f
\circ p(g) = h$ and $f' \circ p(g') = h$ with $f$ and $f'$
irreducible are related by the isomorphism $k:= g' \circ g^{-1}$.

We fix an $x \in \cale$ and let $H=\aut_\cale(x)$, $G=\aut_\calb(p(x))$.
The covering $p$ induces an inclusion of
$H$ into $G$. Consider the following diagram.
$$\xymatrix@C=3pc@R=3pc{K_0(\IC\calb) \ar[r]_S^{\cong} \ar@/^2pc/[rrr]^{\rk^{(2)}_\calb} \ar[d]_{p^\ast} & K_0(\IC G) \ar[r]_-{\rk^{(2)}_G} \ar[d]_{p^\ast} & \IR
\ar[d]^{[G:H]\cdot}  & U(\calb)\otimes \IR \ar[l]_-\cong
\\ K_0(\IC\cale) \ar[r]^S_{\cong} \ar@/_2pc/[rrr]_{\rk^{(2)}_\cale}
& K_0(\IC H) \ar[r]^-{\rk^{(2)}_H} & \IR  & U(\cale)\otimes \IR
\ar[l]^-\cong}$$ The left square commutes by
Theorem~\ref{the_splitting_of_K-theory_for_EI_categories}. The
second square commutes by
Example~\ref{exa:L2_rank_for_modules_over_finite_groups}. The top
and bottom diagrams commute by definition of $\rk^{(2)}$. Beginning
in the upper left-hand corner, we have $o(\calb;\IC) \in K_0(\IC
\calb)$. By Theorem~\ref{the:o(Gamma;R)_and_restriction}, we have
$p^*(o(\calb;\IC))=o(\cale;\IC)$. Two applications of
Theorem~\ref{the:comparing_o_and_chi(2)} combined with the
commutativity of the diagrams leads us to
$\chi^{(2)}(\cale)=[G:H]\cdot \chi^{(2)}(\calb)$. An argument
similar to the one in~\eqref{equ:counting_morphisms} shows that
$[G:H]$ is equal to the number of sheets $n$.
\end{proof}

\begin{example}
Let $\cale=\{0 \leftrightarrow 1\}$ and let $\calb$ be the category
with one object and one nontrivial arrow, which is its own inverse.
By Example~\ref{exa:chi_f(calg)L2_for_finite_groupoids}, the
$L^2$-Euler characteristics are $\chi^{(2)}(\cale)=1$ and
$\chi^{(2)}(\calb)=1/2$.  The unique covering $\cale\to \calb$ is
2-sheeted and we have
$$\chi^{(2)}(\cale)=2\chi^{(2)}(\calb).$$
\end{example}

\begin{corollary}
Any $n$-sheeted covering functor between connected finite groupoids
is equivalent to the inclusion of an index $n$ subgroup into a
finite group. More precisely, if $p\colon \cale \to \calb$ is an
$n$-sheeted covering between connected finite groupoids and $e \in
\cale$, then the diagram
$$\xymatrix{\widehat{\aut(e)} \ar@{^{(}->}[r] \ar[d]_{p\vert_{\widehat{\aut(e)}}}
& \cale \ar[d]^p \\ \widehat{\aut(p(e))} \ar@{^{(}->}[r] & \calb}$$
commutes, the horizontal functors are equivalences of categories,
the left vertical functor is mono, and $[\aut(p(e)):p(\aut(e))]=n$.
\end{corollary}

\begin{remark}
Examples of covering functors are obtained from coverings of
topological spaces: a covering of topological spaces induces a
covering functor between the associated fundamental groupoids.
\end{remark}

We next turn to compatibility of $\chi^{(2)}$ with isofibrations.

\begin{definition}[Isofibration]
A functor $p\colon \cale \to \calb$ is an \emph{isofibration} if for every isomorphism in $\calb$
of the form $g\colon b \cong p(e)$ there is an isomorphism $f$ in $\cale$ such that $p(f)=g$
\end{definition}

We remark that if $\cale$ and $\calb$ are groupoids, then isofibrations and
Grothendieck fibrations coincide (because isomorphisms in the domain category
are always cartesian arrows).

\begin{theorem}[Compatibility of the $L^2$-Euler characteristic with isofibrations of finite groupoids] \label{the:compatibility_with_isofibrations} Let $p\colon \cale
  \to \calb$ be an isofibration between connected finite groupoids. If $b \in
  \calb$ and $p^{-1}(b)$ is connected, then
\begin{equation} \label{equ:compatibility_with_isofibrations}
\chi^{(2)}(\cale)=\chi^{(2)}(p^{-1}(b))\cdot\chi^{(2)}(\calb).
\end{equation}
\end{theorem}
\begin{proof}
  As in the proof of Theorem~\ref{the:compatibility_with_coverings}, we reduce
  to the case where the base groupoid has only one object.  If $b \in \calb$ and
  $\cale_b$ denotes the groupoid $p^{-1}(\widehat{\aut(b)})$, then the diagram
$$\xymatrix{\cale_b \ar@{^{(}->}[r] \ar[d]_{p\vert_{\cale_b}}
& \cale \ar[d]^p \\ \widehat{\aut(b)} \ar@{^{(}->}[r] & \calb}$$
commutes, the horizontal functors are equivalences of categories, and $\cale_b$
is connected. The fiber groupoid $p\vert_{\cale_b}^{-1}(b)$ is the same as the
fiber groupoid $p^{-1}(b)$, so $p\vert_{\cale_b}^{-1}(b)$ is also connected.
Since $\chi^{(2)}(\cale)=\chi^{(2)}(\cale_b)$ and
$\chi^{(2)}(\calb)=\chi^{(2)}(\widehat{\aut(b)})$, we have
\eqref{equ:compatibility_with_isofibrations} if
$\chi^{(2)}(\cale_b)=\chi^{(2)}(p\vert_{\cale_b}^{-1}(b))\cdot\chi^{(2)}(\widehat{\aut(b)})$. We
have reduced to the case where the base groupoid has only one object.

Suppose now that $\calb$ has only one object $b$, so that
$\calb=\widehat{\aut(b)}$.  For $e \in p^{-1}(b)$, we write simply
$p_e$ for the group homomorphism $\aut(e) \to \aut(b)$. Then $p_e$
is surjective.  If $g$ is an automorphism of $b$, there exists an
$f\colon e' \to e$ with $p(f)=g$. The connectivity of the fiber
$p^{-1}(b)$ then gives us an isomorphism $h\colon e \to e'$, and an
automorphism $f \circ h$ of $e$ such that $p_e(f \circ h)=g$.

Finally,
\[\chi^{(2)}(\cale)=\frac{1}{\vert\aut(e)\vert}
=\frac{1}{\vert \ker p_e \vert \cdot \vert \aut(b) \vert}=\chi^{(2)}(p^{-1}(b)) \cdot \chi^{(2)}(\calb).\qedhere\]
\end{proof}


\typeout{--------------------  Section 6: Moebius inversion ---------------------------}
\section{M\"obius inversion}
\label{sec:Moebius-inversion}

We extend the $K$-theoretic M\"obius inversion
of~L\"uck~\cite[Chapter~16]{Lueck(1989)} from finite to quasi-finite
EI-categories and apply it to the finiteness obstruction and the
Euler characteristic of a category. Throughout this section let
$\Gamma$ be an EI-category (see Definition~\ref{def:EI-category}).
We have already introduced the splitting $(S,E)$ of $K_0(R\Gamma)$
in Theorem~\ref{the_splitting_of_K-theory_for_EI_categories}.
Provided that $\Gamma$ is a quasi-finite EI-category, we obtain a
second splitting $(\Res,I)$ in Theorem
\ref{the:second_splitting_of_K_0(RGamma)}. The $K$-theoretic
M\"obius inversion $(\mu, \omega)$ will compare these two splittings
in
Theorem~\ref{the:Two_splittings_and_the_K-theoretic_Moebius_inversion}.
As a consequence, in
Theorem~\ref{the:The_finiteness_obstruction_and_Moebius_inversion_finite_EI-categories}
we obtain explicit formulas for the various Euler characteristics of
finite EI-categories. Important special cases of our $K$-theoretic
M\"obius inversion include Philip Hall's M\"obius inversion formula
for finite posets and Leinster's M\"obius inversion formula for
finite skeletal categories with only trivial endomorphisms. See
Examples~\ref{Finite_partially_ordered_set}~and~\ref{exa:Moebius_Inversion_for_fin_sk_cat_with_triv_endos}.

After treating the second splitting $(\Res,I)$ and the $K$-theoretic
M\"obius inversion $(\mu, \omega)$ in
Subsections~\ref{subsec:A_second_Splitting}~and~\ref{subsec:The-K-theoretic-Moebius-Inversion},
we turn to the relationship between the $K$-theoretic M\"obius
inversion $(\mu, \omega)$ and the $L^2$-rank in
Subsection~\ref{subsec:The-K-theoretic-Moebius-Inversion_and_the_L2-rank}.
There we construct a pair of homomorphisms $\overline{\mu}^{(2)}
\colon U(\Gamma)\otimes_\IZ \IQ \rightleftarrows
U(\Gamma)\otimes_\IZ \IQ \colon \overline{\omega}^{(2)}$ that are
inverse to one another if $\Gamma$ is a quasi-finite, free
EI-category, and commute appropriately with $(\mu,\omega)$ and
$\rk_\Gamma^{(2)}$ as in
Theorem~\ref{the:K_theoretic-Moebius_inversion_and_L2-rank}. All of
these homomorphisms and splittings are illustrated for
$G$-$H$-bisets (viewed as two-object EI-categories) in
Subsection~\ref{subsec:The_example_of_a_biset}.

In general, the finiteness obstruction and Euler characteristics of $\Gamma^{\op}$ are different from those of $\Gamma$,
as we see in Subsection~\ref{subsec:The_passage_to_the_opposite_category} with a biset example.
However, in the case of a finite EI-category $\Gamma$, the groups $K_0(\IQ\Gamma)$ and $K_0(\IQ\Gamma^{\op})$ are isomorphic, and
we say more about the respective splittings in Subsection~\ref{subsec:The_passage_to_the_opposite_category_for_finite_EI-categories}.

In Section \ref{sec:Moebius-inversion} we also introduce the proper orbit category $\uor{G}$, an important quasi-finite, free EI-category to which we shall return in Section \ref{sec:The_proper_orbit_category}.

\subsection{A second splitting}
\label{subsec:A_second_Splitting}

Given an object $x$ in a (small) category $\Gamma$, define the \emph{restriction functor at $x$}
\begin{eqnarray}
& \Res_x \colon \MOD\text{-}R\Gamma \to \MOD\text{-}R[x]&
\label{Res_x}
\end{eqnarray}
by evaluating an $R\Gamma$-module $N$ at the object $x$.  This functor
is exact but does not respect finitely generated projective in general.
Given an EI-category $\Gamma$, the \emph{inclusion functor at $x$}
\begin{eqnarray}
& I_x \colon  \MOD\text{-}R[x] \to \MOD\text{-}R\Gamma&
\label{I_x}
\end{eqnarray}
sends a right $R[x]$-module $M$ to the $R\Gamma$-module given by
$$I_xM(y) := \begin{cases}
M \otimes_{R[x]} R\mor(y,x) & \text{if } \overline{y} =
\overline{x};
\\
0 & \text{if } \overline{y} \not= \overline{x}.
\end{cases}
$$
Notice that we need the EI-condition to ensure that this definition makes sense.
This functor is compatible with direct sums, but does not respect finitely generated
projective in general.

\begin{lemma} \label{lem:adjunctions} Let $\Gamma$ be an EI-category.
Then we obtain for every $x \in \ob(\Gamma)$  adjoint
pairs of functors $(E_x,\Res_x)$ and $(S_x,I_x)$,
where $E_x$, $\Res_x$, $S_x$ and $I_x$ are the
functors defined in~\eqref{E_x},~\eqref{Res_x},~\eqref{S_x} and~\eqref{I_x}.
\end{lemma}
\begin{proof}
See~L\"uck~\cite[Lemma~9.31 on page~171]{Lueck(1989)}.
\end{proof}

The EI-property ensures that we obtain a
well-defined partial ordering on $\iso(\Gamma)$ by
\begin{eqnarray}
  \overline{x} \le \overline{y}
  & \Longleftrightarrow &
  \mor(x,y) \not= \emptyset.
  \label{partial_ordering_on_iso(Gamma)}
\end{eqnarray}

\begin{definition}[Length of an element]
\label{def:length_of_an_element}
  Given an element $x \in \iso(\Gamma)$, define its \emph{length}
  $$l(\overline{x}) \in \{0,1,2, \ldots \} \amalg \{\infty\}$$
  to be the supremum over
  the natural numbers $n$, for which there exists elements
  $\overline{x_n}$, $\overline{x_{n-1}}$, \ldots ,$\overline{x_0}$ in
  $\iso(\Gamma)$ with $\overline{x_n} < \overline{x_{n-1}} < \ldots <
  \overline{x_0}$ and $\overline{x_0} = \overline{x}$.
  \end{definition}

The length of $\overline{x}$ is zero if and only if
every morphism with $x$ as target is an isomorphism.

\begin{definition}[Finite, quasi-finite, and free categories]
  \label{def:finite_quasi-finite_and_free_category}
  Let $\Gamma$ be a (small) category.

  We call $\Gamma$ \emph{quasi-finite} if for every $\overline{x} \in \iso(\Gamma)$
  the set $\{\overline{y} \in \iso(\Gamma) \mid \overline{y} \le \overline{x}\}$ is finite,
  and for every two objects $x,y \in \ob(\Gamma)$ the right $\aut(x)$-set
  $\mor(x,y)$ is proper and cofinite, i.e., every isotropy group under
  the right $\aut(x)$-action is finite and the quotient $\mor(x,y)/\aut(x)$ is finite.

  We call $\Gamma$ \emph{finite} if $\iso(\Gamma)$ is finite and $\mor(x,y)$ is
  finite for every two objects $x,y \in \ob(\Gamma)$. A small
  category is finite if and only if it is equivalent to a category
  with finitely many objects and finitely many morphisms.

  We call $\Gamma$ \emph{free} if the left $\aut(y)$-action on
  $\mor(x,y)$ is free for every two objects $x,y \in \ob(\Gamma)$.
\end{definition}

One of our main examples for $\Gamma$ will be the orbit category.

\begin{definition}[Orbit category and proper orbit category]\label{def:orbit_category}
  Let $G$ be a group. The \emph{orbit category} $\Or(G)$ has as objects homogeneous spaces
  $G/H$ and as morphisms $G$-equivariant maps. The \emph{proper orbit category}
  $$\uor{G} = \OrGF{G}{\calfin},$$
  sometimes also called
  the orbit category associated to the family $\calfin$ of finite subgroups,
  is defined to be the full subcategory of $\Or(G)$
  consisting of objects $G/H$ with finite $H$.
\end{definition}

\begin{lemma}\label{lem:orbit_category_morphisms}
Let $H$ and $K$ be subgroups of a group $G$. If $g \in G$ and
$g^{-1}Hg \subseteq K$, then we get a well-defined $G$-equivariant
map $$\xymatrix{R_g \colon G/H \ar[r] & G/K}$$
$$\xymatrix{g'H \ar@{|->}[r] & g'gK}.$$
Every $G$-equivariant map $G/H \to G/K$ is of the form
$R_g$. We have $R_g=R_{g'}$ if and only if $g^{-1}g' \in K$ holds.
In particular, we have a bijection
\begin{eqnarray} \label{equ:orbit_category_morphisms}
& \xymatrix{\mor(G/H, G/K) \ar[r] & \{gK \mid g^{-1}Hg \subseteq
K\}}
\\ & \xymatrix{f \ar@{|->}[r] & f(1H)}.\nonumber
\end{eqnarray}
We also have $R_{g_2} \circ R_{g_1}=R_{g_1g_2}$.
\end{lemma}
\begin{proof}
See tom~Dieck~\cite[I.1.14]{Dieck(1987)} and L\"uck~\cite[Lemma~1.31
on page~22]{Lueck(1989)}.
\end{proof}

\begin{lemma} \label{lem:orbit_category_free_EI}
The orbit category $\Or(G)$ is a free EI-category.
\end{lemma}
\begin{proof}
A direct consequence of Lemma~\ref{lem:orbit_category_morphisms} is
that the monoid $\map(G/H,G/H)$ is isomorphic to the Weyl group
$N_GH/H$, so every endomorphism of $\Or(G)$ is an automorphism.

If $G/H$ and $G/K$ are two objects in $\Or(G)$, and  $f \colon G/H
\to G/K$ and $a \colon G/K \to G/K$ are $G$-equivariant maps, then
$a \circ f = f$ implies $a = \id$ since $f$ is surjective. Hence
$\Or(G)$ is free.
\end{proof}

\begin{lemma} \label{lem:orbit_category_quasi-finite}
The proper orbit category $\uor{G}$ is a quasi-finite and free
EI-category.
\end{lemma}
\begin{proof}
The proper orbit category $\uor{G}$ is a full subcategory of the
orbit category $\Or(G)$, which is a free EI-category, so $\uor{G}$
is also a free EI-category.

For the quasi-finiteness, we first observe from the bijection
\eqref{equ:orbit_category_morphisms} that $$\mor(G/H,G/K)\neq
\emptyset$$ if and only if $H$ is $G$-conjugate to a subgroup of
$K$. If $H$ and $H'$ are $G$-conjugate, then $G/H$ and $G/H'$ are
isomorphic objects of $\Or(G)$. Thus for a fixed $G/K$, the number
of isomorphism classes $\overline{G/H}$ with $\mor(G/H,G/K)\neq
\emptyset$ is at most the number of $G$-conjugacy classes of
subgroups of $K$. Whenever $K$ is a finite group, this number is
finite. Thus, $\{\overline{G/H} \in \iso(\uor{G})   \mid
\overline{G/H} \le \overline{G/K}\}$ is finite.

Continuing the notation of Lemma~\ref{lem:orbit_category_morphisms},
consider a morphism $R_{g_2} \colon G/H \to G/K$ in $\uor{G}$.
Suppose $R_{g_1} \in \aut(G/H)$ fixes $R_{g_2}$. Then
$R_{g_1g_2}=R_{g_2}$ and $g_2^{-1}g_1g_2 \in K$, so
that $g_1 \in g_2 K g_2^{-1}$. But $g_2 K g_2^{-1}$ is finite, so
there are only finitely many possibilities for $g_1$. Thus every
isotropy group for the right $\aut(G/H)$-action on $\mor(G/H,G/K)$
is finite.

For objects $G/H$ and $G/K$ in $\uor{G}$, the quotient
$\mor(G/H,G/K)/\aut(G/H)$ is in bijective correspondence with
\begin{equation} \label{equ:quotient_space_of_proper_orbit_morphisms}
\{g_2K \mid g_2^{-1}Hg_2 \subseteq K\}/\sim
\end{equation}
by Lemma~\ref{lem:orbit_category_morphisms}, where $g_2K \sim
g_1g_2K$ if $g_1 \in G$ and $g_1^{-1}Hg_1 \subseteq H$. Since $H$ is
finite, $g_1^{-1}Hg_1 \subseteq H$ implies $g_1^{-1}Hg_1 = H$. But
\eqref{equ:quotient_space_of_proper_orbit_morphisms} is in bijective
correspondence with $G$-conjugates of $H$ contained in $K$, of which
there are only finitely many because $K$ is finite. Thus the
quotient $\mor(G/H,G/K)/\aut(G/H)$ is finite.
\end{proof}

%
%
%

\begin{lemma} \label{lem:quasi-finite_(FP)_and_cardinality_of_iso}\
  \begin{enumerate}
  \item \label{lem:quasi-finite_(FP)_and_cardinality_of_iso:general}
    Suppose for the EI-category $\Gamma$ that for every $\overline{x}
    \in \iso(\Gamma)$ the set $\{\overline{y} \in \iso(\Gamma) \mid
    \overline{y} \le \overline{x}\}$ is finite. Let $M$ be a finitely generated
    $R\Gamma$-module $M$. Then
    $$\bigl\{\overline{x} \in \iso(\Gamma)\mid M(x) \not= 0\bigr\}$$
    is finite;

  \item \label{lem:quasi-finite_(FP)_and_cardinality_of_iso:special}
    If $\Gamma$ is a quasi-finite EI-category of type (FP$_R$), then
    $\iso(\Gamma)$ is finite.
  \end{enumerate}
\end{lemma}
\begin{proof}\ref{lem:quasi-finite_(FP)_and_cardinality_of_iso:general}
  Choose a finite subset $I \subseteq \iso(\Gamma)$ and natural
  numbers $n_i \ge 1$ for each $i \in I$ such that there exists an
  epimorphism of $R\Gamma$-modules
  $$\bigoplus_{i \in I} R\mor(?,x_i)^{n_i} \to M.$$
  Then for every $\overline{y} \in \iso(\Gamma)$ with $M(y) \not= 0$ there is
  $i \in I$ with $\overline{y} \le \overline{x_i}$.
  Since $I$ is finite, $\bigl\{\overline{x} \in \iso(\Gamma)\mid M(x) \not= 0\bigr\}$ is
  finite.
  \\[1mm]\ref{lem:quasi-finite_(FP)_and_cardinality_of_iso:special}
  This follows from
  assertion~\ref{lem:quasi-finite_(FP)_and_cardinality_of_iso:general} applied
  to the constant module $\underline{R}$.
\end{proof}

  \begin{definition}[Length of a module]
  \label{def:length_of_a_module}
  The \emph{length}
  $l(M) \in \{-1,0,1,2 \ldots\} \amalg \{\infty\}$ of an
  $R\Gamma$-module $M$  is defined to be $-1$
  if $M$ is zero and otherwise to be the supremum of the length
  of elements $\overline{x} \in \iso(\Gamma)$ with $M(x) \not= 0$.
\end{definition}

  If $\Gamma$ is quasi-finite and hence
  $\{\overline{y} \in \iso(\Gamma) \mid \overline{y} \le \overline{x}\}$
  is finite for every $\overline{x} \in \iso(\Gamma)$, the length of
  $R\mor(?,x)$ is finite for every object $x \in \ob(\Gamma)$ and hence every
  finitely generated $R\Gamma$-module has finite length.

\begin{lemma}\label{lem:finite_homological_dimension}
  Suppose that $\Gamma$ is a quasi-finite EI-category. Suppose for any
  morphism $f \colon x \to y$ in $\Gamma$ that the order of the finite group
  $\{g \in \aut(x) \mid f \circ g = f\}$ is invertible in $R$.

  \begin{enumerate}

  \item \label{lem:finite_homological_dimension:Res_x}
    Consider $x \in \ob(\Gamma)$. Let $M$ be an
    $R\Gamma$-module which is finitely generated projective or which
    possesses a finite projective
    $R\Gamma$-resolution respectively.    Then the
    $R\aut(x)$-module $\Res_xM = M(x)$ is finitely generated projective or
    has a finite projective $R[x]$-resolution respectively;

  \item \label{lem:finite_homological_dimension:criterion} Let $M$ be
    an $R\Gamma$-module such that the set
    \[\bigl\{\overline{x} \in \iso(\Gamma)\mid M(x) \not= 0\bigr\}\]
    is finite. If $\Res_x M$ possesses a finite
    projective $R[x]$-resolution for all $x \in \ob(\Gamma)$, then $M$ possesses a finite
    projective $R\Gamma$-resolution, ;

  \item \label{lem:finite_homological_dimension:I_x} Let $x \in
    \ob(\Gamma)$ and let $N$ be an $R[x]$-module which possesses a
    finite projective $R[x]$-resolution.  Then the $R\Gamma$-module
    $I_xN$ defined in~\eqref{I_x} possesses a finite projective $R\Gamma$-resolution;

   \item \label{lem:finite_homological_dimension:criterion_for_FP}
   $\Gamma$ is of type (FP$_R$) if and only if
   $\iso(\Gamma)$ is finite and for every object $x \in \ob(\Gamma)$
   the trivial $R[x]$-module $R$ is of type (FP$_R$) respectively;

   \item \label{lem:finite_homological_dimension:finite_Gamma}
   Let $\Gamma$ be a finite EI-category. Assume that for every object $x$
   the order of the finite group $\aut(x)$ is invertible in $R$.
   Then an $R\Gamma$-module $M$ possesses a finite projective
   resolution if for every object $x$ the $R$-module $M(x)$ possesses a finite
   projective $R$-resolution. In particular $\Gamma$ is of type (FP$_R$).

  \end{enumerate}
\end{lemma}
\begin{proof}\ref{lem:finite_homological_dimension:Res_x}
  Since $\Res_x$ is exact, it suffices to show that $\Res_xR\mor(?,y)
  = R\mor(x,y)$ is a finitely generated projective $R[x]$-module for every $y \in
  \ob(\Gamma)$. This follows from the assumptions that the right
  $\aut(x)$-set $\mor(x,y)$ is a finite union of homogeneous
  $\aut(x)$-spaces of the form $H \backslash \aut(x)$ for finite
  $H \subseteq \aut(x)$ such that $|H| \cdot 1_R$ is a unit in $R$.
  \\[1mm]\ref{lem:finite_homological_dimension:criterion} Since $\Gamma$ is quasi-finite and $M$ has finite support,
  the $R\Gamma$-module $M$ has finite length. We do
  induction over the length of the $R\Gamma$-module
  $M$.  The induction beginning $l = -1$ is trivial, the induction step
  from $l-1$ to $l \ge 0$ done as follows.

  If $0 \to M_1 \to M_1 \to M_3 \to 0$ is an exact sequence of
  $R\Gamma$-modules such that two of the $R\Gamma$-modules
  $M_1$, $M_2$, and $M_3$ possess finite projective
  $R\Gamma$-resolutions, then all three possess finite projective
  $R\Gamma$-resolutions (see~L\"uck~\cite[Lemma~11.6 on
  page~216]{Lueck(1989)}). Thus, using the Filtration Theorem
  (see~L\"uck~\cite[Theorem~16.8 on page~326]{Lueck(1989)}) and
  the induction hypothesis, it suffices to
  show for any object $x$ of length $l$ and any $R[x]$-module $N$
  which admits a finite projective $R[x]$-resolution that $I_xN$ has a
  finite projective $R\Gamma$-resolution.  Since $I_x$ is exact, it is
  enough to consider the case $N = R[x]$. Consider the epimorphism $f
  \colon R\mor(?,x) \to I_x(R[x])$ sending $\id_x$ to $1_{R[x]} \otimes
  \id_x \in R[x] \otimes_{R[x]} R\mor(x,x) = I_x(R[x])$.  Its kernel
  $\ker(f)$ is an $R\Gamma$-module of length $\le l-1$ and satisfies
  $\Res_y(\ker(f)) = R\mor(y,x) = \Res_y R\mor(?,x)$ for $\overline{y}
  < \overline{x}$ and $\Res_y(\ker(f)) = 0$ otherwise.
  Assertion~\ref{lem:finite_homological_dimension:Res_x}  implies that
  $\Res_y(\ker(f))$ possesses a finite projective $R[y]$-resolution
  for all objects $y \in \ob(\Gamma)$. Hence $\ker(f)$ possesses a finite
  projective $R\Gamma$-resolution  by induction hypothesis.  This
  implies that $I_xR[x]$ possesses a finite projective
  $R\Gamma$-resolution.  This finishes the proof of the induction
  step.
  \\[1mm]\ref{lem:finite_homological_dimension:I_x} This
  follows directly from
  assertion~\ref{lem:finite_homological_dimension:criterion}.
  \\[1mm]~\ref{lem:finite_homological_dimension:criterion_for_FP}
  This follows directly from
  Lemma~\ref{lem:quasi-finite_(FP)_and_cardinality_of_iso}~%
\ref{lem:quasi-finite_(FP)_and_cardinality_of_iso:special} and
  assertions~\ref{lem:finite_homological_dimension:Res_x}
  and~\ref{lem:finite_homological_dimension:criterion}.
  \\[1mm]\ref{lem:finite_homological_dimension:finite_Gamma}
  Since $|\aut(x)|$ is invertible in $R$ and finite, an $R[x]$-module
  possesses a finite projective
  $R[x]$-resolution if and only if it possesses a finite projective   $R$-resolution.
  Now apply assertion~\ref{lem:finite_homological_dimension:criterion}.
\end{proof}

Our main example for $R$ will of course be $\IQ$.

\begin{theorem}[A second splitting of $K_0(R\Gamma)$]
  \label{the:second_splitting_of_K_0(RGamma)}
  Suppose that $\Gamma$ is a quasi-finite EI-category.
  Suppose for any morphism $f \colon x \to y$ in $\Gamma$ that the order of the finite
  group $\{g \in \aut(x) \mid f \circ g = f\}$ is invertible in $R$.

  Then we obtain isomorphisms $\Res$ and $I$ which are inverse to one another.
  \begin{eqnarray*}
    \Res \colon K_0(R\Gamma)  \to  \Split K_0(R\Gamma),
    & & [P] \mapsto  \{[\Res_xP] \mid \overline{x} \in \iso(\Gamma)\}
    \\
    I \colon \Split K_0(R\Gamma)  \to  K_0(R\Gamma),
    & & \{[Q_x] \mid \overline{x} \in \iso(\Gamma)\}  \mapsto
    \sum_{\overline{x} \in \iso(\Gamma)} [I_xQ_x]
  \end{eqnarray*}
\end{theorem}
\begin{proof}
  Consider a finitely generated projective $R\Gamma$-module $P$.  Then
  for any object $x \in \ob(\Gamma)$ the $R[x]$-module $\Res_xP$
  possesses a finite projective $R[x]$-resolution
  (see Lemma~\ref{lem:finite_homological_dimension}~\ref{lem:finite_homological_dimension:Res_x})
  and hence defines an
  element in $K_0(R[x])$, namely its finiteness obstruction in the
  sense of Definition~\ref{def:finiteness_obstruction_of_a_module}.
  Since $\Gamma$ is by assumption quasi-finite and hence
  $\{\overline{y} \in \iso(\Gamma) \mid \overline{y} \le
  \overline{x}\}$ is finite for every object $x \in
  \ob(\Gamma)$, there are only finitely many elements $\overline{x} \in
  \iso(\Gamma)$ with $\Res_x P \not= 0$ by
  Lemma~\ref{lem:quasi-finite_(FP)_and_cardinality_of_iso}~%
\ref{lem:quasi-finite_(FP)_and_cardinality_of_iso:general}. Hence we obtain a
  well-defined element
  $$\Res([P]) := \{[\Res_x P] \mid \overline{x} \in \iso(\Gamma)\}
  \; \in \; \bigoplus_{\overline{x} \in \iso(\Gamma)} K_0(R[x])
  = \Split K_0(R\Gamma).$$
  Thus we obtain a homomorphism
  $$\Res \colon K_0(R\Gamma) \to \Split K_0(R\Gamma).$$
  Define
  $$I \colon \Split K_0(R\Gamma) \to K_0(R\Gamma)$$
  analogously using
  Lemma~\ref{lem:finite_homological_dimension}~\ref{lem:finite_homological_dimension:I_x}.

  One obtains $\Res \circ I = \id$ from the fact
  that the functor $\Res_y \circ I_x \colon \MOD\text{-}R[x] \to
  \MOD\text{-}R[y]$ is naturally isomorphic to the identity functor
  if $x = y$ and is trivial if $\overline{x} \not= \overline{y}$.  It
  remains to show that $I$ is surjective. This is done by induction
  over the length, which is finite by
  Lemma~\ref{lem:quasi-finite_(FP)_and_cardinality_of_iso}~\ref{lem:quasi-finite_(FP)_and_cardinality_of_iso:general},
  of a finitely generated projective $R\Gamma$-module representing
  a class in $K_0(R\Gamma)$  using
  Lemma~\ref{lem:finite_homological_dimension} and the Filtration Theorem
  (see~L\"uck~\cite[Theorem~16.8 on page~326]{Lueck(1989)}).
\end{proof}

\subsection{The $K$-theoretic M\"obius inversion}
\label{subsec:The-K-theoretic-Moebius-Inversion}

\begin{convention}
Suppose for the remainder of this subsection
that $\Gamma$ is a quasi-finite EI-category and that for
every morphism $f \colon x \to y$ in $\Gamma$ the order of the finite
group $\{g \in \aut(x) \mid f \circ g = f\}$ is invertible in $R$.
\end{convention}
We obtain a well-defined homomorphism
$$\omega_{x,y} \colon K_0(R[x]) \to K_0(R[y]),
\quad [P] \mapsto [P \otimes_{R[x]} R\mor(y,x)]$$
since the right $R[y]$-module $R\mor(y,x) = \Res_yR\mor(?,x)$ is
finitely generated projective by
Lemma~\ref{lem:finite_homological_dimension}~\ref{lem:finite_homological_dimension:Res_x}.
Define
\begin{eqnarray}
& \omega \colon \Split K_0(R\Gamma) \to \Split K_0(R\Gamma)&
\label{omega_splitK_0_to_splitK_0}
\end{eqnarray}
by the matrix of homomorphisms
$$\left(\omega_{x,y}\right)_{\overline{x},\overline{y} \in \iso(\Gamma)}
\colon \bigoplus_{\overline{x} \in \iso(\Gamma)} K_0(R[x]) \to
\bigoplus_{\overline{y} \in \iso(\Gamma)}  K_0(R[y]).$$ This
definition makes sense since for a given $\overline{x} \in
\iso(\Gamma)$ there are only finitely many $\overline{y} \in
\iso(\Gamma)$ with $\omega_{x,y} \not= 0$.

\begin{example}
If $R=\IQ$ and $\Gamma$ is a finite skeletal category with trivial
automorphism groups, then $K_0(\IQ[x])=\IZ$ and
$\omega_{x,y}=|\mor_\Gamma(y,x)|$ for all $x,y \in \ob(\Gamma)$. In
this case of $R$ and $\Gamma$, the matrix for $\omega$ is the
transpose of the zeta function considered by Leinster in Section 1
of \cite{Leinster(2008)}. See also
Example~\ref{exa:Moebius_Inversion_for_fin_sk_cat_with_triv_endos}.
\end{example}

\begin{definition}[$l$-chain in $\iso(\Gamma)$]\label{def:chains}
  Let $\Gamma$ be an EI-category.  Given a natural number $l \ge 1$, an
  \emph{$l$-chain in $\iso(\Gamma)$} is a sequence $c = \overline{x_0}
  < \overline{x_1} < \cdots < \overline{x_l}$.  Denote by
  $\ch_l(\Gamma)$ the set of $l$-chains in $\Gamma$.

  Given two objects $x$ and $y$, let
  $\ch_l(y,x)$ be the set of $l$-chains
  $c = \overline{x_0} < \overline{x_1} < \cdots < \overline{x_l}$ with
  $\overline{x_0} = \overline{y}$ and $\overline{x_l} = \overline{x}$.
  Define for an $l$-chain
  $c = \overline{x_0} < \overline{x_1} < \cdots  < \overline{x_l}$
  in $\ch_l(y,x)$ the $\aut(x)$-$\aut(y)$-biset
  $$S(c) = \mor(x_{l-1},x) \times_{\aut(x_{l-1})}  \mor(x_{l-2},x_{l-1}) \times_{\aut(x_{l-2})}
  \cdots \times_{\aut(x_{1})} \mor(y,x_1)$$ for some choice of
  representatives $x_i \in \overline{x_i}$ for $0 < i < l-1$.
  (If $l = 1$ then $S(c)$ is to be understood as the
   $\aut(x)$-$\aut(y)$-biset $\mor(y,x)$.)

  Define $\ch_0(\Gamma)$ to be $\iso(\Gamma)$. Define
  $\ch_0(y,x)$ to be empty if
  $\overline{x} \not= \overline{y}$  and to be $\overline{y}$ if
  $\overline{x} = \overline{y}$. If $\overline{x} = \overline{y}$,
  put $S(c) = \mor(x,x)$ for $c \in \ch_0(y,x)$.
\end{definition}

Notice that the $\aut(x)$-$\aut(y)$-biset $S(c)$ is unique up to
isomorphism of $\aut(x)$-$\aut(y)$-bisets. Since $\Gamma$ is quasi-finite and hence for
every two objects $x,y \in \ob(\Gamma)$ the right $\aut(y)$-set
$\mor(y,x)$ is proper and cofinite, each set $S(c)$ is a proper cofinite
right $\aut(y)$-set, and the $R[y]$-module $RS(c)$ is finitely generated projective.
Hence we obtain a well-defined homomorphism for $c \in \ch_l(y,x)$
$$\mu_{x,y}(c) \colon K_0(R[x]) \to K_0(R[y]), \quad [P] \mapsto [P \otimes_{R[x]} RS(c)].$$
Define a homomorphism
\begin{eqnarray}
& \mu \colon \Split K_0(R\Gamma) \to \Split K_0(R\Gamma)&
\label{mu_splitK_0_to_splitK_0}
\end{eqnarray}
by the matrix of homomorphisms
$$\left(\sum_{l \ge 0} (-1)^l \cdot \sum_{c \in \ch_l(y,x)} \mu_{x,y}(c)
\right)_{\overline{x},\overline{y} \in \iso(\Gamma)} \colon
\bigoplus_{\overline{x} \in \iso(\Gamma)} K_0(R[x]) \to
\bigoplus_{\overline{y} \in \iso(\Gamma)}  K_0(R[y]).$$ This
definition  makes sense since for a given $\overline{x} \in
\iso(\Gamma)$ there are only finitely many $\overline{y} \in
\iso(\Gamma)$ with $\mu_{x,y} \not= 0$.

\begin{theorem}[Two splittings and the $K$-theoretic M\"obius inversion]
\label{the:Two_splittings_and_the_K-theoretic_Moebius_inversion}
Suppose that $\Gamma$ is a quasi-finite EI-category.  Suppose for any
morphism $f \colon x \to y$ in $\Gamma$ that the order of the finite
group $\{g \in \aut(x) \mid f \circ g = f\}$ is invertible in $R$.

\begin{enumerate}
\item \label{the:Two_splittings_and_the_K-theoretic_Moebius_inversion:diagram}
Then we obtain pairs of inverse isomorphisms $(S,E)$
(see Theorem~\ref{the_splitting_of_K-theory_for_EI_categories}),
$(\Res,I)$ (see Theorem~\ref{the:second_splitting_of_K_0(RGamma)}) and
$(\omega,\mu)$ (see~\eqref{omega_splitK_0_to_splitK_0}
and~\eqref{mu_splitK_0_to_splitK_0}).  They are compatible with one another
in the sense that the following diagram commutes
$$\xymatrix{
  & K_0(R\Gamma) \ar@/_/[ddl]_S \ar@/_/[ddr]_{\Res} &
  \\
  & &
  \\
  \Split K_0(R\Gamma) \ar@/_/[uur]_{E} \ar@/_/[rr]_{\omega} && \Split
  K_0(R\Gamma). \ar@/_/[uul]_{I} \ar@/_/[ll]_{\mu} }
$$
\item \label{the:Two_splittings_and_the_K-theoretic_Moebius_inversion:o}
Suppose that $\Gamma$ is of type (FP$_R$), or, equivalently, that
$\iso(\Gamma)$ is finite and for each object $x \in \ob(\Gamma)$ the
trivial $R[x]$-module $R$ possesses a finite projective
$R[x]$-resolution. Let $\eta \in \Split K_0(R\Gamma)$ be the element
whose component at $\overline{x} \in \iso(\Gamma)$ is given by the
class $[R] \in K_0(R[x])$ of the trivial $R[x]$-module $R$. That is,
the component of $\eta$ at each $\overline{x}$ is the finiteness
obstruction $o(\widehat{\aut(x)};R)\in K_0(R\aut(x))$. Then
$$S\left(o(\Gamma;R)\right) = \mu(\eta).$$
\end{enumerate}
\end{theorem}
\begin{proof}\ref{the:Two_splittings_and_the_K-theoretic_Moebius_inversion:diagram}
  We have already shown in
  Theorem~\ref{the_splitting_of_K-theory_for_EI_categories} that $S$
  and $E$ are inverse to one another and in
  Theorem~\ref{the:second_splitting_of_K_0(RGamma)} that $\Res$ and
  $I$ are inverse to one another. Obviously $\omega = \Res \circ E$.
  Hence it remains to show that $\mu \circ \omega = \id$.
  This follows analogously to the argument at the end of the proof
  of~L\"uck~\cite[Theorem~16.27 on page~330]{Lueck(1989)}.
\\[1mm]\ref{the:Two_splittings_and_the_K-theoretic_Moebius_inversion:o}
This follows from
assertion~\ref{the:Two_splittings_and_the_K-theoretic_Moebius_inversion:diagram}
and
Lemma~\ref{lem:finite_homological_dimension}~\ref{lem:finite_homological_dimension:Res_x}~and~\ref{lem:finite_homological_dimension:criterion_for_FP}.
Namely, $\Res_x[R]=[R]$, so $\Res [R] = \eta$, and
$S\left(o(\Gamma;R)\right) = \mu \Res \left(o(\Gamma;R)\right)= \mu
\Res [R] =\mu(\eta).$
\end{proof}

We can now apply M\"obius inversion to calculate the finiteness
obstruction and Euler characteristics of finite EI-categories in
terms of chains.

\begin{theorem}[The finiteness obstruction and Euler characteristics of finite EI-categories]
  \label{the:The_finiteness_obstruction_and_Moebius_inversion_finite_EI-categories}
  Suppose that $\Gamma$ is a finite EI-category.  Suppose that for every object
  $x \in \ob(\Gamma)$ the order of its automorphism group $|\aut(x)|$
  is invertible in $R$. Then $\Gamma$ is of type (FP$_R$) and we have:

  \begin{enumerate}

  \item \label{the:The_finiteness_obstruction_and_Moebius_inversion_finite_EI-categories:o}
  The image of the finiteness obstruction $o(\Gamma;R)$ under the isomorphism
  $$S \colon K_0(R\Gamma) \xrightarrow{\cong}
  \bigoplus_{\overline{y} \in \iso(\Gamma)} K_0(R[y])$$
  has as component for $\overline{y} \in \iso(\Gamma)$ the element in
  $K_0(R[y])$ given by
  $$\sum_{l \ge 0} (-1)^l \cdot \sum_{\overline{x} \in \iso(\Gamma)}\;
  \sum_{c \in \ch_l(y,x)} \bigl[R(\aut(x)\backslash S(c))\bigr],$$
  where $\aut(x)\backslash S(c)$ is the finite right $\aut(y)$-set
  obtained from the $\aut(x)$-$\aut(y)$-biset $S(c)$ (see Definition~\ref{def:chains})
  by dividing out the left $\aut(x)$-action and $R(\aut(x)\backslash S(c))$
  is the associated right $R[y]$-module;

  \item \label{the:The_finiteness_obstruction_and_Moebius_inversion_finite_EI-categories:chi_f}
  The functorial Euler characteristic $\chi_f(\Gamma;R) \in U(\Gamma)$ has at
  $\overline{y}$ the value
  $$\sum_{l \ge 0} (-1)^l \cdot \sum_{\overline{x} \in \iso(\Gamma)}\;
  \sum_{c \in \ch_l(y,x)} \bigl|\aut(x)\backslash S(c)/\aut(y)\bigr|,$$
  where $\bigl|\aut(x)\backslash S(c)/\aut(y)\bigr|$ is the order of
  the set obtained from $S(c)$ by dividing out the $\aut(x)$-action and the
  $\aut(y)$-action;

  \item \label{the:The_finiteness_obstruction_and_Moebius_inversion_finite_EI-categories:chi}
  The Euler characteristic $\chi(\Gamma,R)$ and topological Euler characteristic $\chi(B\Gamma;R)$ are equal and are both
  given by the integer
  $$\sum_{l \ge 0}  (-1)^l \cdot \sum_{\overline{x},\overline{y} \in \iso(\Gamma)}\;
  \sum_{c \in \ch_l(y,x)} \bigl|\aut(x)\backslash S(c)/\aut(y)\bigr|;$$

  \item \label{the:The_finiteness_obstruction_and_Moebius_inversion_finite_EI-categories:chi_f(2)}
  The functorial $L^2$-Euler characteristic
  $\chi_f^{(2)}(\Gamma) \in U^{(1)}(\Gamma)$ has at
  $\overline{y}$ the value
  $$\sum_{l \ge 0} (-1)^l \cdot \sum_{\overline{x} \in \iso(\Gamma)}\; \sum_{c \in \ch_l(y,x)}
  \dim_{\caln(y)}\bigl(\IC(\aut(x)\backslash S(c)) \otimes_{\IC[y]} \caln(y)\bigr),$$
  where $\dim_{\caln(y)}\bigl(\IC(\aut(x)\backslash S(c)) \otimes_{\IC[y]} \caln(y)\bigr)$
  is $\sum_{i \in I,|L_i| < \infty} 1/|L_i|$
  if the cofinite right $\aut(y)$-set $\aut(x)\backslash S(c)$
  is the disjoint union of homogeneous $\aut(y)$-spaces
  $\coprod_{i \in I} L_i\backslash\aut(y)$;

  \item \label{the:The_finiteness_obstruction_and_Moebius_inversion_finite_EI-categories:chi(2)}
  The $L^2$-Euler characteristic $\chi^{(2)}(\Gamma)$ is given by
  $$\sum_{l \ge 0}  (-1)^l \cdot \sum_{\overline{x},\overline{y} \in \iso(\Gamma)}\;
  \sum_{c \in \ch_l(y,x)}
  \dim_{\caln(y)}\bigl(\IC(\aut(x)\backslash S(c)) \otimes_{\IC[y]} \caln(y)\bigr).$$
\end{enumerate}
\end{theorem}
  \begin{proof}
  The category $\Gamma$ is of type (FP$_R$) by
  Lemma~\ref{lem:finite_homological_dimension}~\ref{lem:finite_homological_dimension:finite_Gamma}.
  \\[1mm]\ref{the:The_finiteness_obstruction_and_Moebius_inversion_finite_EI-categories:o}
  This follows from
  Theorem~\ref{the:Two_splittings_and_the_K-theoretic_Moebius_inversion}~%
\ref{the:Two_splittings_and_the_K-theoretic_Moebius_inversion:o}
  since the $R[y]$-modules $R \otimes_{R\aut(x)} RS(c)$ and
  $R(\aut(x)\backslash S(c))$ are isomorphic.
  \\[1mm]\ref{the:The_finiteness_obstruction_and_Moebius_inversion_finite_EI-categories:chi_f}
   and~\ref{the:The_finiteness_obstruction_and_Moebius_inversion_finite_EI-categories:chi}
  follow now from
  assertion~\ref{the:The_finiteness_obstruction_and_Moebius_inversion_finite_EI-categories:o},
  Lemma~\ref{lem:directly_finite_and_idempotents_and_EI}, and
  Theorem~\ref{the:chi_f_determines_chi}.
  \\[1mm]\ref{the:The_finiteness_obstruction_and_Moebius_inversion_finite_EI-categories:chi_f(2)}
  and~\ref{the:The_finiteness_obstruction_and_Moebius_inversion_finite_EI-categories:chi(2)}
  follow from Theorem~\ref{the:comparing_o_and_chi(2)},
  Example~\ref{exa:examples of the dimension}~\ref{exa:von_Neumann_dimension_and_permutation_modules}
  and assertion~\ref{the:The_finiteness_obstruction_and_Moebius_inversion_finite_EI-categories:o}.
\end{proof}

\begin{example}[M\"obius inversion for a finite partially ordered set]
  \label{Finite_partially_ordered_set}
  Let $(I,\le)$ be a partially ordered set. It defines an
  EI-category $\Gamma(I)$ whose set of objects is $I$ and for which
  $\mor(x,y)$ consists of precisely one element if $x \le y$ and is
  empty otherwise.

  Suppose that $I$ is finite.  Take $R = \IQ$. Then
  $$\Split K_0(\IQ\Gamma(I))  = \IZ I = \bigoplus_I \IZ$$
  and the homomorphism $\omega$ is given by the matrix
  $A = \bigl(a_{i,j}\bigr)_{i,j \in I}$ with $a_{i,j} = 1$ if $j \le i$ and $w_{i,j} = 0$
  otherwise. Let $B = \bigl(b_{i,j}\bigr)_{i,j \in I}$ be the matrix given by
  $$b_{i,j} = \sum_{l \ge 0} (-1)^l \cdot |\ch_l(j,i)|,$$
  where $|\ch_0(j,i)|$ is $0$ if $j \not= i$ and $1$ otherwise, and for
  $l \ge 1$, $\ch_l(j,i)$ is the set of chains $j = k_0 < k_1 < \ldots
  < k_{l-1} < k_l = i$.  Then we conclude from
  Theorem~\ref{the:Two_splittings_and_the_K-theoretic_Moebius_inversion}
  that the matrices $A$ and $B$ are inverse to one another. This is the
  classical M\"obius inversion in combinatorics (see for
  instance~Aigner~\cite[IV.2]{Aigner(1979)}).

  We get from
  Theorem~\ref{the:The_finiteness_obstruction_and_Moebius_inversion_finite_EI-categories}~%
\ref{the:The_finiteness_obstruction_and_Moebius_inversion_finite_EI-categories:chi}
and~\ref{the:The_finiteness_obstruction_and_Moebius_inversion_finite_EI-categories:chi(2)}
$$\chi(\Gamma;\IQ) = \chi^{(2)}(\Gamma) = \sum_{i,j \in I} b_{i,j}.$$
\end{example}

\begin{example}[M\"obius inversion for a finite skeletal category with trivial
endomorphisms]\label{exa:Moebius_Inversion_for_fin_sk_cat_with_triv_endos}
Generalizing Example~\ref{Finite_partially_ordered_set}, let
$\Gamma$ be a finite skeletal category in which every endomorphism
is an identity, and take $R = \IQ$. Recall that a category is
\emph{skeletal} if for any two objects $x$ and $y$ with $x \cong y$,
we have $x=y$. Then
$$\Split K_0(\IQ\Gamma) = \IZ \ob(\Gamma) = \bigoplus_{\ob(\Gamma)}
\IZ$$ and the homomorphism $\omega$ is given by the matrix $A =
\bigl(a_{x,y}\bigr)_{x,y \in \ob(\Gamma)}$ with $a_{x,y} =
|\mor(y,x)|$.

The (bi)set $S(c)$ in Definition \ref{def:chains} is simply the set of non-degenerate paths
$x_0 \to x_1 \to \cdots \to x_l$,
and $\mu_{x,y}(c)=|S(c)|$. Let $B = \bigl(b_{x,y}\bigr)_{x,y \in \ob(\Gamma)}$ be the matrix given by
$$b_{x,y} = \sum_{l \ge 0} (-1)^l \cdot \sum_{c \in \ch_l(y,x)} |S(c)|=
\sum_{l \ge 0} (-1)^l \cdot |\{\text{non-degenerate $l$-paths from
$y$ to $x$}\}|.$$ Then we conclude from
Theorem~\ref{the:Two_splittings_and_the_K-theoretic_Moebius_inversion}
that the matrices $A$ and $B$ are inverse to one another. That is to
say, in the terminology of~Leinster~\cite{Leinster(2008)}, the
category $\Gamma$ has M\"obius inversion given by $B$. Thus
Corollary~1.5 of~Leinster~\cite{Leinster(2008)} is a special case of
the $K$-theoretic M\"obius inversion of
Theorem~\ref{the:Two_splittings_and_the_K-theoretic_Moebius_inversion}~\ref{the:Two_splittings_and_the_K-theoretic_Moebius_inversion:diagram}.
See also Example
\ref{exa:Rational_Moebius_Inversion_for_fin_sk_cat_free_EI}, which
illustrates rational M\"obius inversion for a finite, skeletal, free
EI-category. See also the related proof of Lemma
\ref{lem:chi(2)_and_chi}, which shows that the $L^2$-Euler
characteristic coincides with Leinster's Euler characteristic in the
case of a finite, skeletal, free EI-category.
\end{example}

\subsection{The $K$-theoretic M\"obius inversion and the $L^2$-rank}
\label{subsec:The-K-theoretic-Moebius-Inversion_and_the_L2-rank}

In this subsection we investigate when the homomorphisms $\omega$
and $\mu$ factorize over the homomorphism given by the $L^2$-rank.

\begin{condition}[Condition (I) for groups and categories]
\label{con:condition_(I)}
  A group $G$ \emph{satisfies condition (I)} if the map induced by the
  various inclusions of finite subgroups
  $$\bigoplus_{H \subseteq G, |H| < \infty} K_0(\IQ H) \otimes_{\IZ} \IQ  \to
  K_0(\IQ G) \otimes_{\IZ}  \IQ$$
  is surjective. A category $\Gamma$ \emph{satisfies condition (I)} if for every object $x$ its
  automorphism group $\aut_\Gamma(x)$ satisfies condition (I).
\end{condition}

Obviously any finite group and any finite category satisfy condition (I).

\begin{remark}[Condition (I) and the Farrell-Jones Conjecture]
  \label{rem:groups_satisfying_condition_(I)}
  Let $\calf\calj(\IQ)$ be the class of groups for which the
  $K$-theoretic Farrell-Jones Conjecture with coefficients in $\IQ$
  holds.
  By~Bartels--L\"uck--Reich~\cite[Theorem~0.5]{Bartels-Lueck-Reich(2008appl)},
  every group in $\calf\calj(\IQ)$ satisfies condition (I).  This
  class $\calf\calj(\IQ)$ is analyzed for instance
  by Bartels--L\"uck in~\cite{Bartels-Lueck(2009borelhyp)} and Bartels--L\"uck--Reich in~\cite{Bartels-Lueck-Reich(2008hyper)}
  and~\cite{Bartels-Lueck-Reich(2008appl)}.
  It contains for instance subgroups of finite products
  of hyperbolic groups or $\CAT(0)$-groups, directed colimits of hyperbolic groups
  or $\CAT(0)$-groups, and all elementary amenable groups. For a survey article
  on the Farrell-Jones Conjecture we refer for instance
  to~L\"uck--Reich~\cite{Lueck-Reich(2005)}.
\end{remark}

\begin{lemma} \label{lem:L2-rank_and_bisets} Let $G$ and $H$ be
  groups. Suppose that $H$ satisfies condition (I) defined
  in~\eqref{con:condition_(I)}. Let $S$ be an $H$-$G$-biset which is
  cofinite proper as a right $G$-set and free as a left $H$-set.

  \begin{enumerate}

  \item \label{lem:L2-rank_and_bisets:image_of_rank(2)_is_in_Q}
  The image of
  $$\rk_{H}^{(2)} \colon K_0(\IQ H) \to \IR, \quad
  [P] \mapsto \dim_{\caln(H)}\bigl(P \otimes_{\IQ H} \caln(H)\bigr)$$
  lies in $\IQ$;

  \item \label{lem:L2-rank_and_bisets:diagram}
  The following diagram commutes
  $$
  \xymatrix{K_0(\IQ H) \ar[r]^{\omega_S} \ar[d]_{\rk_H^{(2)}} & K_0(\IQ
  G) \ar[d]^{\rk_G^{(2)}}
  \\
  \IR \ar[r]_{\overline{\omega}_S} & \IR}
  $$
  where $\omega_S$ sends $[P]$ to $[P \otimes_{\IQ H} \IQ S]$,
  and $\overline{\omega}_S$
  is multiplication with the rational number
  $\dim_{\caln(G)}\bigl(\IQ S \otimes_{\IQ G} \caln(G)\bigr)$.
 \end{enumerate}
\end{lemma}
\begin{proof}\ref{lem:L2-rank_and_bisets:image_of_rank(2)_is_in_Q}
Because $H$ satisfies condition (I), this follows from
Lemma~\ref{lem:dim_caln(G)_and_induction_and_restriction}~%
\ref{lem:dim_caln(G)_and_induction_and_restriction:dim_and_induction}
and Example~\ref{exa:examples of the dimension}~\ref{exa:The_von_Neumann_dimension_for_finite_groups} .
\\[1mm]\ref{lem:L2-rank_and_bisets:diagram}
  For a finite group $H'$ every element in  $K_0(\IQ H') \otimes_{\IZ} \IQ$
  can be written as a $\IQ$-linear combination of elements of the form
  $\bigl[\IQ [K\backslash H']\bigr]$ (see~Serre~\cite[Theorem~30 in Chapter~13 on page~103]{Serre(1977)}).
  Since $H$ in the claim satisfies condition (I), we can find for every element
  $\eta \in K_0(\IQ H)$ a natural number $k \ge 1$, finitely many
  finite subgroups $K_1$, $K_2$, \ldots, $K_r$ of $H$, and integers
  $n_1$, $n_2$, \ldots , $n_r$ such that we get in $K_0(\IQ H)$
  $$k \cdot \eta = \sum_{i=1}^r n_i \cdot \bigl[\IQ[K_i\backslash H]\bigr].$$
  Hence it suffices
  to show for any finite subgroup $K \subseteq H$
\begin{multline*}
  \dim_{\caln(G)}\bigl(\IQ[K\backslash H] \otimes_{\IQ H} \IQ S \otimes_{\IQ G}
  \caln(G)\bigr)
  \\
  = \dim_{\caln(H)}\bigl(\IQ[K\backslash H] \otimes_{\IQ H} \caln(H)\bigr)
  \cdot \dim_{\caln(G)}\bigl(\IQ
  S \otimes_{\IQ G} \caln(G)\bigr).
\end{multline*}
We get from
Example~\ref{exa:examples of the dimension}~\ref{exa:von_Neumann_dimension_and_permutation_modules}
\begin{eqnarray*}
  \dim_{\caln(H)}\bigl(\IQ[K\backslash H]  \otimes_{\IQ H} \caln(H)\bigr) & = & \frac{1}{|K|};
  \\
  \dim_{\caln(G)}\bigl(\IQ[K\backslash H] \otimes_{\IQ H} \IQ S \otimes_{\IQ G} \caln(G)\bigr)
  & = &
  \dim_{\caln(G)}\bigl(\IQ[K \backslash S] \otimes_{\IQ G} \caln(G)\bigr).
\end{eqnarray*}
Hence it suffices to show for a $K$-$G$-biset $T$ which is proper and
cofinite as a $G$-set and free as a left $K$-set
\begin{eqnarray*}
  |K| \cdot \dim_{\caln(G)}\bigl(\IQ[K \backslash T] \otimes_{\IQ G} \caln(G)\bigr)
  & = &
  \dim_{\caln(G)}\bigl(\IQ T \otimes_{\IQ G} \caln(G)\bigr).
\end{eqnarray*}
  We can interpret the $K$-$G$-biset $T$ as a right $(K \times G)$-set by
  putting $t \cdot (k,g) = k^{-1} t g$ for $k\in K$, $g \in G$ and $t \in T$,
  and vice versa.  Since $K$ is finite, $T$ is free as a left $K$-set, and $T$
  is cofinite and proper as a right
  $G$-set, the $(K \times G)$-set $T$ is a finite union
  of homogeneous spaces of the form $L\backslash (K \times G)$, where $L$ is a
  finite subgroup of $K \times G$ with $\left( K \times \{1\} \right) \cap L = \{1\}$.
  Hence we can assume without loss of
  generality that $T$ is of the form $L\backslash (K \times G)$ for
  finite $L \subseteq K \times G$ with $\left( K \times \{1\} \right) \cap L = \{1\}$

  The projection $\pr \colon K \times G \to G$ induces a
  bijection $L \xrightarrow{\cong} \pr(L)$.  Since the $G$-sets
  $K\backslash \bigl(L\backslash (K \times G)\bigr)$ and
  $\pr(L)\backslash G$ are $G$-isomorphic, we conclude from
  Example~\ref{exa:examples of the dimension}~\ref{exa:von_Neumann_dimension_and_permutation_modules}
  \begin{eqnarray*}
  |K| \cdot \dim_{\caln(G)}\left(\IQ\bigl[K\backslash \bigl(L\backslash (K \times G)\bigr)\bigr]
  \otimes_{\IQ G} \caln(G)\right)
  & = &
  \frac{|K|}{|L|}.
  \end{eqnarray*}
  We conclude from
  Lemma~\ref{lem:dim_caln(G)_and_induction_and_restriction} and
  Example~\ref{exa:examples of the dimension}~\ref{exa:von_Neumann_dimension_and_permutation_modules}
  \begin{eqnarray*}
  \lefteqn{\dim_{\caln(G)}\left(\IQ\bigl[L\backslash (K \times G)\bigr] \otimes_{\IQ G} \caln(G)\right)}
  & &
  \\
  & = &
  \dim_{\caln(G)}\left(\IQ\bigl[L\backslash (K \times G)\bigr] \otimes_{\IQ [K \times G]}
    \IQ [K \times G] \otimes_{\IQ G} \caln(G)\right)
  \\
  & = &
  \dim_{\caln(G)}\left(\IQ\bigl[L\backslash (K \times G)\bigr]
  \otimes_{\IQ[K \times G]} \caln(K \times G)\right)
  \\
  & = &
  |K| \cdot
  \dim_{\caln(K \times G)}\left(\IQ\bigl[L\backslash (K \times G)\bigr]
   \otimes_{\IQ [K \times G]} \caln(K \times G)\right)
  \\
  & = &
  \frac{|K|}{|L|}.
  \end{eqnarray*}
  This finishes the proof of Lemma~\ref{lem:L2-rank_and_bisets}.
\end{proof}

Let $\Gamma$ be a quasi-finite, free EI-category. Define the
$\IQ$-homomorphism
\begin{eqnarray}
& \overline{\omega}^{(2)} \colon U(\Gamma) \otimes_{\IZ} \IQ \to U(\Gamma) \otimes_{\IZ} \IQ
\label{overlineomega(2)}
\end{eqnarray}
by the matrix over the rational numbers
$$\biggl(\dim_{\caln(y)}\bigl(\IQ\mor(y,x) \otimes_{\IQ[y]} \caln(y)\bigr)
\biggr)_{\overline{x},\overline{y} \in \iso(\Gamma)}.$$
Define the $\IQ$-homomorphism
\begin{eqnarray}
& \overline{\mu}^{(2)} \colon U(\Gamma) \otimes_{\IZ} \IQ \to U(\Gamma) \otimes_{\IZ} \IQ
\label{overlinemu(2)}
\end{eqnarray}
by the matrix over the rational numbers
\[
\biggl(\sum_{l \ge 0} (-1)^l \cdot \sum_{c \in \ch_l(y,x)}
\dim_{\caln(y)}\bigl(\IQ S(c) \otimes_{\IQ[y]} \caln(y)\bigr)
\biggr)_{\overline{x},\overline{y} \in \iso(\Gamma)}.
\]
Notice that these homomorphisms are well-defined because of
Example~\ref{exa:examples of the
dimension}~\ref{exa:von_Neumann_dimension_and_permutation_modules}
since the right $\aut(y)$-sets $\mor(y,x)$ and $S(c)$ are proper
cofinite and for a given $\overline{x} \in \iso(\Gamma)$ there are
only finitely many $\overline{y} \in \iso(\Gamma)$ for which the
sets $\mor(y,x)$ and $S(c)$ are non-empty.

\begin{theorem}[Rational M\"obius inversion]
  \label{theorem:rational_Moebius_inversion}
  Let $\Gamma$ be a quasi-finite, free EI-category. Then the
  homomorphisms $\overline{\omega}^{(2)}$ of~\eqref{overlineomega(2)}
  and $\overline{\mu}^{(2)}$ of~\eqref{overlinemu(2)} are isomorphisms and inverse
  to one another.
\end{theorem}
\begin{proof}
  Let
  \begin{eqnarray*}
   \overline{\iota} \colon U(\Gamma) = \bigoplus_{\overline{x} \in \iso(\Gamma)} \IZ
  & \to &
 \Split K_0(\IQ\Gamma) = \bigoplus_{\overline{x} \in \iso(\Gamma)} K_0(\IQ[x])
 \end{eqnarray*}
 be the homomorphism that sends $\{n_{\overline{x}} \mid \overline{x} \in \iso(\Gamma)\}$
  to $\bigl\{n_x \cdot [\IQ[x]] \;\big|\;\overline{x} \in \iso(\Gamma)\bigr\}$.
  A direct computation shows that
  $$\rk_{\Gamma}^{(2)} \circ \omega \circ \overline{\iota} = \overline{\omega}^{(2)}.$$
  The image of $\omega \circ \overline{\iota}$ in $\Split K_0(\IQ \Gamma)$ has
  the property that its value at any $\overline{x} \in \iso(\Gamma)$
  is an element in $K_0(\IQ[x])$ given by a $\IZ$-linear combination
  of classes of the form $\bigl[\IQ[K\backslash \aut(x)]\bigr]$ for finite subgroups $K\subseteq
  \aut(x)$.  Hence the argument in the proof of
  Lemma~\ref{lem:L2-rank_and_bisets}~\ref{lem:L2-rank_and_bisets:diagram}
  shows (without using condition (I)) that $\rk_{\Gamma}^{(2)} \circ
  \mu = \overline{\mu}^{(2)} \circ \rk_{\Gamma}^{(2)}$ is true on the image
  of $\omega \circ \overline{\iota}$. This implies
  $$\overline{\mu}^{(2)} \circ \overline{\omega}^{(2)}
  = \overline{\mu}^{(2)} \circ \rk_{\Gamma}^{(2)}\circ \omega \circ \overline{\iota} =
  \rk_{\Gamma}^{(2)} \circ \mu \circ \omega \circ \overline{\iota}.
  $$
  We conclude $\mu \circ \omega=\id$  from
  Theorem~\ref{the:Two_splittings_and_the_K-theoretic_Moebius_inversion}.
  A direct computation shows $\rk_{\Gamma}^{(2)} \circ \overline{\iota} = \id$.
  Hence
  $$\overline{\mu}^{(2)} \circ \overline{\omega}^{(2)} = \id.$$
  Since the matrix defining $\overline{\omega}^{(2)}$ is a triangular matrix whose
  entries on the diagonal are all $1$, $ \overline{\omega}^{(2)}$ is
  an isomorphism. Hence $\overline{\omega}^{(2)}$
  of~\eqref{overlineomega(2)} and $\overline{\mu}^{(2)}$
  of~\eqref{overlinemu(2)} are isomorphisms and inverse
  to one another.
\end{proof}

\begin{remark} \label{rem:condition_free}
Notice that the condition free is not needed when we want to define
the finiteness obstruction or to compute it as long as we stay on
the $K$-theory level. It does enter, when we want to consider the
rank or $L^2$-rank of the finiteness obstruction, to ensure that
certain comparisons can be done on the level of the Euler
characteristics, or, equivalently, certain maps on the $K_0$-level
factorize over the rank or $L^2$-rank homomorphism from
$K_0(R\Gamma)$ to $U(\Gamma)$ or $U(\Gamma) \otimes_\IZ \IR$.
\end{remark}

\begin{example}[Rational M\"obius inversion for a finite, skeletal, free EI-category]
\label{exa:Rational_Moebius_Inversion_for_fin_sk_cat_free_EI}
Generalizing Example~\ref{Finite_partially_ordered_set}, let
$\Gamma$ be a finite skeletal EI-category which is free in the sense
of Definition~\ref{def:finite_quasi-finite_and_free_category}, and
take $R = \IQ$. Then
$$U(\Gamma) \otimes_\IZ \IQ = \bigoplus_{\ob(\Gamma)} \IQ$$
and the homomorphism $\overline{\omega}^{(2)}$ is given by the
matrix $$\biggl(\dim_{\caln(y)}\bigl(\IQ\mor(y,x) \otimes_{\IQ[y]}
\caln(y)\bigr) \biggr)_{x,y \in \ob(\Gamma)}=
\biggl(\frac{|\mor(y,x)|}{|\aut(y)|}\biggr)_{x,y \in \ob(\Gamma)}.$$
The last equality follows from
Example~\ref{exa:examples of the dimension}~\ref{exa:von_Neumann_dimension_and_permutation_modules}. If
we let $\omega_L$ be the matrix
$$\left(|\mor_{\Gamma}(y,x)|\right)_{x,y \in \ob(\Gamma)}$$ and
$D$ is the diagonal matrix with entry $|\aut(y)|$ at $(y,y)$ for $y
\in \ob(\Gamma)$, then $D \circ \overline{\omega}^{(2)}=\omega_L$.

Then by Theorem \ref{theorem:rational_Moebius_inversion}, the
homomorphism $\overline{\omega}^{(2)}$ is invertible and its inverse
is $\overline{\mu}^{(2)}$. Hence $\omega_L$ admits an inverse
$\mu_L:=(D \circ \overline{\omega}^{(2)})^{-1}=\overline{\mu}^{(2)}
\circ D^{-1}$. We calculate $\mu_L$ by way of the matrix for
$\overline{\mu}^{(2)}$ using the formula just after
equation~\eqref{overlinemu(2)}. For any $l$-chain $c \in \ch_l(y,x)$
with $c=x_0<x_1< \cdots < x_l$ we have
$$|S(c)|=
\frac{|\mor(x_{l-1},x_l)|\cdot|\mor(x_{l-2},x_{l-1})|\cdot
\cdots\cdot |\mor(x_0,x_1)|}{|\aut(x_{l-1})| \cdot |\aut(x_{l-2})|
\cdot \cdots \cdot |\aut(x_{1})| }$$ by freeness. Then,
\begin{multline*}
\dim_{\caln(y)}\bigl(\IQ S(c) \otimes_{\IQ[y]} \caln(y)\bigr)\\=
\frac{|S(c)|}{|\aut(y)|}=\frac{|\mor(x_{l-1},x_l)|\cdot|\mor(x_{l-2},x_{l-1})|\cdot
\cdots\cdot |\mor(x_0,x_1)|}{|\aut(x_{l-1})| \cdot |\aut(x_{l-2})|
\cdot \cdots \cdot |\aut(x_{1})| \cdot |\aut(x_{0})|}
\end{multline*}
by
Example~\ref{exa:examples of the dimension}~\ref{exa:von_Neumann_dimension_and_permutation_modules}.
Summing up, we have
$$\aligned
\mu_L &=\overline{\mu}^{(2)} \circ D^{-1} \\
&= \biggl(\sum_{l \ge 0} (-1)^l \cdot \sum_{c \in \ch_l(y,x)}
\dim_{\caln(y)}\bigl(\IQ S(c) \otimes_{\IQ[y]} \caln(y)\bigr)
\biggr)_{x,y \in \ob(\Gamma)} \circ D^{-1} \\
&=\biggl(\sum_{l \ge 0} (-1)^l \cdot \sum_{c \in \ch_l(y,x)}
\frac{|\mor(x_{l-1},x_l)|\cdot|\mor(x_{l-2},x_{l-1})|\cdot
\cdots\cdot |\mor(x_0,x_1)|}{|\aut(x_{l-1})| \cdot |\aut(x_{l-2})|
\cdot \cdots \cdot |\aut(x_{1})| \cdot |\aut(x_{0})|} \biggr)_{x,y \in \ob(\Gamma)} \circ D^{-1}\\
&=\biggl(\sum_{l \ge 0} (-1)^l \cdot \sum_{c \in \ch_l(y,x)}
\frac{|\mor(x_{l-1},x_l)|\cdot|\mor(x_{l-2},x_{l-1})|\cdot
\cdots\cdot |\mor(x_0,x_1)|}{|\aut(x_{l})| \cdot |\aut(x_{l-1})|
\cdot |\aut(x_{l-2})| \cdot \cdots \cdot |\aut(x_{1})| \cdot
|\aut(x_{0})|}\biggr)_{x,y \in \ob(\Gamma)} \\
&=\biggl(\sum_{l \ge 0} (-1)^l \cdot \sum \frac{1}{|\aut(x_{l})|
\cdot |\aut(x_{l-1})| \cdot |\aut(x_{l-2})| \cdot \cdots \cdot
|\aut(x_{1})| \cdot |\aut(x_{0})|} \biggr)_{x,y \in \ob(\Gamma)}.
\endaligned$$
The final sum is over all $l$-paths $x_0 \to x_1 \to \cdots \to x_l$
from $y$ to $x$ such that $x_0, \ldots, x_l$ are all distinct. Thus,
in the terminology of~Leinster~\cite{Leinster(2008)}, the category
$\Gamma$ has M\"obius inversion given by $\mu_L$, and Leinster's
Euler characteristic $\chi_L(\Gamma)$ is the sum of the entries in
the matrix $\mu_L$ above. The free case of
Leinster~\cite[Theorem~1.4]{Leinster(2008)} is now a special case of
rational M\"obius inversion (Theorem
\ref{theorem:rational_Moebius_inversion}). See also the related
proof of Lemma \ref{lem:chi(2)_and_chi}, which shows that the
$L^2$-Euler characteristic coincides with Leinster's Euler
characteristic in the case of a finite, skeletal, free EI-category.
Thus, the $L^2$-Euler characteristic $\chi^{(2)}(\Gamma)$ is also
given by the sum of the entries in the matrix $\mu_L$ above.
\end{example}

\begin{theorem}[The $K$-theoretic M\"obius inversion and the $L^2$-rank]
\label{the:K_theoretic-Moebius_inversion_and_L2-rank}
 Let $\Gamma$ be a quasi-finite, free EI-category satisfying condition
 (I) defined in~\ref{rem:groups_satisfying_condition_(I)}.
 Then the following diagram commutes
  $$\xymatrix{
  & K_0(\IQ\Gamma) \ar@/_/[ddl]_S \ar@/_/[ddr]_{\Res} &
  \\
  & &
  \\
  \Split K_0(\IQ\Gamma) \ar@/_/[uur]_{E} \ar@/_/[rr]_{\omega}
  \ar[dd]_{\rk_{\Gamma}^{(2)}} &&
  \Split K_0(\IQ\Gamma) \ar@/_/[uul]_{I} \ar@/_/[ll]_{\mu} \ar[dd]^{\rk_{\Gamma}^{(2)}}
  \\
  & &
  \\
  U(\Gamma) \otimes_{\IZ} \IQ \ar@/_/[rr]_{\overline{\omega}^{(2)}}
  &&
  U(\Gamma) \otimes_{\IZ} \IQ   \ar@/_/[ll]_{\overline{\mu}^{(2)}}
}
$$
Here the pairs  $(S,E)$
(see Theorem~\ref{the_splitting_of_K-theory_for_EI_categories}),
$(\Res,I)$ (see Theorem~\ref{the:second_splitting_of_K_0(RGamma)}),
$(\omega,\mu)$ (see~Theorem~\ref{the:Two_splittings_and_the_K-theoretic_Moebius_inversion}),
and $(\overline{\omega}^{(2)},\overline{\mu}^{(2)})$
(see Theorem~\ref{theorem:rational_Moebius_inversion})
are pairs of isomorphisms inverse to one another,
and the map $\rk_{\Gamma}^{(2)}$ comes from the map defined in~\eqref{rk(2)_K_o_to_U}.
\end{theorem}
\begin{proof}
The map $\rk_{\Gamma}^{(2)}$ takes values in $U(\Gamma)
\otimes_{\IZ} \IQ$ by
Lemma~\ref{lem:L2-rank_and_bisets}~\ref{lem:L2-rank_and_bisets:image_of_rank(2)_is_in_Q}.
The other claims follow from
Theorem~\ref{the:Two_splittings_and_the_K-theoretic_Moebius_inversion},
Lemma~\ref{lem:L2-rank_and_bisets}~\ref{lem:L2-rank_and_bisets:diagram},
and Theorem~\ref{theorem:rational_Moebius_inversion}.
\end{proof}

\begin{theorem}[The finiteness obstruction and the (functorial)
  $L^2$-Euler characteristic]
  \label{the:The_finiteness_obstruction_and_the_(functorial)_L2-Euler_characteristic}
  \
  \begin{enumerate}

  \item \label{the:The_finiteness_obstruction_and_the_(functorial)_L2-Euler_characteristic:Res(o)}
    Let $\Gamma$ be a quasi-finite EI-category of type (FP$_\IQ$). Then the
    image of the finiteness obstruction $o(\Gamma;\IQ)$ under the
    homomorphism
    $$\Res \colon K_0(\IQ\Gamma) \to \Split K_0(\IQ \Gamma)$$
    defined in Theorem~\ref{the:second_splitting_of_K_0(RGamma)} has as
    entry at $\overline{x} \in \iso(\Gamma)$ the finiteness obstruction
    $o\bigl(\widehat{\aut(x)};\IQ\bigr)$ of the category $\widehat{\aut(x)}$, i.e., the finiteness obstruction
    $o(\IQ)$ of the $\IQ[x]$-module $\IQ$ with the trivial
    $\aut(x)$-action. This possesses a finite projective $\IQ[x]$-resolution by
    Lemma~\ref{lem:finite_homological_dimension}~\ref{lem:finite_homological_dimension:Res_x}. As usual, we will write $[\IQ]$ for
    $o\bigl(\widehat{\aut(x)};\IQ\bigr)$.

\item \label{the:The_finiteness_obstruction_and_the_(functorial)_L2-Euler_characteristic:chi(2)_(I)}
  Suppose that $\Gamma$ is a quasi-finite, free EI-category of type
  (FP$_\IQ$) satisfying condition (I) or that $\Gamma$ is a
  quasi-finite, free EI-category of type (FF$_\IQ$).

  Then for every object $x$ the $L^2$-Euler characteristic
  $\chi^{(2)}(\aut(x))$ is a rational number and is non-trivial for
  only finitely many $\overline{x} \in \iso(\Gamma)$.  The collection
  $\bigl(\chi^{(2)}(\aut(x)\bigr)_{\overline{x} \in \iso(\Gamma)}$
  defines an element $\eta \in U(\Gamma) \otimes_{\IZ} \IQ$. The
  functorial $L^2$-Euler characteristic $\chi^{(2)}_f(\Gamma)$ lies in
  $U(\Gamma) \otimes_{\IZ} \IQ$. We get
  \begin{eqnarray*}
    \overline{\omega}^{(2)}\bigl(\chi^{(2)}_f(\Gamma)\bigr)  & = & \eta;
    \\
    \overline{\mu}^{(2)}(\eta) & = & \chi^{(2)}_f(\Gamma),
  \end{eqnarray*}
  where $\overline{\omega}^{(2)}$ and $\overline{\mu}^{(2)}$ are the
  homomorphisms defined in~\eqref{overlineomega(2)}
  and~\eqref{overlinemu(2)}.
\end{enumerate}
\end{theorem}
\begin{proof}\ref{the:The_finiteness_obstruction_and_the_(functorial)_L2-Euler_characteristic:Res(o)}
  Since $\Gamma$ is of type (FP$_\IQ$), we conclude from
  Lemma~\ref{lem:finite_homological_dimension}~\ref{lem:finite_homological_dimension:Res_x}
  that the $\IQ[x]$-module $\IQ$ with the trivial $\aut(x)$-action
  possesses a finite projective $\IQ[x]$-resolution and hence defines
  an element in $K_0(\IQ[x])$. Since $\Res_x \colon
  \MOD\text{-}\IQ\Gamma \to \MOD\text{-}\IQ[x]$ is exact, the claim follows from
  Lemma~\ref{lem:finite_homological_dimension}~\ref{lem:finite_homological_dimension:Res_x}.
  \\[1mm]\ref{the:The_finiteness_obstruction_and_the_(functorial)_L2-Euler_characteristic:chi(2)_(I)}
  We begin with the case where $\Gamma$ is a quasi-finite, free
  EI-category of type (FP$_\IQ$) satisfying condition (I).  The map
  $\rk_{\Gamma}^{(2)} \colon \Split K_0(\IQ \Gamma) \to
  \prod_{\overline{x} \in \iso(\Gamma)} \IR$ takes values in $U(\Gamma)
  \otimes_{\IZ} \IQ$ by
  Lemma~\ref{lem:L2-rank_and_bisets}~\ref{lem:L2-rank_and_bisets:image_of_rank(2)_is_in_Q}.
  The image of $o(\Gamma;\IQ)$ under the composite
  $$K_0(\IQ\Gamma) \xrightarrow{S} \Split K_0(\IQ\Gamma)
  \xrightarrow{\rk_{\Gamma}^{(2)}} \prod_{\overline{x} \in \iso(\Gamma)}
  \IR$$
  is by definition $\chi^{(2)}_f(\Gamma)$.  The image of
  $o(\Gamma;\IQ)$ under the composite
  $$K_0(\IQ\Gamma) \xrightarrow{\Res} \Split K_0(\IQ\Gamma)
  \xrightarrow{\rk_{\Gamma}^{(2)}} \prod_{\overline{x} \in \iso(\Gamma)} \IR$$
  is by definition $\eta$. Now the claim follows from
  Theorem~\ref{the:K_theoretic-Moebius_inversion_and_L2-rank}.

  Next we deal with the case where $\Gamma$ is a quasi-finite, free
  EI-category of type (FF$_\IQ$). Since $\Gamma$ is of type (FF$_\IQ$), the image
  of $o(\Gamma;\IQ)$ under the isomorphism $S \colon K_0(\IQ \Gamma)
  \xrightarrow{\cong} \Split K_0(\IQ\Gamma)$ is the image of
  $\chi^{(2)}_f(\Gamma) \in U(\Gamma)$ under the map $\iota \colon
  U(\Gamma) \to \Split K_0(\IQ\Gamma)$ defined in~\eqref{iota}, as
  $\rk^{(2)}_{\Gamma} \circ \iota$ is the inclusion of $U(\Gamma)$, see Lemma~\ref{lem:rank_is_rankL2_for_fin_free_modules}. A
  direct computation shows that $\overline{\omega}^{(2)} =
  \rk_{\Gamma}^{(2)} \circ \omega \circ \iota$.  This implies
  $$\overline{\omega}^{(2)}(\chi^{(2)}(\Gamma)) = \eta.$$
  We get
  \begin{eqnarray*}
  \overline{\mu}^{(2)}(\eta) & = & \chi^{(2)}_f(\Gamma),
  \end{eqnarray*}
  from Theorem~\ref{theorem:rational_Moebius_inversion}.
\end{proof}

\subsection{The example of a biset}
\label{subsec:The_example_of_a_biset}

  Let $H$ and $G$ be groups and let $S$ be a $G$-$H$-biset. They define an
  EI-category $\Gamma(S)$ with two objects $x$ and $y$, where the
  automorphism group of $x$ is $H$, the automorphism group of $y$ is
  $G$, the set of morphisms from $x$ to $y$ is $S$, the set of
  morphisms from $y$ to $x$ is empty and the composition in
  $\Gamma(S)$ comes from the group structure on $H$ and $G$ and the
  $G$-$H$-biset structure on $S$. Any EI-category with precisely two
  objects which are not isomorphic arises as $\Gamma(S)$ for some $S$.  The category
  $\Gamma(S)$ is free if and only if $S$ is free as a left $G$-set.  The
  category $\Gamma(S)$ is quasi-finite if and only if $S$ is proper
  and cofinite as a right $H$-set. The set of isomorphism classes of objects
  contains precisely two elements, namely $x$ and $y$.

  Suppose that $\Gamma(S)$ is quasi-finite. Then $\Gamma(S)$ is of type
  (FP$_\IQ$) if and only if the trivial $\IQ H$-module $\IQ$ has a finite projective
  $\IQ H$-resolution and the trivial $\IQ G$-module $\IQ$ has a finite projective
  $\IQ G$-resolution (see
  Lemma~\ref{lem:finite_homological_dimension}~\ref{lem:finite_homological_dimension:criterion_for_FP}).

  Suppose that $\Gamma(S)$ is quasi-finite and of type (FP$_\IQ$).
  Then the image of the finiteness obstruction
  under the isomorphism
  $$S \colon K_0(\IQ\Gamma(S))  \xrightarrow{\cong} K_0(\IQ H) \oplus K_0(\IQ G)$$
  is the element $\mu([\IQ],[\IQ]\bigr)$ by
  Theorem~\ref{the:Two_splittings_and_the_K-theoretic_Moebius_inversion}~%
\ref{the:Two_splittings_and_the_K-theoretic_Moebius_inversion:o},
where $[\IQ]$ stands, of course, for the finiteness obstruction of
the trivial $\IQ H$-module and trivial $\IQ G$-module $\IQ$,
respectively. That is, $[\IQ]$ means $o(\widehat{H};\IQ)$ or
$o(\widehat{G};\IQ)$ respectively.

  Suppose that $\Gamma(S)$ is quasi-finite, free, and of type (FP$_\IQ$). Then the
  $\IQ H$-module $\IQ G\backslash S$ has a finite projective
  $\IQ H$-resolution and the image of the finiteness obstruction
  under the isomorphism
  $$S \colon K_0(\IQ\Gamma(S))  \xrightarrow{\cong} K_0(\IQ H) \oplus K_0(\IQ G)$$
  is the element
  $$\mu([\IQ],[\IQ])=([\IQ]-[\IQ \otimes_{\IQ G} \IQ S],[\IQ])=\bigl([\IQ] - [\IQ G\backslash S],[\IQ]\bigr)$$
  by Theorem~\ref{the:Two_splittings_and_the_K-theoretic_Moebius_inversion}~%
\ref{the:Two_splittings_and_the_K-theoretic_Moebius_inversion:o}.

  Suppose that  $\Gamma(S)$ is quasi-finite, free, and of type (FP$_\IQ$),
  and that $H$ and $G$ satisfy Condition (I)
  (see~\ref{con:condition_(I)}).
  Then $\Gamma(S)$ satisfies Condition (I) by definition.
  The commutative diagram appearing in
  Theorem~\ref{the:K_theoretic-Moebius_inversion_and_L2-rank}
  $$\xymatrix{
  & K_0(\IQ\Gamma(S)) \ar@/_/[ddl]_S \ar@/_/[ddr]_{\Res} &
  \\
  & &
  \\
  \Split K_0(\IQ\Gamma(S)) \ar@/_/[uur]_{E} \ar@/_/[rr]_{\omega}
  \ar[dd]_{\rk_{\Gamma(S)}^{(2)}} &&
  \Split K_0(\IQ\Gamma(S)) \ar@/_/[uul]_{I} \ar@/_/[ll]_{\mu} \ar[dd]^{\rk_{\Gamma(S)}^{(2)}}
  \\
  & &
  \\
  U(\Gamma(S)) \otimes_{\IZ} \IQ \ar@/_/[rr]_{\overline{\omega}^{(2)}}
  &&
  U(\Gamma(S)) \otimes_{\IZ} \IQ   \ar@/_/[ll]_{\overline{\mu}^{(2)}}
}
$$
  becomes
  $$\xymatrix{
  & K_0(\IQ\Gamma(S)) \ar@/_/[ddl]_S \ar@/_/[ddr]_{\Res} &
  \\
  & &
  \\
  K_0(\IQ H) \oplus K_0(\IQ G)\ar@/_/[uur]_{E} \ar@/_/[rr]_{\omega}
  \ar[dd]_{\rk_{\Gamma(S)}^{(2)}} &&
   K_0(\IQ H) \oplus K_0(\IQ G) \ar@/_/[uul]_{I} \ar@/_/[ll]_{\mu} \ar[dd]^{\rk_{\Gamma(S)}^{(2)}}
  \\
  & &
  \\
  \IQ \oplus  \IQ \ar@/_/[rr]_{\overline{\omega}^{(2)}}
  &&
  \IQ \oplus \IQ   \ar@/_/[ll]_{\overline{\mu}^{(2)}}
}
$$
where $\omega$ sends $\bigl([P],[Q]\bigr)$ to $\bigl([P] +
[Q\otimes_{\IQ G} \IQ S],[Q]\bigr)$ and $\mu$ sends
$\bigl([P],[Q]\bigr)$ to $\bigl([P] - [Q\otimes_{\IQ G} \IQ
S],[Q]\bigr)$. If the proper cofinite right $H$-set $S$ is the
disjoint union $\coprod_{i=1}^r L_i\backslash H$ and $d :=
\sum_{i=1}^r 1/|L_i|$, then the matrices for
$\overline{\omega}^{(2)}$ and $\overline{\mu}^{(2)}$ are
respectively $\squarematrix{1}{0}{d}{1}$ and
$\squarematrix{1}{0}{-d}{1}$ by Example~\ref{exa:examples of the
dimension}~\ref{exa:von_Neumann_dimension_and_permutation_modules}
and
Lemma~\ref{lem:L2-rank_and_bisets}~\ref{lem:L2-rank_and_bisets:diagram}.
We conclude from
Theorem~\ref{the:The_finiteness_obstruction_and_the_(functorial)_L2-Euler_characteristic},
the definitions of $\chi_f(\Gamma(S);\IQ)$ and
$\chi(\Gamma(S);\IQ)$, and Theorem~\ref{the:chi_f_determines_chi}
that
\begin{eqnarray*}
\chi^{(2)}_f(\Gamma(S)) & = & \bigl(\chi^{(2)}(H) - d \cdot \chi^{(2)}(G), \chi^{(2)}(G)\bigr);
\\
\chi^{(2)}(\Gamma(S)) & = & \chi^{(2)}(H) + (1-d) \cdot \chi^{(2)}(G);
\\
\chi_f(\Gamma(S);\IQ) & = &\bigl( 1 - |G\backslash S/H|,1);
\\
\chi(\Gamma(S);\IQ) & = & 2 - |G\backslash S/H|;
\\
\chi(B\Gamma(S);\IQ) & = & 2 - |G\backslash S/H|.
\end{eqnarray*}

The situation above simplifies considerably in the finite case.
\begin{example}[Finite $G$-$H$-biset for finite groups $H$ and $G$]
  \label{exa:H-G-biset}
  Let $H$ and $G$ be finite groups and $S$ a
  finite $G$-$H$-biset. Then the category $\Gamma(S)$
  is a finite EI-category.  We conclude  from
  Theorem~\ref{the:The_finiteness_obstruction_and_Moebius_inversion_finite_EI-categories}
  that $\Gamma(S)$ is of type (FP$_\IQ$). The image of the finiteness
  obstruction under the isomorphism
$$S \colon K_0(\IQ\Gamma(S))  \xrightarrow{\cong} K_0(\IQ H) \oplus K_0(\IQ G)$$
is the element $\mu([\IQ],[\IQ])=\bigl([\IQ] - [\IQ G\backslash S],[\IQ]\bigr)$, and
\begin{eqnarray*}
  \chi^{(2)}_f(\Gamma(S)) & = & \biggl(\frac{1}{|H|} -
  \frac{|G\backslash S|}{|H|}, \frac{1}{|G|}\biggr);
  \\
  \chi^{(2)}(\Gamma(S)) & = &
  \frac{1}{|H|} +   \frac{1}{|G|} - \frac{|G\backslash S|}{|H|};
  \\
  \chi_f(\Gamma(S);\IQ) & = & (1 - |G\backslash S/H|,1);
  \\
 \chi(\Gamma(S);\IQ) & = & 2 - |G\backslash
S/H|;
\\
\chi(B\Gamma(S);\IQ) & = & 2 - |G\backslash S/H|.
\end{eqnarray*}
since
$\dim_{\caln(H)}\bigl(\IC(G\backslash S) \otimes_{\IC H}
\caln(H)\bigr) =  \frac{|G\backslash S|}{|H|}$ by
Example~\ref{exa:examples of the
dimension}~\ref{exa:von_Neumann_dimension_and_permutation_modules}.
If $S$ is free as a left $G$-set, or, equivalently, if $\Gamma(S)$
is free, we obtain
\begin{eqnarray*}
  \chi^{(2)}_f(\Gamma(S))
  & = &
  \biggl(\frac{1}{|H|} -  \frac{|S|}{|G| \cdot |H|}, \frac{1}{|G|}\biggr);
  \\
  \chi^{(2)}(\Gamma(S))
  & = &
  \frac{1}{|H|} +   \frac{1}{|G|} - \frac{|S|}{|G| \cdot |H|},
\end{eqnarray*}
since in this case $ \frac{|G\backslash S|}{|H|} = \frac{|S|}{|G| \cdot |H|}$.
\end{example}

\subsection{The passage to the opposite category}
\label{subsec:The_passage_to_the_opposite_category}

In this subsection we want to compare the invariants of $\Gamma$
with the invariants of the opposite category $\Gamma^{\op}$. The
categories $\Gamma$ and $\Gamma^{\op}$ can be distinguished by $o,
\chi_f, \chi^{(2)}_f,$ and $\chi^{(2)}$.

In general $\Gamma$ and $\Gamma^{\op}$ behave very differently. It
may happen that $\Gamma$ is of type (FP$_R$) but $\Gamma^{\op}$ is
not of type (FP$_R$) or that both $\Gamma$ and $\Gamma^{\op}$ are of
type (FP$_R$), but their finiteness obstructions and functorial
Euler characteristics are very different. This is illustrated by the
following example.

\begin{example}\label{exa:Gamma_versus_Gamma_op}
  Let $G$ be a group. Let $S$ be the $G$-$\{1\}$ biset consisting of precisely one element.
  Let $\Gamma(S)$ be the associated EI-category of
  Subsection~\ref{subsec:The_example_of_a_biset}. It has two objects $x$ and
  $y$. The sets $\mor_{\Gamma(S)}(x,x)$ and $\mor_{\Gamma(S)}(x,y)$ each contain
  precisely one element, the set
  $\mor_{\Gamma(S)}(y,y)$ is equal to $G$, and the set $\mor_{\Gamma(S)}(y,x)$
  is empty.  The category $\Gamma(S)$ is quasi-finite in the sense of
  Definition~\ref{def:finite_quasi-finite_and_free_category} and also directly
  finite in the sense of Definition~\ref{def:directly_finite_category}.  We
  conclude from Lemma~\ref{lem:finite_homological_dimension}~%
~\ref{lem:finite_homological_dimension:criterion_for_FP} that $\Gamma(S)$ is of
  type (FP$_{\IQ}$) if and only if the group $G$ is of type (FP$_{\IQ}$), i.e.,
  the trivial $\IQ G$-module $\IQ$ possesses a finite projective $\IQ
  G$-resolution.

  Now suppose that $G$ is of type (FP$_{\IQ}$). Then the trivial $\IQ G$-module $\IQ$
  has a finite projective $\IQ G$-resolution and defines an element $[\IQ]=o(G;\IQ)
  \in K_0(\IQ G)$.  Let $\alpha \colon K_0(\IQ G) \to K_0(\IQ)$ be the homomorphism which sends
  $[P]$ to $[P\otimes_{\IQ G} \IQ]$.  We conclude from
  Theorem~\ref{the:Two_splittings_and_the_K-theoretic_Moebius_inversion}~%
\ref{the:Two_splittings_and_the_K-theoretic_Moebius_inversion:o}
  that the finiteness obstruction $o(\Gamma;\IQ)$ is sent under the
  isomorphism of~\eqref{S_colon_SplitK_0_toK_0}
  $$S_{\IQ\Gamma(S)} \colon K_0(\IQ \Gamma(S))
  \xrightarrow{\cong} K_0(\IQ) \oplus K_0(\IQ G)$$
  to $\mu([\IQ],[\IQ])=([\IQ] - \alpha([\IQ]),[\IQ])$,

  This implies
  $$\begin{array}{lccl}
  \chi^{(2)}_f(\Gamma(S))
  & = &
  \bigl(1 - \chi(BG), \chi^{(2)}(G)\bigr)
  & \in U(\Gamma(S)) \otimes_{\IZ} \IQ = \IQ \oplus \IQ;
  \\
  \chi^{(2)}(\Gamma(S))
  & = &
  1 - \chi(BG) + \chi^{(2)}(G)
  & \in \IQ;
  \\
  \chi_f(\Gamma(S);\IQ)
  & = &
  \bigl(1 - \chi(BG), \chi(BG)\bigr)
  & \in U(\Gamma(S)) = \IZ \oplus \IZ;
  \\
   \chi(\Gamma(S);\IQ)
  & = &
  1 & \in \IZ;
  \\
  \chi(B\Gamma(S);\IQ) & =  & 1 & \in \IZ.
  \end{array}$$
  If $G$ satisfies condition (I) of~\eqref{con:condition_(I)} or
  $G$ is of type (FF$_\IQ$), then we conclude from
  Lemma~\ref{lem:rank_is_rankL2_for_fin_free_modules}~\ref{lem:rank_is_rankL2_for_fin_free_modules:rank_is_L2-rank}
  $$\chi^{(2)}(\Gamma(S)) = 1.$$

  The opposite category $\Gamma(S)^{\op} = \Gamma(S^{\op})$ has a terminal object,
  namely $x$. Hence it is always of type (FP$_\IQ$)
  and its finiteness obstruction $o(\Gamma(S)^{\op};\IQ)$ is sent under the isomorphism
  of~\eqref{S_colon_SplitK_0_toK_0}
  $$S_{\IQ\Gamma(S)^{\op}} \colon K_0(\IQ\Gamma(S)^{\op})
  \xrightarrow{\cong} K_0(\IQ) \oplus K_0(\IQ G)$$
  to $\mu([\IQ],[\IQ])=([\IQ],0)$.

  This implies
  $$\begin{array}{lccl}
  \chi^{(2)}_f(\Gamma(S)^{\op})
  & = &
  \bigl(1,0\bigr)
  & \in U(\Gamma(S)^{\op}) \otimes_{\IZ} \IQ = \IQ \oplus \IQ;
  \\
  \chi^{(2)}(\Gamma(S)^{\op})
  & = &
  1
  & \in \IQ;
  \\
  \chi_f(\Gamma(S)^{\op};\IQ)
  & = &
  (1,0)
  & \in U(\Gamma(S)^{\op}) = \IZ \oplus \IZ;
  \\
  \chi(\Gamma(S)^{\op};\IQ)
  & = &
  1 & \in \IZ;
  \\
  \chi(B\Gamma(S)^{\op};\IQ) & = & 1 & \in \IZ.
  \end{array}$$
Notice that all the results for $\Gamma(S)$ depend on $G$, whereas
the results for $\Gamma(S)^{\op}$ are all independent of $G$. So for
example, if $G$ is not of type (FP$_\IQ$), then $\Gamma(S)$ is not
of type (FP$_\IQ$), while $\Gamma(S)^{\op}$ is of type (FP$_\IQ$).
\end{example}

\subsection{The passage to the opposite category for finite EI-categories}
\label{subsec:The_passage_to_the_opposite_category_for_finite_EI-categories}

One can say more about the passage from $\Gamma$ to $\Gamma^{\op}$ in the special
case where $\Gamma$ is a finite EI-category.  Let $R$ be a commutative ring.
Given an $R$-module $M$, denote by $M^* := \hom_R(M,R)$ its dual
$R$-module.  Notice that $M^*$ is again an $R$-module since $R$ is
commutative.  This defines a contravariant functor
$$\ast_R \colon \MOD\text{-}R \to \MOD\text{-}R.$$
There is a natural $R$-homomorphism $I(M) \colon M \to (M^*)^*$ which
sends $m \in M$ to $M^* \to R, \; \phi \mapsto \phi(m)$. It is an
isomorphism if $M$ is a finitely generated projective $R$-module.

We obtain a functor
$$\ast_{R\Gamma} \colon \MOD\text{-}R\Gamma \to \MOD\text{-}R\Gamma^{\op}$$
which sends a contravariant $R\Gamma$-module $P$ to the contravariant
$R\Gamma^{\op}$-module, or equivalently, covariant $R\Gamma$-module
$P^*$ given by the composite $\Gamma \xrightarrow{P} \MOD\text{-}R
\xrightarrow{\ast} \MOD\text{-}R$.  The functor $\ast_{R\Gamma}$ is
exact when restricted to $R\Gamma$-modules $M$ for which $M(x)$ is a
finitely generated projective $R$-module for every object $x \in \ob(\Gamma)$. Let
$M$ be an $R\Gamma$-module such that $M(x)$ is a finitely generated
projective $R$-module for every object $x \in \ob(\Gamma)$. Then $M^*$
is an $R\Gamma^{\op}$-module such that $M(x)$ is a finitely generated
projective $R$-module for every object $x \in \ob(\Gamma^{\op})$ and
there is a natural isomorphism of $R\Gamma$-modules $M
\xrightarrow{\cong} (M^*)^*$.

Now assume that the order of the automorphism group of every object in
$\Gamma$ is invertible in $R$. Then an $R\Gamma$-module $M$, for which
the $R$-module $M(x)$ possesses a finite projective $R$-resolution for every object
$x\in \ob(\Gamma)$, possesses a finite projective $R\Gamma$-resolution by
Lemma~\ref{lem:finite_homological_dimension}~%
\ref{lem:finite_homological_dimension:finite_Gamma}.  Hence we obtain a
well-defined homomorphism
\begin{eqnarray}
  & \ast_{R\Gamma} \colon K_0(R\Gamma) \to K_0(R\Gamma^{\op}), \quad [P] \mapsto [P^*] &
  \label{ast_colon_K_0(RGamma)_to_K_0(RGamma_op)}
\end{eqnarray}
The functor $\ast_{R\Gamma}$ sends the constant $R\Gamma$-module
$\underline{R}$ to the constant $R\Gamma^{\op}$-module
$\underline{R}$. We conclude:

\begin{lemma} \label{ast_RGamma_on_K_0}
  Let $\Gamma$ be a finite $EI$-category. Let $R$ be a commutative ring such that the order of
  the automorphism group of every object in $\Gamma$ is invertible in
  $R$.

  \begin{enumerate}
  \item \label{ast_RGamma_on_K_0:ast_iso} The map
  of~\eqref{ast_colon_K_0(RGamma)_to_K_0(RGamma_op)}
  $$\ast_{R\Gamma} \colon K_0(R\Gamma) \to K_0(R\Gamma^{\op})$$
  is bijective, an inverse is
  $$\ast_{R\Gamma^{\op}} \colon K_0(R\Gamma^{\op}) \to K_0(R\Gamma);$$

  \item \label{ast_RGamma_on_K_0:ast_and_o} Both $\Gamma$
  and $\Gamma^{\op}$ are of type (FP$_R$) and
  $$\ast_{R\Gamma}\bigl(o(\Gamma;R)\bigr) = o(\Gamma^{\op};R).$$

\end{enumerate}
\end{lemma}

The map $\ast_{R\Gamma}$ is rather complicated as the next result shows.
\begin{lemma}\label{lem:ast_and_split}
  Let $\Gamma$ be a finite EI-category. Let $R$ be a commutative ring
  such that the order of
  the automorphism group of every object in $\Gamma$ is invertible in
  $R$. Then the following diagram commutes
  $$\xymatrix{
  K_0(R\Gamma) \ar[r]^{\ast_{R\Gamma}}_{\cong} \ar[d]_{S_{R\Gamma}}^{\cong}
  &
   K_0(R\Gamma^{\op})  \ar[d]^{S_{R\Gamma^{\op}}}_{\cong}
  \\
  \Split K_0(R\Gamma) \ar[r]_{\nu}^{\cong}
  &
  \Split K_0(R\Gamma^{\op})
  }$$
  Here $S_{R\Gamma}$ and $S_{R\Gamma^{\op}}$ are the homomorphisms defined
  in~\eqref{S_colon_SplitK_0_toK_0} which are isomorphisms by
  Theorem~\ref{the_splitting_of_K-theory_for_EI_categories},
  the isomorphism $\ast_{R\Gamma}$ has been defined
  in~\eqref{ast_colon_K_0(RGamma)_to_K_0(RGamma_op)} and the isomorphism
  $\nu$ is the composite
  $$\nu \colon \Split K_0(R\Gamma) \xrightarrow{\omega_{R\Gamma}}
  \Split K_0(R\Gamma) \xrightarrow{D} \Split K_0(R\Gamma^{\op}) \xrightarrow{\mu_{R\Gamma^{\op}}}
  \Split K_0(R\Gamma^{\op}),$$
  where $\omega_{R\Gamma}$ is the isomorphism defined
in~\eqref{omega_splitK_0_to_splitK_0} for
  $\Gamma$, $\mu_{R\Gamma^{\op}}$ is the isomorphism defined
in~\eqref{mu_splitK_0_to_splitK_0} for
  $\Gamma^{\op}$ and $D$ is given by the direct sum of the isomorphisms
  $K_0(R\aut_{\Gamma}(x)) \xrightarrow{\cong} K_0(R\aut_{\Gamma^{\op}}(x))$
  sending the class of the finitely generated projective $R\aut_{\Gamma}(x)$-module $P$ to
  the class of the finitely generated projective $R\aut_{\Gamma^{\op}}(x)$-module $P^*$.
\end{lemma}
\begin{proof}
Consider the following diagram.
$$\xymatrix{
  & K_0(R\Gamma) \ar[r]^{\ast_{R\Gamma}} \ar[dl]_{S_{R\Gamma}}
  \ar[d]^{\Res_{R\Gamma}}
  &
   K_0(R\Gamma^{\op})  \ar[dr]^{S_{R\Gamma^{\op}}} \ar[d]_{\Res_{R\Gamma^{\op}}}
  &
  \\
  \Split K_0(R\Gamma) \ar[r]_{\omega_{R\Gamma}}
  &
  \Split K_0(R\Gamma) \ar[r]_{D}
  & \Split K_0(R\Gamma^{\op}) \ar[r]_{\mu_{R\Gamma^{\op}}}
  & \Split K_0(R\Gamma^{\op})
  }$$
The left and right triangles commute by
Theorem~\ref{the:Two_splittings_and_the_K-theoretic_Moebius_inversion}
and the middle square commutes from the definitions, so the entire diagram commutes.
\end{proof}

\begin{lemma}\label{lem:ast_and_split_and_U}
  Let $\Gamma$ be a finite EI-category. Suppose that both $\Gamma$ and
  $\Gamma^{\op}$ are free in the sense of Definition~\ref{def:finite_quasi-finite_and_free_category}.
  Then the following diagram commutes.
  $$\xymatrix{
  K_0(\IQ\Gamma) \ar[r]^{\ast_{\IQ\Gamma}}_{\cong} \ar[d]_{S_{\IQ\Gamma}}^{\cong}
  &
  K_0(\IQ\Gamma^{\op})  \ar[d]^{S_{\IQ\Gamma^{\op}}}_{\cong}
  \\
  \Split K_0(\IQ\Gamma) \ar[r]^{\nu}_{\cong} \ar[d]_{\rk^{(2)}_{\Gamma}}
  &
  \Split K_0(\IQ\Gamma^{\op}) \ar[d]^{\rk^{(2)}_{\Gamma^{\op}}}
  \\
  U(\Gamma) \otimes_{\IZ} \IQ \ar[r]^{\overline{\nu}^{(2)}}_{\cong}
  &
  U(\Gamma) \otimes_{\IZ} \IQ
  }$$
  Here the upper square is taken from Lemma~\ref{lem:ast_and_split},
  the maps $\rk^{(2)}_{\Gamma}$ and $\rk^{(2)}_{\Gamma^{\op}}$ have been defined
  in~\eqref{rk(2)_K_o_to_U}, and the isomorphism
  $\overline{\nu}^{(2)}$ is defined to be
  $\overline{\mu}^{(2)}_{\Gamma^{\op}} \circ \overline{\omega}^{(2)}_{\Gamma}$,
  where $ \overline{\omega}^{(2)}_{\Gamma}$ is the isomorphism defined
  in~\eqref{overlineomega(2)} for $\Gamma$ and $\overline{\mu}^{(2)}_{\Gamma^{\op}}$ is the
  isomorphism defined
  in~\eqref{overlinemu(2)} for $\Gamma^{\op}$.
\end{lemma}
\begin{proof}
  This follows from
  Theorem~\ref{the:K_theoretic-Moebius_inversion_and_L2-rank},
  Lemma~\ref{lem:ast_and_split}, and the easy to verify fact
  that the following diagram commutes
  for the homomorphism $D$ appearing in Lemma~\ref{lem:ast_and_split}.
\[\xymatrix{
  \Split K_0(\IQ\Gamma) \ar[r]^{D} \ar[d]_{\rk^{(2)}_{\Gamma}}
  &
  \Split K_0(\IQ\Gamma^{\op}) \ar[d]^{\rk^{(2)}_{\Gamma^{\op}}}
  \\
  U(\Gamma) \otimes_{\IZ} \IQ \ar[r]_{\id}
  &
  U(\Gamma) \otimes_{\IZ} \IQ
}\]
\end{proof}

\begin{example}[The isomorphism $\ast$ for a
finite $G$-$H$-biset for finite groups $H$ and $G$]
  \label{exa:H-G-biset_and_ast}
  Let $H$ and $G$ be finite groups and $S$ a finite $G$-$H$-biset. We have
  defined a finite EI-category $\Gamma(S)$ in
  Subsection~\ref{subsec:The_example_of_a_biset} and
  Example~\ref{exa:H-G-biset}. We conclude from
  Subsection~\ref{subsec:The_example_of_a_biset} that the commutative diagram
  appearing in Lemma~\ref{lem:ast_and_split} can be identified for $\Gamma(S)$
  with
$$\xymatrix{
  K_0(\IQ\Gamma(S)) \ar[r]^{\ast_{\IQ\Gamma(S)}}_{\cong}
  \ar[d]_{S_{\IQ\Gamma(S)}}^{\cong} & K_0(\IQ\Gamma(S)^{\op})
  \ar[d]^{S_{\IQ\Gamma(S)^{\op}}}_{\cong}
  \\
  K_0(\IQ H) \oplus K_0(\IQ G) \ar[r]_-{\nu}^-{\cong} & K_0(\IQ H^{\op}) \oplus
  K_0(\IQ G^{\op}).  }$$ By the calculation for $\omega$ and $\mu$ in
Subsection~\ref{subsec:The_example_of_a_biset}, the homomorphism $\nu$ sends
$\bigl([P],[Q]\bigr)$ to
$$\bigl([P^*] + [(Q \otimes_{\IQ G} \IQ S)^*],
[Q^*] - [P^* \otimes_{\IQ H^{\op}} \IQ S^{\op}] - [(Q \otimes_{\IQ G} \IQ S)^*
\otimes_{\IQ H^{\op}}\IQ S^{\op}]\bigr)$$ (recall that the roles of $G^{\op}$
and $H^{\op}$ are switched in the formula for $\mu_{\IQ\Gamma^{\op}}$).

Now suppose that both $\Gamma(S)$ and $\Gamma(S)^{\op}$ are free, or, equivalently,
that $G$ acts freely from the left on $S$ and $H$ acts freely from the right on $S$.
Then the commutative diagram appearing in Lemma~\ref{lem:ast_and_split_and_U} can be
identified with
$$\xymatrix{
  K_0(\IQ\Gamma(S)) \ar[r]^{\ast_{\IQ\Gamma(S)}}_{\cong} \ar[d]_{S_{\IQ\Gamma(S)}}^{\cong}
  &
   K_0(\IQ\Gamma(S)^{\op})  \ar[d]^{S_{\IQ\Gamma(S)^{\op}}}_{\cong}
  \\
  K_0(\IQ H) \oplus K_0(\IQ G) \ar[r]^-{\nu}_-{\cong} \ar[d]_{\rk^{(2)}_{\Gamma(S)}}
  &
  K_0(\IQ H^{\op}) \oplus K_0(\IQ G^{\op}) \ar[d]^{\rk^{(2)}_{\Gamma(S)^{\op}}}
  \\
  \IQ \oplus \IQ \ar[r]_{\overline{\nu}^{(2)}} &
   \IQ \oplus \IQ
  }$$
where $\overline{\nu}^{(2)}$ is given by the matrix
$\left(\begin{matrix} 1 & \frac{|S|}{|H|}
\\ - \frac{|S|}{|G|} & 1 - \frac{|S|^2}{|H|\cdot |G|}
   \end{matrix}\right)$.
\end{example}


\typeout{---------- Section 10: Comparison with the invariants of
Baez--Dolan and Leinster
-------------------}
\section{Comparison with the invariants of Baez--Dolan and Leinster}
\label{sec:Comparison_with_Leinsters_invariant}

In this section we compare our invariants with the groupoid
cardinality of Baez--Dolan~\cite{Baez-Dolan(2001)} and the Euler
characteristic of Leinster~\cite{Leinster(2008)}. If $\Gamma$ is a
skeletal, finite, free EI-category, then $\Gamma$ is of type
(FP$_\IQ$) and of type ($L^2$), and
Leinster's Euler characteristic coincides with the $L^2$-Euler
characteristic. However, if we leave out the freeness hypothesis,
then Leinster's Euler characteristic can very well be different from
the $L^2$-Euler characteristic, see Remark
\ref{rem_free_is_necessary}.

\subsection{Comparison with the groupoid cardinality of Baez--Dolan} Baez--Dolan define in~\cite{Baez-Dolan(2001)}
the \emph{groupoid cardinality} of a groupoid $\Gamma$ to be
$$\sum_{\overline{x} \in \iso(\Gamma)} \frac{1}{|\aut(x)|},$$
provided this sum converges. In other words, the groupoid
cardinality is the count of the isomorphism classes of objects
inversely weighted by the size of their symmetry groups. This agrees
with the $L^2$-Euler characteristic of such groupoids as seen in
Example~\ref{exa:chi_f(calg)L2_for_finite_groupoids}.

\subsection{Review of Leinster's Euler characteristic}
\label{subsec:Review_of_Leinster_Euler_characteristic}

We briefly review the Euler characteristic due to
Leinster~\cite{Leinster(2008)}. Let $\Gamma$ be a finite category
(see Definition~\ref{def:finite_quasi-finite_and_free_category}). A
\emph{weighting} on $\Gamma$ is a function $k^{\bullet} \colon
\ob(\Gamma) \to \IQ$ such that for all objects $x \in \iso(\Gamma)$
we have $\sum_{y \in \ob(\Gamma)} |\mor(x,y)| \cdot k^y = 1$. A
\emph{coweighting} $k_{\bullet}$ on $\Gamma$ is a weighting on
$\Gamma^{\op}$.

\begin{definition} \label{def:Leinsters_Euler_characteristic}
A finite category $\Gamma$ \emph{has an Euler characteristic in the
sense of Leinster} if it has a weighting and a coweighting. Its
\emph{Euler characteristic in the sense of Leinster} is then defined
as
$$\chi_L(\Gamma) := \sum_{y \in \ob(\Gamma)} k^y = \sum_{x \in \ob(\Gamma)} k_x$$
for any choice of weighting $k^{\bullet}$ or coweighting
$k_{\bullet}$.
\end{definition}

This is indeed independent of the choice of the weighting and the
coweighting. In particular we get $\chi_L(\Gamma) =
\chi_L(\Gamma^{\op})$.

\begin{remark}\label{remark:what_Leinsters-Euler-characteristic_sees} Leinster's
    Euler characteristic can only be defined if the category $\Gamma$ is
    finite and depends only on the set of objects $\ob(\Gamma)$ and the
    orders $|\mor(x,y)|$ for $x,y \in \ob(\Gamma)$. This is different
    from the other invariants such as the finiteness obstruction. For
    instance $\chi_L$ does not distinguish between the category $\Gamma$
    appearing in
    Example~\ref{exa:finiteness_obstruction_of_projection_category} and
    the groupoid $\widehat{\IZ/2}$, whereas the finiteness obstructions
    and the $L^2$-Euler characteristic do.
\end{remark}

\subsection{Finite, free, skeletal, EI-categories and comparison of
$\chi^{(2)}$ with $\chi_L$}
\label{subsec:Finite_free_skeletal_EI-categories}

\begin{lemma} \label{lem:chi(2)_and_chi}
Let $\Gamma$ be a finite, free, EI-category which is skeletal, i.e.,
two isomorphic objects are already equal.

Then $\Gamma$ is of type (FP$_\IC$)
and of type ($L^2$), and has an Euler characteristic in the sense of
Leinster. We get for the $L^2$-Euler characteristic
$\chi^{(2)}(\Gamma)$ of
Definition~\ref{def:L2-Euler_characteristic_of_a_category} and
Leinster's Euler characteristic $\chi_L(\Gamma)$ of
Definition~\ref{def:Leinsters_Euler_characteristic}
$$\chi^{(2)}(\Gamma) = \chi_L(\Gamma).$$
\end{lemma}
\begin{proof}
By \cite[Lemma~1.3 and Theorem~1.4]{Leinster(2008)} the category
$\Gamma^{\op}$ has a M\"obius inversion, i.e., the homomorphism
$$\omega_L \colon U(\Gamma) \otimes_{\IZ} \IQ \to U(\Gamma)\otimes_{\IZ} \IQ$$
given by the matrix
$$\left(|\mor_{\Gamma}(y,x)|\right)_{x,y \in \ob(\Gamma)}$$
is bijective, and has an Euler characteristic in the sense of
Leinster. Then by definition
$$\chi_L(\Gamma) = \chi_L(\Gamma^{\op}) = \sum_{x \in \ob(\Gamma)} k_x$$
for any element $k_{\bullet} \in U(\Gamma) \otimes_{\IZ}\IQ$ such
that $\omega_L(k_{\bullet}) $ is the element $\overline{1} \in
U(\Gamma)$ which assigns $1$ to every element in $\ob(\Gamma)$.

We conclude from
Theorem~\ref{the:The_finiteness_obstruction_and_Moebius_inversion_finite_EI-categories}
that $\Gamma$ is of type (FP$_\IC$)
and hence of type ($L^2$). Hence it remains to show
$$\omega_L\bigl(\chi^{(2)}_f(\Gamma)\bigr) = \overline{1} \quad \in U(\Gamma),$$
since by definition $\chi^{(2)}(\Gamma) = \sum_{x \in \ob(\Gamma)}
\chi^{(2)}_f(\Gamma)(x)$.

Since $\aut(y)$ is finite,
Example~\ref{exa:examples of the dimension}~\ref{exa:von_Neumann_dimension_and_permutation_modules}
implies
$$\dim_{\caln(y)}\bigl(\IC\mor(y,x) \otimes_{\IC[y]} \caln(y)\bigr)
= \frac{|\mor(y,x)|}{|\aut(y)|}$$ for every $x,y \in \ob(\Gamma)$.
Hence the homomorphism $\omega_L$ agrees with the composite $D \circ
\overline{\omega}^{(2)}$, where $\overline{\omega}^{(2)}$ is defined
in~\eqref{overlineomega(2)} and $D$ is the isomorphism given by the
diagonal matrix with entry $|\aut(y)|$ at $(y,y)$ for $y \in
\ob(\Gamma)$. Since $D \circ \overline{\omega}^{(2)}$ maps
$\chi^{(2)}_f(\Gamma)$ to $\overline{1}$ because of
Theorem~\ref{the:The_finiteness_obstruction_and_the_(functorial)_L2-Euler_characteristic}~%
\ref{the:The_finiteness_obstruction_and_the_(functorial)_L2-Euler_characteristic:chi(2)_(I)}
and because of $\chi^{(2)}(\aut(x)) = 1/|\aut(x)|$,
Lemma~\ref{lem:chi(2)_and_chi} follows. We need $\Gamma$ to be free in the sense of Definition \ref{def:finite_quasi-finite_and_free_category} in order to apply Theorem ~\ref{the:The_finiteness_obstruction_and_the_(functorial)_L2-Euler_characteristic}~%
\ref{the:The_finiteness_obstruction_and_the_(functorial)_L2-Euler_characteristic:chi(2)_(I)}.
\end{proof}

\begin{remark} \label{rem_free_is_necessary} The condition in
  Lemma~\ref{lem:chi(2)_and_chi} that $\Gamma$ is free is necessary as
  the following example shows. Let $H$ and $G$ be finite groups and
  $S$ be a finite $G$-$H$-biset. Let $\Gamma(S)$ be the associated
  finite EI-category of Example~\ref{exa:H-G-biset}. We conclude from
  Example~\ref{exa:H-G-biset} and the definition of
  $\chi_L(\Gamma(S))$ that
  \begin{eqnarray*}
    \chi^{(2)}(\Gamma(S)) & = &
    \frac{1}{|H|} +   \frac{1}{|G|} - \frac{|G\backslash S|}{|H|};
    \\
    \chi(\Gamma(S)) & = & 2 - |G\backslash S/H|;
    \\
    \chi(B\Gamma(S)) & = & 2 - |G\backslash S/H|;
    \\
    \chi_L(\Gamma(S)) & = & \frac{1}{|H|} +   \frac{1}{|G|} - \frac{|S|}{|G|\cdot |H|}.
  \end{eqnarray*}
  Hence $\chi^{(2)}(\Gamma(S)) = \chi_L(\Gamma(S))$ holds if and only
  if $|G\backslash S| = \frac{|S|}{|G|}$. The latter is equivalent to
  the condition that $\Gamma(S)$ is free.

  Notice that $\chi(\Gamma(S))$ and $\chi(B\Gamma(S))$ are always integers and are in general
  different from both $\chi^{(2)}(\Gamma(S))$ and
  $\chi_L(\Gamma(S))$.
\end{remark}

\begin{remark}[Homotopy colimit formula]
In~\cite{FioreLueckSauerHoColim(2009)}, we prove the compatability
of various Euler characteristics of categories with homotopy
colimits. There we compare our homotopy colimit results with
Leinster's results on Grothendieck fibrations.
\end{remark}

\subsection{Passage to the opposite category and initial and terminal objects}
\label{subsec:Passage_from_Gamma_to_Gammaop}

Leinster's Euler characteristic $\chi_L(\Gamma)$ and the topological
Euler characteristic $\chi(B\Gamma)$ do not see a difference between
$\Gamma$ and $\Gamma^{\op}$. We have discussed in detail in
Subsection~\ref{subsec:The_passage_to_the_opposite_category} that
$\Gamma$ and $\Gamma^{\op}$ can be distinguished by the finiteness
obstruction $o(\Gamma;R)$, the functorial Euler characteristic
$\chi_f(\Gamma;R)$, the functorial $L^2$-Euler characteristic
$\chi_f^{(2)}(\Gamma)$, and the $L^2$-Euler characteristic
$\chi^{(2)}(\Gamma)$.

Suppose that $\Gamma$ has a terminal object $x$. Let $i \colon
\{\ast\} \to \Gamma$ be the inclusion of the trivial category with
value $x$. Then the finiteness obstruction is the image of $[R]$
under $i_* \colon K_0(R) \to K_0(R\Gamma)$ by
Example~\ref{exa:terminal_object}. The  functorial Euler
characteristic $\chi_f(\Gamma;R) \in U(\Gamma)$ and the functorial
$L^{2}$-Euler characteristic $\chi^{(2)}_f(\Gamma) \in
U^{(1)}(\Gamma)$ agree and are given by the element $1 \cdot
\overline{x}$. The Euler characteristic $\chi(\Gamma;R)$, the
$L^{2}$-Euler characteristic $\chi^{(2)}(\Gamma) \in
U^{(1)}(\Gamma)$, and topological Euler characteristic
$\chi(B\Gamma;R)$ are all equal to $1$. Since $\Gamma$ has a
terminal object, it admits a weighting, see Leinster~\cite[Example
1.11.c]{Leinster(2008)}. If $\Gamma$ additionally admits a
coweighting, then Leinster's Euler characteristic $\chi_L(\Gamma)$
is equal to $1$.

If $\Gamma$ has an initial object, we cannot predict the values of
$o(\Gamma;R)$, $\chi_f(\Gamma;R)$, $\chi^{(2)}_f(\Gamma)$, and
$\chi^{(2)}(\Gamma)$ in general, as the results in
Subsections~\ref{subsec:The_example_of_a_biset}
and~\ref{subsec:The_passage_to_the_opposite_category} illustrate. In
particular, $\chi^{(2)}(\Gamma)$ is not necessarily $1$ if $\Gamma$
has an initial object. For instance, Example~\ref{exa:H-G-biset}
yields for $H = 1$, $S = \{\ast\}$, and $G$ any finite group
$\chi^{(2)}(\Gamma(S)) = 1/|G|$. If $\Gamma$ has an initial object,
then $\Gamma$ admits a coweighting. If $\Gamma$ additionally admits
a weighting, then Leinster's Euler characteristic $\chi_L(\Gamma)$
is equal to $1$.

The topological Euler characteristic $\chi(B\Gamma;R)$ is of course
equal to $1$ if $\Gamma$ has an initial or a terminal object.


\subsection{Relationship between weightings and free resolutions}

\begin{theorem}[Weighting from a free resolution]
\label{the:weighting_from_free_resolution} Let $\Gamma$ be a finite
category. Suppose that the constant $R \Gamma$-module
$\underline{R}$ admits a finite free resolution $P_*$. If $P_n$ is
free on the finite $\ob(\Gamma)$-set $C_n$, that is
\begin{equation}\label{the:weighting_from_free_resolution:equation}
P_n=B(C_n)=\bigoplus_{y \in \ob(\Gamma)} \bigoplus_{C^y_n}
R\mor(?,y),
\end{equation}
then the function $k^{\bullet} \colon \ob(\Gamma) \to \IQ$ defined
by
$$k^y:=\sum_{n \geq 0}(-1)^n \cdot |C^y_n|$$
is a weighting on $\Gamma$.
\end{theorem}
\begin{proof}
At each object $x$ of $\Gamma$, the $R$-chain complex $P_*(x)$ has
Euler characteristic 1, since it is a resolution of $R$. Further,
calculating the Euler characteristic of $P_*(x)$ using
equation~\eqref{the:weighting_from_free_resolution:equation} yields
\begin{align*} 1  = \chi(P_*(x))  &= \sum_{n \geq 0}(-1)^n \rk_R P_n(x)
\\ &=\sum_{n \geq 0}(-1)^n \left(\sum_{y \in \ob(\Gamma)} |C_n^y|\cdot|\mor(x,y)|\right)
\\ &= \sum_{y\in \ob(\Gamma)}|\mor(x,y)| \left(\sum_{n \geq 0}(-1)^n |C_n^y|\right)
\\ &= \sum_{y\in \ob(\Gamma)}|\mor(x,y)|k^y.\qedhere
\end{align*}
\end{proof}

In~\cite{FioreLueckSauerHoColim(2009)}, we recall the
$\Gamma$-$CW$-complexes of~Davis--L\"uck~\cite{Davis-Lueck(1998)} in
the context of Euler characteristics and homotopy colimits.

\begin{corollary}[Construction of a weighting from a
finite $\Gamma$-$CW$-model for the classifying $\Gamma$-space]
\label{cor:weighting_from_finite_model} Let $\Gamma$ be a finite
category. Suppose that $\Gamma$ admits a finite $\Gamma$-$CW$-model
$X$ for the classifying $\Gamma$-space $E\Gamma$. Then the function
$k^{\bullet} \colon \ob(\Gamma) \to \IQ$ defined by
$$k^y:=\sum_{n \geq 0}(-1)^n(\text{number of $n$-cells of $X$ based at $y$})$$
is a weighting on $\Gamma$.
\end{corollary}
\begin{proof}
The composite of the cellular $R$-chain complex functor with $X$
is a finite free resolution of the constant $R \Gamma$-module
$\underline{R}$. The number of $n$-cells of $X$ based at $y$ is
$|C^y_n|$.
\end{proof}

\begin{remark}
We may think of $k^\bullet$ in
Corollary~\ref{cor:weighting_from_finite_model} as the
\emph{$\Gamma$-Euler characteristic of the $\Gamma$-$CW$-space $X$}.
If $R=\IC$ and $\Gamma$ is skeletal and directly finite, then the
function $k^{\bullet}$ is just $\chi_f(\Gamma;\IC) =
\chi^{(2)}_f(\Gamma)$ by
Lemma~\ref{lem:rank_RGamma_and_free_modules}~\ref{lem:rank_RGamma_and_free_modules_rank_determines_f.f.-modules}
and Theorem~\ref{the:coincidence}. The role of direct finiteness is
to guarantee that the splitting functors $S_x$ are defined.
\end{remark}

\begin{example}
Let $\Gamma=\{1 \leftarrow 0 \rightarrow 2\}$ be the category with
objects $0$, $1$, and $2$ and only two nontrivial morphisms, one
from $0$ to $1$ and one from $0$ to $2$. A finite
$\Gamma$-$CW$-model for $E\Gamma$ has two zero-cells $\mor(?,1)$ and
$\mor(?,2)$ and one $1$-cell $\mor(?,0) \times D^1$ whose attaching
map $\mor(?,0) \times S^0 \to \mor(?,1) \amalg\mor(?,2)$ is the
disjoint union of the canonical maps $\mor(?,0) \to \mor(?,1)$ and
$\mor(?,0)\to \mor(?,2)$. This finite model produces the weighting
$(k^0,k^1,k^2)=(-1,1,1)$ by
Corollary~\ref{cor:weighting_from_finite_model}. This is the same
weighting as Leinster~\cite[1.11.a]{Leinster(2008)}.
\end{example}

\begin{example}
Let $\Gamma=\{a \rightrightarrows b\}$ be the category consisting of
two objects and a single pair of parallel arrows between them. A
finite $\Gamma$-$CW$-model for $E\Gamma$ has a single 0-cell based
at $b$ and a single 1-cell based at $a$. The gluing map $\mor(-,a)
\times S^0 \to \mor(-,b)$ is induced by the two parallel arrows $a
\rightrightarrows b$.
Corollary~\ref{cor:weighting_from_finite_model} then produces the
weighting $(k^a,k^b)=(-1,1)$, the same weighting as
Leinster~\cite[3.4.b]{Leinster(2008)}.
\end{example}

\begin{example}
Let $\Gamma$ be the category with objects the non-empty subsets of
$[q]=\{0,1, \dots, q\}$ and a unique arrow $J \to K$ if and only if
$K \subseteq J$. In~\cite{FioreLueckSauerHoColim(2009)}, we
construct a finite $\Gamma$-$CW$-model with precisely one $|J|-1$
cell based at $J$ for each nonempty $J \subseteq [q]$. By
Corollary~\ref{cor:weighting_from_finite_model}, we obtain a
weighting $k^\bullet$ on $\Gamma$ by defining $k^J:=(-1)^{|J|-1}$.
This is the same weighting as Leinster~\cite[3.4.d]{Leinster(2008)}.
\end{example}

\begin{remark}
For a finite group $G$, there is no finite model. So it appears the
above method of finding the weighting does not work.  However, if we
use the $L^2$-rank, something similar does. Every finite group $G$
has a finite projective resolution of $\underline{\IQ}$, namely
$\underline{\IQ}$ itself. Then we obtain for the weighting
$$k^*=\sum_{n \geq 0}(-1)^n\dim_{\caln(G)}\underline{\IQ}_*=\dim_{\caln(G)}\underline{\IQ}=1/|G|,$$
precisely as by Leinster.
\end{remark}


\typeout{-- Section 8: The proper orbit category}
\section{The proper orbit category}
\label{sec:The_proper_orbit_category}

The principal virtue of the finiteness-obstruction approach to Euler
characteristics is the wide variety of examples and familiar notions
it encompasses. We have already seen the topological Euler
characteristic of a category and the classical $L^2$-Euler
characteristic of a group \cite[Chapter~7]{Lueck(2002)} as special
cases. We turn now to another special case: the equivariant Euler
characteristic of the classifying space $\eub{G}$ for proper
$G$-actions. Recall from Definition~\ref{def:orbit_category} that
the \emph{proper orbit category} $\uor{G}$ has as objects the
homogeneous spaces $G/H$ with $H$ a finite subgroup of $G$, and as
morphisms the $G$-equivariant maps. We have shown in
Lemma~\ref{lem:orbit_category_quasi-finite} that $\uor{G}$ is a
quasi-finite and free EI-category. We will explain in this section
that the finiteness obstructions and Euler characteristic notions
for $\Gamma=\uor{G}$ correspond to established notions in
equivariant topology for the classifying space $\eub{G}$ for proper
$G$-actions. This gives in particular the possibility to compute and
relate the invariants for $\uor{G}$ to more geometric notions.

In
Subsection~\ref{subsec:The_classifying_space_for_proper_G-actions}
we recall $G$-$CW$-complexes, the classifying space for proper
$G$-actions, and the relationship between equivariant invariants of
$\eub{G}$ and our category-theoretic invariants of $\uor{G}$. In
Subsection~\ref{subsec:The_Moebius_inversion_for_uor(G)} we discuss
M\"obius inversion for $\uor{G}$ in the case where $\eub{G}$ admits
a finite model. If $G_0$ is a subgroup of $G_1$ and $G_2$, then the
Euler characteristics of $\uor{G_1 \ast_{G_0} G_2}$ are computed
additively from those of $\uor{G_0}$, $\uor{G_1}$, and $\uor{G_2}$
in Subsection \ref{subsec:Additivity_of_o_and_chi_for_uor(G)}. In
Subsection
\ref{subsec:The_Burnside_integrality_relations_and_the_classical_Burnside_congruences}
we derive the Burnside congruences from an integrality condition
involving $(\overline{\mu}^{(2)}, \overline{\omega}^{(2)})$. We work
everything out explicitly for $G$ the infinite dihedral group in
Subsection~\ref{subsec:infinite_dihedral_group}. Fundamental
groupoids are considered in
Subsection~\ref{subsec:The_fundamental_category}.

\subsection{The classifying space for proper $G$-actions}
\label{subsec:The_classifying_space_for_proper_G-actions}

\begin{definition}[$G$-$CW$-complex] \label{def:G-CW-complex} A
  \emph{$G$-$CW$-complex $X$} is a $G$-space $X$ together with a
 filtration by $G$-spaces
  $X_{-1} = \emptyset \subseteq X_0 \subseteq X_1 \subseteq
  \ldots \subseteq X = \bigcup_{n \ge 0} X_n$
  such that $X = \colim_{n \to \infty} X_n$ and for each $n$ there is a $G$-pushout, that is,
a pushout in the category of $G$-spaces
\[
\comsquare{\coprod_{i \in I_n} G/H_i \times S^{n-1}}{\coprod_{i \in
    I_n} q^n_i}{X_{n-1}} {}{}{\coprod_{i \in I_n} G/H_i \times
  D^n}{\coprod_{i \in I_n} Q^n_i}{X_n.}
\]
\end{definition}

For more information about $G$-$CW$-complexes we refer
to~L\"uck~\cite[Chapters~1 and~2]{Lueck(1989)}. A $G$-$CW$-complex
is \emph{proper} if and only if all its isotropy groups are finite
(see~L\"uck~\cite[Theorem~1.23 on page~18]{Lueck(1989)}).

A $G$-$CW$-complex is \emph{finite}, i.e., is built out of finitely
many equivariant cells $G/H_i \times D^n$ if and only if it is
\emph{cocompact}, i.e., $G \backslash X$ is compact.  A
$G$-$CW$-complex $X$ is \emph{finitely dominated} if and only if
there exists a finite $G$-$CW$-complex $Y$ and $G$ maps $i \colon X
\to Y$ and $r \colon Y \to X$ with $r \circ i \simeq_G \id_X$.

\begin{definition}[Classifying space for proper $G$-actions]
  \label{def:classifying_space_for_proper_G-actions}
  A model for the \emph{classifying space for proper $G$-actions} is a
  $G$-$CW$-complex $\eub{G}$ such that the subspace of $H$-fixed points $\eub{G}^H$ is  contractible
  for every finite subgroup $H \subseteq G$ and is empty for every infinite subgroup $H
  \subseteq G$.
\end{definition}

For much more information about $\eub{G}$ than presented here we
refer the reader to the survey article \cite{Lueck(2005s)} of
L\"uck. We have $EG = \eub{G}$ if and only if $G$ is torsion-free.
We can choose $G/G$ as a model for $\eub{G}$ if and only if $G$ is
finite.

\begin{remark} \label{rem:universal_property_of_eub} The classifying
  space for proper $G$-actions has the following universal property.
  If $X$ is a proper $G$-$CW$-complex, then there is up to
  $G$-homotopy precisely one $G$-map from $X$ to $\eub{G}$.  In other
  words, a model for $\eub{G}$ is a terminal object in the
  $G$-homotopy category of proper $G$-$CW$-complexes. In particular,
  two models for $\eub{G}$ are $G$-homotopy equivalent.
\end{remark}

Recall from Notation~\ref{U(Gamma)_and_augmentation} that
$U(\Gamma):=\IZ \iso(\Gamma)$ for any category $\Gamma$.

\begin{definition}[Equivariant Euler characteristic]
  \label{def:equivariant_Euler_characteristic}
  Let $X$ be a finite $G$-$CW$-complex
  (see~Definition~\ref{def:G-CW-complex}).  The \emph{equivariant
  Euler characteristic of X}
  $$\chi^G(X) \in U(\uor{G}) $$
  is
  $$\chi^G(X) := \sum_{n \ge 0} (-1)^n \cdot \sum_{i \in I_n} \overline{G/H_i}$$
  for any choice of $G$-pushout appearing in Definition~\ref{def:G-CW-complex}.
\end{definition}

\begin{theorem}[The relation between $\eub{G}$ and $\uor{G}$]
  \label{the:uor(G)_and_eub(G)} \
  \begin{enumerate}
  \item \label{the:uor(G)_and_eub(G):(FF)} If there exists a finite
    $G$-$CW$-model for $\eub{G}$, then the EI-category $\uor{G}$ is of type (FF$_R$)
    for any ring $R$;

  \item \label{the:uor(G)_and_eub(G):(FP)} If there exists a finitely
    dominated $G$-$CW$-model for $\eub{G}$, then $\uor{G}$ is of type
    (FP$_R$) for any ring $R$;

  \item \label{the:uor(G)_and_eub(G):converse} Suppose that $G$
    contains only finitely many conjugacy classes of finite subgroups
    and for every finite subgroup $H \subset G$ its Weyl group $W_GH
    := N_GH/H$ is finitely presented. Suppose that $R = \IZ$.  Then
    the converses of assertions~\ref{the:uor(G)_and_eub(G):(FF)}
    and~\ref{the:uor(G)_and_eub(G):(FP)} are true;

  \item \label{the:uor(G)_and_eub(G):o_for_finitely_dominated} If
    $\eub{G}$ is a finitely dominated $G$-$CW$-complex, then the
    equivariant finiteness obstruction of L\"uck~\cite[Definition~14.4 on
    page~278]{Lueck(1989)} agrees with the finiteness obstruction
    $o(\uor{G};\IZ)$ of Definition~\ref{def:finiteness_obstruction_of_a_category};

  \item \label{the:uor(G)_and_eub(G):o_for_finite} Suppose that there is a
    finite $G$-$CW$-complex model for $\eub{G}$. Then its equivariant Euler characteristic
    $\chi^G(\eub{G}) \in U(\uor{G})$ agrees with the functorial Euler
    characteristic $\chi_f(\uor{G};\IZ)$ and
    the functorial $L^{(2)}$-Euler characteristic
    $\chi^{(2)}_f(\uor{G})$. Moreover, its
    finiteness obstruction $o(\uor{G};R)$ is the
    image of $\chi_f(\uor{G};\IZ)$ under the composite
    $$U(\uor{G}) \xrightarrow{\iota} K_0(\IZ \uor{G})
    \xrightarrow{c} K_0(R\uor{G})$$
    where $\iota$ has been defined in~\eqref{iota} and $c$ is the
    obvious change of coefficients homomorphism.
    \end{enumerate}
  \end{theorem}
\begin{proof}\ref{the:uor(G)_and_eub(G):(FF)} The cellular $\IZ \uor{G}$-chain
    complex $C_*(X)$ of a proper $G$-$CW$-complex $X$ sends $G/H$ to
    the cellular chain complex of the $CW$-complex $\map_G(G/H,X) = X^H$.
    It is always free, and it is finite free if and only if $X$ is
    finite (see~L\"uck~\cite[Section~18A]{Lueck(1989)}.

    Since $\eub{G}^H$ is contractible, the cellular $\IZ
    \uor{G}$-chain complex $C_*(\eub{G})$ is a free and hence
    projective resolution of the constant $\IZ \uor{G}$-module
    $\underline{R}$.
    \\[1mm]~\ref{the:uor(G)_and_eub(G):(FP)} This follows
    from~L\"uck~\cite[Proposition~11.11 on page~222]{Lueck(1989)}.
    \\[1mm]~\ref{the:uor(G)_and_eub(G):converse}
    This follows from L\"uck--Meintrup~\cite[Theorem~0.1]{Lueck-Meintrup(2000)}.
    \\[1mm]~\ref{the:uor(G)_and_eub(G):o_for_finitely_dominated}
    This follows now from the definitions.
    \\[1mm]~\ref{the:uor(G)_and_eub(G):o_for_finite} This
     follows for $\chi_f(\uor{G};\IZ)$ from the definitions.
    For  $\chi^{(2)}_f(\uor{G})$  apply
    Theorem~\ref{the:coincidence}.
\end{proof}

\begin{remark}\label{rem:about_eub(G)}
      The classifying spaces for proper $G$-actions $\eub{G}$ play a
      prominent role in the Baum-Connes Conjecture
      (see~Baum--Connes--Higson~\cite[Conjecture 3.15 on page
      254]{Baum-Connes-Higson(1994)}) and they have been intensively
      studied in their own right.

      Given a group $G$, there are often nice geometric models for
      $\eub{G}$ which are finite. If there is a finitely dominated
      model for $BG$, then $G$ must be torsion-free.  This is not the
      case for $\eub{G}$.
\end{remark}

\begin{example}[Groups with finite $\eub{G}$]
        \label{exa:groups_with_finite_eub(G)}
        If $G$ is a hyperbolic group in the sense of Gromov, then its
        Rips complex (for an appropriate parameter) is a finite model
        for $\eub{G}$ (see~Meintrup--Schick~\cite{Meintrup-Schick(2002)}).

        If the group $G$ acts simplicially cocompactly and properly by
        isometries on a $\CAT(0)$-space $X$, i.e., a complete
        Riemannian manifold with non-positive sectional curvature or a
        tree, then $X$ is a finite $G$-$CW$-model for $\eub{G}$.  This
        follows from Bridson--Haefliger~\cite[Corollary II.2.8 on page
        179]{Bridson-Haefliger(1999)}.

        Further groups admitting finite models for $\eub{G}$ are
        mapping class groups, the group of outer
        automorphisms of a finitely generated free group,
        finitely generated one-relator groups,
        and cocompact lattices in connected Lie groups.
\end{example}

\subsection{The M\"obius inversion for the proper orbit category}
\label{subsec:The_Moebius_inversion_for_uor(G)}

Next we take a closer look at
Theorem~\ref{the:K_theoretic-Moebius_inversion_and_L2-rank} in the
case of $\Gamma = \uor{G}$ for a group $G$ with a finite model for
$\eub{G}$.

Given an object $G/H$, we obtain by Lemma~\ref{lem:orbit_category_morphisms}
an isomorphism of groups
\begin{eqnarray} W_GH :=N_GH/H & \xrightarrow{\cong} & \aut(G/H)
  \label{W_GH_is_aut(G/H)}
\end{eqnarray}
by sending the class $gH \in N_GH/H$ to the $G$-automorphism $G/H \to
G/H, g'H \mapsto g'g^{-1}H$.

We obtain a bijection
\begin{eqnarray}
  \{(H) \mid H\subseteq G, |H| < \infty\} \xrightarrow{\cong} \iso(\uor{G}),
  \quad (H) \mapsto \overline{G/H}
  \label{iso(uor(G)_and_conjugacy_classes}
\end{eqnarray}
where $(H)$ denotes the conjugacy class of the subgroup $H$. Define a partial
ordering on $\{(H) \mid H\subseteq G, |H| < \infty\}$ by
\begin{eqnarray}
  (H) \le (K) & \Leftrightarrow & H \;\text{is conjugate to a subgroup of}\; K
  \label{le_on_conjugacy_classes}
\end{eqnarray}
Then the bijection~\eqref{iso(uor(G)_and_conjugacy_classes} is
compatible with the partial orderings
of~\eqref{partial_ordering_on_iso(Gamma)}
and~\eqref{le_on_conjugacy_classes}.

Given two elements $\overline{G/H}, \overline{G/K} \in \iso(\uor{G})$,
an $l$-chain $c \in \ch_l(\overline{G/K},\overline{G/H})$ in the sense
of Definition~\ref{def:chains} is, under the
bijection~\eqref{iso(uor(G)_and_conjugacy_classes}, the same as a
sequence of conjugacy classes of subgroups $(H_0) < (H_1) < \ldots < (H_l)$ with $(H_0) = (K)$ and
$(H_l) = (H)$. The $\aut(G/H)$-$\aut(G/K)$-biset $S(c)$ becomes under
this identification and the identification~\eqref{W_GH_is_aut(G/H)}
the $W_GH$-$W_GK$-biset
\begin{eqnarray*}
  S(c) & = &
  \map_G(G/H_{l-1},G/H) \times_{W_GH_{l-1}} \map_G(G/H_{l-2},G/H_{l-1})
  \times_{W_GH_{l-2}}
  \\
  & & \quad \quad \ldots \times_{W_GH_1} \map_G(G/K,G/H_1)
  \\
  & = &
  (G/H)^{H_{l-1}} \times_{W_GH_{l-1}} (G/H_{l-1})^{H_{l-2}}  \times_{W_GH_{l-2}}
  \ldots \times_{W_GH_1} (G/H_1)^K
\end{eqnarray*}
where we can arrange $K \subsetneq H_1 \subsetneq H_2 \subsetneq
\ldots \subsetneq H_{l-1} \subsetneq H$.

The commutative diagram appearing in
Theorem~\ref{the:K_theoretic-Moebius_inversion_and_L2-rank} becomes
the following diagram
  $$\xymatrix{
    & K_0(\IQ\uor{G}) \ar@/_/[ddl]_S \ar@/_/[ddr]_{\Res} &
    \\
    & &
    \\
    \bigoplus_{(H), |H| < \infty} K_0(\IQ W_GH) \ar@/_/[uur]_{E}
    \ar@/_/[rr]_{\omega} \ar[dd]_{\bigoplus_{(H), |H| <
        \infty}\rk_{W_GH}^{(2)}} && \bigoplus_{(H), |H| < \infty}
    K_0(\IQ W_GH) \ar@/_/[uul]_{I} \ar@/_/[ll]_{\mu}
    \ar[dd]^{\bigoplus_{(H), |H| < \infty}\rk_{W_GH}^{(2)}}
    \\
    & &
    \\
     \bigoplus_{(H), |H| < \infty} \IQ
    \ar@/_/[rr]_{\overline{\omega}^{(2)}} &&
    \bigoplus_{(H), |H| <
      \infty} \IQ \ar@/_/[ll]_{\overline{\mu}^{(2)}} }
$$
where $\rk_{W_G H}^{(2)} \colon K_0(\IQ W_GH) \to \IQ$ sends $[P]$ to
$\dim_{\caln(W_GH)}\bigl(P \otimes_{\IQ W_GH}
\caln(W_GH)\bigr)$, the map $\omega$ is given by the collection of
homomorphisms
$$\omega_{(H),(K)} \colon K_0(\IQ W_GH) \to K_0(\IQ W_GK), \quad
[P] \mapsto \bigl[P \otimes_{\IQ W_GH} \IQ \map_G(G/K,G/H)\bigr],$$
the map $\mu$
is given by the collection of homomorphisms
\begin{multline*}
  \mu_{(H),(K)} \colon K_0(\IQ W_GH) \to K_0(\IQ W_GK),
  \\
  [P] \mapsto \sum_{l \ge 0} (-1)^l \cdot \sum_{c \in \ch_l((K),(H))}
  [P \otimes_{\IQ W_GH} \IQ S(c)],
\end{multline*}
the map $\overline{\omega}^{(2)}$ is given by the matrix
$\left(\overline{\omega}^{(2)}_{(H),(K)}\right)$ over $\IQ$, where
$$\overline{\omega}^{(2)}_{(H),(K)} = \sum_{i=1}^r \frac{1}{|L_i|}$$
if the right $W_GK$-set $\map_G(G/K,G/H) = \left( G/H \right)^K$ is
the disjoint union $\sum_{i=1}^r L_i\backslash W_GK$, and the map
$\overline{\mu}^{(2)}$ is given by the matrix
$\left(\overline{\mu}^{(2)}_{(H),(K)}\right)$ over $\IQ$, where
$$\overline{\mu}^{(2)}_{(H),(K)} = \sum_{l \ge 0} (-1)^l \cdot \sum_{c \in \ch_l((K),(H))}\;
\sum_{i=1}^r \frac{1}{|L_i(c)|}$$ if the right $W_GK$-set
$$S(c)  = (G/K)^{H_{l-1}} \times_{W_GH_{l-1}} (G/H_{l-1})^{H_{l-2}}  \times_{W_GH_{l-2}}
\ldots \times_{W_GH_1} (G/H_1)^H$$ is the disjoint union $\sum_{i=1}^r
L_i(c)\backslash W_GK$.

\subsection{Additivity of the finiteness obstruction and the Euler characteristic for
the proper orbit category}
\label{subsec:Additivity_of_o_and_chi_for_uor(G)}

\begin{theorem}[Additivity of the finiteness obstruction and the Euler
  characteristic for the proper orbit category]
  \label{the:Additivity_of_for_uor(G)}
  Consider two groups $G_1$ and $G_2$ with a common subgroup
  $G_0$. Let $G$ be the amalgamated product $G = G_1 \ast_{G_0} G_2$.
  Then:

\begin{enumerate}
\item \label{the:Additivity_of_for_uor(G):G-pushout} We obtain a
  $G$-pushout of $G$-$CW$-complexes
  $$\xycomsquareminus{G \times_{G_0} \eub{G_0}}{j_1}{G \times_{G_1} \eub{G_1}}
  {j_2}{} {G \times_{G_2} \eub{G_2}}{}{\eub{G}}
  $$
  where $j_1$ and $j_2$ are inclusions of $G$-$CW$-complexes;
\item \label{the:Additivity_of_for_uor(G):Additivity_for_o} If
  $\uor{G_k}$ is of type (FP$_R$) for $k = 0,1,2$, then $\uor{G}$ is of
  type (FP$_R$) and we get for the finiteness obstruction
  \begin{multline*}
  \quad \quad
  o(\uor{G};R) = (i_1)_*\bigl(o(\uor{G_1};R)\bigr)
  + (i_2)_*\bigl(o(\uor{G_2};R)\bigr)
   \\ - (i_0)_*\bigl(o(\uor{G_0};R)\bigr) \quad \in K_0(R\uor{G}),
  \end{multline*}
  where $(i_k)_* \colon K_0(R\uor{G_k}) \to K_0(R\uor{G})$ is the homomorphism
  induced by the functor $(i_k)_* \colon \uor{G_k} \to \uor{G}$ coming from
  induction associated to the inclusion $i_k \colon G_k \to G$ for
  $k = 0,1,2$;

\item \label{the:Additivity_of_for_uor(G):Additivity_for_chi_f_and_chi}
  If $\uor{G_k}$ is of type (FP$_R$) for $k = 0,1,2$, then $\uor{G}$
  is of type (FP$_R$) and we get for the functorial Euler
  characteristic
  \begin{multline*}
  \quad \quad \chi_f(\uor{G};R) = (i_1)_*\bigl(\chi_f(\uor{G_1};R)\bigr) +
  (i_2)_*\bigl(\chi_f(\uor{G_2};R)\bigr)
  \\ - (i_0)_*\bigl(\chi_f(\uor{G_0};R)\bigr) \quad \in U(\uor{G}),
  \end{multline*}
  where $(i_k)_* \colon U(\uor{G_k}) \to U(\uor{G})$ is the homomorphism
  induced by the functor $(i_k)_* \colon \uor{G_k} \to \uor{G}$ coming from
  induction associated to the inclusion $i_k \colon G_k \to G$ for
  $k = 0,1,2$, and we get for the Euler characteristic
  $$\chi(\uor{G};R)
  = \chi(\uor{G_1};R) + \chi(\uor{G_2};R) - \chi(\uor{G_0};R) \quad \in \IZ.$$
  If $R$ is additionally Noetherian, then
  $\chi(B\uor{G_k};R)=\chi(\uor{G_k};R)$ and we get for the topological Euler
  characteristic
$$\chi(B\uor{G};R)
  = \chi(B\uor{G_1};R) + \chi(B\uor{G_2};R) - \chi(B\uor{G_0};R) \quad \in \IZ.$$

\item \label{the:Additivity_of_for_uor(G):Additivity_for_chi_f(2)_and_chi(2)}
  If $\uor{G_k}$ is of type ($L^2$) for $k = 0,1,2$, then $\uor{G}$
  is of type ($L^2$) and we get for the functorial $L^2$-Euler
  characteristic
  \begin{multline*}
  \quad \quad \chi_f^{(2)}(\uor{G}) = (i_1)_*\bigl(\chi_f^{(2)}(\uor{G_1})\bigr) +
  (i_2)_*\bigl(\chi_f^{(2)}(\uor{G_2})\bigr)
  \\ - (i_0)_*\bigl(\chi_f^{(2)}(\uor{G_0})\bigr) \quad \in U^{(1)}(\uor{G}),
  \end{multline*}
  where $(i_k)_* \colon U^{(1)}(\uor{G_k}) \to U^{(1)}(\uor{G})$ is the homomorphism
  induced by the functor $(i_k)_* \colon \uor{G_k} \to \uor{G}$ coming from
  induction associated to the inclusion $i_k \colon G_k \to G$ for
  $k = 0,1,2$, and we get for the $L^2$-Euler characteristic
  $$\chi^{(2)}(\uor{G})
  = \chi^{(2)}(\uor{G_1}) + \chi^{(2)}(\uor{G_2}) - \chi^{(2)}(\uor{G_0})\quad \in \IR.$$
\end{enumerate}
\end{theorem}
\begin{proof}~\ref{the:Additivity_of_for_uor(G):G-pushout} Associated
  to $G = G_1 \ast_{G_0} G_2$ there is a $1$-dimensional contractible
  $G$-$CW$-complex $T$ which is obtained as a $G$-pushout
  $$\comsquare{G/G_0 \times S^0}{\pr_1 \amalg \pr_2}{G/G_1 \amalg G/G_2}
  {}{} {G/G_0 \times D^1}{}{T}$$ where $\pr_k \colon G/G_0 \to G/G_k$ is the
  projection (see~Serre~\cite[Theorem~7 in~I.4 on page~32]{Serre(1980)}).

  Since for every finite
  subgroup $H \subseteq G$ the $H$-fixed point set $T^H$ is a non-empty
  subtree, by Serre~\cite[Proposition~19 in~I.4 on p.~36]{Serre(1980)},
  and thus contractible, the product with the diagonal $G$-action $T
  \times \eub{G}$ is again a model for $\eub{G}$. Note that
  $\res^{G_k}_G\eub{G}$ is a model for $\eub{G_k}$ and
  \[G/G_k\times\eub{G}\xrightarrow{\cong_G} G\times_{G_k}\res^{G_k}_G\eub{G},~~(gG_k,x)\mapsto (g,g^{-1}x)\]
  is a $G$-equivariant homeomorphism. Combining everything,
  we obtain
  the following $G$-pushout by crossing the $G$-pushout for $T$ above with
  $\eub{G}$
  \[\xycomsquareminus{G \times_{G_0} \eub{G_0}  \times S^0}
  {} {G\times_{G_1} \eub{G_1} \amalg G\times_{G_2} \eub{G_2}} {}{} {G
  \times_{G_0} \eub{G_0} \times D^1} {} {\eub{G}.}\]
  We can write the preceding $G$-pushout equivalently as
  \[\xycomsquareminus{G \times_{G_0} \eub{G_0}\times D^1}{j_1}{G \times_{G_1} \eub{G_1}}
  {j_2}{} {G \times_{G_2} \eub{G_2}}{}{\eub{G}}\]
  where $j_1$ and $j_2$ are inclusions of $G$-$CW$-complexes. Furthermore,
  $\eub{G_0}\times D^1$ is just another model of $\eub{G_0}$.
  \\[1mm]~\ref{the:Additivity_of_for_uor(G):Additivity_for_o} For $k = 0,1,2$ we get
  $$\ind_{i_k} C_*(\eub{G_k}) \cong C_*\bigl(G \times_{G_k} \eub{G_k}\bigr)$$
  where $C_*(\eub{G_k})$ is the cellular $\IZ\uor{G_k}$-chain complex of
  the $G_k$-$CW$-complex $\eub{G_k}$, $C_*\bigl(G \times_{G_k} \eub{G_k}\bigr)$
  is the cellular $\IZ\uor{G}$-chain complex of the
  $G$-$CW$-complex $G \times_{G_k} \eub{G_k}$.  From the $G$-pushout
  of assertion~\ref{the:Additivity_of_for_uor(G):G-pushout} we obtain
  a short exact sequence of $\IZ\uor{G}$-chain complexes
  $$0 \to \ind_{i_0} C_*(\eub{G_0}) \to
  \ind_{i_1} C_*(\eub{G_1}) \oplus \ind_{i_2} C_*(\eub{G_2}) \to C_*(\eub{G}) \to 0.$$
  Now apply~L\"uck~\cite[Theorem~11.2 on page~212]{Lueck(1989)},
  Theorem~\ref{the:relating_o_and_chi}, and Theorem~\ref{the:uor(G)_and_eub(G)}.
  \\[1mm]~\ref{the:Additivity_of_for_uor(G):Additivity_for_chi_f_and_chi}
  This follows from the definition of $\chi_f(\uor{G};R)$ since
  $\rk_{R\Gamma} \colon K_0(R\uor{G}) \to U(\uor{G})$ is compatible with
  induction homomorphisms induced from group homomorphisms. The
  category $\uor{G}$ is directly finite by
  Lemma~\ref{lem:directly_finite_and_idempotents_and_EI}, so
  Theorem~\ref{the:chi_f_determines_chi} applies.
  \\[1mm]~\ref{the:Additivity_of_for_uor(G):Additivity_for_chi_f(2)_and_chi(2)}
  We obtain for any object $G/H$ in $\uor{G}$ a short exact sequence of
  $\IZ\uor{G}$-chain complexes
  \begin{multline*}0 \to S_{G/H}\bigl(\ind_{i_0} C_*(\eub{G_0})\bigr) \to
  S_{G/H}\bigl(\ind_{i_1} C_*(\eub{G_1})\bigr) \oplus
  S_{G/H}\bigl(\ind_{i_2} C_*(\eub{G_2})\bigr)
  \\ \to   S_{G/H}\bigl(C_*(\eub{G})\bigr) \to 0.
  \end{multline*}
  For every finite subgroup
  $H \subset G_k$ and $k = 0,1,2$ the inclusion $G_k \to G$ induces an
  injection $W_{G_k}H \to W_GH$. The splitting functor is compatible
  with induction.  Now apply
  Theorem~\ref{lem:basic_properties_of_L2-Euler_characteristic}.
\end{proof}

\subsection{The Burnside integrality relations and the classical Burnside congruences}
\label{subsec:The_Burnside_integrality_relations_and_the_classical_Burnside_congruences}

Let $G$ be a group and let $X$ be a finite proper $G$-$CW$-complex. We
have defined its equivariant Euler characteristic $\chi^G(X) \in
U(\uor{G})$ in Definition~\ref{def:equivariant_Euler_characteristic}. The map
$$\overline{\omega}^{(2)} \colon \bigoplus_{(H), |H| < \infty} \IQ \to
\bigoplus_{(H), |H| < \infty} \IQ$$
defined in Subsection~\ref{subsec:The_Moebius_inversion_for_uor(G)} sends
$$\chi^G(X) \in U(\uor{G}) \subseteq U(\uor{G}) \otimes_{\IZ} \IQ
= \bigoplus_{(H), |H| < \infty} \IQ$$
to the collection
$\bigl(\chi^{(2)}(X^H;\caln(W_GH))\bigr)_{(H), |H| < \infty}$ of the
$L^2$-Euler characteristics of the $\caln(W_GH)$-chain complexes
$C_*(X^H) \otimes_{\IZ W_GH} \caln(W_GH)$. If $X = \eub{G}$, then
$\chi^{(2)}(X^H;\caln(W_GH)) = \chi^{(2)}(W_GH)$. Notice that we get
for the map
$$\overline{\mu}^{(2)} \colon U(\uor{G}) \otimes_{\IZ} \IQ
\to U(\uor{G}) \otimes_{\IZ} \IQ$$
defined in
Subsection~\ref{subsec:The_Moebius_inversion_for_uor(G)}
$$\overline{\mu}^{(2)}
\biggl(\bigl(\chi^{(2)}(X^H;\caln(W_GH))\bigr)_{(H), |H| <
\infty}\biggr) = \chi^G(X).$$

\begin{lemma}
  \label{lem:integrality_conditions}
  Consider $\eta = \bigl(\eta_{(H)}\bigr)_{(H), |H| < \infty} \in
  \prod_{(H), |H| < \infty} \IR$.  Then there is a finite proper
  $G$-$CW$-complex $X$ with $\chi^{(2)}(X^H;\caln(W_GH)) = \eta_{(H)}$
  for every finite subgroup $H \subseteq G$ if and only if
  $\eta \in  U(\uor{G}) \otimes_{\IZ} \IQ = \bigoplus_{(H), |H| < \infty} \IQ$ and
  $\overline{\mu}^{(2)}(\eta)$ lies in $U(\uor{G})$.
\end{lemma}
\begin{proof}
The direction ``$\Rightarrow$'' was proved
in the sentences preceding the Lemma. For the direction
``$\Leftarrow$'', we first note that every element of $U(\uor{G})$
can be realized as $\chi^G(X)$ for some $G$-$CW$-complex $X$.
Namely, $\overline{G/H}$ is realized by the 0-dimensional
$G$-$CW$-complex $G/H$, and $-\overline{G/H}$ is realized by the
1-dimensional $G$-$CW$-complex given by two $G$-1-cells $G/H \times
D^1$ attached to a single $G$-0-cell $G/H$. All other elements of
$U(\uor{G})$ arise from finite disjoint unions of $G$-$CW$-complexes
of these two forms. If $\eta \in U(\uor{G}) \otimes_{\IZ} \IQ$ and
$\overline{\mu}^{(2)}(\eta) \in U(\uor{G})$, then we realize
$\overline{\mu}^{(2)}(\eta)$ as $\chi^G(X)$ and apply
$\overline{\omega}^{(2)}$ with
Theorem~\ref{theorem:rational_Moebius_inversion} to obtain
$\chi^{(2)}(X^H;\caln(W_GH)) = \eta_{(H)}$
  for every finite subgroup $H \subseteq G$.
\end{proof}

\begin{lemma} \label{lem:condition_for_finite_eub(G)}
Let $G$ be a group such that $\uor{G}$ is of type
(FP$_\IQ$).

\begin{enumerate}
\item \label{lem:condition_for_finite_eub(G):chi(2)_f}
If $\uor{G}$ satisfies condition (I), then
$$\chi_f^{(2)}(\uor{G}) =
\overline{\mu}^{(2)}\biggl(\bigl(\chi^{(2)}(W_GH)\bigr)_{(H), |H| < \infty}\biggr);$$

\item \label{lem:condition_for_finite_eub(G):type_(FF)}
If there is a finite model for $\eub{G}$, then
the following integrality condition is satisfied
$$\overline{\mu}^{(2)}\biggl(\bigl(\chi^{(2)}(W_GH)\bigr)_{(H), |H| < \infty}\biggr)
\in U(\uor{G}).$$
\end{enumerate}
\end{lemma}
\begin{proof}
This follows from
Theorem~\ref{the:The_finiteness_obstruction_and_the_(functorial)_L2-Euler_characteristic},
Theorem~\ref{the:uor(G)_and_eub(G)},
and Lemma~\ref{lem:integrality_conditions}.
\end{proof}

\begin{example}[Burnside congruences]\label{lem:Burnside_congruences}
  These considerations are already interesting in the case of a finite
  group $G$. Since we assume $G$ is finite in this example, we refrain here from writing $|H| < \infty$ when
  summing over conjugacy classes $(H)$ of subgroups of $G$.
  For every finite $G$-$CW$-complex $X$, the map
  $$\overline{\omega}^{(2)} \colon \bigoplus_{(H)} \IQ \to
  \bigoplus_{(H)} \IQ$$ sends the equivariant Euler characteristic $\chi^G(X)$ to
  the collection $\bigl(\chi(X^H)/|W_GH|\bigr)_{(H)}$, where
  $\chi(X^H)$ is the classical Euler characteristic of the $H$-fixed
  point set. We conclude from Lemma~\ref{lem:integrality_conditions} that
  for an element  $\eta = (\eta_{(H)})_{(H)} \in\bigoplus_{(H)} \IQ$
  there exists a finite $G$-$CW$-complex $X$ such that
  $\chi(X^H)/|W_GH| = \chi^{(2)}(X^H;\caln(W_GH))$ agrees with $\eta_{(H)}$
  for any subgroup $H \subseteq G$, if and only if
  $\overline{\mu}^{(2)}(\eta) \in U(\uor{G})$. The latter is a kind of
  integrality condition. In the case of a finite group $G$ it can be transformed
  into equivalent congruence conditions for integers.

  Let $$\ch = \ch^G \colon U(\uor{G})  \to \bigoplus_{(H)} \IZ$$
  be the map uniquely determined by the property that it sends
  $\chi^G(X)$ to the collection $\bigl(\chi(X^H)\bigr)_{(H)}$ for every
  finite $G$-$CW$-complex $X$.  Under the obvious identification of
  $U(\uor{G})$ with the Burnside ring $A(G)$, the map $\ch$ corresponds
  to the character map which sends a finite $G$-set $S$ to the
  collection $\bigl(|S^H|\bigr)_{(H)}$. We have
  $$i \circ \ch = D \circ \overline{\omega}^{(2)} \circ i,$$
  if $i \colon U(G) \to U(G) \otimes_{\IZ} \IQ$ is the obvious inclusion and
  the map $D \colon U(G) \otimes_{\IZ} \IQ \to U(G) \otimes_{\IZ} \IQ$ is given by the
  diagonal matrix whose entry at $(H)$ is $|W_GH|$.  Let
  $$\nu \colon \bigoplus_{(H)} \IZ \to   \bigoplus_{(H)} \IZ$$
  be the map uniquely determined by
  $i \circ \nu = D \circ \overline{\mu}^{(2)} \circ D^{-1} \circ i$.
  One easily checks that it is given by the
  integer matrix whose entry at $\bigl((H),(K)\bigr)$ is
  $$\sum_{l \ge 0} (-1)^l \cdot \sum_{(H_0) < \cdots < (H_l) \in \ch_l\bigl((K),(H)\bigr)}
  \prod_{i=1}^l\;\bigl|W_GH_{i+1}\backslash \map_G(G/H_i,G/H_{i+1})\bigr|.$$
  Notice that $i \circ \nu \circ \chi = D \circ i$.
  We conclude that an element $\xi \in \bigoplus_{(H)} \IZ$ lies in the image of
  $\ch$ if and only if, for every conjugacy class $(H)$ of subgroups of the finite group $G$, the
  following congruence of integers holds:
  $$\nu(\xi)_{(H)} \equiv 0 \mod |W_GH|.$$
  These are the \emph{Burnside ring congruences}. For more information about
  the Burnside ring we refer for instance to~tom~Dieck~\cite[Chapter~1]{Dieck(1979)}.

  If $G$ is the cyclic group $\IZ/p$ of order $p$ for a prime $p$,
  then $U(\uor{\IZ/p}) = \IZ^2$,
  $$\ch = \left(\begin{array}{cc} p & 1 \\ 0 & 1 \end{array}\right) \colon
  U(\uor{\IZ/p}) = \IZ^2 \to U(\uor{\IZ/p}) = \IZ^2,$$
  and
  $$\nu = \left(\begin{array}{cc} 1 & -1 \\ 0 & 1 \end{array}\right) \colon
  U(\uor{\IZ/p}) = \IZ^2 \to U(\uor{\IZ/p}) = \IZ^2.$$
  The Burnside ring congruences reduce to one congruence, namely
  $$\eta_{(\IZ/p)/\{1\}} - \eta_{(\IZ/p)/(\IZ/p)} \equiv 0 \mod p.$$
  The latter reflects the fact that the cardinality of $S -S^{\IZ/p}$  is a multiple of $p$
  for a finite $\IZ/p$-set $S$.
\end{example}

\begin{example}[Amenable $G$]\label{exa:amenable_G}
  Let $G$ be an amenable group. Suppose that $\uor{G}$ is of type
  (FP$_\IQ$). Then
  $\chi^{(2)}(\uor{G})$ is the image of $\eta =
  \bigl(\eta_{(H)}\bigr)_{(H), |H| < \infty}$ under
  $$\overline{\mu}^{(2)} \colon U(\uor{G}) \otimes_{\IZ} \IQ
  \to U(\uor{G}) \otimes_{\IZ} \IQ$$
  where $\eta_{(H)} = 0$ if $W_GH$ is infinite and
  $\eta_{(H)} = 1/|W_GH|$ if $W_GH$ is finite.

  In particular, if $W_GH$ is infinite for every finite subgroup
  $H \subseteq G$, then $\chi^{(2)}(\uor{G})$  vanishes.

  This follows from
  Theorem~\ref{the:The_finiteness_obstruction_and_the_(functorial)_L2-Euler_characteristic},
  Lemma~\ref{lem:condition_for_finite_eub(G)}, and the result of
  Cheeger and Gromov that all the $L^2$-Betti numbers of any infinite
  amenable group $G$ vanish (see~Cheeger--Gromov~\cite{Cheeger-Gromov(1986)}
  and L\"uck~\cite[Theorem~7.2 on page~294]{Lueck(2002)}).
\end{example}

\subsection{The infinite dihedral group}
\label{subsec:infinite_dihedral_group}
Consider the infinite dihedral group
$$D_{\infty} = \langle t,s \mid s^2 = 1, sts = t^{-1} \rangle
\cong \IZ \rtimes \IZ/2
\cong \IZ/2 \ast \IZ/2.$$
As an illustration we want to make all the material of this section explicit for this
easy special case.

The infinite dihedral group $D_{\infty}$
has three conjugacy classes of finite subgroups $(C_1)$, $(C_2)$,
and $(T)$, where $C_1 = \langle s \rangle$ and $C_2 = \langle ts
\rangle$ have order two and $T$ is the trivial group.

One easily checks that $W_{D_{\infty}}C_i$ is trivial for $i = 1,2$
and $W_{D_{\infty}} T = D_{\infty}$. Hence we get
$$\Split K_0(\IQ\uor{D_{\infty}}) = K_0(\IQ D_{\infty}) \oplus K_0(\IQ) \oplus K_0(\IQ)
=  K_0(\IQ D_{\infty}) \oplus \IZ \oplus \IZ$$
by the discussion in Subsection~\ref{subsec:The_Moebius_inversion_for_uor(G)}.

The $W_{D_{\infty}}C_i$-$W_{D_{\infty}}T$-biset
$\map_{D_{\infty}}({D_{\infty}}/T,{D_{\infty}}/C_i)$ is given by the
right ${D_{\infty}}$-set ${C_i\backslash D_{\infty}}$ for $i = 1,2$.
The $W_{D_{\infty}}T$-$W_{D_{\infty}}T$-biset
$\map_{D_{\infty}}(D_{\infty}/T,D_{\infty}/T)$ is $D_{\infty}$
regarded as $D_{\infty}$-$D_{\infty}$-biset.  The
$W_{D_{\infty}}C_j$-$W_{D_{\infty}}C_i$-biset
$\map_{D_{\infty}}(D_{\infty}/C_i,D_{\infty}/C_j)$ is empty for $i
\not= j$ and is the $\{1\}$-$\{1\}$-biset consisting of one point
for $i = j$.  The $W_{D_{\infty}}T$-$W_{D_{\infty}}C_i$-biset
$\map_{D_{\infty}}(D_{\infty}/C_i,D_{\infty}/T)$ is empty for $i
=1,2$.  There are exactly two $1$-chains in $\uor{D_{\infty}}$,
namely $(T) < (C_1)$ and $(T) < (C_2)$.

Hence we get
\begin{multline*}
  \omega \colon K_0(\IQ D_{\infty}) \oplus \IZ \oplus \IZ \to K_0(\IQ
  D_{\infty}) \oplus \IZ \oplus \IZ,
  \\
  (x,n_1,n_2) \mapsto \bigl(x + n_1 \cdot [\IQ C_1\backslash D_{\infty}] + n_2
  \cdot [\IQ C_2\backslash D_{\infty}], n_1, n_2\bigr),
\end{multline*}
\begin{multline*}
  \mu \colon K_0(\IQ D_{\infty}) \oplus \IZ \oplus \IZ \to K_0(\IQ
  D_{\infty}) \oplus \IZ \oplus \IZ,
  \\
  (x,n_1,n_2) \mapsto \bigl(x - n_1 \cdot [\IQ C_1\backslash D_{\infty}] - n_2
  \cdot [\IQ C_2\backslash D_{\infty}], n_1, n_2\bigr),
\end{multline*}
\begin{multline*}
  \overline{\omega}^{(2)} = \left(\begin{array}{ccc} 1 & 1/2 & 1/2
      \\
      0 & 1 & 0
      \\
      0 & 0 & 1
    \end{array}\right)\colon \IZ^3 \to \IZ^3,
  \quad (n_0,n_1,n_2) \mapsto (n_0 +n_1/2 +n_2/2,n_1,n_2),
\end{multline*}
and

\begin{multline*}
  \overline{\mu}^{(2)} = \left(\begin{array}{ccc} 1 & -1/2 & - 1/2
      \\
      0 & 1 & 0
      \\
      0 & 0 & 1
    \end{array}\right)\colon \IZ^3 \to \IZ^3,
  \quad (n_0,n_1,n_2) \mapsto (n_0 - n_1/2 - n_2/2,n_1,n_2).
\end{multline*}

The map
$$\rk^{(2)}_{\uor{D_{\infty}}} \colon K_0(\IQ D_{\infty}) \oplus \IZ \oplus \IZ
\to \IZ \oplus \IZ \oplus \IZ$$
sends $\bigl([P],n_1,n_2\bigr)$ to
$\bigl(\dim_{\caln(D_{\infty})}(P \otimes_{\IQ D_{\infty}} \caln(D_{\infty})),n_1,n_2\bigr)$.

There is the isomorphism
$$\IZ \oplus \IZ \oplus \IZ \xrightarrow{\cong} K_0(\IQ D_{\infty}),
\quad \bigl(n_0,n_1,n_2\bigr) \mapsto n_0 \cdot [\IQ D_{\infty}] +
n_1 \cdot [\IQ C_1\backslash D_{\infty}] + n_2 \cdot [\IQ
C_2\backslash D_{\infty}]$$ (see for example the Mayer-Vietoris
sequence for amalgated products in Waldhausen~\cite[Corollary 2.15
on page 216]{Waldhausen1978a} and the subsequent remarks there).
Under this identification
$$\rk^{(2)}_{\uor{D_{\infty}}} =
\left(\begin{array}{ccccc} 1 & 1/2 & 1/2 & 0 & 0
    \\
    0 & 0 & 0 & 1 & 0
    \\
    0 & 0 & 0 & 0 & 1
  \end{array}\right)
\colon \IZ^5 \to \IZ^3,
$$
$$\omega = \left(\begin{array}{ccccc}
    1 & 0 & 0 & 0 & 0
    \\
    0 & 1 & 0 & 1 & 0
    \\
    0 & 0 & 1 & 0 & 1
    \\
    0 & 0 & 0 & 1 & 0
    \\
    0 & 0 & 0 & 0 & 1
  \end{array}\right)\colon \IZ^5 \to \IZ^5,$$
and
$$\mu = \left(\begin{array}{ccccc}
    1 & 0 & 0 & 0 & 0
    \\
    0 & 1 & 0 & - 1 & 0
    \\
    0 & 0 & 1 & 0 & - 1
    \\
    0 & 0 & 0 & 1 & 0
    \\
    0 & 0 & 0 & 0 & 1
  \end{array}\right)\colon \IZ^5 \to \IZ^5.$$

The infinite dihedral group $D_{\infty} = \IZ \rtimes \IZ/2$ acts on
$\IR$ by the  action of $\IZ$ on $\IR$ given by addition and
the action of $\IZ/2$ in $\IR$ given by multiplication with
$(-1)$.  There is a $D_{\infty}$-$CW$-structure on $\IR$ such that
there are three equivariant cells of the type $D_{\infty}/C_1 \times D^0$,
$D_{\infty}/C_2 \times D^0$, and $D_{\infty}/T \times D^1$. One easily
checks that this is a model for $\eub{D_{\infty}}$. Hence we get for
the equivariant Euler characteristic of $\eub{D_{\infty}}$
$$\chi^{D_{\infty}}(\eub{D_{\infty}})
= \overline{D_{\infty}/C_1} + \overline{D_{\infty}/C_2} - \overline{D_{\infty}/T} \in U(\uor{D_{\infty}}).$$
By Theorem~\ref{the:Two_splittings_and_the_K-theoretic_Moebius_inversion}~%
\ref{the:Two_splittings_and_the_K-theoretic_Moebius_inversion:o}
and Theorem~\ref{the:uor(G)_and_eub(G)}~\ref{the:uor(G)_and_eub(G):o_for_finite}
the image of the finiteness obstruction $o(\uor{D_{\infty}})$ under the isomorphism
$$S \colon K_0(\IQ \uor{D_{\infty}})
\xrightarrow{\cong} \Split K_0(\IQ \uor{D_{\infty}}) = K_0(\IQ
D_{\infty}) \oplus \IZ \oplus \IZ = \IZ^5$$ is $(-1,0,0,1,1)$.
The image of this element under $\omega$ is $(-1,1,1,1,1)$.
All this is consistent with Theorem~\ref{the:Additivity_of_for_uor(G)}
applied to $D_{\infty} = \IZ/2 \ast \IZ/2$.

The trivial $\IQ D_{\infty}$-module $\IQ$ has a finite projective
$\IQ D_{\infty}$-resolution of the form $0 \to \IQ D_{\infty} \to
\IQ D_{\infty}/C_1 \oplus \IQ D_{\infty}/C_1 \to \IQ \to 0$ coming
from the $\IQ D_{\infty}$-chain complex of $\IR$. This implies that
the homomorphism
$$\Res \colon K_0(\IQ \uor{D_{\infty}})
\xrightarrow{\cong} \Split K_0(\IQ \uor{D_{\infty}}) = K_0(\IQ
D_{\infty}) \oplus \IZ \oplus \IZ = \IZ^5$$ sends
$o(\uor{D_{\infty}};\IQ)$ to $(-1,1,1,1,1)$ (see
Theorem~\ref{the:The_finiteness_obstruction_and_the_(functorial)_L2-Euler_characteristic}~\ref{the:The_finiteness_obstruction_and_the_(functorial)_L2-Euler_characteristic:Res(o)}).
This is consistent with the fact that $\omega$ sends the image of
the finiteness obstruction $o(\uor{D_{\infty}})$ under $S$, which is
given by $(-1,0,0,1,1) \in \IZ^5$, to the element $(-1,1,1,1,1) \in
\IZ^5$ (see
Theorem~\ref{the:Two_splittings_and_the_K-theoretic_Moebius_inversion}).

We have $\chi^{(2)}_f(\uor{D_{\infty}};\IQ) = (-1,1,1) \in U(\uor{D_{\infty}}) =
\IZ^3$.  The composite
$$\rk_{\uor{D_{\infty}}}^{(2)} \circ \Res \colon
K_0(\IQ \uor{D_{\infty}}) \to U(\uor{D_{\infty}}) = \IZ^3$$ sends
$o(\uor{D_{\infty}};\IQ)$ to
$\bigl(\chi^{(2)}(D_{\infty}),\chi^{(2)}(\{1\}),\chi^{(2)}(\{1\})\bigr)$.
Since the $L^2$-Euler characteristic of an infinite amenable group
vanishes (see~Cheeger--Gromov~\cite{Cheeger-Gromov(1986)}) and the
$L^2$-Euler characteristic of the trivial group is $1$, we get
$\bigl(\chi^{(2)}(D_{\infty}),
\chi^{(2)}(\{1\}),\chi^{(2)}(\{1\})\bigr) = (0,1,1)$. This is
consistent with the fact that $\overline{\omega}^{(2)}$ sends
$(-1,1,1)$ to $(0,1,1)$ and with Example~\ref{exa:amenable_G}.

\subsection{The fundamental category}
\label{subsec:The_fundamental_category}

Let $X$ be a $G$-space. Consider the functor
$$F \colon \Or(G) \to \GROUPOIDS, \quad G/H \mapsto \Pi\bigl(\map_G(G/H,X)\bigr),$$
which sends $G/H$ to the fundamental groupoid of $X^H =
\map_G(G/H,X)$.  Its homotopy colimit is by definition the
\emph{fundamental groupoid} $\Pi(G,X)$ which plays an important role
in transformation groups (see~L\"uck~\cite[Definition~8.13 on
page~144]{Lueck(1989)}).

Denote by $\underline{\Pi}(G,X)$ the homotopy colimit of the functor
$F$ above restricted to $\uor{G}$. If all isotropy groups of $X$ are
finite, then $\Pi(G,X)$ and $\underline{\Pi}(G,X)$ agree.

Suppose that there is a finite $G$-$CW$-model for $\eub{G}$. Let
$I_n$ be the set of equivariant $n$-cells $c = G/H_c \times (D^n -
S^{n-1})$. Consider a $G$-$CW$-complex $X$.  Suppose that for every
finite subgroup $H \subseteq G$ each groupoid $\Pi(X^H)$ is of type
(FP$_{\IQ}$).  This is equivalent to requiring that for every finite
subgroup $H \subseteq G$ the set $\pi_0(X^H)$ is finite and at each
base point $x \in X^H$ the fundamental group $\pi_1(X^H,x)$ is of
type (FP$_{\IQ}$).  This follows from~Brown~\cite[Exercise~8 in
VIII.6 on page 205]{Brown(1982)} using the facts that  $W_GH$ is of
type (FP$_{\IQ}$) because $\uor{G}$ is of type (FP$_{\IQ}$) (see
Theorem~\ref{the:uor(G)_and_eub(G)}~\ref{the:uor(G)_and_eub(G):(FF)}
and Lemma~\ref{lem:finite_homological_dimension}~%
\ref{lem:finite_homological_dimension:Res_x}) and for every object
$x \colon G/H \to X$ in $\underline{\Pi}(G,X)$ there exists an exact
sequence
\begin{eqnarray}
  & 1 \to \pi_1(X^H,x) \to \aut(x \colon G/H \to X) \to W_GH(x) \to 1
  & \label{exact_sequence_for_aut(x)}
\end{eqnarray}
for the subgroup $W_GH(x) \subseteq W_GH$ of finite index which is
the isotropy group of the component in $X^H$ determined by $x$ under
the $W_GH$-action on $\pi_0(X^H)$ (see~L\"uck~\cite[Proposition~8.33
on page~150]{Lueck(1989)}). Hence the homotopy colimit formula
of~Fiore--Sauer--L\"uck~\cite{FioreLueckSauerHoColim(2009)} applies.
For instance we get
\begin{eqnarray*}
  \chi^{(2)}(\underline{\Pi}(G;X))
  & = &
  \sum_{n \ge 0} (-1)^n \cdot  \sum_{c \in I_n} \; \sum_{C \in \pi_0(X^{H_c})/W_GH_c}
  \chi^{(2)}\bigl(\aut(x(C))\bigr);
  \\
  \chi(\underline{\Pi}(G;X);\IQ)
  & = &
  \sum_{n \ge 0} (-1)^n \cdot  \sum_{c \in I_n} \; \sum_{C \in \pi_0(X^{H_c})/W_GH_c}
  \chi\bigl(B\aut(x(C));\IQ\bigr),
\end{eqnarray*}
where for a component $C \in \pi_0(X^{H_c})$ we denote by $x(C) \colon G/H_c
\to X$ an object in $\underline{\Pi}(G,X)$ such that $x(C)(eH_c)$ lies
in the component $C$ and $\aut(x(C))$ is its automorphism group in
$\underline{\Pi}(G,X)$ which fits into the exact
sequence~\eqref{exact_sequence_for_aut(x)}.

If we take $X = \pt$ itself, we get back
Theorem~\ref{the:uor(G)_and_eub(G)}~\ref{the:uor(G)_and_eub(G):o_for_finite}.

One can define for a functor $\mu \colon \Or(G) \to \GROUPOIDS$ its
\emph{equivariant Eilenberg Mac Lane space $E(\mu,1)$} which is a
$G$-$CW$-complex such that $\mu$ can be identified with the functor
$\Or(G) \to \GROUPOIDS$ sending $G/H$ to $\Pi(E(\mu,1)^H)$ and we
have $\pi_n(E(\mu,1)^H,x)$ is trivial for all $n \ge 2$, $H
\subseteq G$ and $x \in E(\mu,1)^H$
(see~L\"uck~\cite{Lueck(1987a)}). There is a natural equivalence
$\hocolim_{\Or(G)} \mu \to \Pi(G;E(\mu,1))$ which induces an
isomorphism
$$K_0\bigl(\IZ\hocolim_{\Or(G)} \mu\bigr) \to K_0\bigl(\IZ\Pi(G;E(\mu,1))\bigl).$$
Under this isomorphism the finiteness obstruction of
$\hocolim_{\Or(G)} \mu$ in the sense of
Definition~\ref{def:finiteness_obstruction_of_a_category}
corresponds to the finiteness obstruction of $E(\mu,1)$ in the sense
of~L\"uck~\cite[Definition~14.4 on page~278]{Lueck(1989)}.


\typeout{-----------  Section 9: An example of a finite category without property EI ---}
\section{An example of a finite category without property EI}
\label{sec:An_example_of_a_finite_category_without_property_EI}

For the remainder of this section we will consider the following
category $\Gamma$.  It has precisely two objects $x$ and $y$.  There
is precisely one morphism $u \colon x \to y$ and precisely one
morphism $v \colon y \to x$.  There are precisely two endomorphisms
of $x$, namely, $v \circ u$ and $\id_x$.  There are precisely two
endomorphisms of $y$, namely, $u \circ v$ and $\id_x$.  We have $vuv
= v$ and $uvu= u$.  Obviously $\Gamma$ is a free finite category. It
has two idempotents which are not the identity, namely, $vu$ and
$uv$. It is directly finite but it is not Cauchy complete and not an
EI-category. In this section we compute the homomorphisms $S$, $E$,
and $\Res$ for $K_0(R\Gamma)$ and determine the finiteness
obstruction.

Given an $R$-module $M$, we define three $R\Gamma$-modules $I_xM$,
$I_yM$, and $I_cM$ as follows. The contravariant functor $I_xM$ sends
$x$ to $M$ and $y$ to $\{0\}$ and every morphism except $\id_x$ to
the zero homomorphism. The contravariant functor $I_yM$ sends $y$ to
$M$ and $x$ to $\{0\}$ and every morphism except $\id_y$  in
$\Gamma$ to the zero homomorphism. The contravariant functor $I_cM$
sends both $x$ and $y$ to $M$ and every morphism in $\Gamma$ to the
identity $\id_M$.

\begin{lemma} \label{lem_structure_of_Rcalc-modules} Let $M$ be an
  $R\Gamma$-module. Then there is an isomorphism of $R\Gamma$-modules,
  natural in $M$
  $$f  \colon
  I_x\bigl(\ker(M(vu))\bigr) \oplus I_y\bigl(\ker(M(uv))\bigr) \oplus
  I_c(\im(vu)) \xrightarrow{\cong} M.$$
\end{lemma}
\begin{proof}
  The transformation $f$ is given at the object $x$ by the direct sum
  of the obvious inclusions
  $$i_x \oplus j_x \colon \ker(M(vu)) \oplus \im(M(vu)) \xrightarrow{\cong} M(x).$$
  This is an isomorphism since $M(vu)^2 = M((vu)^2) = M(vu)$.  The
  transformation $f$ is given at the object $y$ by the direct sum of the
  inclusion $i_y$ and the map induced by $M(v)$
  $$i_y \oplus M(v)|_{\im(M(uv))} \colon \ker(M(uv)) \oplus \im(M(vu))
  \xrightarrow{\cong} M(y).$$
  This is an isomorphism of $R$-modules, an inverse is given by
  $$(\id - M(uv)) \times M(u) \colon M(y) \to \ker(M(uv)) \oplus \im(M(vu)).$$
  It remains to check that $f$ is a transformation. We check this for
  the morphism $v$, the proof for $u$ is analogous. We have to show that
  the following diagram is commutative
  $$\xymatrix{\ker(M(vu)) \oplus \im(M(vu)) \ar[r]^{0 \oplus \id}\ar[d]_{i_x \oplus j_x}
  & \ker(M(uv)) \oplus \im(M(vu)) \ar[d]^{i_y \oplus
    M(v)|_{\im(M(uv))}}
  \\
  M(x) \ar[r]_{M(v)} & M(y) }$$ This is equivalent to showing that
 $M(v)|_{\ker(M(vu))} = 0$. This follows from $M(v) = M(vuv) = M(v) \circ M(vu)$.
\end{proof}

\begin{lemma} \label{lem:fin:gen._projective_Rcalc_modules} Let $M$ be
  an $R$-module.
  \begin{enumerate}
  \item \label{lem:fin:gen._projective_Rcalc_modules:restriction} The
    functors $\Res_x$ and $\Res_y$ respectively from
    $\MOD\text{-}R\Gamma$ to $\MOD\text{-}R$, which are given by
    evaluation at $x$ and $y$ respectively, are exact and send
    finitely generated projective $R\Gamma$-modules to finitely
    generated projective $R$-modules;

  \item \label{lem:fin:gen._projective_Rcalc_modules:I} The following
    assertions are equivalent:
    \begin{enumerate}
    \item $M$ is a finitely generated projective $R$-module;
    \item $I_xM$ is a finitely generated projective $R\Gamma$-module;
    \item $I_yM$ is a finitely generated projective $R\Gamma$-module;
    \item $I_cM$ is a finitely generated projective $R\Gamma$-module.
    \end{enumerate}
  \end{enumerate}
\end{lemma}
\begin{proof}\ref{lem:fin:gen._projective_Rcalc_modules:restriction}
  Obviously $\Res_x$ and $\Res_y$ are exact. Hence it remains to show that
  they send both $R\mor(?,x)$ and $R\mor(?,y)$ to a finitely generated projective
  $R$-module. This is obviously true.
  \\[1mm]\ref{lem:fin:gen._projective_Rcalc_modules:I}
  Suppose that $I_xM$ is a finitely generated projective
  $R\Gamma$-module.  Then $M$ is a finitely generated $R$-module
  because of
  assertion~\ref{lem:fin:gen._projective_Rcalc_modules:restriction}
  since $I_x(M)(x) = M$. Analogously one shows that $M$ is finitely
  generated projective if $I_yM$ or $I_cM$ is a finitely generated
  projective $R\Gamma$-module.

  Suppose that $M$ is a finitely generated projective $R$-module.
  We want to show that $I_xM$, $I_yM$, and $I_cM$ are finitely generated
  projective $R\Gamma$-modules. Since the functors $I_x$, $I_y$, and $I_c$
  are exact, it suffices to check this in the special case $M =
  R$. This follows from Lemma~\ref{lem_structure_of_Rcalc-modules}
  since $R\mor(?,x)$ and $R\mor(?,y)$ are free $R\Gamma$-modules and
  $I_xR$, $I_yR$, and $I_cR$ are direct summands in $R\mor(?,x)$ or
  $R\mor(?,y)$.
\end{proof}

\begin{corollary}
The constant functor $\underline{R} \colon \Gamma^{\op} \to
R\text{-}\MOD$ with value $R$ defines a projective $R\Gamma$-module.
In particular, $\underline{R}$ admits a finite projective resolution
and $\Gamma$ is of type (FP$_R$).
\end{corollary}

\begin{lemma}\label{lem:K_0(Rcalc)}
  We obtain isomorphisms, inverse to one another,
  \begin{multline*}
  \alpha \colon K_0(R) \oplus K_0(R) \oplus K_0(R) \xrightarrow{\cong} K_0(R\Gamma),
  \\
  \bigl([P_1],[P_2],[P_3]\bigr) \mapsto [I_x(P_1)] + [I_y(P_2)] + [I_c(P_3)]
  \end{multline*}
  and
  $$\beta \colon K_0(R\Gamma) \xrightarrow{\cong} K_0(R) \oplus K_0(R) \oplus K_0(R),
  \quad [P] \mapsto \bigl([S_xP], [S_yP], [\Res_xP] -[S_xP]\bigr),$$
  where the functors $S_x$ and $S_y$ are the splitting functors defined
  in~\eqref{S_x}.
\end{lemma}
\begin{proof}
  This follows from Lemma~\ref{lem_structure_of_Rcalc-modules}
  and~Lemma~\ref{lem:fin:gen._projective_Rcalc_modules}.
\end{proof}

Consider the following commutative diagram
$$\xymatrix{
  & & \Split K_0(R\Gamma)
  \\
  & K_0(R\Gamma) \ar[ru]^{S} \ar[ld]_{\Res} &
  \\
  \Split K_0(R\Gamma) \ar[d]_{\rk_R} & & \Split K_0(R\Gamma)
  \ar[uu]_{\id}\ar[ll]^{\omega} \ar[lu]_E \ar[d]^{\rk_R}
  \\
  U(\Gamma) & & U(\Gamma) \ar[ll]^{\overline{\omega}} }$$ where the
homomorphisms $S$ and $E$ have been defined
in~\eqref{S_colon_SplitK_0_toK_0} and
in~\eqref{E_colon_SplitK_0_to_K_0} and satisfy $S \circ E = \id$ by
Lemma~\ref{S_circ_E_is_id}, the homomorphism $\Res$ sends $[P]$ to
$\bigl([\Res_x P], [\Res_y P]\bigr)$, the homomorphism $\omega$ has
been defined in~\eqref{omega_splitK_0_to_splitK_0}, the map $\rk_R$
is given by the direct sum of the homomorphisms $K_0(R) \to \IZ$
sending $[P]$ to $\rk_R(P)$ and $\overline{\omega}$ is given by the
matrix $\squarematrix{2}{1}{1}{2}$. Under the identification
$\alpha$ of Lemma~\ref{lem:K_0(Rcalc)} and the definitions $\Split
K_0(R\Gamma) := K_0(R) \oplus K_0(R)$ and $U(\Gamma) = \IZ \oplus
\IZ$, where the first summand corresponds to $x$ and the second to
$y$, this diagram becomes
$$\xymatrix{
& &  K_0(R) \oplus K_0(R)
\\
& &
\\
& K_0(R) \oplus K_0(R) \oplus K_0(R)
\ar[ruu]^{\begin{array}{c}S\end{array}=\left(\begin{array}{ccc}\id &
0 & 0
\\0 & \id & 0
\end{array}\right)\quad\quad} \ar[ldd]_{\left(\begin{array}{ccc}\id & 0 &
\id \\0 & \id & \id \end{array}\right)\quad\quad}^{\begin{array}{c}\Res\end{array}} &
\\
& &
\\
K_0(R) \oplus K_0(R)
\ar[dd]_{\begin{array}{c}\rk_R\end{array}=\squarematrix{\rk_R}{0}{0}{\rk_R}}
& & K_0(R) \oplus K_0(R)
\ar[uuuu]_{\squarematrix{\id}{0}{0}{\id}=\begin{array}{c}\id\end{array}}
\ar[ll]^{\squarematrix{2 \cdot \id}{\id}{\id}{2 \cdot
\id}}_{\begin{array}{c}\omega\end{array}}
\ar[luu]_{\quad\left(\begin{array}{cc}\id & 0  \\ 0 & \id \\ \id &
\id \end{array}\right)}^{\begin{array}{c}E\end{array}}
 \ar[dd]^{\squarematrix{\rk_R}{0}{0}{\rk_R}=\begin{array}{c}\rk_R\end{array}}
\\
& &
\\
\IZ \oplus \IZ & & \IZ \oplus \IZ.
\ar[ll]^{\squarematrix{2}{1}{1}{2}}_{\begin{array}{c}\overline{\omega}\end{array}}
}$$

The finiteness obstruction $o(\Gamma;R) \in K_0(R\Gamma)$ of
Definition~\ref{def:finiteness_obstruction_of_a_category}
corresponds under the identification $\alpha$ of
Lemma~\ref{lem:K_0(Rcalc)} to the element $(0,0,[R]) \in K_0(R)
\oplus K_0(R) \oplus K_0(R)$. Its image under $S \colon K_0(R\Gamma)
\to \Split K_0(R\Gamma) = K_0(R) \oplus K_0(R)$ is $(0,0)$.  Its
image under $\Res \colon K_0(R\Gamma) \to \Split K_0(R\Gamma) =
K_0(R) \oplus K_0(R)$ is $([R],[R])$.  Its image under the composite
$\rk_R \circ \Res \colon K_0(R\Gamma) \to U(\Gamma) = \IZ \oplus
\IZ$ is $(1,1)$. An inverse $\overline{\mu}$ of the isomorphism
induced by $\overline{\omega} \colon U(\Gamma) \otimes_{\IZ} \IQ \to
U(\Gamma) \otimes_{\IZ} \IQ$ is given by
$$\squarematrix{2/3}{-1/3}{-1/3}{2/3} \colon \IQ \oplus \IQ \to \IQ \oplus \IQ.$$
The Euler characteristic in the sense of
Leinster~\cite{Leinster(2008)} is $2/3 + (-1/3) + (-1/3) + 2/3 = 2/3$.
We see that the Euler characteristic in the sense of
Leinster~\cite{Leinster(2008)} is the image of the finiteness
obstruction under the composite
$$K_0(R\Gamma) \xrightarrow{\Res} \Split K_0(R\Gamma) \xrightarrow{\rk_R}
U(\Gamma) \xrightarrow{i} U(\Gamma) \otimes_{\IZ} \IQ
\xrightarrow{\overline{\mu}} U(\Gamma) \otimes_{\IZ} \IQ
\xrightarrow{\epsilon} \IQ$$
where $i$ is the obvious inclusion and $\epsilon$
is the augmentation homomorphism.


\typeout{------------ Section 14: A finite category without property
(FP)  -------------}
\section{A finite category without property (FP$_R$)}
\label{sec:A_finite_category_without_property_(FP)}

In this section we investigate the finite category $\IA$ appearing
in~Leinster~\cite[Example~1.11.d]{Leinster(2008)}, recalled below.
Leinster showed that $\IA$ has no weighting. Obviously $\IA$ is
Cauchy complete but not directly-finite and in particular not an
EI-category.  We will show that it is not of type (FP$_R$), give a
full classification of the finitely generated projective
$R\IA$-modules, and compute $K_0(R\IA)$, $G_0(R\IA)$, and
$H_n(B\IA;R) = H_n(\IA;R)$.

The nontrivial morphisms of Leinster's example $\IA$ are drawn in
the diagram below.
$$
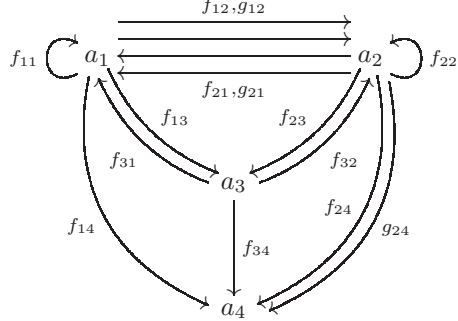
\xymatrix@C=3pc@R=3pc{a_1 \ar@<3.2ex>[rr]^{f_{12}, g_{12}} \ar@<1.7ex>[rr] \ar@(dl,ul)[]^{f_{11}}
\ar@/_.75pc/@<.6ex>[dr]^{f_{13}} \ar@<-.5ex>@/_2pc/[drd]_{f_{14}} &
& a_2 \ar@<-.2ex>[ll] \ar@<1.3ex>[ll]^{f_{21},g_{21}}
\ar@(dr,ur)[]_{f_{22}} \ar@/^.75pc/@<-.6ex>[dl]_{f_{23}}
\ar@<.3ex>@/^2pc/[ddl]_{f_{24}} \ar@<1.5ex>@/^2pc/[ddl]^{g_{24}} \\
& a_3 \ar@/^.75pc/@<.6ex>[ul]^{f_{31}} \ar@/_.75pc/@<-.6ex>[ur]_{f_{32}} \ar[d]^{f_{34}} & \\
& a_4 & }$$ He also defines $f_{33}:=\id_{a_3}$ and
$f_{44}:=\id_{a_4}$. Composition in the category $\IA$ is:
for any composable pair $\xymatrix@1{a_i \ar[r]^p & a_j
\ar[r]^q & a_k}$ in $\IA$ for which neither $p$ nor $q$ is an
identity we have $q \circ p = f_{ik}$.

\begin{lemma} \label{lem:Leinster's_Example_is_h.e._to_a_point}
The space $\vert N \IA \vert$ is homotopy equivalent to a point.
\end{lemma}
\begin{proof}
We consider the subcategory $\IU$ of $\IA$ which does not contain
$g_{24}$, but otherwise is the same as $\IA$. The object $a_4$ is a
terminal object for $\IU$, so $\vert N \IU \vert \simeq \ast$.

But $\vert N \IU \vert \simeq \vert N \IA \vert$. We have the
inclusion $i\colon \IU \to \IA$. The functor $r \colon \IA \to \IU$
is the identity functor, except on $g_{24}$, which $r$ maps to
$f_{24}$. Then $ri = \id_{\IU}$ and we also have a natural
transformation $\alpha\colon ir \Rightarrow \id_{\IA}$ defined by
$$\aligned \alpha(a_1) &= \id_{a_1} \\
\alpha(a_2) &= f_{22} \\
\alpha(a_3) &= \id_{a_3} \\
\alpha(a_4) &= \id_{a_4}. \endaligned$$ The continuous maps $\vert
Nr \vert$ and $\vert N i\vert$ are homotopy inverses.
\end{proof}

Although $\IA$ has the homotopy type of a point, $\IA$ is not
equivalent to the trivial category, for the unique functor $\IA \to
\ast$ is not fully faithful. Alternatively, we note that the trivial
category is of type (FP$_R$) while $\IA$ is not of type (FP$_R$), as
we now show.


\subsection{Property (FP$_R$)}
\label{subsec:Property_(P)}

\begin{theorem}\label{the:Leinsters_example_is_not_FP}
  The above finite category $\IA$ appearing
  in~Leinster~\cite[Examples~1.11.d]{Leinster(2008)} is not of type $(FP_R)$ for any associative, commutative ring $R$ with identity.
\end{theorem}

\begin{proof} In the sequel we use the notation in $\IA$ appearing
  in~Leinster~\cite[Examples~1.11.d]{Leinster(2008)}, recalled above.  Let $M$ be the
  $R\IA$-module $M$ which is uniquely determined by $M(a_i) = \{0\}$
  for $i = 1,3,4$, $M(a_2) = R$, and $M(f_{22}) = 0$. Such an
  $R\IA$-module $M$ exists since $\id_{a_2} = a \circ b$ implies $a =
  b = \id_{a_2}$. Let $u_0\colon R\mor(?,a_4) \to \underline{R}$ be
  the $R\IA$-homomorphism uniquely defined by the property that it
  sends $\id_{a_4}$ to $1 \in R$.  Let $u_1 \colon M \to
  R\mor(?,a_4)$ be the $R\IA$-homomorphism uniquely determined by the
  property that its evaluation at $a_2$ sends $1 \in R = M(a_2)$ to
  $f_{24} - g_{24}$.  Let $v_1 \colon R\mor(?,a_2) \to M$ be the
  $R\IA$-homomorphism uniquely determined by the property that it
  sends $\id_{a_2}$ to $1 \in R = M(a_2)$. Let $v_2 \colon
  R\mor(?,a_1) \to R\mor(?,a_2)$ be the $R\IA$-homomorphism uniquely
  determined by the property that it sends $\id_{a_1}$ to $g_{12} \in
  R\mor(a_1,a_2)$. Let $v_3 \colon M \to R\mor(?,a_1)$ be the
  $R\IA$-homomorphism uniquely determined by the property that its
  evaluation at $a_2$ sends $1 \in R = M(a_2)$ to $f_{21} -
  g_{21}$. Then we obtain exact sequences of $R\IA$-modules
  \begin{eqnarray}
  & 0 \to M \xrightarrow{u_1} R\mor(?,a_4) \xrightarrow{u_0}  \underline{R} \to 0,&
  \label{exact_sequence_M_to-Ra_4_to_underlineR}
  \end{eqnarray}
  and
  \begin{eqnarray}\label{eq:periodic sequence}
  & 0 \to M \xrightarrow{v_3} R\mor(?,a_1) \xrightarrow{v_2} R\mor(?,a_2)
  \xrightarrow{v_1} M \to 0.
  &
  \label{exact_sequence_M-to_ra_1_toRa_2_to_M}
  \end{eqnarray}
  The first exact sequence and~L\"uck~\cite[Lemma~11.6 on
  page~216]{Lueck(1989)} imply that $\underline{R}$ has a finite-dimensional projective
  $R\IA$-resolution if and only if $M$ has. By concatenating copies
  of~\ref{eq:periodic sequence} we obtain an exact sequence
  \[
    0\rightarrow M\rightarrow F_n\rightarrow\dots\rightarrow F_0\rightarrow M\rightarrow 0
  \]
  with free $R\IA$-modules $F_i$ of arbitrarily long length $n$.
  Thus, using~Brown~\cite[Lemma (2.1) on p.~184]{Brown(1982)}, $M$ has a finite-dimensional
  projective $R\IA$-resolution if and only if $M$ is projective.
  Hence $\underline{R}$ has a finite-dimensional projective
  $R\IA$-resolution if and only if $M$ is projective. Since $v_1$ is
  surjective, $M$ is projective only if $v_1$ has a section.  Hence it
  suffices to show that $v_1$ has no section.

  Let $s \colon M \to R\mor(?,a_2)$ be any $R\IA$-homomorphism.
  Consider the homomorphism $g_{12}^* \colon R\mor(a_2,a_2) \to
  R\mor(a_1,a_2)$ given by composition with $g_{12}$. It sends the
  $R$-basis $\{\id_{a_2},f_{22}\}$ bijectively to the $R$-basis
  $\{g_{12},f_{12}\}$ and is hence an isomorphism.  The composite
  $g_{12}^* \circ s(a_2) \colon M(a_2) \to R\mor(a_1,a_2)$ factorizes
  through $M(a_1)$ and hence is trivial since $M(a_1) = \{0\}$.  Hence
  the $R\IA$-morphism $s \colon M \to R\mor(?,a_2)$ is trivial and
  cannot be a section of $v_1$.
\end{proof}


\subsection{Finitely generated projective modules}
\label{subsec:Finitely_generated_projective_modules}

We want to classify all finitely generated projective
$R\IA$-modules.  Let $P$ be a finitely generated projective
$R$-module.  For $i = 1,2$ let $K_1(P)$ be the $R\IA$-module whose
evaluation at both $a_1$ and $a_2$ is $P$ and whose evaluation at
both $a_3$ and $a_4$ is $\{0\}$. We require that $g_{21}$ for $i =
1$ and that $g_{12}$ for $i = 2$ induces the identity $\id \colon P
\to P$, whereas all other morphisms in $\IA$ besides the identity
morphisms of the objects $a_1$ and $a_2$ induce the zero
homomorphism. Then

\begin{theorem} \label{the:fin.-gen.proj.RA-mod} Let $P$ be an
  $R\IA$-module.
  \begin{enumerate}

  \item \label{the:fin.-gen.proj.RA-mod:decomposition}
  $P$ is finitely generated projective if and only if
  there exists finitely generated projective $R$-modules $P_1$, $P_2$,
  $P_3$, and $P_4$ such that
  $$P \cong K_1(P_1) \oplus K_2(P_2) \oplus E_{a_3}(P_3) \oplus E_{a_4}(P_4),$$
  where $E_{a_3}$ and $E_{a_4}$ denote the extension functors defined in
 \eqref{E_x};

  \item \label{the:fin.-gen.proj.RA-mod:uniqueness}
  Suppose that there exists finitely generated projective $R$-modules $P_1$, $P_2$,
  $P_3$, and $p_4$ such that
  $$P \cong K_1(P_1) \oplus K_2(P_2) \oplus E_{a_3}(P_3) \oplus E_{a_4}(P_4).$$
  Then
  \begin{eqnarray*}
  P_1 & \cong & S_{a_1}P;
  \\
  P_2 & \cong & S_{a_1}P;
  \\
  P_3 & \cong & \coker\bigl(P(f_{34}) \colon P(a_4) \to P(a_3)\bigr);
  \\
  P_4 & \cong & P(a_4),
  \end{eqnarray*}
  where $S_{a_i}$ is the splitting functor defined in~\eqref{S_x};

   \item \label{the:fin.-gen.proj.RA-mod:free}
  $P$ is finitely generated free if and only if
  there exists finitely generated free $R$-modules $F_1$, $F_2$,
  $F_3$, and $F_4$ such that
  $$P \cong K_1(F_1) \oplus K_2(F_2) \oplus E_{a_3}(F_1 \oplus F_2 \oplus F_3)
  \oplus E_{a_4}(F_4).$$
\end{enumerate}
\end{theorem}
\begin{proof}\ref{the:fin.-gen.proj.RA-mod:decomposition}
  Recall that the extension functor $E_{a_j}$ satisfies $E_{a_j}(R) =
  R\mor(?,a_j)$, is compatible with direct sums, and sends finitely
  generated projective modules to finitely generated projective
  modules (see Lemma~\ref{lem:basic_properties_of-E_x_and_S_x}~%
\ref{lem:basic_properties_of-E_x_and_S_x:E_x}). In particular
  $E_{a_3}(P_3)$ and $E_{a_4}(P_4)$ are finitely generated projective
  $R\IA$-modules if $P$ is a finitely generated projective $R$-module.

  Given a category $\Gamma$ and an endomorphism $u \colon x \to x$ of an
  object in $\Gamma$ and an $R[x]$-module $Q$, we obtain a morphism of
  $R\Gamma$-modules $u_* \colon E_xQ \to E_xQ$ as follows. Its
  evaluation at an object $y$ is given by
  $$Q \otimes_{R[x]} R\mor(y,x) \to Q \otimes_{R[x]} R\mor(y,x),
  \quad q \otimes v \mapsto q \otimes uv.$$
  Obviously $(\id_x)_* =  \id_{E_xQ}$ and
  $(u_1)_* \circ (u_2)_* = (u_1 \circ u_2)_*$ for two
  endomorphisms $u_1$ and $u_2$.

  Consider a finitely generated projective $R$-module $P$. Consider
  $i \in \{1,2\}$. The construction above applied to the idempotent
  $f_{ii} \colon a_i \to a_i$ yields an idempotent endomorphism of
  $R\IA$-modules $(f_{ii})_* \colon E_{a_i}P \to E_{a_i}P$.  We obtain
  a direct sum decomposition of finitely generated projective
  $R\IA$-modules
  \begin{eqnarray}
  E_{a_i}P & \cong & \im\bigl((f_{ii})_*) \oplus \ker\bigl((f_{ii})_*).
  \label{dec_E_a_3}
  \end{eqnarray}
  Next we show for $i = 1,2$
  \begin{eqnarray}
  \im\bigl((f_{ii})_*\bigr) & \cong & E_{a_3}P;
  \label{im(f_ii)}
  \\
  \ker\bigl((f_{ii})_*\bigr) & \cong & K_i(P).
  \label{ker(f_ii)}
  \end{eqnarray}
  We only treat the case $i = 1$, the case $i = 2$ is completely
  analogous.  Let
  \begin{eqnarray}
  & \alpha \colon E_{a_3}P \to E_{a_1}P &
  \label{map_alpha}
  \end{eqnarray}
  be the $R\Gamma$-homomorphism which is the adjoint under the adjunction
  of~L\"uck~\cite[Lemma~9.31~a) on page~171]{Lueck(1989)} of the
  $R$-homomorphism $P \to E_{a_1}P(a_3) = P \otimes_R R\mor(a_3,a_1)$
  sending $p$ to $p \otimes f_{31}$. Explicitly the evaluation of
  $\alpha$ at an object $a_j$ is given by
  $$P \otimes_R R\mor(a_j,a_3) \to P \otimes_R R\mor(a_j,a_1), \quad
  p \otimes u \mapsto p \otimes (f_{31} \circ u).$$ One easily checks
  that $\alpha$ is injective. The image of $\alpha(a_j)$ is $\{0\}$
  for $j = 4$ and is $\{p \otimes f_{j1} \mid p \in P\}$ for
  $j =1,2,3$. This is the same as the image of
  $(f_{11})_* \colon  E_{a_1}P \to E_{a_1}P$ and~\eqref{im(f_ii)} follows. The
  cokernel of $\alpha$ is isomorphic to $\ker\bigl((f_{11})_*\bigr)$ since
  $(f_{11})_*$ is an idempotent.  Obviously the cokernel evaluated at
  $a_4$ and $a_3$ is $\{0\}$. The cokernel evaluated at the objects
  $a_1$ and $a_2$ is isomorphic to $R$. The element $\id_{a_1}$
  projects down to a generator in $\coker(\alpha)(a_1)$ and the element $g_{21}$
  projects down to a generator in $\coker(\alpha)(a_2)$. Hence the
  morphism $g_{21}$ induces a map $\coker(\alpha)(a_1)$ to
  $\coker(\alpha)(a_2)$ that respects these generators.  The morphisms
  $f_{11}$, $f_{12}$, $f_{22}$ and $g_{12}$ induce the trivial
  homomorphism on the cokernel of $\alpha$. Now~\eqref{ker(f_ii)}
  follows.

  In particular we see that $K_i(P)$ is a finitely generated
  projective $R\IA$-module if $P$ is a finitely generated projective
  $R$-module.

  Now consider a finitely generated projective $R\IA$-module
  $P$. Choose a finitely generated free $R\Gamma$-module $F$ together
  with $R\Gamma$-maps $i \colon P \to F$ and $r \colon F \to P$. Let
  $$\sigma_{a_4}(P) \colon E_{a_4}P(a_4) \to P$$ be the adjoint of the
  adjunction of~L\"uck~\cite[Lemma~9.31 on page~171]{Lueck(1989)} of the
  $R$-homomorphism $\id_{a_4} \colon P(a_4) \to P(a_4)$. Explicitly
  its evaluation at $a_j$ is given by
  $$P(a_4) \otimes_R R\mor(a_j,a_4) \to P(a_j), \quad p \otimes u \mapsto P(u)(p).$$
  The map $\sigma_{a_4}(P)$ is natural in $P$. Let $\overline{P}$ and
  $\overline{F}$ respectively be the cokernel of $\sigma_{a_4}(P)$ and
  $\sigma_{a_4}(F)$ respectively. Denote by $\pr(P) \colon P \to
  \overline{P}$ and $\pr(F) \colon P \to \overline{F}$ the canonical
  projections.

  Choose non-negative integers $m_1$, $m_2$, $m_3$, and $m_4$ such that
  $$F \cong \bigoplus_{j=1}^4 R\mor(?,a_j)^{m_j}.$$
  Since the are no morphisms from $a_4$ to the other objects $a_1$, $a_2$ and $a_3$,
  one easily checks that the  sequence
  $$E_{a_4}F(a_4) \xrightarrow{\sigma_{a_4}} F \xrightarrow{\pr(F)} \overline{F}$$
  can be identified with the obvious  split exact sequence
  $$R\mor(?,a_4)^{m_4} \to \bigoplus_{j=1}^4 R\mor(?,a_j)^{m_j} \to
  \bigoplus_{j=1}^3 R\mor(?,a_j)^{m_j}.$$
  We obtain a commutative diagram
  $$\xymatrix@C=18mm{
  0 \ar[d] & 0 \ar[d] & 0 \ar[d]
  \\
  E_{a_4}(P(a_4)) \ar[r]^{E_{a_4}(i(a_4))} \ar[d]^{\sigma_{a_4}(P)}
  &
  E_{a_4}(F(a_4)) \ar[r]^{E_{a_4}(r(a_4))} \ar[d]^{\sigma_{a_4}(F)}
  &
  E_{a_4}(P(a_4)) \ar[d]^{\sigma_{a_4}(P)}
  \\
  P \ar[r]^{i} \ar[d]^{\pr(P)}
  &
  F \ar[r]^{r} \ar[d]^{\pr(F)}
  &
  P \ar[d]^{\pr(P)}
  \\
  \overline{P} \ar[r]^{\overline{i}} \ar[d]
  &
  \overline{F} \ar[r]^{\overline{r}} \ar[d]
  &
  P\ar[d]
  \\
  0 & 0 & 0}
$$
where $\overline{i}$ and $\overline{r}$ are the maps induced by $i$
and $r$.  We know already that the middle row is exact. We conclude
$E_{a_4}(r(a_4)) \circ E_{a_4}(i(a_4)) = \id$ and $\overline{r}
\circ \overline{i} = \id$ from $r \circ i =\id$. An easy diagram
shows that all rows are exact.

Hence $\overline{P}$ is a finitely generated projective
$R\IA$-module, we have the isomorphisms
\begin{eqnarray*}
  P & \cong & E_{a_4}(P(a_4)) \oplus \overline{P};
  \label{reducing_P_to_overlineP}
  \\
  \overline{F} & \cong & \bigoplus_{j=1}^3 R\mor(?,a_j)^{m_j},
  \label{overlineF}
\end{eqnarray*}
and $R\IA$-homomorphisms $\overline{i} \colon \overline{P} \to
\overline{F}$ and $\overline{r} \colon \overline{F} \to
\overline{P}$ with $\overline{r} \circ \overline{i} = \id$. The
$R$-module $P(a_4)$ is a finitely generated projective $R$-module
since it is a direct summand in the finitely generated free
$R$-module $F(a_4) = R^{m_4}$. Hence it suffices to proof the claim
for $\overline{P}$.

Now we more or less repeat the argument above, but nor replacing
$a_4$ by $a_3$.  So we define
\begin{eqnarray*}
  & \sigma_{a_3}(\overline{P}) \colon E_{a_3}\overline{P}(a_3) \to \overline{P} &
  \\
  & \sigma_{a_3}(\overline{P})  \colon E_{a_3}\overline{F}(a_3) \to \overline{F} &
\end{eqnarray*}
as above. Denote by $\overline{\overline{P}}$ and
$\overline{\overline{F}}$ respectively the cokernel of
$\sigma_{a_3}(\overline{P})$ and $\sigma_{a_3}(\overline{F})$
respectively. Let $\pr(\overline{P}) \colon \overline{P} \to
\overline{\overline{P}}$ and $\pr(\overline{F}) \colon \overline{F}
\to \overline{\overline{F}}$ be the canonical projections. Denote by
$\overline{\overline{i}}\colon \overline{\overline{P}} \to
\overline{\overline{F}}$ and $\overline{\overline{r}}\colon
\overline{\overline{F}} \to \overline{\overline{P}}$ the maps
induced by $\overline{i}$ and $\overline{r}$.  The maps
$\sigma_{a_3}(P)$ are natural in $P$ and compatible with direct
sums.  One easily checks that the $R\IA$-homomorphism
$\sigma_{a_3}(R\mor(?,a_3))$ is an isomorphism. Hence also the
$R\IA$-homomorphism
$$\sigma_{a_3}(R\mor(?,a_3)^{m_3}) \colon E_{a_3}R\mor(a_3,a_3)^{m_3}
\to R\mor(?,a_3)^{m_3}$$ is an isomorphism. The map
$\sigma_{a_3}(R\mor(?,a_1)) \colon E_{a_3}R\mor(a_3,a_1) \to
R\mor(?,a_1)$ is the same as the map $\alpha$ defined in~\eqref{map_alpha}.
Hence it is injective and its cokernel is
$K_1(R)$. This implies that
$$\sigma_{a_3}(R\mor(?,a_1)^{m_1}) \colon E_{a_3}R\mor(a_3,a_1)^{m_1}
\to R\mor(?,a_1)^{m_1}$$ is injective with the finitely generated
projective $R\IA$-module $K_1(R^{m_1})$ as cokernel. Analogously one
shows that
$$\sigma_{a_3}(R\mor(?,a_2)^{m_2}) \colon E_{a_3}R\mor(a_3,a_2)^{m_2}
\to R\mor(?,a_2)^{m_2}$$ is injective with the finitely generated
projective $R\IA$-module $K_2(R^{m_2})$ as kernel. This implies
$$\overline{F} \cong K_1(R^{m_1}) \oplus K_2(R^{m_2}).$$
As above we obtain a commutative diagram with exact rows

$$\xymatrix@C=18mm{
  0 \ar[d] & 0 \ar[d] & 0 \ar[d]
  \\
  E_{a_3}(\overline{P}(a_3)) \ar[r]^{E_{a_3}(\overline{i}(a_3))}
  \ar[d]^{\sigma_{a_3}(\overline{P})} & E_{a_3}(\overline{F}(a_3))
  \ar[r]^{E_{a_3}(\overline{r}(a_3))}
  \ar[d]^{\sigma_{a_3}(\overline{F})} & E_{a_3}(\overline{P}(a_3))
  \ar[d]^{\sigma_{a_3}(\overline{P})}
  \\
  \overline{P} \ar[r]^{\overline{i}} \ar[d]^{\pr(\overline{P})} &
  \overline{F} \ar[r]^{\overline{r}} \ar[d]^{\pr(\overline{F})} &
  \overline{P} \ar[d]^{\pr(\overline{P})}
  \\
  \overline{\overline{P}} \ar[r]^{\overline{\overline{i}}} \ar[d] &
  \overline{\overline{F}} \ar[r]^{\overline{\overline{r}}} \ar[d] &
  \overline{\overline{P}}\ar[d]
  \\
  0 & 0 & 0}
$$
Hence $\overline{\overline{P}}$ is a finitely generated projective
$R\IA$-module with is a direct summand in $\overline{F} \cong
K_1(R^{m_1}) \oplus K_2(R^{m_2})$ and we obtain an isomorphism
$$\overline{P} \cong
E_{a_3}(\overline{P}(a_3)) \oplus \overline{\overline{P}}.$$ Since
$\overline{P}(a_3)$ is a direct summand in the finitely generated
free $R$-module $\overline{F}(a_3) \cong R^{m_1 + m_2 + m_3}$, it is
finitely generated projective $R$-module.  Hence it remains to prove
the claim for $\overline{\overline{P}}$.

Since $\overline{\overline{P}}$ is a direct summand in $K_1(R^{m_1})
\oplus K_2(R^{m_2})$, one easily checks that we have exact sequences
of finitely generated projective $R$-modules
$$0 \to \im\bigl(\overline{\overline{P}}(g_{12})\bigr) \xrightarrow{i_1}
\overline{\overline{P}}(a_1)
\xrightarrow{\overline{\overline{P}}(g_{21})}
\im\bigl(\overline{\overline{P}}(g_{21})\bigr) \to 0,$$ and
$$0 \to \im\bigl(\overline{\overline{P}}(g_{21})\bigr) \xrightarrow{i_2}
\overline{\overline{P}}(a_2)
\xrightarrow{\overline{\overline{P}}(g_{12})}
\im\bigl(\overline{\overline{P}}(g_{12})\bigr) \to 0,$$ where $i_1$
and $i_2$ are the inclusions. Choose $R$-maps
\begin{eqnarray*}
r_1 \colon \overline{\overline{P}}(a_1) & \to &
\im\bigl(\overline{\overline{P}}(g_{12})\bigr),
\\
r_2 \colon \overline{\overline{P}}(a_2) & \to &
\im\bigl(\overline{\overline{P}}(g_{21})\bigr),
\end{eqnarray*}
satisfying $r_1 \circ i_1 = \id$ and $r_2 \circ i_2 = \id$. Next we
define an $R\IA$-isomorphism
$$\beta \colon \overline{\overline{P}} \xrightarrow{\cong}
K_1\left(\im\bigl(\overline{\overline{P}}(g_{21})\bigr)\right)
\oplus
K_2\left(\im\bigl(\overline{\overline{P}}(g_{12})\bigr)\right).
$$
Its evaluation at $a_1$ is given by the $R$-isomorphism
$$\overline{\overline{P}}(g_{21}) \oplus r_1 \colon
\overline{\overline{P}}(a_1) \xrightarrow{\cong}
\im\bigl(\overline{\overline{P}}(g_{21})\bigr) \oplus
\im\bigl(\overline{\overline{P}}(g_{12})\bigr)
$$
and its evaluation at $a_2$ by the  $R$-isomorphism
$$r_2 \oplus \overline{\overline{P}}(g_{12}) \colon
\overline{\overline{P}}(a_2) \xrightarrow{\cong}
\im\bigl(\overline{\overline{P}}(g_{21})\bigr) \oplus
\im\bigl(\overline{\overline{P}}(g_{12})\bigr).
$$
This finishes the proof of
assertion~\ref{the:fin.-gen.proj.RA-mod:decomposition} of
Theorem~\ref{the:fin.-gen.proj.RA-mod}.
\\[1mm]~\ref{the:fin.-gen.proj.RA-mod:uniqueness}
Recall that $K_i(P_i)$ is a direct summand in $E_{a_i}(P_i)$ for $i
= 1,2$ (see~\eqref{ker(f_ii)}).
Using Lemma~\ref{lem:basic_properties_of-E_x_and_S_x}~%
\ref{lem:basic_properties_of-E_x_and_S_x:S_x_circE_x_of_free} one
easily checks
$$\begin{array}{lclcl}
S_{a_i}(P) & \cong & S_{a_i}(K_i(P_i)) &\cong & P_i \quad
\text{for}\; i = 1,2;
\\
P(a_4) & \cong & P_4. & &
\end{array}
$$
A direct computation shows
\begin{eqnarray*}
\lefteqn{\coker\bigl(P(f_{34}) \colon P(a_4) \to P(a_3)\bigr)} & &
\\
& \cong & \bigoplus_{i=1}^2\coker\bigl(K_i(p_i)(f_{34})\bigr) \oplus
\coker\bigl(E_{a_3}(P_3)(f_{34})\bigr) \oplus
\coker\bigl(E_{a_4}(P_4)(f_{34})\bigr)
\\
& \cong & \coker\bigl(E_{a_3}(P_3)(f_{34})\bigr)
\\
& \cong & P_3.
\end{eqnarray*}
This finishes the proof of
assertion~\ref{the:fin.-gen.proj.RA-mod:uniqueness}.
\\[1mm]~\ref{the:fin.-gen.proj.RA-mod:free}
This follows from
assertions~\ref{the:fin.-gen.proj.RA-mod:decomposition}
and~\ref{the:fin.-gen.proj.RA-mod:uniqueness} and the isomorphism
for $i = 1,2$ (see~\eqref{dec_E_a_3}, \eqref{im(f_ii)}
and~\eqref{ker(f_ii)})
$$R\mor(?,a_i) \cong R\mor(?,a_3) \oplus K_1(R).$$
This finishes the proof of Theorem~\ref{the:fin.-gen.proj.RA-mod}.
\end{proof}

\begin{remark}\label{rem:not_natural}
Notice that the decomposition of
Theorem~\ref{the:fin.-gen.proj.RA-mod}~\ref{the:fin.-gen.proj.RA-mod:decomposition}
is not natural in $P$. However, the cofiltration by epimorphisms
$$P \to \overline{P} \to \overline{\overline{P}}$$
and the identifications
\begin{eqnarray*}
\overline{\overline{P}} & \cong & K_1(S_{a_1}(P)) \oplus
K_2(S_{a_2}(P));
\\
\ker\bigl(\overline{P} \to
\overline{P}/\overline{\overline{P}}\bigr) & \cong &
E_{a_3}\left(\coker\bigl(P(f_{34}) \colon P(a_4) \to
P(a_3)\bigr)\right);
\\
\ker\bigl(P \to \overline{P}\bigr) & \cong & E_{a_4}(P(a_4)),
\end{eqnarray*}
are natural in $P$.
\end{remark}

Let $K_0^f(R\IA)$ be the Grothendieck group of finitely generated
free $R\IA$-modules.  Let
$$\iota \colon U(\Gamma) \to K_0^f(R\IA)$$
be the homomorphism which sends a basis element $\overline{x} \in
\iso(\IA)$ to the class of $R\mor(?,x)$.

\begin{theorem}[$K_0(R\IA)$]
\label{the:K_0(RIA)}
\begin{enumerate}
\item\label{the:K_0(RIA):K_0}
The maps
\begin{eqnarray*}
\xi \colon K_0(R)^4 &\xrightarrow{\cong} &  K_0(R\IA),
\\
\eta \colon K_0(R\IA) & \xrightarrow{\cong} & K_0(R)^4,
\end{eqnarray*}
defined by
\begin{eqnarray*}
\xi\bigl([P_1], [P_2], [P_3], [P_4]\bigr) & = & [K_1(P_1)] +
[K_2(P_2)] + [E_{a_3}(P_3)] + [E_{a_4}(P_4)],
\\
\eta([P]) & = & \left([S_{a_1}P], [S_{a_2}P],
\bigl[\coker\bigl(P(f_{34}) \colon P(a_4) \to P(a_3)\bigr)\bigr] ,
[S_{a_4}P]\right),
\end{eqnarray*}
are isomorphisms, inverse to another.
\item\label{the:K_0(RIA):U(IA)_and_K_of(RIA)}
The map
$$\iota \colon U(\IA) \xrightarrow{\cong}  K_0^f(R\IA)$$
is bijective. If $R$ is  a principal domain, then the forgetful map
$$F^f \colon K_0^f(R\IA) \xrightarrow{\cong} K_0(R\IA)$$
is bijective.
\end{enumerate}
\end{theorem}
\begin{proof}~\ref{the:K_0(RIA):K_0}
This follows from
Theorem~\ref{the:fin.-gen.proj.RA-mod}~\ref{the:fin.-gen.proj.RA-mod:decomposition}
and~\ref{the:fin.-gen.proj.RA-mod:uniqueness}.
\\[1mm]~\ref{the:K_0(RIA):U(IA)_and_K_of(RIA)}
The map $\iota$ is obviously surjective. The composite
$$U(\Gamma) \xrightarrow{\iota} K_0^f(R\IA) \xrightarrow{F^f} K_0(R\IA) \xrightarrow{\eta}
K_0(R)^4 \xrightarrow{\rk_R} \IZ^4$$ can be identified with the
injection
$$\IZ^4 \xrightarrow{\cong} \IZ^4, \quad \bigl(m_1, m_2, m_3, m_4\bigr)
\mapsto \bigl(m_1, m_2, m_1 + m_2 +m_3, m_4\bigr)$$ by
Theorem~\ref{the:fin.-gen.proj.RA-mod}~\ref{the:fin.-gen.proj.RA-mod:free}.
The forgetful map $F^f \colon K_0^f(R\IA) \to K_0(R\IA)$ is
surjective by
Theorem~\ref{the:fin.-gen.proj.RA-mod}~\ref{the:fin.-gen.proj.RA-mod:free}
provided that $R$ is an integral domain and hence $\IZ \to K_0(R),\;
n \mapsto [R^n]$ is an isomorphism. This finishes the proof of
Theorem~\ref{the:K_0(RIA)}.
\end{proof}


\subsection{$K_0$ versus $G_0$}
\label{subsec:K_0_versus_G_0}

    Let $R$ be a commutative Noetherian ring and let $\Gamma$ be a finite category
    (see Definition~\ref{def:finite_quasi-finite_and_free_category}).
    Denote by $G_0(\IQ\Gamma)$ the  Grothendieck group of
    finitely generated $\IQ\Gamma$-modules. Since
    $\Gamma$ is finite, an $R\Gamma$-module is finitely generated if and only if
    for every object $x$ the $\IQ$-module $M(x)$ is finitely generated as an $R$-module.
    In particular the category of $R\Gamma$-modules is Noetherian, i.e., a submodule
    of a finitely generated $R\Gamma$-module is finitely generated
    (see~L\"uck~\cite[Lemma~16.10 on page~327]{Lueck(1989)}).

    \begin{remark} \label{rem:element_in_G_0(R)}
    Notice that the constant
    $R$-module $\underline{R}$ defines an element $[\underline{R}]$ in $G_0(R\Gamma)$
    which may be viewed as a kind of analogue of the finiteness obstruction.
    Only if $\Gamma$ is of type (FP$_R$), then we get also an element
    $o(\Gamma;R) := [\underline{R}]$ in $K_0(R\Gamma)$
    which is mapped under the forgetful homomorphism
    \begin{eqnarray*}
    F_{R\Gamma} \colon K_0(R\Gamma) & \to & G_0(R\Gamma).
    \label{K_0(RGamma)_to_G_0(RGamma)}
    \end{eqnarray*}
    to $[\underline{R}] \in G_0(R\Gamma)$.
    \end{remark}

    Notice that $F_{R\Gamma}$ is bijective if $\Gamma$ is a finite
    EI-category and the order $\aut(x)$ is invertible in $R$ for every object
    $x$ in $\Gamma$ (see~L\"uck~\cite[Proposition~16.28 on page~332]{Lueck(1989)}).
    This is not true in general as the following example shows.

    \begin{example}\label{exa:K_0_versus_G_0_for_IA}
    We conclude from~\eqref{exact_sequence_M-to_ra_1_toRa_2_to_M} that
    \begin{eqnarray}
    [R\mor(?,a_1)] = [R\mor(?,a_2)] & & \in G_0(R\IA).
    \label{Qa_1_and_Qa_2_agree_in_G_0}
    \end{eqnarray}
    This together with Theorem~\ref{the:K_0(RIA)}~\ref{the:K_0(RIA):U(IA)_and_K_of(RIA)}
    implies that
    $$F \colon K_0 (R\IA) \to G_0(R\IA)$$
    is not injective.
    \end{example}

Define a map
\begin{eqnarray}
&& \Res \colon G_0(R\Gamma) \to \bigoplus_{\overline{x} \in
\iso(\Gamma)} G_0(R[x]), \quad [P] \mapsto \{[P(x)] \mid
\overline{x} \in \iso(\Gamma)\}. \label{res_for_G_0}
\end{eqnarray}

Provided that the order $\aut(x)$ is invertible in $R$ for every
object $x$ in $\Gamma$, we also obtain a map
\begin{eqnarray}
& & \Res \colon K_0(R\Gamma) \to \bigoplus_{\overline{x} \in
\iso(\Gamma)} K_0(R[x]), \quad [P] \mapsto \{[P(x)] \mid
\overline{x} \in \iso(\Gamma)\}, \label{res_for_K_0}
\end{eqnarray}
and we get a commutative diagram
  $$\xymatrix@C=20mm{
  K_0(R\Gamma) \ar[r]^{F_{R\Gamma}} \ar[d]^{\Res}
  & G_0(R\Gamma) \ar[d]^{\Res}
  \\
  \bigoplus_{\overline{x} \in \iso(\Gamma)} K_0(R[x])
  \ar[r]^{\bigoplus_{\overline{x} \in \iso(\Gamma)} F_{R[x]}}_{\cong}
  &
  \bigoplus_{\overline{x} \in \iso(\Gamma)} G_0(R[x])
  }$$
whose lower horizontal arrow is an isomorphism.

Now we consider the special case $\Gamma = \IA$ and $R = \IQ$. For a
$\IQ$-module $P$ and $k \in \{1,2,4\}$ denote by $I_k(P)$ the
$\IQ\IA$-module for which $I_k(\IQ)(a_k) = \IQ$, $I_k(\IQ)(a_j) =
\{0\}$ for $j \not= k$ and all morphisms except $\id_{a_k}$ induce
the trivial homomorphism. One easily checks that this is a
well-defined $\IQ$-module. (Notice that this definition does not
make sense for the object $a_3$).

\begin{theorem}[$G_0(\IQ\IA)$] \label{the:G_0(QA)}
    The homomorphisms
    \begin{multline*}
    \omega \colon \IZ^4 \to G_0(R\IA),
    \\
    (n_1,n_2,n_3.n_4) \mapsto n_1 \cdot [I_1(\IQ)] + n_2 \cdot
    [I_2(\IQ)] + n_3 \cdot [R\mor(?,a_3)] + n_4 \cdot [I_4(\IQ)]
  \end{multline*}
    and the composite
    $$G_0(\IQ\IA) \xrightarrow{\Res} \bigoplus_{i=1}^4 G_0(\IQ)
    \xrightarrow{\bigoplus_{i=1}^4 \rk_{\IQ}} \IZ^4$$
    are isomorphisms.
\end{theorem}
\begin{proof}
  The composite
$$\IZ^4 \xrightarrow{\omega} G_0(R\IA) \xrightarrow{\Res} \bigoplus_{i=1}^4 G_0(\IQ)
\xrightarrow{\bigoplus_{i=1}^4 \rk_{\IQ}} \IZ^4$$ sends
$(m_1,m_2,m_3,m_4)$ to $(m_1 + m_3,m_2 + m_3, m_3, m_4)$ and is
hence an isomorphism.  Therefore it suffices to show that $\omega$
is surjective.

Consider a finitely generated $\IQ\IA$-module $M$. There is the
epimorphism of $\IQ\IA$-modules $M \to I_4(M(a_4))$ whose evaluation
at $a_4$ is the identity.  Let $N$ be its kernel. Then we get $[M] =
[N] + [I_4(M(a_4)]$ in $G_0(\IQ\IA)$ and $N(a_4) = \{0\}$. Hence it
suffices to prove that $[N]$ lies in the image of $\omega$.

Consider the $\IQ\IA$-homomorphism $f \colon E_3(N(a_3)) \to N$
uniquely determined by the property that its evaluation at $a_3$ is
the isomorphism $N(a_3) \otimes_{\IQ} \IQ\mor(a_3,a_3)
\xrightarrow{\cong} N(a_3)$ sending $x \otimes \id_{a_3}$ to $x$.
Let $K$ be its kernel and $L$ be its cokernel. We get in $[N] =
[E_3(N(a_3))] + [L] - [K]$ in $G_0(\IQ\IA)$ and $K(a_3) = K(a_4) =
L(a_3) = L(a_4) = \{0\}$. Hence it suffices to show that $K$ lies in
the image of $\omega$ if $K$ is a finitely generated $\IQ\IA$-module
with $K(a_3) = K(a_4) = 0$.

Notice that the all morphisms in $\IA$ possibly except $g_{12}$ and
$g_{21}$ and the identity morphisms for $a_1$ and $a_2$ induce the
trivial homomorphism on $K$ since they factor through the object
$a_3$ or $a_4$ and $K(a_3) = K(a_4) = 0$. Consider the
$\IQ\IA$-homomorphism
$$g \colon I_1\bigl(\ker(N(g_{21})\bigr) \to K$$
given by the inclusion $\ker(N(g_{21})) \to N(a_1)$. Let $P$ be its
cokernel. By construction the map $P(g_{21}) \colon P(a_1) \to
P(a_2)$ is injective. Since $P(a_3) = 0$, we get
$$P(g_{21}) \circ P(g_{12}) = P(g_{12} \circ g_{21}) = P(f_{11})
= P(f_{31} \circ f_{13}) = P(f_{13}) \circ P(f_{31}) = 0.$$ Since
$P(g_{21})$ is injective, $P(g_{12}) = 0$. Hence the identity on
$P(a_2)$ induces an injection $I_{a_2}(P(a_2)) \to P$. Let $Q$ be
its cokernel. Then $Q(a_2) = Q(a_3) = Q(a_4)$. This implies $Q =
I_{a_1}(Q(a_1))$. Hence we get in $G_0(\IQ\IA)$
$$[K] =   [I_{a_1}(\ker(N(g_{21}))] + [I_{a_2}(P(a_2))] + [I_{a_1}(Q(a_1))].$$
This finishes the proof of Theorem~\ref{the:G_0(QA)}.
\end{proof}

\begin{example} Put $R = \IQ$ and $\Gamma = \IA$. Then the following diagram
commutes
$$\xymatrix@C=20mm{
  U(\IA) = \IZ^4 \ar[r]^{\iota}_{\cong}  \ar[rdd]^A
  & K_0(\IQ\IA) \ar[r]^{F_{\IQ\IA}} \ar[d]^{\Res}
  & G_0(\IQ\IA) \ar[d]^{\Res}_{\cong}
  \\
  & \bigoplus_{i=1}^4 K_0(\IQ) \ar[d]^{\bigoplus_{i=1}^4\rk_{\IQ}}_{\cong}
  \ar[r]^{\bigoplus_{\overline{x} \in \iso(\IA)} F_{\IQ}}_{\cong}
  &
  \bigoplus_{i=1}^4 G_0(\IQ) \ar[d]^{\bigoplus_{i=1}^4\rk_{\IQ}}_{\cong}
  \\
  & \IZ^4 \ar[r]^{\id}_{\cong} & \IZ^4
  }$$
where $A$ is given by the matrix
$$\left(\begin{matrix}
2 & 2 & 1 & 1
\\
2 & 2 & 1 & 2
\\
1 & 1 & 1 & 1
\\
0 & 0 & 0 & 1
\end{matrix}\right)$$
Notice that $(1,1,1,1)$ is not in the image of $A \colon \IZ^4 \to
\IZ^4$. Obviously $[\underline{\IQ}] \in G_0(\IQ\IA)$ is sent under
the composite
$$\bigoplus_{i=1}^4 \rk_{\IQ} \circ \Res \colon G_0(\IQ\IA) \to \IZ^4$$
to $(1,1,1,1)$. Hence we see again that $\IA$ is not of type
(FP$_R$), since otherwise $[\underline{R}] \in G_0(\IQ\IA)$ lies in
the image of $F_{\IQ\IA}$ and hence  $(1,1,1,1)$ lies in the image
of $A \colon \IZ^4 \to \IZ^4$.

\end{example}


\subsection{Homology of $\IA$}
\label{subsec:homology_of_IA}

We obtain from the short exact
sequence~\eqref{exact_sequence_M-to_ra_1_toRa_2_to_M} the following
periodic projective resolution $P_*$ of the $R\IA$-module $M$
\begin{multline*}
\cdots \xrightarrow{v_3 \circ v_1} R\mor(?,a_1) \xrightarrow{v_2}
R\mor(?,a_2) \xrightarrow{v_3 \circ v_1} R\mor(?,a_1)
\xrightarrow{v_2} R\mor(?,a_2) \xrightarrow{v_1} M.
\end{multline*}
Recall that $v_2$ sends $\id_{a_1}$ to $g_{12}$ and $v_3 \circ v_1$
sends $\id_{a_2}$ to $f_{21}-g_{21}$. The $R$-chain complex
$P_*\otimes_{R\IA} \underline{R}$ looks like
$$\cdots \xrightarrow{0} R \xrightarrow{\id} R
\xrightarrow{0} R \xrightarrow{\id} R$$ Hence we get for $n \ge 0$
\begin{equation} \label{equ:homology}
H_n^R(\IA;M) := H_n\bigl(P_* \otimes_{R\IA} M\bigr) =\{0\}.
\end{equation}
We conclude from $\underline{R} \otimes_{R\IA} \underline{R} \cong
R$, from~\eqref{equ:homology}, and the short exact
sequence~\eqref{exact_sequence_M_to-Ra_4_to_underlineR} that
$$H_n(B\IA;R) = H_n(\IA;R) = H_n^R(\IA;\underline{R}) =
\begin{cases}
R & \text{if} \; n = 0
\\
\{0\} & \text{if} \; n > 0,
\end{cases}$$
as we may expect from the contractibility of $B\IA$.


\typeout{---------------------------- References ---------------------------------------}

\bibliographystyle{abbrv}


\end{document}